\documentclass[titlepage,twoside,headsepline,dvips]{scrartcl}
\usepackage{latexsym, amssymb,amsmath,german}
\usepackage[arrow, matrix, curve]{xy}
\usepackage{makeidx}
\def\GLn{\operatorname{GL}_{n}(K)}
\def\Sp{\operatorname{Sp}}

\def\On{\operatorname{O}}
\def\SL2{\operatorname{SL}_{2}(K)}
\def\SO{\operatorname{SO}}
\def\GL2{\operatorname{GL}_{2}(K)}
\def\Ga{{\mathbb G}_{a}}
\def\Gm{{\mathbb G}_{m}}
\def\INVSL2{$K[V]^{operatorname{SL}_{2}(K)}$}
\def\INVSO2{$K[V]^{operatorname{SO}_{2}(K)}$}
\def\INVGL2{$K[V]^{operatorname{GL}_{2}(K)}$}

\def\Bew{\noindent \textit{Beweis. }}
\def\qed{\hfill $\Box$}
\def\Magma{{\sc Magma }} 

\def\depth{\operatorname{depth}}
\def\height{\operatorname{height}}
\def\cmdef{\operatorname{cmdef}}
\def\trdeg{\operatorname{trdeg}}
\def\Hom{\operatorname{Hom}}
\def\GL{\operatorname{GL}}
\def\SL{\operatorname{SL}}
\def\id{\operatorname{id}}
\def\projdim{\operatorname{pd}}
\def\chr{\operatorname{char}}
\def\Quot{\operatorname{Quot}}
\def\bigP{\mathcal{P}}
\def\diag{\operatorname{diag}}
\def\rang{\operatorname{rang}}
\def\sgn{\operatorname{sgn}}
\def\LM{\operatorname{LM}}
\def\spol{\operatorname{spol}}
\def\Spur{\operatorname{Spur}}
\def\Ass{\operatorname{Ass}}
\def\Ann{\operatorname{Ann}}
\def\rang{\operatorname{rang}}
\def\im{\operatorname{im}}
\def\Ext{\operatorname{Ext}}
\def\myforall{\textrm{ f"ur alle }}

\newtheorem{Lemma}{Lemma}[section]
\newtheorem{Satz}[Lemma]{Satz}
\newtheorem{Def}[Lemma]{Definition}
\newtheorem{Korollar}[Lemma]{Korollar}
\newtheorem{Bemerkung}[Lemma]{Bemerkung}
\newtheorem{SatzDef}[Lemma]{Satz und Definition}
\newtheorem{Prop}[Lemma]{Proposition}
\newtheorem{Hauptsatz}[Lemma]{Hauptsatz}
\newtheorem{Bsp}[Lemma]{Beispiel}

\newtheorem{BemDef}[Lemma]{Bemerkung und Definition}

\newenvironment{Algorithmus}[1][Algorithmus \arabic{section}.\arabic{Lemma}]{\refstepcounter{Lemma}\begin{trivlist}
\item[\hskip \labelsep {\bfseries #1}]}{\end{trivlist}}

\newenvironment{BemRoman}[1][Bemerkung \arabic{section}.\arabic{Lemma}]{\refstepcounter{Lemma}\begin{trivlist}
\item[\hskip \labelsep {\bfseries #1}]}{\end{trivlist}}

\newenvironment{BspRoman}[1][Beispiel \arabic{section}.\arabic{Lemma}]{\refstepcounter{Lemma}\begin{trivlist}
\item[\hskip \labelsep {\bfseries #1}]}{\end{trivlist}}

\makeindex
\begin{document}
\pagenumbering{roman}
\begin{titlepage}
\pagestyle{empty}
\vspace{-7cm}
\begin{center}
{\Large
Technische Universit\"at M\"unchen\\
Zentrum Mathematik}
\end{center}
\vspace{2cm}
\begin{center}
\bf\LARGE  "Uber die Tiefe von Invariantenringen unendlicher Gruppen\\

\end{center}
\vspace{0,8cm}
\begin{center}
{\Large Martin Wilhelm Kohls\\}
\end{center}
\vspace{2.3cm}
\noindent Vollst\"andiger Abdruck der von der Fakult\"at f\"ur
Mathematik der Technischen Universit\"at M\"unchen zur Erlangung des
akademischen Grades eines
\begin{center}
Doktors der Naturwissenschaften (Dr.\,rer.\,nat.)
\end{center}
genehmigten Dissertation.\vspace{3,3cm}\\
\begin{tabular}{lrl}
$\qquad$ Vorsitzender: && Univ.-Prof. Dr. Peter Rentrop\\
\vspace{-.4cm}\\
$\qquad$ Pr\"ufer der Dissertation:
& 1. & Univ.-Prof. Dr. Gregor Kemper\\ \vspace{-,4cm}\\
& 2. & Prof. Dr. Harm Derksen\\ \vspace{-.4cm}\\
&& University of Michigan, USA\\ \vspace{-.4cm}\\
&&(schriftliche Beurteilung)\\ \vspace{-.4cm}\\
& 3. & Univ.-Prof. Dr. Werner Heise\\ \vspace{-.4cm}\\
&&(m"undliche Pr"ufung)\\
&&\\
\end{tabular}
\vspace{2cm}\\
Die Dissertation wurde am 15.05.2007 bei der Technischen
Universit\"at eingereicht und durch die Fakult\"at f\"ur Mathematik am
17.10.2007 angenommen.
\end{titlepage}

\clearpage{\thispagestyle{empty}}
\begin{verbatim}

\end{verbatim}
\newpage
\tableofcontents

\newpage
\pagenumbering{arabic}
\addsec{Abstract}
The topic of this thesis is the depth of invariant rings of infinite
groups. The main result is a lower bound for the Cohen-Macaulay defect \[\cmdef
K[V]^{G}:=\dim  K[V]^{G}-\depth K[V]^{G}\]  of the
invariant ring $K[V]^{G}$ of a rational representation $V$ (also called a
$G$-module) of a reductive
group $G$. With $K$ we denote an algebraically closed field.
So far, such bounds have only been known for finite groups. In
particular, we will show the following result:\\ 

For each reductive group $G$
that is not linearly reductive, there exists a faithful rational representation
$V$ of $G$  such that
\[
\cmdef K\left[\bigoplus_{k=1}^{n}V\right]^{G}\ge n-2 \quad \textrm{for all } n \in {\mathbb N},
\]
in particular we have
\[
\lim_{n\to\infty}\cmdef K\left[\bigoplus_{k=1}^{n}V\right]^{G}=\infty.
\]\\

In the proof, such a $G$-module $V$ will be given explicitely. We show a
similar result for additive and unipotent groups of positive characteristic
(these groups are not reductive).

In another
chapter, we concentrate on the group $\SL_{2}$, and refining our methods that
lead to the above result, we receive the following:\\

Let  $\langle X,Y \rangle$
 denote the natural representation of the group $\SL_{2}$, $\chr K=p>0$. 
Let further denote $\langle X^{p},Y^{p}\rangle$  the submodule of the
 $p$th symmetric power $S^{p}(\langle X,Y \rangle)$ spanned by the $p$th power
 of the variables. Then
 we have
\[
\cmdef K\left[\langle X^{p},Y^{p}\rangle \oplus\bigoplus_{i=1}^{k}\langle X,Y
  \rangle\right]^{\SL_{2}}\ge k-3.\]
It is worth noting that the underlying module is self-dual and completely
  reducible (as a direct sum of self-dual, irreducible modules).

In the final chapter, we develop an algorithm to compute these special
invariant rings. With an implementation in {\sc Magma}, this allowed us to compute
the Cohen-Macaulay defect explicitely for the cases $(p,k)=(2,4),(2,5),(3,4)$
to $1,2,1$, so the lower bound is sharp in these cases.\\

This thesis also contains a chapter on the cohomolgy of groups, chapter
\ref{Kohomologievongruppen}, which is almost
independent of the others. In this chapter we investigate annihilators of
$n$-cocycles of infinite groups. In particular, we show the following: Let $G$
be any (possibly infinite) group, $d_{0}: KG\rightarrow K$ the augmentation map (which sends each
$\sigma\in G$ to $1\in K$). Then for each $KG$-module $V$ (not necessary
finite-dimensional) and each $g\in H^{n}(G,V)$ we have $d_{0}\cdot g=0\in H^{n}(G,\Hom_{K}(KG,V))$.

\newpage
\addsec{Einleitung}
Die Tiefe von Invarianten \emph{endlicher} Gruppen ist ein in den letzten 25
Jahren viel beachtetes und gut untersuchtes Gebiet. Als eine (sicher unvollst"andige)
Liste von Arbeiten zu diesem Thema seien hier
\cite{Bourguiba,CampEtAl,Ellingsrud,FleischmannCohomology,FleischmannCohomConnectivity,FleischmannDepthTrace,Gordeev,Henn,KemperFiniteGroups,KemperOnCM,KemperDepthCohom,KemperLoci,LandweberStongLCM,LorenzPathak,NeuselDepth,ShankWehlau,SmithKoszul}
genannt. Als bahnbrechend muss hier sicher die Arbeit von Ellingsrud und Skjelbred
\cite{Ellingsrud} hervorgehoben werden.

Dagegen ist die vorliegende Dissertation die erste Arbeit, die sich mit der Tiefe von Invariantenringen
\emph{unendlicher} (genauer: zusammenh"angender) Gruppen befasst. Genauer betrachten
wir folgende Situation: $K$ ist ein algebraisch abgeschlossener
K"orper und $G$ eine "uber $K$ definierte lineare algebraische Gruppe, die
linear und rational auf einem endlich-dimensionalen Vektorraum $V$ operiert. Dadurch wird auch eine
Operation auf dem Ring der Polynomfunktionen $K[V]$ auf $V$ induziert, und die
Menge aller invarianten Polynomfunktionen bezeichnen wir mit $K[V]^{G}$. Wir
interessieren uns dabei nur f"ur den Fall, dass
$K[V]^{G}$ endlich erzeugt ist (z.B. wenn $G$ eine reduktive Gruppe ist), es
sich also um eine graduierte affine Algebra handelt. Wichtige Kenngr"o"sen
einer solchen Algebra $R$ sind die \emph{Tiefe} $\depth
R$, die \emph{Dimension} $\dim R$ und der Cohen-Macaulay-Defekt
$\cmdef R:=\dim R-\depth R$, welcher stets gr"o"ser gleich $0$ ist. Je
gr"o"ser der Cohen-Macaulay-Defekt, desto "`komplizierter"' ist die Struktur
von $R$. Ringe mit $\cmdef R=0$ (die "`sch"onsten"' also) hei"sen
\emph{Cohen-Macaulay}. Der Satz von Hochster und Roberts besagt, dass f"ur
\emph{linear reduktive} Gruppen $G$ und jeden $G$-Modul $V$ der
Invariantenring $K[V]^{G}$ stets Cohen-Macaulay ist, also stets $\cmdef
K[V]^{G}=0$ gilt. Ziel dieser Arbeit ist es, \emph{untere} Schranken f"ur $\cmdef
K[V]^{G}$ zu finden (da $\dim K[V]^{G}$ meist bekannt, ist dies "aquivalent
zum finden oberer Schranken f"ur die Tiefe $\depth K[V]^{G}$), wenn $G$ nicht
linear reduktiv ist. Das Ergebnis
ist unser (technischer) Hauptsatz \ref{BigMainTheorem}. Ein Gro"steil der
restlichen Arbeit besch"aftigt sich dann mit der expliziten Konstruktion von
$G$-Moduln $V$, auf die der Hauptsatz anwendbar ist. Als "`griffiges"'
Ergebnis erhalten wir (siehe Korollar \ref{depthInfVects}, Bemerkung
\ref{BemerkungNachDepthInvVects} und Abschnitt \ref{ExpliziterKonstruktion}):\\

F"ur jede reduktive, nicht linear reduktive Gruppe $G$ existiert ein
  treuer $G$-Modul $V$ (den wir \emph{explizit} angeben k"onnen) mit
\[
\cmdef K\left[\bigoplus_{k=1}^{n}V\right]^{G}\ge n-2 \quad \textrm{f"ur alle } n \in {\mathbb N},
\]
insbesondere also
\[
\lim_{n\to\infty}\cmdef K\left[\bigoplus_{k=1}^{n}V\right]^{G}=\infty.
\]\\

Dies ist eine Versch"arfung des Resultats von Kemper \cite{KemperLinRed},
welches besagt, dass f"ur jede reduktive, nicht linear reduktive Gruppe $G$
ein $G$-Modul $V$ mit $\cmdef K[V]^{G}>0$ existiert.\\ 

Die nach diesem allgemeinen Resultat konstruierten $G$-Moduln $V$ haben eine
relativ gro"se Dimension und komplizierte Struktur. In Satz \ref{ZweitesHauptResultat}, dem zweiten
Hauptresultat dieser Arbeit, wird eine sehr einfach gebaute Folge von
$\SL_{2}$- bzw. $\Ga$-Moduln angegeben, deren Cohen-Macaulay-Defekt der
zugeh"origen Invariantenringe gegen unendlich geht. Unter anderem zeigen wir dort (Satz
\ref{ZweitesHauptResultat}):\\

 Ist $\langle X,Y \rangle$
die nat"urliche Darstellung der $\SL_{2}$, $\chr K=p>0$, so gilt
\[
\cmdef K\left[\langle X^{p},Y^{p}\rangle \oplus\bigoplus_{i=1}^{k}\langle X,Y
  \rangle\right]^{\SL_{2}}\ge k-3.\]\\
Alle hier dargestellten Ergebnisse spielen sich in positiver Charakteristik
ab, denn in Charakteristik $0$ ist jede reduktive Gruppe auch linear reduktiv,
und damit Invariantenringe reduktiver Gruppen stets Cohen-Macaulay. Bevor ich mich dem Thema
der vorliegenden Arbeit gewidmet habe, habe ich ein halbes Jahr lang ohne
  nennenswertes Ergebnis versucht, einen
nicht Cohen-Macaulay Invariantenring in Charakteristik $0$ zu
konstruieren. Dabei war mein Ansatz haupts"achlich experimenteller Natur, d.h. ich habe
f"ur verschiedene Darstellungen $V$ nicht reduktiver Gruppen $G$ den Invariantenring
heuristisch auf die Cohen-Macaulay-Eigenschaft getestet (mit Hilfe von {\sc
  Magma}). Zwar kann man auch f"ur nicht reduktive Gruppen mit Lemma
\ref{Hauplemma} leicht drei Elemente $a_{1},a_{2},a_{3}\in K[V]^{G}$
konstruieren, die als Annullatoren eines Kozyklus $g\in H^{1}(G,K[V])$ keine regul"are Sequenz in $K[V]^{G}$ sind und ein phsop in
$K[V]$ bilden. (Dies ging sogar leichter als f"ur reduktive Gruppen, da die
Kozyklen einfacher waren). Doch leider ist f"ur reduktive Gruppen Lemma
\ref{redphsop} falsch (dies ist das Haupthindernis), also hat man kein a
priori Kriterium f"ur das Vorliegen eines phsops im Invariantenring. Um zu
testen, ob nun doch ein phsop vorliegt, bleibt also nichts "ubrig, als
(heuristisch) Generatoren f"ur den Invariantenring $K[V]^{G}$ und deren Relationenideal zu
berechnen, und so den Test im Restklassenring eines Polynomrings mit
\Magma durchzuf"uhren.

Da kein Algorithmus
bekannt ist, um f"ur nicht reduktive Gruppen den Invariantenring zu berechnen,
habe ich einfach jeweils gen"ugend viele Invarianten berechnet, so dass ich von der
davon erzeugten Unteralgebra $A$ vermutete, dass sie gleich $K[V]^{G}$
ist. Dabei ist im Allgemeinen nicht einmal garantiert, dass $K[V]^{G}$
"uberhaupt endlich erzeugt ist. In
den F"allen, in denen ich das Relationenideal von $A$ berechnen konnte, hat
sich $A$ immer als Cohen-Macaulay herausgestellt. Insbesondere haben die
drei Annullatoren eines Kozyklus kein phsop in $A$ gebildet.

F"ur einige Darstellungen,
die mir als "`heisse"' Kandidaten vorkamen, lies sich das Relationenideal von
$A$ nicht in akzeptabler Zeit berechnen, so dass ich hier keine Aussage machen
konnte. 

\subsection*{Aufbau der Arbeit}
{\bf Abschnitt \ref{Grundlagen}} gibt die ben"otigten Grundlagen wieder. Obwohl sich
hier keine wesentlich neuen Resultate finden, ist der Abschnitt sehr
ausf"uhrlich gehalten. Zum einen sollte die Darstellung nat"urlich so
selbsttragend wie m"oglich sein, zum anderen werden etwa die Resultate, die wir
f"ur noethersche \emph{graduierte} Ringe ben"otigen, in der Literatur (etwa Eisenbud
\cite{Eisenbud} oder Bruns und Herzog \cite{BrunsHerzog}) fast ausschliesslich
f"ur noethersche \emph{lokale} Ringe formuliert. Die dortigen Beweise gehen
zwar meist mit kleinen Modifikationen auch f"ur den graduierten Fall durch,
aber um nicht durch "`sinngem"a"ses Zitieren"' eine L"ucke entstehen zu lassen
habe ich wo n"otig einen vollst"andigen Beweis f"ur den graduierten Fall
gegeben. 

Hierzu noch eine Bemerkung, die auch f"ur die anderen Abschnitte gilt:
Nat"urlich kann auch diese Arbeit nicht bei Null
beginnen. Vorausgesetzt werden die Grundlagen der kommutativen Algebra und der
algebraischen Geometrie, wie sie sich etwa in (der ersten H"alfte von) Standardlehrb"uchern wie
Eisenbud \cite{Eisenbud} oder Kunz \cite{Kunz} finden. Der Leser sollte also
mit Begriffen wie noethersche Ringe und Moduln, Ringerweiterungen und
Ganzheit, Lying-over, Going-up und Going-down, Dimension und H"ohe, 
assoziierte Primideale, Gr"obner Basen, affine
Variet"aten und Morphismen etwas anfangen
k"onnen. Die meisten der hier erw"ahnten Begriffe werden wir aber zumindest
erkl"aren und Referenzen geben. Wenn m"oglich verweisen wir f"ur die
Grundlagen der Theorie der linearen algebraischen Gruppen immer auf die erste
Auflage des 
vergleichsweise elementar gehaltenen Lehrbuchs Springer
\cite{SpringerLin}. Teilweise m"ussen wir auf dessen tiefer gehende zweite
Auflage \cite{SpringerNew} bzw. auf Humphreys \cite{Humphreys} verweisen. Die
Grundlagen der "`modernen"' (im Gegensatz zur "`klassischen"')
Invariantentheorie finden sich etwa in Derksen und Kemper \cite{DerksenKemper}.  
Unter diesen Voraussetzungen sind die Beweise dieser Arbeit vollst"andig. 
Ein Leser ohne die genannten Kenntnisse wird zwar nicht jedes Detail auf
Korrektheit verifizieren, aber wie ich hoffe doch mit Hilfe der
angegebenen Referenzen nachvollziehen k"onnen.

Im Unterabschnitt "uber lineare algebraische Gruppen geben wir ein Kriterium f"ur die lineare Reduktivit"at von
Gruppen, Satz \ref{NagataLinRedCrit}, welches sich innerhalb eines Beweises von Nagata \cite{NagataComplete}
findet und wenig bekannt zu sein scheint. Wir werden von diesem Kriterium
intensiv Gebrauch machen, und geben daher einen (bei den genannten Vorkenntnissen) vollst"andigen, gegl"atteten Beweis.

 Es folgt eine Einf"uhrung in die erste Kohomologie
algebraischer Gruppen, welche unser wichtigstes technisches Hilfsmittel
darstellt. Wir schliessen mit der Bereitstellung der von uns ben"otigten Resultate der
Invariantentheorie.\\

{\bf Abschnitt \ref{Kohomologievongruppen}} hat als wesentliches Ziel, 
ein Resultat von Kemper "uber die Annullation von $1$-Kozyklen auf
$n$-Kozyklen zu erweitern, siehe Satz \ref{AnnulatorDurchSeq}. Insbesondere erhalten wir hieraus, dass die
Augmentationsabbildung $KG\rightarrow K$ ($G$ eine beliebige, nicht notwendig
endliche Gruppe) jeden $n$-Kozyklus annulliert, Korollar
\ref{augmentationannulliert}, was das entsprechende Resultat f"ur endliche
Gruppen verallgemeinert.

Dieses Ergebnis steht in dieser Arbeit f"ur sich
allein und wird hier nicht mehr verwendet. Daher ist es m"oglich, dieses
  Kapitel beim ersten Lesen zu "uberspringen und nur bei Interesse an dem
  dargestellten Resultat darauf zur"uckzukommen. Das Verst"andnis der
restlichen Kapitel wird dabei nicht gest"ort.\\ 

{\bf Abschnitt \ref{Hauptabschnitt}} enth"alt den technischen Hauptsatz, die
bereits angesprochene untere Schranke f"ur den Cohen-Macaulay Defekt. Wir vergleichen
unser Resultat mit den entsprechenden Ergebnissen f"ur endliche Gruppen.\\

Im ankn"upfenden {\bf Abschnitt \ref{Anwendungen}} flie"sen die Ergebnisse der
vorherigen Abschnitte zusammen, und wir erhalten so explizit Darstellungen,
deren Cohen-Macaulay-Defekt der Vektorinvarianten gegen unendlich geht. Wir
untersuchen die explizite Konstruktion f"ur einige F"alle konkret und kommen
so f"ur diese F"alle noch zu Vereinfachungen. Der Abschnitt schlie"st mit
der bereits erw"ahnten einfach gebauten Serie von $\SL_{2}$-Invariantenringen f"ur beliebige
Charakteristik $p>0$, deren Cohen-Macaulay-Defekt gegen unendlich geht.\\ 

Im letzten {\bf Abschnitt \ref{Algorithmik}} untersuchen wir die Ringe
dieser Serie mit algorithmischen Methoden (mit Hilfe von {\sc Magma})
genauer. Insbesondere entwickeln wir einen Algorithmus, der in kurzer Rechenzeit
und mit wenig Speicherbedarf (hieran scheitert der Standardalgorithmus)
Generatoren dieser speziellen Invariantenringe liefert. F"ur drei F"alle
k"onnen wir so den Cohen-Macaulay-Defekt exakt bestimmen (wir erhalten als
Defekte zweimal $1$ und einmal $2$).

\subsection*{Dank}
An dieser Stelle m"ochte ich mich ganz herzlich bei meinem Betreuer, Prof. Dr. Gregor
Kemper bedanken. Zum einen weil er es immer geschafft hat, mich wenn n"otig
f"ur eine bestimmte Richtung zu
motivieren, zum anderen aber auch, weil er mir immer gen"ugend Freiraum
gelassen hat, um eigenverantwortlich zu forschen.

Bei Dr. Frank Himstedt bedanke ich mich f"ur die immer kompetente Hilfe bei meinen Fragen zur
Darstellungstheorie von Gruppen.

Weiter bedanke ich mich beim Graduiertenkolleg "`Angewandte Algorithmische
Mathematik"', von dem ich im ersten Jahr dieser Arbeit finanziell unterst"utzt wurde. 

Zuguterletzt gilt mein Dank den restlichen Mitgliedern der Lehr- und Forschungseinheit
M11 - f"ur die vielen Kaffee-Pausen.

\newpage

\section{Grundlagen der kommutativen Algebra und der Invariantentheorie}\label{Grundlagen}
Wir geben zun"achst die ben"otigten Grundlagen und Begriffe der
kommutativen Algebra und Invariantentheorie
wieder. Dabei soll auch jeweils der algorithmische Aspekt ber"ucksichtigt
werden.\\

\noindent {\bf Standardvoraussetzung.} In dieser Arbeit bezeichnen wir mit $K$
stets einen algebraisch
abgeschlossenen K"orper der Charakteristik $p\ge 0$. 
\subsection{Graduierte affine Algebren} \label{gradAffineAlgs}
Eine \emph{$K$-Algebra}\index{Algebra} ist (in dieser Arbeit, es gibt inkompatible Definitionen)
ein $K$-Vektorraum $R$, der zus"atzlich ein kommutativer Ring mit $1\ne 0$ ist, so dass
stets $\lambda (ab)=(\lambda a)b=a(\lambda b)$ f"ur $\lambda \in K, a,b \in R$
gelten (insbesondere ist $K\cong K\cdot 1 \subseteq R$). $R$ hei"st
\emph{graduiert},\index{Ring!graduierter} falls es eine direkte Zerlegung in
Vektorr"aume $R=\bigoplus_{i=0}^{\infty}R_{i}$ gibt, wobei zus"atzlich
$R_{i}R_{j}\subseteq R_{i+j}$ f"ur alle $i,j \in \mathbb{N}_{0}$ gelten soll.
Dabei soll in dieser Arbeit dann auch immer $R_{0}=K=K\cdot 1$ gelten, d.h. $R$ soll
dann zus"atzlich \emph{zusammenh"angend}\index{Ring!zusammenh\"angend} sein. Dann ist das Ideal 
$R_{+}:=\bigoplus_{i=1}^{\infty}R_{i}$ wegen $R/R_{+}\cong R_{0}=K$ maximal, es
hei"st das \emph{maximale homogene Ideal von R}\index{Ideal!maximales homogenes}.

 Die Elemente aus $R_{i}
\setminus \left\{0\right\}$ hei"sen \emph{homogen vom Grad
  $i$}.\index{homogenes Element} In der
"ublichen Weise ist dann auch eine Gradfunktion $\deg$ auf $R\setminus\{0\}$
definiert (oft setzen wir auch $\deg 0:=-\infty$). Ein Ideal $I \unlhd R$
hei"st \emph{homogen}, \index{Ideal!homogenes} wenn es von homogenen Elementen erzeugt
wird, was gleichbedeutend damit ist, dass es mit jedem $f\in I$ auch alle
homogenen Komponenten von $f$ enth"alt.
Im wesentlichen interessieren
uns in dieser Arbeit nur \emph{affine} \index{Algebra!affine}(d.h. endlich erzeugte) graduierte
Algebren. Dann gibt es homogene Elemente $a_{1},\ldots,a_{n} \in R$ mit
$R=K[a_{1},\ldots,a_{n}]$. (Eine zusammenh"angende graduierte Algebra ist genau
dann affin, wenn sie noethersch ist. In unserem Standardfall f"ur Algebren
sind also affin und noethersch synonym.)
Man hat dann einen Epimorphismus $\phi$ vom Polynomring
$P=K[X_{1},\ldots,X_{n}]$ auf $R$, gegeben durch $X_{i}\mapsto a_{i}$, welcher
bei Verwendung der Graduierung $\deg X_{i}:=\deg a_{i}$ auf $P$ homogen wird,
d.h. homogene Elemente auf homogene Element gleichen Grades (oder die 0)
abbildet. Insbesondere ist dann $I:=\ker \phi$ ein \emph{homogenes} Ideal in
$P$, das \emph{Relationenideal}, \index{Ideal!Relationen}
und es gilt $R \cong P/I$. Eine nullteilerfreie, endlich erzeugte $K$-Algebra
hei"st ein \emph{$K$-Bereich}. \index{Bereich}\\

{\bf Standardvoraussetzung. }Wenn nicht anders vermerkt, bezeichnen wir mit $R$
von nun an stets einen zusammenh"angenden, graduierten, endlich erzeugten
$K$-Bereich.\\

 Um mit affinen graduierten Algebren rechnen zu
k"onnen, ben"otigt man praktisch immer eine solche Darstellung als
Restklassenring eines Polynomrings modulo des Relationenideals (au"ser falls
zus"atzliche Eigenschaften bekannt sind). Ist $R$
endlich erzeugte Unteralgebra eines Polynomrings, so kann man das Relationenideal $I$ mit Standard-Methoden (Eliminationstheorie) berechnen, siehe etwa Eisenbud
\cite[Proposition 15.30]{Eisenbud}. In \Magma \cite{Magma} ist ein solcher
Algorithmus bereits implementiert und kann mit dem Befehl {\tt RelationIdeal}
aufgerufen werden. Oft scheitert hier bereits das weitere Vorgehen aufgrund der
zu langen Rechenzeit.

\subsubsection{Dimension, H"ohe und Noether-Normalisierung}
Sei $R$ eine affine Algebra, $\wp$ ein Primideal von $R$. Dann ist die
\emph{H"ohe}\index{H\"ohe} von $\wp$, bezeichnet mit $\height(\wp)$, die maximale L"ange $n$
einer echt aufsteigenden Kette von Primidealen $\wp_{0} \subset \wp_{1} \subset
\ldots \subset \wp_{n}$ mit $\wp_{n}=\wp$. Ist $I\ne R$ ein beliebiges Ideal
von $R$, so ist $\height(I)$ das Minimum der H"ohen der $I$ umfassenden
Primideale. (Dabei haben echte Ideale in noetherschen Ringen nach dem
Krullschen Hauptidealsatz stets endliche H"ohe.)
Die \emph{Krulldimension}\index{Krulldimension}\index{Dimension} von $R$ ist das Maximum der H"ohen aller Primideale
von $R$, d.h. die maximale L"ange $n$
einer echt aufsteigenden Kette von Primidealen $\wp_{0} \subset \wp_{1} \subset
\ldots \subset \wp_{n}$ von $R$, und wird mit $\dim R$ bezeichnet. Ist $R$
ein $K$-Bereich, so gilt $\dim R=\trdeg \Quot(R)/K$ (\cite[Theorem 13.A]{Eisenbud}).
Man setzt $\dim I:=\dim R/I$. Zwischen H"ohe und Dimension eines Ideals besteht
folgende fundamentale Relation:

\begin{Prop}\label{dimRdimIheightI}
Ist $R$ eine affine Algebra und $I\ne R$ ein Ideal, so gilt
\[
\dim R \ge \dim I +\height I.
\]
Ist $R$ ein affiner Bereich, so gilt sogar Gleichheit.
\end{Prop}

\Bew Sei $\wp_{0}\supset \wp_{1}\supset\ldots\supset\wp_{m}\supseteq I$ eine Kette
von Primidealen von $R$ mit $m=\dim R/I$. Offenbar gilt dann $\dim R\ge
m+\height \wp_{m}\ge \dim R/I + \height I$. F"ur das Gleichheitszeichen im
Fall eines affinen Bereichs siehe \cite[Corollary 13.4]{Eisenbud}. \qed\\

Ist $R$ ein Polynomring, so l"asst sich
$\dim I$ mit Gr"obner-Basis Methoden berechnen \cite[Algorithm
1.2.4]{DerksenKemper}. Insbesondere l"asst sich die Dimension eines affinen
Rings berechnen, wenn das Relationenideal bekannt ist. Zentral ist der folgende 

\begin{Satz}[Noether-Normalisierung] \label{NoetherNorm}
Sei $R$ eine graduierte affine $K$-Algebra. Dann gibt es homogene,
"uber $K$ algebraisch unabh"angige Elemente
$a_{1},\ldots,a_{n} \in R_{+}$, so dass R ganz "uber $A:=K[a_{1},\ldots,a_{n}]$
ist (oder "aquivalent: R ist endlich erzeugt als Modul "uber A). Dabei ist $n=\dim R$ eindeutig bestimmt.
\end{Satz}

Ein Beweis folgt nach Korollar \ref{ganzPhsop}.

\begin{Def}
Eine Menge $\left\{a_{1},\ldots,a_{n}\right\}$ mit den im Noetherschen Normalisierungssatz genannten
Eigenschaften hei"st ein \emph{homogenes Parametersystem
  (hsop)}.\index{hsop}\index{homogenes Parametersystem}\index{phsop} Eine
Menge hei"st \emph{partielles homogenes Parametersystem (phsop)}, wenn sie
Teilmenge eines hsops ist. 
\end{Def}

F"ur eine alternative Charakterisierung von
Parametersystemen (in affinen \emph{Bereichen}) ben"otigen wir zun"achst ein allgemeines Lemma.

\begin{Lemma} \label{htIS}
Sei $S/R$ eine ganze Erweiterung noetherscher Integrit"atsringe und zus"atzlich $R$
\emph{normal}. Ist $I \ne R$ ein Ideal von $R$, so ist $\height(I)=\height(SI)$.\end{Lemma}

\Bew Unter den gemachten Voraussetzungen gilt "`going-down"' \cite[Theorem 13.9]{Eisenbud}. Insbesondere
gilt f"ur jedes Primideal $\bigP \lhd S$, dass 
\[
\height(\bigP)=\height(\bigP\cap R)
\]
 (siehe Kunz \cite[Korollar VI.2.9]{Kunz}). Sei nun $\wp \lhd R$ ein
$I$ umfassendes Primideal mit $\height(\wp)=\height(I)$. Nach "`lying-over"'
(\cite[Satz VI.2.3]{Kunz}) existiert ein Primideal $\bigP \lhd S$ mit
$\wp=\bigP \cap R$. Da $I\subseteq\wp\subseteq\bigP$, gilt auch $SI \subseteq\bigP$, und
damit
\[
\height(SI) \le \height(\bigP)=\height(\bigP\cap R)=\height(\wp)=\height(I).
\]
Sei umgekehrt $\bigP \lhd S$ ein $SI$ umfassendes Primideal mit
$\height(\bigP)=\height(SI)$. Dann ist $I\subseteq \bigP \cap R$, und damit
\[
\height(I) \le \height(\bigP \cap R)=\height(\bigP)=\height(SI).
\]
Insgesamt gilt also $\height(I)=\height(SI)$. \qed\\

(Der in \cite[Lemma 1.5.e]{KemperOnCM} gegebene Beweis f"ur diesen Satz ohne
die gemachten Voraussetzungen an $R$ und $S$ ($R$ normal und $S$
nullteilerfrei fehlen) enth"alt eine L"ucke.) Wie
"ublich bezeichnen wir f"ur einen $R$-Modul $M$ mit $\Ass_{R}M$ die Menge der
assoziierten Primideale von $M$ in $R$ (siehe Eisenbud \cite[Chapter 3.1]{Eisenbud}).

\begin{Lemma} \label{phsopHeight}\index{phsop}\index{hsop}
Sei $R$ eine graduierte affine $K$-Algebra.  Es
gilt stets $\height R_{+}=\dim R$. Seien weiter $a_{1},\ldots,a_{k}\in R$
($k=0$ m"oglich) homogene
Elemente positiven Grades, sowie $I \ne R$ ein homogenes Ideal von $R$.
\begin{enumerate}
\renewcommand{\labelenumi}{(\alph{enumi})}

\item Gilt $\height(a_{1},\ldots,a_{k})=k$, so ist $a_{1},\ldots,a_{k}$ ein
phsop von $R$. Ist hierbei $k=\dim R$, so handelt es
sich sogar um ein hsop. Ist $R$ nullteilerfrei und $a_{1},\ldots,a_{k}$ ein
phsop von $R$, so gilt umgekehrt $\height(a_{1},\ldots,a_{k})=k$.

\item Ist $a_{1},\ldots,a_{k}\in I$ mit $\height (a_{1},\ldots,a_{k})=k$ und
  $r=\height I$, so
  gibt es homogene Elemente positiven Grades $a_{k+1},\ldots,a_{r}\in I$ mit
  $\height (a_{1},\ldots,a_{r})=r$. Ist $R$ nullteilerfrei, so kann man also
  ein in $I$ liegendes phsop  zu einem in $I$ liegenden phsop mit $\height(I)$ Elementen erg"anzen.
\end{enumerate}
F"ur einen affinen Bereich ist also $a_{1},\ldots,a_{k}$ genau dann ein phsop,
wenn $\height (a_{1},\ldots,a_{k})=k$ gilt.
\end{Lemma}
 
\Bew (b) Sei $\height(a_{1},\ldots,a_{k})=k$. Ist $k=\height(I)$, so sind
wir fertig, sei also $k < \height(I)$. Es gen"ugt offenbar f"ur den Beweis,
ein weiteres homogenes Element $a_{k+1}\in I$ zu finden, so
dass $\height (a_{1},\ldots,a_{k+1})=k+1$ ist (da $I\ne R$ hat $a_{k+1}$ dann
automatisch positiven Grad).
Sei $\{\wp_{1},\ldots,\wp_{s}\}$ die Menge der minimalen,
$(a_{1},\ldots,a_{k})$ umfassenden Primideale. Diese Menge ist endlich und
alle Primideale sind homogen, da sie Teilmenge der
endlichen Menge homogener Primideale $\Ass_{R}(R/(a_{1},\ldots,a_{k}))$ ist (\cite[Theorem
3.1.a]{Eisenbud}, \cite[Satz C.23]{Kunz}), und es gilt $\height(\wp_{i}) \le k
\, \myforall i$ nach dem
Krullschen Hauptidealsatz (\cite[Theorem 10.2]{Eisenbud}). Es gilt sogar
Gleichheit, da $\height(\wp_{i}) \ge \height(a_{1},\ldots,a_{k})=k$.
 W"are $I\subseteq \bigcup_{i=1}^{s}\wp_{i}$, so g"abe es
ein $i$ mit $I\subseteq \wp_{i}$ (Lemma "uber das Vermeiden von Primidealen, \cite[Lemma 3.3]{Eisenbud} oder \cite[Lemma III.3.6]{Kunz}). Dann w"are
$\height(I)\le \height(\wp_{i})=k$, im Widerspruch zur Voraussetzung. Also
gibt es ein nach dem zitierten Lemma sogar homogenes $a_{k+1}\in
I\setminus\bigcup_{i=1}^{s}\wp_{i}$. Ein minimales Primideal $\wp$ mit
$(a_{1},\ldots,a_{k+1})\subseteq\wp$ umfasst dann eins der minimalen
Primideale $\wp_{i}$ (\cite[Satz III.1.10.b]{Kunz}), und da $a_{k+1}\notin \wp_{i}$,
ist die Inklusion echt. Es folgt $\height(\wp) \ge \height(\wp_{i})+1=k+1$. Nach
Krulls Hauptidealsatz gilt auch $\height(\wp) \le k+1$, also Gleichheit. Da
dies f"ur jedes solche minimale Primideal $\wp$ gilt, folgt also
$\height(a_{1},\ldots,a_{k+1})=k+1$, was zu zeigen war. 

(Ist $R$ nullteilerfrei, so bilden  $a_{1},\ldots,a_{k}\in I$ nach (a) genau dann
ein phsop, wenn $\height (a_{1},\ldots,a_{k})=k$. Daher kann man ein in $I$
liegendes phsop zu einem in $I$ liegenden phsop mit $\height(I)$ Elementen
erg"anzen. Dies beweist dann den Zusatz in (b), den wir nat"urlich nicht f"ur
den folgenden Beweis von (a) brauchen.) 

(a) (i) Sei $\height(a_{1},\ldots,a_{k})=k$. Nach  (b) mit $I=R_{+}$ und
$n:=\height R_{+}$ lassen sich
die $k$ Elemente aus $R_{+}$ zu einer in $R_{+}$ liegenden Menge homogener Elemente
$\{a_{1},\ldots,a_{n}\}$ mit
$\height(a_{1},\ldots,a_{n})=n=\height R_{+}$ erg"anzen. Wir zeigen,
dass eine solche
Menge ein hsop ist, und das $\height R_{+}=\dim R$ gilt, was die ersten beiden
Teile der Behauptung (a) sowie die Behauptung im Vorspann zeigt.
 Ein "uber $(a_{1},\ldots,a_{n})$ liegendes minimales Primideal $\wp$ ist
 homogen (s.o.), liegt also im maximalen homogenen Ideal $R_{+}$, und aus $n=\height(a_{1},\ldots,a_{n})\le \height(\wp) \le \height(R_{+})=n$ folgt
$\wp=R_{+}$. Also ist $\wp$ maximal, und damit ist $\dim R/(a_{1},\ldots,a_{n})=0$. Ist $A:=K[a_{1},\ldots,a_{n}]$, so bedeutet dies $\dim R/A_{+}R=0$, d.h.
$\dim_{K} R/A_{+}R < \infty$. Nach dem graduierten Nakayama Lemma \cite[Lemma
3.5.1]{DerksenKemper} ist dann 
$R$ endlich erzeugt als $A$-Modul, also ist $R$ ganz "uber $A$. Damit ist
$\dim A=\dim R$.
Ist $J\lhd K[X_{1},\ldots,X_{n}]=:K[X]$ das
Relationenideal von $a_{1},\ldots,a_{n}$, so ist $A\cong K[X]/J$. W"aren
$a_{1},\ldots,a_{n}$ nicht algebraisch unabh"angig, also $J\ne (0)$, so w"are
$\dim R=\dim A=\dim K[X_{1},\ldots,X_{n}]/J < n=\height (R_{+})$. Da aber
offenbar $\height R_{+}\le \dim R$ gilt, ist dies ein Widerspruch und $J=(0)$. Also ist
$a_{1},\ldots,a_{n}$ ein hsop, und wir haben $\dim R=\dim A=n=\height R_{+}$ gezeigt.

(ii) Sei nun $R$ zus"atzlich nullteilerfrei. Man erg"anze das phsop $a_{1},\ldots,a_{k}$ zu einem hsop
$a_{1},\ldots,a_{n}$ von $R$. Dann ist $R$ ganz "uber dem (normalen, da
faktoriellen) Polynomring
$A=K[a_{1},\ldots,a_{n}]$. Mit Lemma \ref{htIS}, das wir wegen der
Nullteilerfreiheit von $R$ und der Normalit"at von $A$ auf die Erweiterung $R/A$
 anwenden k"onnen, folgt daher sofort
$\height(a_{1},\ldots,a_{k})_{R}=\height(a_{1},\ldots,a_{k})_{A}=k$. 

Der Zusatz ist lediglich ein Spezialfall von (a). 
\qed\\

\begin{BemRoman}\label{phsopZweiGradus}
Da die H"ohe eines Ideals nicht von der Graduierung abh"angt, folgt im
nullteilerfreien Fall aus Teil
$(a)$, dass Elemente die bez"uglich zweier verschiedener Graduierungen von $R$
homogen und positiven Grades sind und ein phsop bez"uglich einer dieser Graduierungen bilden, auch
bez"uglich der anderen ein phsop bilden.
\end{BemRoman}

\begin{Korollar} \label{ganzPhsop}
Sei $S/R$ eine ganze Erweiterung noetherscher graduierter affiner Bereiche und zus"atzlich $R$
\emph{normal}. Homogene Elemente $a_{1},\ldots,a_{k}\in R_{+}$ bilden genau
dann ein phsop in $R$, wenn sie eines in $S$ bilden.
\end{Korollar}

\Bew Dies folgt sofort mit Lemma \ref{phsopHeight} (a) und
\[
\height(a_{1},\ldots,a_{k})_{R}=\height(a_{1},\ldots,a_{k})_{S}
\]
nach Lemma \ref{htIS}. \qed\\

\noindent \textit{Beweis des Noetherschen Normalisierungssatzes
  \ref{NoetherNorm}.} Man erg"anze die leere Menge gem"a"s Lemma
  \ref{phsopHeight} (b) mit $I=R_{+}$ zu einer Menge homogener Elemente
  positiven Grades $a_{1},\ldots,a_{n}$ mit
  $\height(a_{1},\ldots,a_{n})=n=\height R_{+}=\dim R$. Nach Teil (a) des
  Lemmas ist diese Menge ein hsop von $R$. Ist umgekehrt $R$ ganz "uber einem
  Polynomring $A$, so ist $n=\dim A=\dim R$ eindeutig bestimmt.\qed\\

\begin{BemRoman}
Ist $a_{1},\ldots,a_{n}$ ein hsop von $R$, (also $n=\dim R$), so gilt umgekehrt stets
$\height(a_{1},\ldots,a_{n})=n$ (auch wenn es Nullteiler gibt). Sei n"amlich $A:=K[a_{1},\ldots,a_{n}]$ und
$\wp\supseteq (a_{1},\ldots,a_{n})$ ein minimales Primideal. Dann ist $\wp$
homogen, also $\wp\subseteq R_{+}$. Damit gilt $A_{+}\subseteq \wp \cap
A\subseteq R_{+}\cap A=A_{+}$ (Definition der Graduierung von $A$), also
Gleichheit. "Uber dem Primideal $A_{+}$ von $A$ (wg. $A/A_{+}=K$) liegen also
die Primideale $\wp\subseteq R_{+}$ von $R$, und da $R/A$ ganz, darf es keine
echte Inklusion geben \cite[Satz VI.2.3]{Kunz}, d.h. $\wp=R_{+}$. Also ist $\height
(a_{1},\ldots,a_{n})=\height \wp=\height R_{+}=\dim R=n$.
\end{BemRoman}

Damit wir die dargestellten Resultate stets voll verwenden k"onnen, haben wir
"`$R$ nullteilerfrei"' in unsere Standardvoraussetzung mit aufgenommen.\\

Aus der Definition der H"ohe folgt sofort, dass ein Ideal $I$ gleiche H"ohe wie
sein Radikalideal $\sqrt{I}$ hat, also
\begin{equation} \label{heightOfRadikal}
\height (I)=\height(\sqrt{I}).
\end{equation}
Dies ergibt

\begin{Lemma} \label{phsopPower}
Seien $f_{1},\ldots,f_{k}\in R_{+}$ homogen und $i_{1},\ldots,i_{k}\ge
1$. Dann ist $f_{1},\ldots,f_{k}$ genau dann ein phsop in $R$, wenn
$f_{1}^{i_{1}},\ldots,f_{k}^{i_{k}}$ eines ist.
\end{Lemma}\index{phsop}

\Bew Dies folgt sofort aus Lemma \ref{phsopHeight} (a) und
\[
\height (f_{1},\ldots,f_{k})\stackrel{~\eqref{heightOfRadikal}}{=}\height
\sqrt{(f_{1},\ldots,f_{k})}=\height
\sqrt{(f_{1}^{i_{1}},\ldots,f_{k}^{i_{k}})}\stackrel{~\eqref{heightOfRadikal}}{=}\height
(f_{1}^{i_{1}},\ldots,f_{k}^{i_{k}}).
\]
\qed
\begin{BemRoman}
Wir haben den Begriff "`phsop"' nur f"ur graduierte affine
  Algebren definiert. Ist $R$ allgemeiner ein Noetherscher Ring, so hei"st
  etwa in Kemper \cite[vor Lemma 1.5]{KemperOnCM}
  $a_{1},\ldots,a_{k}\in R$ ein \emph{partielles Parametersystem (psop)}, wenn $(a_{1},\ldots,a_{k})\ne R$ und
  $\height(a_{1},\ldots,a_{i})=i$ f"ur \emph{alle} $i=1,\ldots,k$ gilt. F"ur graduierte affine Bereiche stimmt diese
  Definition im Fall homogener $a_{i}$ mit unserer "uberein, denn wenn die letzte Gleichung f"ur $i=k$
  gilt, also nach Lemma \ref{phsopHeight} ein phsop vorliegt, so ist aufgrund
  unserer Definition erst Recht $a_{1},\ldots,a_{i}$ ein phsop f"ur jedes
  $i\le k$, und
  wieder nach Lemma \ref{phsopHeight} gilt dann auch
  $\height(a_{1},\ldots,a_{i})=i$.
\end{BemRoman}
Ist $R=P/I$ als Faktorring eines Polynomrings gegeben und nullteilerfrei, so
kann man mit Lemma \ref{phsopHeight} leicht
entscheiden, ob gegebene homogene Elemente $(a_{1},\ldots,a_{k})$ ein
phsop bilden. Dann ist n"amlich $\height(a_{1},\ldots,a_{k})_{R}=\dim R - \dim
R/(a_{1},\ldots,a_{k})_{R}=\dim P/I - \dim P/(I+(a_{1},\ldots,a_{k})_{P})$,
und in Polynomringen l"asst sich die Dimension von Idealen wie bereits bemerkt
berechnen.\\

Um eine Algebra $R$ genauer zu untersuchen, ist die
Konstruktion eines (p)hsops oft unerl"asslich. Der Algorithmus von Kemper
\cite{KemperOptHSOP} 
berechnet ein "`optimales"' hsop, wobei (u.a.) das obige Kriterium verwendet
wird. Hier muss also  das Relationenideal bekannt sein. 

\subsection*{Exkurs: Berechnung eine hsops}
{\it In diesem Exkurs verlassen wir kurz den systematischen Aufbau der
  Grundlagen. Die Resultate dieses Abschnitts werden nicht weiter verwendet, so
  dass dieser Abschnitt bedenkenlos "ubersprungen werden kann.}\\

Die meisten Beweise des
Noetherschen Normalisierungssatzes sind in irgendeiner Form
ebenfalls konstruktiv, siehe z.B. Eisenbud \cite[Theorem 13.3]{Eisenbud} oder Fogarty
\cite[Theorem 5.44]{Fogarty}. Wir wollen hier kurz ein Verfahren zur Bestimmung eines hsops angeben, das im Prinzip im
explizit machen der Beweise aus Eisenbud/Fogarty besteht, und ohne die
Berechnung des vollst"andigen Relationenideals auskommt. Leider l"asst sich
der angegebene Algorithmus zwar theoretisch leicht
umsetzen, scheitert in der Praxis aber daran, dass die in jedem
Schritt auftretenden Relationenhauptideale von Polynomen von immer gr"o"ser
werdendem Grad erzeugt wurden, bis diese nicht mehr (in akzeptabler Zeit)
gefunden werden konnten. Hier also zun"achst die Rohfassung des Verfahrens:
\begin{Algorithmus}\label{makephsop}
Bestimmung eines hsops einer Unteralgebra $R$
eines Polynomrings~$K[X]$\\
{\bf Eingabe:} Erzeuger $f_{1},\ldots,f_{n} \in K[X]$ von $R$, also
$R=K[f_{1},\ldots,f_{n}] \subseteq K[X]$.\\
{\bf Ausgabe:} Ein hsop f"ur $R$.\\
\begin{samepage}
{\bf BEGIN}
\begin{enumerate}
\item $hsop_{0}:=\emptyset$
\item For $i:=1$  to $n$ do
\begin{enumerate}
\item Falls $hsop_{i-1}\cup \left\{f_{i}\right\}$ algebraisch unabh"angig ist, so
  setze $hsop_{i}:=hsop_{i-1}\cup\left\{f_{i}\right\}$.
\item Sonst bestimme ein $hsop$ f"ur $K[hsop_{i-1}\cup\{f_{i}\}]$ und setze $hsop_{i}:=hsop$.
\end{enumerate}
\item Return($hsop_{n}$)
\end{enumerate}
{\bf END}
\end{samepage}
\end{Algorithmus}

Bevor wir erkl"aren, wie die Schritte (a) und (b) umzusetzen sind, machen wir
zun"achst folgende {\it Beobachtung}: Es ist stets $hsop_{i}$ ein hsop von
$K[f_{1},\ldots,f_{i}]$ (insbesondere also $hsop_{n}$ ein hsop von $R$).\\
\Bew Offenbar ist in jedem Schritt $hsop_{i}$ algebraisch unabh"angig. 
Wir machen Induktion nach $i$, um die Ganzheitsrelation zu zeigen.\\
 $i=0:$ Die leere Menge $hsop_{0}$ ist ein hsop f"ur
$K$.\\
$(i-1)\rightarrow i:$ Nach Voraussetzung sind $f_{1},\ldots,f_{i-1}$ ganz
"uber $K[hsop_{i-1}]$. Da diese Elemente dann auch ganz "uber
$K[hsop_{i-1},f_{i}]$ sind, und $f_{i}$ dies nat"urlich ebenfalls ist, ist
also $K[f_{1},\ldots,f_{i}]$ ganz "uber $K[hsop_{i-1},f_{i}]$. Damit  folgt
die Aussage im Fall, dass $hsop_{i}$ nach (a) berechnet wird. Im Fall (b) ist
nun $K[hsop_{i-1},f_{i}]$ nach Konstruktion ganz "uber $K[hsop_{i}]$, und
aufgrund der Transitivit"at der Relation "`ganz"' ist also
$K[f_{1},\ldots,f_{i}]$ ganz "uber $K[hsop_{i}]$. \qed\\

\noindent Umsetzung der Schritte (a) und (b):

{\it Zu (a).} Es sei zun"achst bemerkt, dass der Algorithmus nicht falsch wird,
wenn man in jedem Schritt nach (b) vorgeht, z.B. weil aufgrund des folgenden
Korollars keine Aussage gemacht werden kann.

\begin{Satz}[Jacobi-Kriterium]\index{Jacobi Kriterium}
Seien $g_{1},\ldots,g_{n}\in K[X_{1},\ldots,X_{n}]$ Elemente eines
Polynomrings, und $J:=\left(\frac{\partial g_{i}}{\partial
    X_{j}}\right)_{i,j=1,\ldots,n}$ die zugeh"orige \emph{Jacobi-Matrix}. Dann
gilt: 

$K(X_{1},\ldots,X_{n})/K(g_{1},\ldots,g_{n})$ ist genau dann eine
endliche und separable K"orpererweiterung, falls $\det J \ne 0$ ist.
\end{Satz}

\Bew Siehe \cite{Benson}, Proposition 5.4.2. \qed

\begin{Korollar}
Seien $g_{1},\ldots,g_{k}\in K[X_{1},\ldots,X_{n}]$ Elemente eines
Polynomrings, und $J:=\left(\frac{\partial g_{i}}{\partial
    X_{j}}\right)_{i=1,\ldots,k,j=1,\ldots,n}$ die zugeh"orige \emph{Jacobi-Matrix}. Dann
gilt: Gibt es eine
$k \times k$ Teilmatrix von $J$ mit Determinante ungleich $0$, so sind
$g_{1},\ldots,g_{k}$ algebraisch unabh"angig. In Charakteristik $0$ gilt sogar
die Umkehrung dieser Aussage.
\end{Korollar}

\Bew Falls es eine solche Teilmatrix mit Determinante ungleich $0$ gibt und
$j_{1},\ldots,j_{k}$ die zugeh"origen Spalten sind, so hat die zu
$\{g_{1},\ldots,g_{k}\} \cup \left\{X_{i}: i \in \{1,\ldots,n\}\setminus
  \{j_{1},\ldots,j_{k}\}\right\}$ geh"orige Jacobi-Determinante Wert ungleich
$0$. Nach dem Jacobi-Kriterium ist die betrachtete Menge also eine
Transzendenzbasis, und die Teilmenge $\{g_{1},\ldots,g_{k}\}$ damit algebraisch unabh"angig.

Sind umgekehrt $\{g_{1},\ldots,g_{k}\}$ algebraisch unabh"angig, so kann man
diese Menge mit $n-k$ Elementen aus der Transzendenzbasis
$\{X_{1},\ldots,X_{n}\}$ zu einer Transzendenzbasis erg"anzen (Austauschsatz
f"ur Transzendenzbasen). Sei also
$\{g_{1},\ldots,g_{k},X_{j_{1}},\ldots,X_{j_{n-k}}\}$ eine Transzendenzbasis,
mit zugeh"origer Jacobi-Determinante $J$. Dann ist die K"orpererweiterung
$K(X_{1},\ldots,X_{n})/K(g_{1},\ldots,g_{k},X_{j_{1}},\ldots,X_{j_{n-k}})$
  endlich und in Charakteristik $0$ automatisch separabel, also $\det J
  \ne 0$. Dies bedeutet, dass die zu den Zeilen/Spalten $\{1,\ldots,n\}\setminus
\{j_{1},\ldots,j_{n-k}\}$ geh"orige Teilmatrix von $J$ Determinante ungleich
$0$ hat. \qed\\

Der Beweis zeigt auch, wie man ggf. die gegebenen Polynome mit den Variablen
zu einer Transzendenzbasis erg"anzen kann. Im Falle positiver Charakteristik $p$ wird die
Umkehrung des Kriteriums falsch; etwa ist $X^{p}$ eine Transzendenzbasis von $K(X)$, aber
$J=(pX^{p-1})=(0)$ - die Erweiterung $K(X)/K(X^{p})$ ist zwar endlich, aber
nicht separabel.

In positiver Charakteristik wird man bei der Durchf"uhrung des Algorithmus
also so vorgehen: Man testet im Schritt (a) zun"achst auf die Existenz einer
nicht-singul"aren $i\times i$ Teilmatrix. Hat man eine solche, so ist die betrachtete
Menge algebraisch unabh"angig und man kann wie beschrieben vorgehen. Ansonsten
kann man keine Aussage machen, und man muss Schritt (b) durchf"uhren.\\

Die Voraussetzung des Korollars bedeutet nat"urlich nichts anderes, als das
$\operatorname{rang} J = k$ "uber dem rationalen Funktionenk"orper
$K(X_{1},\ldots,X_{n})$, und man kann dies mit dem Gauss-Algorithmus
entscheiden, und erh"alt so ggf. auch die Spalten einer nicht-singul"aren
Teilmatrix. Die folgende Variante um zu testen, ob $J$ eine solche Teilmatrix
besitzt, lies sich mit den in \Magma vorhandenen Funktionen jedoch etwas
schneller implementieren, und lieferte in den betrachteten F"allen das
Ergebnis praktisch in Nullzeit.

{\it Test ob eine nichtsingul"are $k \times k$ Teilmatrix von $J$ existiert,
wobei gleich die zugeh"origen Spalten $j_{1},\ldots,j_{n}$ dieser Matrix,
falls sie existiert, mitbestimmt werden:} W"ahle zun"achst $j_{1}$ als den
kleinsten Index $j$ mit $\frac{\partial g_{1}}{\partial X_{j}} \ne 0$. Falls es
keinen solchen Eintrag gibt, so gibt es auch keine nichtsingul"are $k \times
k$ Teilmatrix. Seien nun $j_{1},\ldots,j_{i}$ bereits konstruiert, so dass die
zu diesen Spaltenindizes und zu den ersten $i$ Zeilen geh"orige Teilmatrix
Determinante ungleich $0$ hat. F"ur $j_{i+1}$ w"ahle man nun den kleinsten
Index $1,\ldots,n$, so dass die zugeh"orige $(i+1) \times (i+1)$ Determinante
einen Wert ungleich $0$ hat. Gibt es keinen solchen Index, so gibt es auch
keine nichtsingul"are $k \times k$ Teilmatrix.\\

{\it Zu (b) in Algorithmus \ref{makephsop}, Schritt 2.} Als erstes m"ussen wir das Relationenideal von
$hsop_{i-1}\cup\left\{f_{i}\right\}$ bestimmen. Wir gehen hier davon aus, dass man im Schritt (b) ankommt, weil bekannt
ist, dass die in
(a) betrachtete Menge $hsop_{i-1}\cup\left\{f_{i}\right\}$ algebraisch abh"angig ist. Ist dies nicht der Fall, etwa
in positiver Charakteristik, so muss man eines der
Standard-(Gr"obner-Basis)-Verfahren verwenden, um das Relationenideal von
$hsop_{i-1}\cup\left\{f_{i}\right\}$ zu bestimmen. Ist es das Nullideal, so zeigt
das nachtr"aglich die algebraische Unabh"angigkeit dieser Menge. 
Das folgende Lemma zeigt, dass das Relationenideal jedenfalls stets ein
Hauptideal ist (nach Konstruktion ist $hsop_{i-1}$ stets algebraisch
unabh"angig, und daher ist das folgende Lemma anwendbar).

\begin{Lemma}
Sei $\{f_{1},\ldots,f_{k+1}\}\subseteq K[Y]$ eine algebraisch abh"angige Menge von
Polynomen, so dass $\{f_{1},\ldots,f_{k}\}$ algebraisch unabh"angig sind. 
Sei $\tilde{p}(X_{k+1})\in K[f_{1},\ldots,f_{k}][X_{k+1}]$ dasjenige Polynom, das aus dem
Minimalpolynom $p$ von $f_{k+1}$ "uber dem Quotientenk"orper
$K(f_{1},\ldots,f_{k})$ durch Multiplikation mit dem Hauptnenner der
Koeffizienten hervorgeht.

Dann
ist das Relationenideal von $\{f_{1},\ldots,f_{k+1}\}$ in
$K[X_{1},\ldots,X_{k+1}]$ ein Hauptideal, erzeugt von dem Polynom, das aus $\tilde{p}$
durch Ersetzen von (den algebraisch unabh"angigen Elementen)
$f_{1},\ldots,f_{k}$ durch $X_{1},\ldots,X_{k}$ hervorgeht. 
\end{Lemma} 
\Bew Nach Voraussetzung ist $R:=K[f_{1},\ldots,f_{k}]\cong
K[X_{1},\ldots,X_{k}]$, insbesondere ist $R$ faktoriell. Sei
$F:=\Quot(R)=K(f_{1},\ldots,f_{k})$ und $p(X_{k+1})$ das Minimalpolynom von
$f_{k+1}$ "uber $F$. Mit $c$ bezeichnen wir den \emph{Inhalt} eines Polynoms
aus $F[X_{k+1}]$, der mit Hilfe eines fest gew"ahlten Repr"asentantensystems
der Primelemente von $R$ berechnet wird. (Erinnerung: Der Inhalt $c(f)$ eines
Polynoms $f \in R[X]$ ist der $ggT$ seiner Koeffizienten, berechnet mit dem
Repr"asentantensystem ($c(0):=0$). Ist $f \in \Quot(R)[X]$ und $a\in R\setminus\{0\}$ mit $af \in R[X]$,
so ist $c(f):=c(af)/c(a)$. F"ur $f,g \in \Quot(R)[X]$ gilt $c(fg)=c(f)c(g)$,
und aus $c(f)=1$ folgt $f \in R[X]$. Insbesondere ist f"ur $f\ne 0$ stets
$f/c(f)\in R[X]$.) Dann ist $\tilde{p}=p/c(p)$.

Sei nun $0\ne g\in K[X_{1},\ldots,X_{k+1}]\cong R[X_{k+1}]$ eine Relation,
d.h. $g(f_{1},\ldots,f_{k+1})=0$. Fassen wir $g$ als Element von $R[X_{k+1}]$
auf, so bedeutet dies $g(f_{k+1})=0$. Es ist zu zeigen, dass $g$ polynomiales
Vielfaches von $p/c(p) \in R[X_{k+1}]$ ist. Da $p$ das Minimalpolynom von
$f_{k+1}$ "uber $F$ ist, gibt es ein $h\in F[X_{k+1}]$ mit $ph=g$, also
$c(p)c(h)=c(g)$. Damit ist $g=\frac{p}{c(p)}\cdot\frac{h}{c(h)}\cdot c(g)$, und
da jeder der drei Faktoren in $R[X_{k+1}]$ liegt, folgt die Behauptung.~\qed\\

Da in unserem Fall alle Polynome homogen sind, ist die das Relationenideal
erzeugende Relation ebenfalls homogen bei Gewichtung $\deg X_{i}:=\deg
f_{i}$. Falls man weis, dass die in (a) betrachtete Menge algebraisch
abh"angig ist, es also Relationen positiven Grades gibt, kann man nacheinander
alle Grade $d=1,2,3,\ldots$ durchlaufen, in jedem Grad $d$ die Menge der
Monome vom Grad $d$ in $hsop_{i-1}\cup\left\{f_{i}\right\}$ bilden und hiervon
dann eine nichttriviale Linearkombination der $0$ suchen - in jedem Schritt
wird also ein homogenes lineares Gleichungssystem gel"ost. Sobald man eine
L"osung hat, ist diese dann somit die minimale, erzeugende Relation, und das
Verfahren terminiert. Wie bereits erw"ahnt wurden die minimalen Grade $d$ in
der Praxis meist bald so gro"s, dass die erzeugende Relation nicht mehr
gefunden werden konnte.\\

Hat man die Relation erstmal gefunden, ist die Berechnung des hsops
vergleichsweise leicht. Sind n"amlich $f=f_{1},f_{2},\ldots,f_{n} \in K[X_{1},\ldots,X_{n}]=:K[X]$ und
interpretiert man $f$ als Relation, so sind $f_{2},\ldots,f_{n}$ genau dann
ein hsop f"ur $K[X]/(f)$, falls $f,f_{2},\ldots,f_{n}$ ein hsop f"ur $K[X]$
ist. Nach Eisenbud \cite[Lemma 13.2.c]{Eisenbud} kann man dazu f"ur $i\ge 2$ z.B.  
$f_{i}=X_{i}-a_{i}X_{1}$ (ggf. vorher homogenisiert) w"ahlen mit $a_{i} \in K$ geeignet. In der
Praxis reicht es meist, zwei der $a_{i}$ ungleich $0$ zu w"ahlen. H"aufig ist
z.B. (nach umnumerieren) $f\in X_{1}^{i}X_{2}^{j}+(X_{3},\ldots,X_{n})$, und
f"ur $f_{2}$ bietet sich dann eine Homogenisierung von $X_{1}+X_{2}$ an, und
 f"ur die restlichen Elemente $f_{i}=X_{i}, \, i>2$. Dies ist dann
 tats"achlich ein hsop, da die gemeinsame Nullstellenmenge der $n$ Polynome
 offenbar nur aus $0 \in K^{n}$ besteht (aus $X_{n}=\ldots=X_{3}=0$ folgt mit
 $f=0$ dann $X_{1}X_{2}=0$; zusammen mit
 $X_{1}+X_{2}=0$ dann auch $X_{1}=X_{2}=0$). Hilberts Nullstellensatz liefert
 die Ganzheit von $K[X]/K[f_{1},\ldots,f_{n}]$.
 Man sieht bereits, dass
das Relationenideal des n"achsten hsops entsprechend komplizierter wird.

\subsubsection{Tiefe und regul"are Sequenzen}
\begin{Def}\label{DefOfDepth}\index{regul\"are
    Sequenz}\index{Tiefe}\index{Tiefe!messende Sequenz}
Sei $R$ kommutativer Ring mit $1$, und $M$ ein $R$-Modul.
Eine Folge von Elementen $a_{1},\ldots,a_{n}\in R$ hei"st eine \emph{$M$-regul"are
  Sequenz (der L"ange $n$)}, falls gelten:

(a) $(a_{1},\ldots,a_{n})M \ne M$.

(b) $a_{i}$ ist kein Nullteiler von $M/(a_{1},\ldots,a_{i-1})M$  $\myforall
i=1,\ldots,n$. (D.h.: Ist $a_{1}m_{1}+\ldots+a_{i}m_{i}=0$ mit
$m_{1},\ldots,m_{i}\in M$, so ist $m_{i}\in(a_{1},\ldots,a_{i-1})M$.)

Ist $I$ ein Ideal von $R$ mit $IM \ne M$, so wird mit $\depth
(I,M)$ die \emph{Tiefe} von $I$ auf $M$ bezeichnet, die gr"osste auftretende L"ange einer
in $I$ liegenden, $M$-regul"aren Sequenz (diese Zahl ist endlich f"ur $R$
noethersch und $M$ endlich erzeugter $R$-Modul (\cite[Proposition 18.2]{Eisenbud})). Man schreibt auch $\depth(I):=\depth(I,R)$. Eine solche regul"are Sequenz
gr"o"stm"oglicher L"ange hei"st dann eine \emph{die ($I$-)Tiefe messende}
regul"are Sequenz.

Eine in $I$ liegende regul"are Sequenz $(a_{1},\ldots,a_{n})$ hei"st
\emph{maximal},\index{regul\"are Sequenz!maximale} wenn sie nicht zu einer l"angeren regul"aren Sequenz in $I$
erg"anzt werden kann, d.h. $I$ besteht nur aus Nullteilern von
$M/(a_{1},\ldots,a_{n})M$.

Ist $R$ sogar eine eine graduierte Algebra, und $M$ ein graduierter $R$-Modul
(d.h. $M=\bigoplus_{i=b}^{\infty}M_{i}$ als $K$-Vektorr"aume, $b\in {\mathbb Z}$,
und $R_{i}M_{j}\subseteq M_{i+j}\,\myforall i,j$), und sind alle Elemente einer regul"aren
Sequenz zus"atzlich homogen, so spricht man von einer \emph{homogenen
  regul"aren Sequenz}.\index{regul\"are Sequenz!homogene} Ist $R_{+}$ das maximale homogene Ideal von $R$, so hei"st
$\depth(M):=\depth(R_{+},M)$ die \emph{Tiefe von $M(\ne 0)$}, insbesondere ist
$\depth(R)=\depth(R_{+},R)$ die \emph{Tiefe von $R$}. \index{Tiefe}
\end{Def}

Ein Ideal $I$ von $R$ ist auch ein $R$-Modul. Mit $\depth(I)$ meint man dann
aber immer $\depth(I,R)$ und nie $\depth(R_{+},I)$.\\

\label{MagmaTestAufRegSeq}
Ist $R=P/I$ als Quotient eines Polynomrings $P$ nach einem Ideal $I$ gegeben,
so l"asst sich algorithmisch entscheiden, ob eine Folge $a_{1},\ldots,a_{n}\in
R$
$R$-regul"ar ist. Offenbar ist n"amlich $a_{1}$ genau dann kein Nullteiler von
$R$, wenn f"ur das Quotientenideal $I:(a_{1})=I$ gilt. Weiter ist genau dann $a_{2}$
kein Nullteiler auf $R/(a_{1}) \cong  P/(I+(a_{1}))$, wenn
$(I+(a_{1})):(a_{2})=I+(a_{1})$ gilt u.s.w. Quotientenideale lassen sich aber
mit Gr"obner-Basen berechnen, siehe z.B. \cite[section
1.2.4]{DerksenKemper}. Ein entsprechender Algorithmus ist in \Magma
mit dem Befehl {\tt IsZeroDivisor} implementiert.\\

\begin{Satz}[Rees]\label{Rees}
Sei $R$ ein noetherscher Ring, $M$ ein endlich erzeugter $R$-Modul und $I$ ein
Ideal von $R$ mit $IM\ne M$. Dann haben je zwei maximale $M$-regul"are Sequenzen
in $I$ die gleiche L"ange, n"amlich $\depth (I,M)$. Eine $M$-regul"are Sequenz
in $I$ ist also genau dann maximal, wenn sie die Tiefe misst. Insbesondere
kann jede in $I$ liegende $M$-regul"are Sequenz zu einer die Tiefe von $I$ messenden $M$-regul"aren Sequenz
erg"anzt werden. Ist daher $(a_{1},\ldots,a_{r})$ eine $M$-regul"are Sequenz in $I$, so ist
\begin{equation} \label{DepthErg}
\depth (I,M/(a_{1},\ldots,a_{r})M)=\depth (I,M) - r.
\end{equation}
\end{Satz}

\Bew Siehe Bruns und Herzog \cite[Theorem 1.2.5]{BrunsHerzog} oder \cite[Satz
E.7]{Kunz}. \qed\\

Der folgende Satz besagt, dass wir uns zur Bestimmung der Tiefe in den uns
interessierenden F"allen auf homogene
regul"are Sequenzen beschr"anken k"onnen.

\begin{Satz} \label{homDepth}
Ist $R$ eine graduierte \emph{affine} Algebra, $M$ ein \emph{endlich
  erzeugter}, graduierter  $R$-Modul, und $I$ ein \emph{homogenes} Ideal von
  $R$ mit $IM\ne M$, so existiert eine \emph{homogene} $M$-regul"are Sequenz der L"ange
  $\depth (I,M)$ in $I$. Genauer l"asst sich jede in $I$ liegende
  homogene $M$-regul"are Sequenz zu einer die Tiefe $\depth(I,M)$ messenden, in
  $I$ liegenden homogenen $M$-regul"aren Sequenz erg"anzen.
\end{Satz}

\Bew (Vgl. Bruns und
 Herzog \cite[Propositon 1.5.11]{BrunsHerzog}.) Sei $a_{1},\ldots,a_{k}\in I$
 eine homogene $M$-regul"are Sequenz ($k=0$ erlaubt). Ist $k=\depth(I,M)$, so sind wir
 fertig. Sei also $k<\depth(I,M)$. Nach Satz \ref{Rees} ist die $M$-Sequenz
 $a_{1},\ldots,a_{k}$ dann nicht maximal, d.h. $I$ besteht nicht nur aus
 Nullteilern von $N:=M/(a_{1},\ldots,a_{k})M\ne 0$. Da die Menge der Nullteiler von
 $N$ in $R$ (und der $0$) durch $\bigcup \Ass_{R}N$ gegeben ist (Eisenbud \cite[Theorem
 3.1.b]{Eisenbud}), gilt also $I\not\subseteq \bigcup \Ass_{R}N$. Da die
 $a_{i}$ homogen sind, ist $N$ ein graduierter $R$-Modul, und damit sind die
 (endlich vielen) Elemente von $\Ass_{R}N$ homogene Primideale (Eisenbud \cite[Theorem
3.1.a, 3.12]{Eisenbud}). Nach dem (graduierten) Lemma "uber das Vermeiden von
 Primidealen (\cite[Lemma III.3.6]{Kunz} oder \cite[Lemma
 1.5.10]{BrunsHerzog}) existiert dann ein homogenes $a_{k+1}\in I$ mit
 $a_{k+1}\not\in\bigcup \Ass_{R}N$. Dann ist $a_{k+1}$ kein Nullteiler von
 $N=M/(a_{1},\ldots,a_{k})$ und damit $a_{1},\ldots,a_{k+1}$ eine homogene
 $M$-regul"are Sequenz in $I$. Durch Induktion folgt die Behauptung. \qed\\

Tr"agt $R$ also zwei verschiedene Graduierungen, die beide das gleiche
maximale homogene Ideal $R_{+}$ haben, so gibt es f"ur jede der Graduierungen
jeweils eine homogene regul"are Sequenz der gleichen L"ange $\depth R$.



 Die
G"ultigkeit dieser S"atze zu gew"ahrleisten ist einer der Gr"unde f"ur die
Forderung $IM \ne M$ in der Definition der Tiefe eines Ideals. In der graduierten Situation bedeutet die Standardbedingung $IM\ne M$ "ubrigens
lediglich $M\ne 0$ und $I\subseteq R_{+}$.

Homogene
regul"are Sequenzen haben - im Gegensatz zu den nicht homogenen - sehr
angenehme Eigenschaften. Zum Beispiel gilt

\begin{Satz} \label{PermReg}\index{regul\"are Sequenz!Permutationen}
Sei $R$ eine graduierte Algebra, $M\ne 0$ ein graduierter $R$-Modul. Ist dann
$(a_{1},\ldots,a_{n})$ eine \emph{homogene} $M$-regul"are Sequenz, so ist f"ur jede
Permutation $\pi \in S_{n}$ auch $(a_{\pi(1)},\ldots,a_{\pi(n)})$ eine
homogene $M$-regul"are Sequenz.
\end{Satz}

\Bew Da $S_{n}$ von den Transpositionen $(k,k+1)$ mit $1\le k\le n-1$ ezeugt
wird, gen"ugt es den Satz f"ur eine solche Permutation zu beweisen. Hierf"ur
ist lediglich noch zu zeigen, dass $a_{k+1}$ kein Nullteiler von
$N:=M/(a_{1},\ldots,a_{k-1})$ und $a_{k}$ kein Nullteiler von $N/(a_{k+1})N$
ist, wobei $a_{k},a_{k+1}$ eine $N$-regul"are Sequenz ist. Es gen"ugt also den
Satz f"ur $n=2$ und $\pi=(1,2)$ zu beweisen. 

Wir zeigen zuerst, dass $a_{2}$
kein Nullteiler auf $M$ ist. Sei n"amlich $m\in M$ mit $a_{2}m=0$. Da dies
dann auch f"ur jede homogene Komponente von $m$ gilt (weil $a_{2}$ homogen
ist), k"onnen wir $m$ als homogen annehmen.
Es gilt dann
erst recht $a_{2}m=0\in M/(a_{1})M$, und wegen der $M$-Regularit"at von
$a_{1},a_{2}$ folgt $m\in (a_{1})M$. Es gibt also $m'\in M$ mit $m=a_{1}m'$. Da
$m$ und $a_{1}$ homogen sind, k"onnen wir dann auch $m'$ homogen w"ahlen, und
es ist wegen $\deg a_{1}>0$ dann $m=0$ (dann sind wir fertig) oder $\deg m'<\deg m$.
Aus $0=a_{2}m=a_{1}a_{2}m'$ und weil $a_{1}$ kein Nullteiler ist, folgt
$a_{2}m'=0$. Da die Graduierung von $M$ nach Definition nach unten beschr"ankt
ist, folgt durch Induktion nach $\deg m$ (genauer: "`Jagd nach dem kleinsten Verbrecher"')
 dann $m'=0$ und damit dann doch $m=a_{1}m'=0$.

Sei nun $m_{1}\in M$ mit $a_{1}m_{1}=0\in M/(a_{2})M$, d.h. es gibt $m_{2}\in
M$ mit $a_{1}m_{1}+a_{2}m_{2}=0$. Wir m"ussen $m_{1}\in (a_{2})M$ zeigen. Da $a_{1},a_{2}$ $M$-regul"ar ist, gilt 
$m_{2}\in (a_{1})M$, d.h. es gilt $m_{2}=a_{1}m'$ mit $m'\in M$. Also ist
$0=a_{1}m_{1}+a_{2}m_{2}=a_{1}m_{1}+a_{1}a_{2}m'=a_{1}(m_{1}+a_{2}m')$. Da
$a_{1}$ kein Nullteiler auf $M$ ist, folgt also $m_{1}=-a_{2}m'\in (a_{2})M$,
was zu zeigen war. (Vgl. auch \cite[Korollar E.16 und Satz E.17]{Kunz} f"ur
eine Verallgemeinerung (mit anderem Beweis).) \qed\\

\begin{Satz}
Sei $R$ eine graduierte affine $K$-Algebra, $M$ ein endlich erzeugter
graduierter $R$-Modul, $I$ ein homogenes Ideal von $R$ und $a$ ein homogenes
Element von $R$. Ist dann $(I+(a))M \ne M$, so gilt
\[
\depth(I+(a),M) \le \depth (I,M) + 1.
\]
\end{Satz}

Wir geben zun"achst einen allgemeinen Beweis, und dann einen elementaren
Beweis unter einer Zusatzannahme, die in dieser Arbeit stets eintreffen
wird, da bei uns immer  $M=S$ eine  nullteilerfreie Oberalgebra von $R$ sein wird. Der 1. Beweis kann daher ggf. "ubersprungen werden.\\

\noindent {\it 1. Beweis (mit Homologie).} Wir "ubersetzen den in Eisenbud
\cite[Lemma 18.3]{Eisenbud} gegebenen Beweis (f"ur einen lokalen Ring $R$) in
unsere graduierte Situation. Sei $x_{1},\ldots,x_{n}$ ein homogenes Erzeugendensystem
von $I$, und sei 
\[
k+1:=\depth(I+(a),M).
\]
 Sei dann $K(x_{1},\ldots,x_{n},a)$
der zu diesem Erzeugendensystem von $I+(a)$ geh"orende \emph{Koszul-Komplex},\index{Koszul-Komplex} siehe
\cite[Section 17.2]{Eisenbud}.
Nach der Charakterisierung der Tiefe durch das Verschwinden der Homologie des
Koszul-Komplexes \cite[Theorem 17.4]{Eisenbud} ist dann $H^{i}(M\otimes
K(x_{1},\ldots,x_{n},a))=0$ f"ur $i \le k$. Nun hat man nach \cite[Corollary
17.11]{Eisenbud} eine exakte Sequenz
\[
H^{i}(M\otimes K(x_{1},\ldots,x_{n})) \stackrel{a}{\rightarrow} H^{i}(M\otimes
K(x_{1},\ldots,x_{n})) \rightarrow H^{i+1}(M\otimes
K(x_{1},\ldots,x_{n},a)),
\]
wobei die erste Abbildung durch Multiplikation mit $a$ gegeben ist. F"ur $i \le
k-1$ sind die rechten Homologiemoduln gleich $0$, d.h. die Linksmultiplikation
mit $a$ ist surjektiv. Also gilt
\[
H^{i}(M\otimes K(x_{1},\ldots,x_{n}))=(a)_{R}\cdot H^{i}(M\otimes K(x_{1},\ldots,x_{n})).
\]
Nun sind die Homologien des Koszul-Komplexes graduierte
$R$-Moduln, und $(a)_{R}$ ist ein homogenes Ideal in $R$. Aufgrund obiger
Gleichung folgt dann aus dem graduierten Nakayama-Lemma \cite[Exercise
  4.6]{Eisenbud} oder \cite[Lemma A.9]{Kunz}, dass $H^{i}(M\otimes
K(x_{1},\ldots,x_{n}))=0$ f"ur $i \le k-1$. Wieder aufgrund der
Charakterisierung der Tiefe mittels Homologie folgt $\depth(I,M) \ge k$, und
mit der Definition von $k$ die Behauptung. \qed\\ 

Die Idee f"ur den folgenden Beweis hat
mir Gregor Kemper beim Kaffee-trinken mitgeteilt.\\

\noindent {\it 2. Beweis (elementar)} unter der Zusatzvoraussetzung, dass $a$ kein
Nullteiler auf $M$ ist, $(a)$ also selber eine
regul"are Sequenz der L"ange $1$ in $I+(a)$ ist. Nach Satz \ref{homDepth} kann
man $a$ zu einer maximalen homogenen regul"aren Sequenz
$a,b_{1}+r_{1}a,b_{2}+r_{2}a,\ldots,b_{k}+r_{k}a$, mit $b_{i} \in I, r_{i}\in
R$ homogen, in $I+(a)$ erg"anzen. Insbesondere gilt dann 
\begin{equation} \label{depthk}
\depth (I+(a),M)=k+1.
\end{equation}
Nun gilt aber f"ur alle $i=1,\ldots,k$, dass
\[
(a,b_{1}+r_{1}a,b_{2}+r_{2}a,\ldots,b_{i-1}+r_{i-1}a)M=(a,b_{1},b_{2},\ldots,b_{i-1})M,
\]
und f"ur $m\in M$ ist dann 
\begin{equation}\label{regSeqBew1}
(b_{i}+r_{i}a)m \in
(a,b_{1}+r_{1}a,b_{2}+r_{2}a,\ldots,b_{i-1}+r_{i-1}a)M
\Leftrightarrow b_{i}m \in (a,b_{1},b_{2},\ldots,b_{i-1})M,
\end{equation}
und ebenso
\begin{equation}\label{regSeqBew2}
m \in
(a,b_{1}+r_{1}a,b_{2}+r_{2}a,\ldots,b_{i-1}+r_{i-1}a)M
\Leftrightarrow  m \in (a,b_{1},b_{2},\ldots,b_{i-1})M.
\end{equation}
Da $a,b_{1}+r_{1}a,b_{2}+r_{2}a,\ldots,b_{k}+r_{k}a$ regul"ar ist, ist dann
also auch $(a,b_{1},b_{2},\ldots,b_{k})$ $M$-regul"ar - man lese dazu
\eqref{regSeqBew1} von rechts nach links, verwende die Regularit"at und lese
dann \eqref{regSeqBew2} von links nach rechts. Zus"atzlich ist die
Sequenz homogen, und damit ist nach Satz \ref{PermReg} auch die Permutation
$(b_{1},b_{2},\ldots,b_{k},a)$ $M$-regul"ar. Insbesondere ist dann
$(b_{1},b_{2},\ldots,b_{k})$ eine regul"are, in $I$ liegende Sequenz. Es folgt 
$k \le \depth(I)$, mit \eqref{depthk} also $\depth(I+(a),M) \le
\depth(I,M)+1$. \qed\\

\begin{Korollar} \label{depthkKorollar}
Sind in obiger Situation $a_{1},\ldots,a_{n}\in R$ homogen mit
$(I+(a_{1},\ldots,a_{n}))M \ne M$, so gilt
\[
\depth(I+(a_{1},\ldots,a_{n}),M) \le \depth(I,M)+n
\]
\end{Korollar}

\Bew Induktiv folgt $\depth(I+(a_{1},\ldots,a_{n}),M) \le
\depth(I+(a_{1},\ldots,a_{n-1}),M)+1 \le \depth(I+(a_{1},\ldots,a_{n-2}),M)+2
\le \ldots \le \depth(I,M) + n$ \qed\\

Mit Hilfe des folgenden Satzes l"asst sich die Tiefe eines Ideals oft
bestimmen. Er ist inspiriert von Shank und Wehlau \cite[Theorem
2.1]{ShankWehlau}, welche sich nach eigenen Angaben von Eisenbud
\cite[Corollary 17.12]{Eisenbud} inspirieren lie"sen.

\begin{Satz} \label{SatzShankWehlau}
Sei $R$ eine graduierte affine $K$-Algebra, $M\ne 0$ ein endlich erzeugter
graduierter $R$-Modul. Seien $a_{1},\ldots,a_{n}\in R_{+}$ homogen und
$a_{1},\ldots,a_{k}$ eine $M$-regul"are Sequenz ($k\le n$). Gibt es dann ein
$m\in M$ mit $m\notin (a_{1},\ldots,a_{k})M$, aber $a_{i}m\in
(a_{1},\ldots,a_{k})M \,\, \myforall i=1,\ldots,n$, so gilt
\[
\depth((a_{1},\ldots,a_{n})_{R},M)=k.
\]
\end{Satz}

\noindent{\bf "Aquivalente Formulierung:} Ist $I\ne R$ ein homogenes Ideal,
$a_{1},\ldots,a_{k}\in I$ eine homogene $M$-regul"are Sequenz, und gibt es ein $m\in M$ mit
$m\notin (a_{1},\ldots,a_{k})M$, aber $rm \in
(a_{1},\ldots,a_{k})M\,\, \myforall r\in I$, so gilt
\[
\depth (I,M)=k.
\]

\Bew Nach Satz \ref{Rees} gen"ugt es zu zeigen, dass $a_{1},\ldots,a_{k}$
eine maximale regul"are Sequenz in $I:=(a_{1},\ldots,a_{n})_{R}$ ist, denn
dann misst sie die Tiefe. W"are
dem nicht so, so g"abe es ein $a=r_{1}a_{1}+\ldots+r_{n}a_{n}\in I$ (mit
$r_{i}\in R\,\myforall i$), so dass $a_{1},\ldots,a_{k},a$ ebenfalls $M$-regul"ar
ist. Nach Voraussetzung ist aber $am=r_{1}a_{1}m+\ldots+r_{n}a_{n}m \in
(a_{1},\ldots,a_{k})M$, jedoch $m\notin(a_{1},\ldots,a_{k})M$, was im
Widerspruch zur Regularit"at von $a_{1},\ldots,a_{k},a$ steht. 

F"ur die "aquivalente Formulierung w"ahle man einfach homogene
$a_{k+1},\ldots,a_{n}\in I$ mit $I=(a_{1},\ldots,a_{n})$. Genau dann ist $rm
\in (a_{1},\ldots,a_{k})M$ f"ur alle $r\in I$, wenn $a_{i}m
\in (a_{1},\ldots,a_{k})M$ f"ur $i=1,\ldots,n$. 

Man kann auch so schlie"sen: Nach Voraussetzung besteht $I$ nur aus
Nullteilern von $M/(a_{1},\ldots,a_{k})M$. Daher ist $a_{1},\ldots,a_{k}$ eine
maximale $M$-regul"are Sequenz in $I$ und misst damit die Tiefe.
\qed\\

{\noindent} Es gilt auch die {\bf Umkehrung:} {\it Ist $a_{1},\ldots,a_{k}$
  eine maximale homogene $M$-regul"are Sequenz in $I$, so gibt es ein $m\in M$
  mit $m\not\in(a_{1},\ldots,a_{k})M$ aber $Im\subseteq
  (a_{1},\ldots,a_{k})M$.}\\

\Bew Nach Voraussetzung besteht $I$ nur aus Nullteilern von
$N:=M/(a_{1},\ldots,a_{k})M$. Nach \cite[Theorem 3.1.b]{Eisenbud} ist also
$I\subseteq \bigcup_{\wp \in \Ass_{R}(N)}\wp$, und nach \cite[Lemma
3.3]{Eisenbud} gibt es dann $\wp\in\Ass_{R}N$ mit $I\subseteq \wp$. Zu $\wp$
gibt es dann ein $n\in N\setminus\{0\}$ mit $\wp=\Ann_{R}n$, und f"ur $m\in M$
mit $n=m+(a_{1},\ldots,a_{k})M$ gilt dann $m\not\in(a_{1},\ldots,a_{k})M$ aber $Im\subseteq
  (a_{1},\ldots,a_{k})M$. \qed\\

\subsubsection{Die Cohen-Macaulay Eigenschaft}
In diesem Abschnitt f"uhren wir die Cohen-Macaulay Eigenschaft ein und bringen
den Zusammenhang zwischen phsops und regul"aren Sequenzen.

Ab jetzt sei immer $R$ eine graduierte affine $K$-Algebra und $M$ ein endlich
erzeugter graduierter $R$-Modul.

F"ur jedes Ideal $I\ne R$ gilt
\begin{equation} \label{depthheight}
\depth(I) \le \height(I),
\end{equation}
siehe \cite[Proposition 18.2]{Eisenbud}. Daher ist die im Folgenden definierte Zahl stets gr"o"ser gleich~$0$.

\begin{Def} \label{DefCMDef}\index{Cohen-Macaulay!Defekt}
Der \emph{Cohen-Macaulay-Defekt} eines echten Ideals $I$ von $R$ ist die Differenz
von H"ohe und Tiefe,
\[
\cmdef (I):=\height(I)-\depth(I).
\]
Der \emph{Cohen-Macaulay-Defekt} von $R$ ist die Differenz von Krulldimension
und Tiefe, also
\[
\cmdef (R):=\dim(R)-\depth(R)=\cmdef(R_{+}).
\]
R hei"st \emph{Cohen-Macaulay},\index{Cohen-Macaulay} falls $\cmdef(R) = 0$, also wenn es eine
homogene regul"are Sequenz der L"ange $\dim(R)$ gibt. 
\end{Def}

Beispielsweise sind Polynomringe $R=K[X_{1},\ldots,X_{n}]$ Cohen-Macaulay, da
offenbar die Folge der Variablen $X_{1},\ldots,X_{n}$ eine regul"are Sequenz der L"ange $n=\dim R$ bildet.

\begin{Satz} \label{regPhsop}
(a) Jede homogene $R$-regul"are Sequenz $a_{1},\ldots,a_{k}\in R_{+}$
erzeugt ein Ideal der H"ohe $k$, $\height(a_{1},\ldots,a_{k})=k$, ist also
insbesondere ein phsop. 

(b) Ist $R$ Cohen-Macaulay, so ist
auch umgekehrt jedes phsop (der L"ange $k$) eine $R$-regul"are Sequenz (und erzeugt damit nach
(a) ein Ideal der H"ohe $k$).

(c) $R$ ist genau dann Cohen-Macaulay, wenn f"ur ein hsop $a_{1},\ldots,a_{n}$
von $R$ und
$A:=K[a_{1},\ldots,a_{n}]$ der dann "uber $A$ endlich erzeugte $A$-Modul $R$ sogar
frei ist. Wenn diese Eigenschaft f"ur ein hsop gilt (also $R$ Cohen-Macaulay
ist), dann gilt sie sogar f"ur jedes hsop.
\end{Satz}

\Bew (a) Sei $a_{1},\ldots,a_{k}$ eine homogene regul"are Sequenz. Diese ist dann
nat"urlich maximal regul"ar in dem Ideal $(a_{1},\ldots,a_{k})_{R}$, und daher ist
\[
k=\depth(a_{1},\ldots,a_{k})_{R}\stackrel{~\eqref{depthheight}}{\le}
\height(a_{1},\ldots,a_{k})_{R}\le k,
\]
wobei im letzten Schritt Krulls Hauptidealsatz verwendet wurde. Also gilt
Gleichheit, und mit Lemma \ref{phsopHeight} ist $a_{1},\ldots,a_{k}$ ein
phsop.

(b) und (c). (i) Sei zun"achst $a_{1},\ldots,a_{n}$ ein hsop von $R$, $n=\dim R$, so dass $R$ frei "uber
$A:=K[a_{1},\ldots,a_{n}]$ ist. Dann gibt es $g_{1},\ldots,g_{m}\in R$ mit
\[
R=\bigoplus_{j=1}^{m}Ag_{j}
\quad \textrm{ und }\quad A\rightarrow Ag_{j},\,\,\, a\mapsto ag_{j}\, \textrm { injektiv
  f"ur } j=1,\ldots,m.
\]
(Die zweite Bedingung, die f"ur ein
  nichtnullteilerfreies $R$ n"otig ist um "`$R$ frei "uber $A$"' zu formulieren, wird in der Literatur oft vergessen.)
Wir zeigen, dass dann $a_{1},\ldots,a_{n}$ eine
regul"are Sequenz ist - dann ist insbesondere $R$ Cohen-Macaulay. Sei also
$k\le n$ und $r_{1},\ldots,r_{k}\in R$ mit

\[
r_{1}a_{1}+\ldots+r_{k}a_{k}=0.
\]
Zu $r_{i}\in\bigoplus_{j=1}^{m}Ag_{j}, \, i=1,\ldots,k$ gibt es dann $p_{ij} \in A,\, j=1,\ldots,m$ mit
$r_{i}=\sum_{j=1}^{m}p_{ij}g_{j}$. Es folgt
\[
\sum_{i=1}^{k}\sum_{j=1}^{m}p_{ij}g_{j}a_{i}=0.
\]
 Aufgrund der Direktheit der
Summe $R=\bigoplus_{j=1}^{m}Ag_{j}$ und der Injektivit"at von $A\rightarrow Ag_{j}$ folgt $\sum_{i=1}^{k}p_{ij}a_{i}=0$ f"ur
$j=1,\ldots,m$. Aus 
\[
p_{kj}a_{k}=-\sum_{i=1}^{k-1}p_{ij}a_{i},\quad \quad j=1,\ldots,m,
\]
 und weil $A$ ein
Polynomring ist, folgt, dass die rechte
Seite (nach Zusammenfassen) nur aus Monomen besteht, von denen jedes durch
$a_{k}$ und wenigstens durch ein $a_{i},\, i=1,\ldots,k-1$ teilbar ist. Es folgt
\[
p_{kj}\in Aa_{1}+\ldots+Aa_{k-1} \quad j=1,\ldots,m,
\]
und daher
\[
r_{k}=\sum_{j=1}^{m}p_{kj}g_{j}\in Ra_{1}+\ldots+Ra_{k-1}.
\]
Dies zeigt die Regularit"at der Folge $a_{1},\ldots,a_{n}$, und da $n=\dim R$
ist $R$ Cohen-Macaulay.

(ii) Sei nun umgekehrt $R$ Cohen-Macaulay, also $\depth R=\dim R=:n$, und
$a_{1},\ldots,a_{n}$ ein hsop sowie $A=K[a_{1},\ldots,a_{n}]$. Aufgrund des
nachfolgenden Lemmas \ref{depthARM} ist dann ebenfalls
$\depth(A_{+},R)=n$. Da $R$ endlich erzeugter $A$-Modul ist, besitzt $R$ nach
dem Hilbertschen Syzygien Satz (\cite[Cororllary 19.8]{Eisenbud}) eine
endliche freie Aufl"osung "uber dem Polynomring $A$. Insbesondere ist damit
die \emph{projektive Dimension}\index{Dimension!projektive} von $R$ als $A$-Modul endlich, $\projdim_{A} R < \infty$.
Damit sind die Voraussetzungen zur Anwendung der graduierten Auslander-Buchsbaum Formel (Eisenbud
\cite[Exercise 19.8]{Eisenbud}) erf"ullt, und nach dieser gilt dann
\[
\projdim_{A} R=\depth(A_{+},A)-\depth(A_{+},R)=n-n=0.
\]
Daher gibt es eine projektive Aufl"osung von $R$ der L"ange $0$, d.h. $R$ ist
selbst projektiv als $A$-Modul. Endlich erzeugte projektive graduierte Moduln
"uber graduierten noetherschen
Ringen sind aber frei (\cite[Theorem 19.2]{Eisenbud}). Also ist $R$ frei als
$A$-Modul. 

Zusammenfassend haben wir also gesehen: (i) Wenn es ein hsop gibt, so dass $R$
frei "uber $A=K[hsop]$ ist, so ist dieses hsop eine regul"are Sequenz. Insbesondere
ist $R$ Cohen-Macaulay. (ii) Wenn $R$ Cohen-Macaulay ist, so ist $R$ frei
"uber $A$. Insbesondere ist also nach (i) das zu $A$ geh"orige hsop eine
regul"are Sequenz. Dies zeigt (b) und (c). \qed\\

So wie $\cmdef R=0$ bedeutet, dass $R$ frei "uber $A=K[hsop]$ ist, so bedeutet
$\cmdef R=1$, dass der erste Syzygien-Modul (der Generatoren von $R$ als
$A$-Modul) frei "uber $A$ ist. Allgemein folgt aus der Auslander-Buchsbaum-Formel
\[
\projdim_{A} R=\depth(A_{+},A)-\depth(A_{+},R)=\dim R-\depth R=\cmdef R,
\]
dass $\cmdef R=k$ "aquivalent ist zur Existenz einer minimalen freien
Aufl"osung
\[
0\rightarrow A^{n_{k}}\rightarrow\ldots\rightarrow A^{n_{1}}\rightarrow
A^{n_{0}}\rightarrow R\rightarrow 0.
\]
Also bedeutet $\cmdef R=k$, dass der $k$-te Syzygien-Modul $\im
(A^{n_{k}}\rightarrow A^{n_{k-1}})=\ker (A^{n_{k-1}}\rightarrow A^{n_{k-2}})$
(dabei "`$A^{n_{-1}}:=R, \,\,A^{n_{-2}}:=0$"') frei "uber $A$ ist.\\

Ist $M$ ein $R$-Modul und $A$ ein Unterring von $R$ (mit $1_{R}\in A$), so ist
$M$ auch ein $A$-Modul. Sind $a_{1},\ldots,a_{n}\in A$, so gilt
\[
(a_{1},\ldots,a_{n})M:=(a_{1},\ldots,a_{n})_{R}M=(a_{1},\ldots,a_{n})_{A}M.
\]
\Bew Die Inklusion "`$\supseteq$"' ist wegen $R\supseteq A$ offensichtlich. Ist
$b=\left(\sum_{i=1}^{n}a_{i}r_{i}\right)m\in(a_{1},\ldots,a_{n})_{R}M$ mit
$r_{1},\ldots,r_{n}\in R$, $m\in M$, so ist wegen $r_{i}m\in
M$ auch $b=\sum_{i=1}^{n}a_{i}(r_{i}m)\in (a_{1},\ldots,a_{n})_{A}M$. Da
die betrachteten Elemente $b$ den linken Modul erzeugen,
zeigt dies die Inklusion "`$\subseteq$"'. \qed\\

Aufgrund dieser Eigenschaft kann man also den Index $R$ bzw. $A$ bei der
betrachteten Menge weglassen.

\begin{Lemma} \label{depthARM}
Sei $A$ eine graduierte affine Unteralgebra von $R$, so dass $R$ ganz
"uber $A$ ist (z.B. $A$ die von einem homogenen Parametersystem erzeugte
Unteralgebra), und $M\ne 0$ ein endlich erzeugter graduierter $R$-Modul (z.B. $M=R$). Dann
gilt
\[
\depth(R_{+},M)=\depth(A_{+},M),
\]
d.h. die Tiefe von $M$ als $R$-Modul ist gleich der Tiefe von $M$ als $A$-Modul.
\end{Lemma}

\Bew (vgl. \cite[Lemma 3.7.2]{DerksenKemper}.) Da $R$ endlich erzeugt als
 $A$-Modul, ist auch $M$
 endlich erzeugt als $A$-Modul.
 Sei $a_{1},\ldots,a_{k}\in A_{+}$ eine maximale $M$-regul"are Sequenz in
$A_{+}$. Es gen"ugt zu zeigen, dass diese auch in $R_{+}$ maximal ist (wegen Satz \ref{Rees}), also
dass $R_{+}$ nur aus Nullteilern von $N:=M/(a_{1},\ldots,a_{k})M$
besteht (wir nehmen hier die $0$ auch als Nullteiler). Nach Voraussetzung besteht jedenfalls $A_{+}$ nur aus Nullteilern von
$N$. Dann liegt $A_{+}$ erst recht in der Menge aller Nullteiler von $N$ in $R$,
welche nach \cite[Theorem 3.1]{Eisenbud} gleich der Vereinigung der assoziierten
Primideale von $N$ in $R$, also gleich $\bigcup\Ass_{R}N$ ist. 
Dann gilt also $A_{+}\subseteq\bigcup_{\bigP\in \Ass_{R}N}\bigP\cap A$, und
 rechts steht eine endliche (\cite[Theorem 3.1]{Eisenbud}) Vereinigung von
 Primidealen von $A$. Nach dem Lemma "uber das Vermeiden von Primidealen
 \cite[Lemma 3.3]{Eisenbud} liegt $A_{+}$ also in einem dieser Primideale,
 d.h. es gibt $\bigP\in \Ass_{R}N$ mit
$A_{+}\subseteq A\cap \bigP \subset A$ (da $1\not\in \bigP$). Da $A_{+}$
maximales Ideal ist (wg. $A/A_{+}=K$), folgt also $A_{+}=A\cap \bigP$. Nach
Definition der Graduierung auf $A$ gilt auch $A_{+}=A\cap R_{+}$. Nach \cite[Proposition 3.12]{Eisenbud} ist $\bigP$ ein homogenes Primideal, also
$\bigP\subseteq R_{+}$. Da $R/A$ ganz ist, gibt es zwischen den "uber $A_{+}$
liegendenen Primidealen $\bigP$ und $R_{+}$ aber keine echte Inklusion
\cite[Corollary 4.18]{Eisenbud}.
 Also gilt
$\bigP=R_{+}$, d.h. als assoziiertes Primideal besteht $R_{+}$ nur aus
Nullteilern von $N$. Dies war zu zeigen. \qed\\
\begin{Korollar} \label{KorDepthR}
Ist $a_{1},\ldots,a_{n}$ ein hsop von $R$, so ist
\[
\depth R = \depth (a_{1},\ldots,a_{n})_{R}.
\]
\end{Korollar}

\Bew Nach dem vorigen Lemma gilt mit $A:=K[a_{1},\ldots,a_{n}]$
\begin{eqnarray*}
\depth (R)&=&\depth(R_{+},R)=\depth(A_{+},R)\\&\le& \depth(A_{+}R,R)= \depth
((a_{1},\ldots,a_{n})_{R},R)\\
 &\le& \depth (R_{+},R) = \depth(R).
\end{eqnarray*}
Also gilt "uberall Gleichheit. \qed\\
\begin{Satz} \label{cmdefAbsch}
Sei $R$ eine graduierte affine $K$-Algebra, und $I\ne R$ ein homogenes Ideal. Dann
gilt
\[
\cmdef R \ge \cmdef I =\height I - \depth I.
\] 
Ist insbesondere $R$ Cohen-Macaulay, so gilt f"ur jedes homogene Ideal $I\ne
R$
\[
\depth I=\height I.
\]
\end{Satz}

\Bew Sei $a_{1},\ldots,a_{r}$ eine maximale homogene regul"are Sequenz in $I$, also
$\depth I=r$ (vgl. S"atze \ref{Rees}, \ref{homDepth}). Dann gilt
 $\height(a_{1},\ldots,a_{r})=r$
(Satz \ref{regPhsop} (a)). Nach Lemma \ref{phsopHeight}
(b) gibt es homogene $a_{r+1},\ldots,a_{r+k} \in I$ mit
$\height(a_{1},\ldots,a_{r+k})=r+k=\height(I)$. Dann ist 
\[
k=\height(I)-\depth(I).
\]
Nochmals nach Lemma \ref{phsopHeight} (b) (mit $I=R_{+}$)
erg"anze man $a_{1},\ldots,a_{r+k}$ zu einer Menge homogener Elemente
$a_{1},\ldots,a_{n}\in R_{+}$ mit
$
n=\height(a_{1},\ldots,a_{n})=\height R_{+}=\dim R
$. Nach Lemma \ref{phsopHeight} (a) ist diese Menge insbesondere ein hsop von $R$.
Es folgt
\begin{eqnarray*}
\depth(R)&=&\depth((a_{1},\ldots,a_{n})_{R},R) \quad \textrm{(Korollar \ref{KorDepthR})}\\
&\le& \depth((a_{1},\ldots,a_{r+k})_{R},R)+(n-r-k) \quad \textrm{(Korollar
  \ref{depthkKorollar})}.
\end{eqnarray*}
Nun ist $(a_{1},\ldots,a_{r+k})_{R} \subseteq I$, und $a_{1},\ldots,a_{r}$ ist maximale regul"are Sequenz in
$I$. Also ist es erst recht eine maximale regul"are Sequenz in
$(a_{1},\ldots,a_{r+k})_{R}$,
d.h. 
\[
\depth((a_{1},\ldots,a_{r+k})_{R},R)=r
\]
(Satz \ref{Rees}). Damit folgt aus obiger
Ungleichung also $\depth R \le r + (n-r-k)=n-k$, oder $k=\height I - \depth I
\le n - \depth R = \dim R-\depth R=\cmdef R$.

Ist nun $R$ Cohen-Macaulay, also $\cmdef R=0$, so gilt nach der gerade
bewiesenen Ungleichung $0 \ge \height(I)-\depth(I)$, und wegen $\depth(I)\le
\height (I)$ (siehe ~\eqref{depthheight}) gilt sogar Gleichheit.\qed\\

An der Stelle in obigem Beweis, bei der Korollar \ref{KorDepthR} verwendet
wird, k"onnte man stattdessen auch nochmal Lemma \ref{depthARM} verwenden und
dann mehrmals zwischen der Betrachtung von $R$ als $R$- oder $A$-Modul hin und
her wechseln. Die Betrachtung als $A$-Modul kann auch deshalb n"utzlich sein, weil Ideale in (dem Polynomring!) A oft besser
"uberblickbar sind als in $R$. 

Wir wollen f"ur diesen Satz, der unser wichtigstes Hilfsmittel sein wird, noch
den "`Lehrbuch"'-Beweis angeben (in den Lehrb"uchern findet sich - wie fast
immer - nur der Fall noetherscher lokaler Ringe aufgeschrieben). Die folgende Proposition ist dabei auch von
unabh"angigem Interesse.
\begin{Prop}
Sei $R$ eine graduierte affine Algebra und $M\ne 0$ ein endlich erzeugter graduierter
$R$-Modul. Dann gilt f"ur jedes assoziierte Primideal $\wp\in\Ass_{R}M$
\[
\depth(R_{+},M)\le \dim (R/\wp).
\]
\end{Prop}
\noindent {\bf Bemerkung. }Hieraus folgt mit $\dim M:=\dim
R/\Ann_{R}M=\max\{\dim R/\wp\,|\,\,\wp\in\Ass_{R}M\}$ auch $\depth M\le\dim
M$.\\

\Bew (Bruns und Herzog \cite[Proposition 1.2.13]{BrunsHerzog} oder \cite[Satz
 VII.2.3]{Kunz} f"ur den lokalen Fall). Wir machen Induktion nach der Tiefe $\depth(R_{+},M)$. Ist
 diese gleich $0$, so ist die Aussage klar. Sei also $\depth(R_{+},M)>0$. Dann
 gibt es eine homogene $M$-regul"are Sequenz $a\in R_{+}$ der L"ange $1$, und es gilt
 $\depth(R_{+},M/(a)M)=\depth(R_{+},M)-1$ (Satz \ref{Rees}). Nach
 Induktion gilt also
\begin{equation}\label{IndDepthRp}
\depth(R_{+},M/(a)M)\le \dim (R/{\mathcal P}) \quad \myforall{\mathcal P}\in \Ass_{R}M/(a)M.
\end{equation}
Da $\wp\in \Ass_{R}M$ und $M$ noethersch ist, gibt es ein nach
\cite[Proposition 3.12]{Eisenbud} \emph{homogenes} $0\ne m\in M$, so dass
$Rm$ maximaler von einem homogenen Element erzeugter Modul mit der Eigenschaft $\wp\cdot Rm=0$ ist. Dann
ist die Restklasse $\bar{m}$ von $m$ in $M/(a)M$ ungleich Null: 

Sei n"amlich stattdessen
$m\in (a)M$, also $m=am'$ mit $m'\in M$. Da $m$ und $a$ homogen, und $a$ kein
Nullteiler ist, ist dann auch $m'$ homogen.
Wegen $0=\wp \cdot m=a\cdot
\wp m'$ und der $M$-Regularit"at von $a$ w"are dann auch $\wp\cdot
m'=0$, also $\wp\cdot Rm'=0$ und $Rm\subseteq Rm'$. Da $\deg a >0, m\ne 0$ und $m=am'$
ist jedoch $\deg m'<\deg m$, also $Rm \subset Rm'$ im Widerspruch zur
Maximalit"at von $Rm$.

Da also $\bar{m}\ne 0$ und $\wp \cdot \bar{m}=0$, besteht  $\wp$ nur aus
Nullteilern von $M/(a)M$ und liegt daher in einem assoziierten Primideal
${\mathcal P}\in \Ass_{R} M/(a)M$, $\wp \subseteq {\mathcal P}$. Weiter ist
$a\not\in \wp$ (da $a$ $M$-regul"ar und $\wp\in \Ass_{R}M$) aber
$a\in{\mathcal P}$ (da $a\in \Ann_{R}M/(a)M\subseteq {\mathcal P}$), also $\wp
\subset {\mathcal P}$. Mit Gleichung \eqref{IndDepthRp} folgt also
\[
\depth (R_{+},M)-1=\depth(R_{+},M/(a)M)\le \dim R/{\mathcal P}\le \dim R/\wp-1
\]
und damit die Behauptung. \qed\\

\begin{Korollar}
Sei $R$ eine graduierte affine Algebra und $I\ne R$ ein homogenes Ideal. Dann
gilt
\[
\depth R\le \depth I +\dim I.
\]
Insbesondere gilt $\cmdef R\ge\cmdef I$.
\end{Korollar}
\Bew (Bruns und Herzog \cite[Exercise 1.2.23]{BrunsHerzog}). Sei $a_{1},\ldots,a_{k}\in I$ eine maximale homogene regul"are Sequenz in
$I$, also $k=\depth (I,R)$. Wir w"ahlen $M:=R/(a_{1},\ldots,a_{k})$ in der Proposition. Da $I$ nur
aus Nullteilern von $M$ besteht, also $I\subseteq\bigcup_{\wp\in\Ass_{R}M}\wp$, gibt es ein $\wp\in \Ass_{R} M$ mit
$I\subseteq \wp$. Mit der Proposition und Satz \ref{Rees} gilt dann
\[
\depth (R_{+},M)=\depth R - k \le \dim R/\wp \le \dim R/I,
\]
also $\depth R\le \depth I +\dim I$. Da in jeder affinen Algebra $\dim I +
\height I \le \dim R$ gilt (Proposition \ref{dimRdimIheightI}), folgt hieraus auch
$\depth R\le \depth I +\dim R-\height I$. \qed
\subsection{Lineare algebraische Gruppen und $G$-Moduln}
In diesem Abschnitt beziehen sich topologische Begriffe immer auf die
Zariski-Topologie. \index{Zariski-Topologie} Dabei ist eine Teilmenge von $K^{n}$ genau dann abgeschlossen,
wenn sie eine affine Variet"at\index{affine Variet\"at} ist, d.h. Nullstellenmenge eines Systems von
Polynomen in $K[X_{1},\ldots,X_{n}]$. Morphismen sind Abbildungen zwischen
affinen Variet"aten, die durch Polynome gegeben sind.\index{Morphismus}

\begin{Def}\index{Gruppe!lineare algebraische}
Eine \emph{lineare algebraische Gruppe} ist eine affine Variet"at $G\subseteq
K^{r}$ zusammen mit Morphismen $\cdot:G \times G \rightarrow G$ und $^{-1}: G
\rightarrow G$, so dass $(G,\cdot)$ zusammen mit der durch $^{-1}$ gegebenen
Inversenbildung eine Gruppe wird.
\end{Def}

Wenn nicht anders vermerkt, bezeichnen wir  das Einselement
einer Gruppe immer mit~$\iota$.\\

\noindent {\bf Standardvoraussetzung.} Ab jetzt bezeichnen wir mit $G$ stets
eine lineare algebraische Gruppe, und mit $V$ einen $G$-Modul (siehe die
folgende Definition).

\begin{Def}\index{$G$-Modul}\index{rationale Darstellung}
Ein \emph{G-Modul} oder eine \emph{rationale Darstellung von $G$} ist ein
endlich-dimensionaler $K$-Vektorraum $V$ zusammen mit einer linearen
Operation von $G$ auf $V$, die durch einen Morphismus $G \rightarrow \GL(V)$
gegeben ist.
\end{Def}

In der "ublichen Weise wird $V$ dann Linksmodul "uber dem (im Allgemeinen
nichtkommutativen) \emph{Gruppenring}\index{Gruppenring} $KG$. Dabei handelt
es sich um den $K$-Vektorraum mit Basis $G$, der durch distributive
 Fortsetzung der Multiplikation von $G$ zu einem Ring wird.

\begin{BemRoman}\label{WasIstDarstellung}
Ist $V=K^{n}$, so ist $\GL(V)=\GL_{n}$ als affine Variet"at realisiert durch
die Menge \[X:=\left\{(A,e)\in K^{n\times n}\times K:e\cdot\det A-1=0\right\}.\] Sind
$f_{ij}\in K[G], \, i,j=1,\ldots,n$ so, dass durch $f: G\rightarrow \GL_{n},\,\,
\sigma \mapsto (f_{ij}(\sigma))_{i,j}$ ein Gruppenhomomorphismus gegeben ist,
so ist durch $F: G\rightarrow X,\,\, \sigma\mapsto
\left(f(\sigma),\det(f(\sigma^{-1}))\right)$ ein Morphismus (und damit eine rationale
Darstellung von $G$) gegeben; Denn $\sigma\mapsto\sigma^{-1}$ ist ein
Morphismus (also durch Polynome gegeben), und damit auch
$\sigma\mapsto\det(f(\sigma^{-1}))$, und es ist $\det
f(\sigma^{-1})\det f(\sigma)=\det f(\sigma^{-1}\sigma)=1$.
\end{BemRoman}

Die \emph{$d$-te symmetrische Potenz} \index{symmetrische Potenz} $S^{d}(V)$  eines $G$-Moduls $V$ besteht aus den
"`homogenen Polynomen vom Grad $d$ in den Basiselementen von $V$ und der Null"' - genauer
handelt es sich um den Faktormodul des $d$-fachen Tensorproduktes von
$V$ mit sich selbst modulo allen Relationen, die Polynome erf"ullen
(Kommutativit"at). Insbesondere setzt man $S^{0}(V):=K$. $S^{d}(V)$ ist
in kanonischer Weise ebenfalls ein $G$-Modul.\\ 

\index{Gruppe!Homomorphismus}
\noindent Ein Homomorphismus linearer algebraischer Gruppen $G,H$ ist ein
Morphismus $f: G\rightarrow H$, der zugleich ein Gruppenhomomorphismus
ist. (Insbesondere ist ein $G$-Modul $V$ also nichts anderes als ein algebraischer
Homomorphismus $G\rightarrow \GL(V)$). Dann ist $f(G)\subseteq H$ automatisch eine \emph{abgeschlossene}
Untergruppe von $H$ (\cite[Proposition 2.2.5 (ii)]{SpringerLin}).
Zwei lineare algebraische Gruppen $G,H$ hei"sen (algebraisch) isomorph, wenn
es einen Isomorphismus von Gruppen $f: G\rightarrow H$ gibt, so dass $f$ und
$f^{-1}$ Morphismen von Variet"aten sind.

\begin{BemRoman}\label{LineareFaktorgruppe} \index{Normalteiler}\index{Faktorgruppen}
{\bf (Normalteiler und Faktorgruppen)}
Ist $N$ ein abgeschlossener Normalteiler von $G$, so kann man $G/N$ die
Struktur einer linearen algebraischen Gruppe geben, derart dass die kanonische
Abbildung $G\rightarrow G/N$ ein Homomorphismus linearer algebraischer Gruppen
wird (\cite[Proposition 5.2.5]{SpringerLin} oder \cite[Theorem
11.5]{Humphreys} zusammen mit \cite[Proposition 2.2.5 (ii)]{SpringerLin}).

Ist $X$ eine affine Variet"at und $f: G\rightarrow X$ ein Morphismus von
affinen Variet"aten mit $f(gn)=f(g)\,\,\myforall g\in G,n\in N$, so induziert
$f$ einen \emph{Morphismus} von affinen Variet"aten $G/N\rightarrow X,\,\,
gN\mapsto f(g)$ (\cite[Exercise 5.2.6 (2)]{SpringerLin} oder \cite[Section
12.3]{Humphreys}).

Ist $V$ ein $G$-Modul, so ist $V^{N}$ ein $G/N$-Modul. Ist n"amlich $v\in
V^{N}$ und $g\in G$, so gibt es zu $n\in N$ ein $n'\in N$ mit $ng=gn'$, also
ist $n\cdot gv=g\cdot n'v\stackrel{v\in V^{N}}{=}gv$, also $gv\in
V^{N}$, und damit ist $V^{N}$ ein $G$-Untermodul. Weiter ist f"ur $g\in G,\, n\in N$ stets $gn\cdot v\stackrel{v\in
  V^{N}}{=}g\cdot v$. Damit ist der Morphismus $G\rightarrow \GL(V^{N})$ konstant
auf den Nebenklassen von $N$ und induziert damit nach oben einen Morphismus
$G/N\rightarrow \GL(V^{N})$, so dass $V^{N}$ also ein $G/N$-Modul ist.
\end{BemRoman}

\begin{Satz} \label{treueDarstellung}
Jede lineare algebraische Gruppe ist (algebraisch) isomorph zu einer
abgeschlossenen Untergruppe einer geeigneten $\GL_{n}(K)$. Dann ist $V=K^{n}$
ein $G$-Modul mit einer treuen Darstellung $G\rightarrow \GL(V)$.
\end{Satz}

\Bew Siehe \cite[Theorem 2.3.6]{SpringerLin} \qed\\
 
Eigenschaften, die Matrixgruppen bzw. ihren Elementen zukommen, k"onnen so auch
linearen algebraischen Gruppen zugeordnet werden, indem man sie mittels dieses
Satzes als Untergruppe einer $\GL_{n}(K)$ auffasst. So hei"st ein
$\sigma \in G \subseteq \GL_{n}(K)$ \emph{unipotent},\index{unipotentes Element} wenn es ein $k\in
{\mathbb N}$ gibt mit $(\sigma-\iota)^{k}=0\in K^{n\times n}$ (hier ist
$\iota=I_{n\times n}\in K^{n\times n}$ die Einheitsmatrix).  $G$ hei"st
\emph{unipotent}, \index{Gruppe!unipotente} wenn jedes ihrer
Elemente unipotent ist.
Ein Element $\sigma\in G\subseteq \GL_{n}(K)$ hei"st \emph{halbeinfach},
\index{halbeinfaches Element} wenn $\sigma$
"ahnlich zu einer Diagonalmatrix ist. (Achtung: Eine Gruppe hei"st
\emph{nicht} halbeinfach, wenn jedes ihrer Elemente halbeinfach ist!)

 Die Eigenschaften "`unipotent"' und "`halbeinfach"'
 sind dabei wohldefiniert, denn wenn sie f"ur das Bild eines $\sigma \in G$
 unter \emph{einer} treuen rationalen Darstellung $G\rightarrow \GL_{n}(K)$ gelten, so auch f"ur das Bild unter \emph{jeder} solchen Darstellung
(siehe \cite[Theorem 2.4.8]{SpringerLin}).

\begin{SatzDef}\label{zushgkomp}
 \index{Gruppe!Zusammenhangskomponente}\index{Zusammenhangskomponente einer Gruppe}
Mit $G^{0}$ bezeichnen wir die Zusammenhangskomponente von $G$, die das
neutrale Element $\iota$ enth"alt. Sie ist die eindeutig bestimmte irreduzible
Komponente von $G$, die $\iota$ enth"alt, und ein abgeschlossener Normalteiler
von $G$ von endlichem Index. $G$ hei"st \emph{zusammenh"angend}, wenn $G=G^{0}$. (Siehe
\cite[Proposition 2.2.1]{SpringerLin}). \index{Gruppe!zusammenh\"angende}
\end{SatzDef}

\begin{Def}\label{defred}\index{Gruppe!reduktive}\index{reduktive Gruppe}
$G$ hei"st \emph{reduktiv}, wenn jeder abgeschlossene, zusammenh"angende und
unipotente Normalteiler von $G$ trivial ist, d.h. nur aus dem Einselement $\iota$ besteht.
\end{Def}

 Alle klassischen \index{klassische Gruppen}\index{Gruppe!klassische}
Gruppen $\SL_{n}, \GL_{n}, \Sp_{n}, \SO_{n}, \On_{n}$ sind in allen Charakteristiken
reduktiv (siehe etwa \cite[Chapter 5.9]{SantosRittatore}). Auch alle endlichen Gruppen sind reduktiv, da die einzige
zusammenh"angende Untergruppe die Einsgruppe $\{\iota\}$ ist. 
Dagegen sind die additiven Gruppen $\Ga=(K,+)$ nicht reduktiv, da sie
selbst zusammenh"angend und unipotent sind.\index{Gruppe!additive}

\begin{Def}\index{linear reduktiv}\index{Gruppe!linear reduktive}
$G$ hei"st \emph{linear reduktiv}, wenn jeder $G$-Modul vollst"andig
reduzibel\index{vollst\"andig reduzibel}
ist. D.h. f"ur jeden $G$-Modul $V$ und jeden
$G$-Untermodul $U \le V$ existiert ein Komplement, d.h. ein $G$-Untermodul $W\le V$ mit
$V=U\oplus W$.
\end{Def}

Wir werden noch sehen, dass jede linear reduktive Gruppe auch reduktiv ist.
Die klassischen Gruppen sind nur in Charakteristik $0$ linear reduktiv. In
positiver Charakteristik gibt es nur sehr wenige linear reduktive Gruppen. Es
gilt n"amlich

\begin{Satz}[Nagata \cite{NagataComplete}] \label{NagataTorus}
Sei $\chr K=p > 0$. $G$ ist genau dann linear reduktiv, wenn $G^{0}$ ein
Torus ist und $(G:G^{0})$ nicht durch $p$ teilbar ist.
\end{Satz}

Diesen Satz k"onnen wir erst auf Seite \pageref{BewNagataTorus} beweisen. Deshalb
werden wir Teile von ihm in den folgenden S"atzen nochmals formulieren und ihn
dann aus diesen folgern.

Dabei ist ein Torus \index{Torus} eine zu einer $\Gm^{k}$ isomorphe Gruppe ($k\ge 0$), wobei
$\Gm=\GL_{1}=\{(a,b)\in K^{2}: ab-1=0\}\cong (K\setminus\{0\},\cdot)$ die
multiplikative Gruppe des K"orpers ist. Tori sind in allen Charakteristiken
linear reduktiv \cite[Theorem 2.5.2]{SpringerLin}, in positiver
Charakteristik sind sie nach obigem Satz sogar die
einzigen zusamenh"angenden und linear reduktiven Gruppen.\\

Interessanterweise gen"ugt f"ur die lineare Reduktivit"at sogar die Forderung,
dass ein \emph{spezieller} $G$-Modul
vollst"andig reduzibel ist. Zun"achst noch eine

\begin{BemDef} \label{FpV}
Sei $G$ eine lineare algebraische Gruppe, $\chr K=p>0$, und $V$ ein $G$-Modul
mit Basis $\{X_{1},\ldots,X_{n}\}$. Dann ist der Untervektorraum $\langle
X_{1}^{p},\ldots,X_{n}^{p} \rangle_{K}$ von $S^{p}(V)$ wegen des
Frobenius-Ho\-momorphismus sogar ein $G$-Untermodul und unabh"angig von der speziellen
Wahl der Basis von $V$. 
Daher kann er basisunabh"angig bezeichnet werden als
die \emph{$p$-te Frobenius-Potenz von $V$}, \index{Frobenius Potenz}
\[
F^{p}(V):=\langle
X_{1}^{p},\ldots,X_{n}^{p} \rangle_{K} \subseteq S^{p}(V).
\] 
Offenbar besteht $F^{p}(V)$ genau aus den $p$-ten Potenzen
aller Elemente aus $V$, also \[F^{p}(V)=\left\{f\in S^{p}(V): \textrm{ es gibt
  } v\in V \textrm{ mit }f=v^{p}   \right\}.\]
\end{BemDef}

\begin{Satz}[Nagata \cite{NagataComplete}] \label{NagataLinRedCrit}
Sei $G$ eine lineare algebraische Gruppe, $V$ ein \emph{treuer} $G$-Modul (existiert
stets nach Satz \ref{treueDarstellung})  und $p=\chr K$.

(a) Ist $p=0$, so ist
$G$ genau dann linear reduktiv, wenn $V$ vollst"andig reduzibel ist.

(b) Ist $p>0$ und $G$ zus"atzlich \emph{zusammenh"angend}, so ist $G$
genau dann linear reduktiv, wenn der Untermodul $F^{p}(V)$ von $S^{p}(V)$ ein Komplement in
$S^{p}(V)$ hat. Dies ist weiter genau dann der Fall, wenn $G$ ein Torus ist.
\end{Satz}

\Bew (a) Siehe \cite[Theorem 3]{NagataComplete}.

(b) ist im
Beweis von \cite[Theorem 1]{NagataComplete} versteckt. Da Nagatas Beweis etwas
kryptisch ist und wir das Resultat an entscheidender Stelle verwenden werden,
bringen wir der Vollst"andigkeit halber einen gegl"atteten Beweis im folgenden
Unterabschnitt. Diesen Teil werden wir auch zum Beweis von Satz
\ref{NagataTorus} verwenden.\qed\\

\begin{Korollar} \label{NagatasNoComplementCorollar}
Sei $\chr K=p>0$, $G$ eine lineare algebraische Gruppe, so dass die
Zusammenhangskomponente $G^{0}$ kein Torus ist, und $V$ ein \emph{treuer} $G$-Modul. Dann hat der
Untermodul  $F^{p}(V)$ von $S^{p}(V)$ kein Komplement in
$S^{p}(V)$. Insbesondere ist $G$ nicht linear reduktiv.
\end{Korollar}

\Bew Da $G^{0}$ kein Torus ist, hat nach Satz \ref{NagataLinRedCrit} (b) dann $F^{p}(V)$
(welches sowohl $G^{0}$-, $G$- als auch $\GL(V)$-Untermodul ist) kein
$G^{0}$-invariantes Komplement in $S^{p}(V)$. Dann gibt es erst
recht kein $G$-invariantes Komplement, und damit ist $G$ nicht linear reduktiv. \qed\\

Unter einer speziellen Zusatzvoraussetzung, die f"ur viele
Darstellungen klassischer Gruppen erf"ullt ist, wollen wir kurz noch einen sehr
{\it anschaulichen Beweis} daf"ur geben, dass $F^{p}(V)$ kein Komplement in
$S^{p}(V)$ hat. 
 Wir nehmen 
an, dass $G\subseteq \GL(V)$ eine nichttriviale abgeschlossene unipotente Untergruppe $U$ enth"alt, so dass
$S(V)^{U}=K[f]$ mit einem $f\in V$ gilt. (Da $U$ nichttrivial folgt $\dim V\ge
2$.
Weiter ist $U^{0}\subseteq G^{0}$ eine abgeschlossene, zusammenh"angende nichttriviale unipotente
Untergruppe, so dass $G^{0}$ kein Torus ist - es handelt sich also
tats"achlich um eine Zusatzvoraussetzung).
 Insbesondere ist dann $\dim
S^{p}(V)^{U}=\dim K[f]_{p}=1$. W"are nun $S^{p}(V)=F^{p}(V)\oplus W$ mit einem $G$- (also
erst recht $U$-) invarianten Komplement $W$, so enthielten (etwa nach Springer
\cite[Theorem 2.4.11]{SpringerLin}) sowohl $F^{p}(V)$ als auch $W\ne 0$ (wegen
$\dim V\ge 2$) eine von
$0$ verschiedene
$U$-Invariante, also $\dim S^{p}(V)^{U}\ge 2$, Widerspruch. 

Hier nun eine Situation, unter der die gemachte Voraussetzung gilt.
Wir betrachten $V=\langle X_{1},\ldots,X_{n} \rangle$ als Dual von
$V^{*}:=K^{n}=\langle e_{1},\ldots,e_{n}\rangle$ (mit Standardbasis), also
$X_{i}(e_{j})=\delta_{ij}$ f"ur alle $i,j$. Sei weiter $U\subseteq \GL(V)$ eine
abgeschlossene Gruppe oberer Dreiecksmatrizen mit $\overline{U\cdot
  e_{1}}=\{e_{1}+\sum_{i=2}^{n}\lambda_{i}e_{i}: \lambda_{i}\in K\}$ (Zariski-Abschluss) $(*)$. Dann
gilt $S(V)^{U}=K[V^{*}]^{U}=K[X_{1}]$. Offenbar ist n"amlich $X_{1}\in
K[V^{*}]^{U}$. Sei umgekehrt $f=f(X_{1},\ldots,X_{n})\in K[V^{*}]^{U}$. F"ur
ein festes $x_{1}\in K\setminus\{0\}$ gilt dann wegen $f\in K[V^{*}]^{U}$
stets $f((x_{1},0,\ldots,0))=f(\sigma(x_{1},0,\ldots,0))$ f"ur alle $\sigma\in U$. Wegen $(*)$
und weil Multiplikation mit $x_{1}\ne 0$ ein Hom"oomorphismus ist, gilt dann
$f((x_{1},0,\ldots,0))=f((x_{1},x_{2},\ldots,x_{n}))$ f"ur alle $x_{i}\in
K$ mit $x_{1}\ne 0$. Weil
die Menge der $x_{1}\ne 0$ Zariski-dicht in $K$ liegt, gilt dies dann auch
f"ur $x_{1}\in K$ beliebig.  Insbesondere erhalten wir damit
\[
 f((x_{1},x_{2},\ldots,x_{n}))=f((x_{1},0,\ldots,0))=f(X_{1},0,\ldots,0)(x_{1},x_{2},\ldots,x_{n})
\]
f"ur alle $x_{i}\in K$, und weil $|K|=\infty$ dann $f=f(X_{1},0,\ldots,0)\in K[X_{1}]$. \qed\\

Dagegen ist etwa f"ur $\chr K=2$, $G=Z_{2}=\{\iota,\sigma\}$ (also
$G^{0}=\{\iota\}$ ein Torus) und $V=\langle X,Y\rangle$ die regul"are
Darstellung (mit $\sigma X=Y,\,\,
\sigma Y=X$) durch $K\cdot XY$ ein Komplement zu $F^{2}(V)$ in $S^{2}(V)$
gegeben. Dennoch gilt die Aussage des Satzes in vielen weiteren F"allen, etwa
wenn $p\ge 3$ und $G\subseteq\GL(V)$ eine \emph{Transvektion}
\index{Transvektion} enth"alt, siehe
Satz \ref{PKoz} und die anschliessenden Bemerkungen nach dessen Beweis.\\

Wir f"uhren nun f"ur den Rest der Arbeit noch etwas {\bf Notation}\label{mynotation} ein: Ist $V$ ein
$G$-Modul, so meinen wir mit $V=\langle X_{1},\ldots,X_{n}\rangle$, dass
$\{ X_{1},\ldots,X_{n}\}$ eine Basis von $V$ als $K$-Vektorraum ist,
und eine eventuelle Darstellung $G\rightarrow \GLn,\, \sigma \mapsto A_{\sigma}$
bez"uglich dieser Basis berechnet wird. Ist $G \subseteq \GLn$, so bezeichnen
wir mit $\langle X_{1},\ldots,X_{n}\rangle$ auch die \emph{nat"urliche
  Darstellung} \index{nat\"urliche Darstellung}
von $G$, also f"ur $\sigma=(a_{ij})\in \GLn$ soll $\sigma
X_{j}=\sum_{i=1}^{n}a_{ij}X_{i}$ gelten. Im Fall $n=2$ schreiben wir meist $X$
und $Y$ statt $X_{1}$ und $X_{2}$. Symmetrische Potenzen schreiben wir
entsprechend als homogene Polynome in den Basisvektoren, etwa $S^{2}(\langle X,Y\rangle)=\langle X^{2},Y^{2},XY\rangle$.

Wir werden auch h"aufig den bekannten {\bf Kalk"ul} zum Rechnen mit
Darstellungen verwenden: Hat $V=\langle X_{1},\ldots,X_{n}\rangle$ die
Darstellung $\sigma \mapsto A_{\sigma}=(a_{ij}^{\sigma})$, so hat der \emph{Dual}\index{Dual}
$V^{*}=\Hom_{K}(V,K)=\langle X_{1}^{*},\ldots,X_{n}^{*}\rangle$ mit Operation
$\sigma \cdot \varphi:=\varphi \circ \sigma^{-1}$ (mit $\varphi\in V^{*},
\sigma\in G$) bez"uglich der angegebenen
\emph{Dualbasis}\index{Dual!basis} (also $X_{i}^{*}(X_{j})=\delta_{ij}$) die Darstellung $\sigma
\mapsto A_{\sigma^{-1}}^{T}$.

Ist $W=\langle Y_{1},\ldots,Y_{m}\rangle$ ein weiterer $G$-Modul mit
Darstellung $\sigma\mapsto B_{\sigma}$, so hat $V\otimes W=\langle
X_{1}\otimes Y_{1},\ldots,X_{1}\otimes Y_{m},\ldots,X_{n}\otimes Y_{1},\ldots,X_{n}\otimes Y_{m}\rangle$ die
Darstellung $\sigma\mapsto A_{\sigma}\otimes B_{\sigma}:=(a_{ij}^{\sigma}B_{\sigma})_{i,j=1,\ldots,n}$ (Block-Matrix,
Kronecker-Produkt von Matrizen). Ist $T=(t_{ij})$ Koordinatenmatrix eines
Elements $t=\sum_{i,j}t_{ij}X_{i}\otimes Y_{j}\in V\otimes W$, so hat $\sigma\cdot t$ die
Koordinatenmatrix $A_{\sigma}TB_{\sigma}^{T}$.

\subsubsection{Beweis von Satz \ref{NagataLinRedCrit} (b)}
Wir ben"otigen zun"achst noch zwei Lemmata "uber zusammenh"angende lineare algebraische Gruppen.

\begin{Lemma}\label{LemmaEins}
Falls jedes Element einer \emph{zusammenh"angenden} linearen algebraischen
Gruppe $G$ halbeinfach ist, dann ist $G$ ein Torus.
\end{Lemma}

\Bew Da das einzige halbeinfache und unipotente Element einer Gruppe das
Einselement $\iota$ ist, also insbesondere $\{\iota\}$ die einzige unipotente
Untergruppe von $G$ ist, ist $G$ nach Definition reduktiv. Als
zusammenh"angende reduktive Gruppe wird $G$ nach Humphreys \cite[Theorem 26.3(d)]{Humphreys} von einem maximalen Torus und den "`Wurzeluntergruppen"'
erzeugt. Da die Wurzeluntergruppen nur aus unipotenten Elementen bestehen
(ebenfalls \cite[Abschnitt 26.3]{Humphreys}), sind sie also trivial, und damit
ist $G$ ein Torus.  \qed\\

Ich danke Frank Himstedt f"ur die Hilfe bei diesem Beweis; In
 Nagatas Arbeit \cite[Lemma 3]{NagataComplete} steht hierzu einfach: "`the following was proved by Borel"'.

\begin{Lemma}\label{LemmaZwei}
Ist $G$ zusammenh"angend und $u\in G$ unipotent, so existiert eine
abgeschlossene, zusammenh"angende, unipotente Untergruppe $U$ von $G$ mit
$u\in U$.
\end{Lemma}

\Bew Nach \cite[Theorem 7.3.3 (i)]{SpringerLin} liegt $u$ in einer
zusammenh"angenden, abgeschlossenen aufl"osbaren Untergruppe $B$ von $G$,
einer sog. "`Boreluntergruppe"' (vgl. \cite[vor Theorem
7.2.6]{SpringerLin}). F"ur die aufl"osbare, zusammenh"angende Gruppe $B$ bildet
die Menge all ihrer unipotenten Elemente $U$ nach \cite[Corollary
6.9 (ii)]{SpringerLin} ebenfalls eine abgeschlossene,
zusammenh"angende, unipotente Untergruppe, und offenbar ist $u\in U$. \qed\\

Schliesslich noch ein einfaches Lemma "uber das Transformationsverhalten
gewisser Koordinatenringe.

\begin{Lemma} \label{Koordring}
Sei $V\subseteq K^{m\times n}$ eine affine Variet"at, $\sigma \in \GL_{m},
\tau\in\GL_{n}$ und $W:=\sigma V \tau \subseteq K^{m \times n}$. Dann ist $\varphi: V
\rightarrow W, \,\,v\mapsto \sigma v \tau$ ein Isomorphismus von Variet"aten. Seien
$K[V]=K[X_{ij}: i=1,\ldots,m, \,\, j=1,\ldots,n]$ und $K[W]:=K[Y_{ij}: i=1,\ldots,m, \,\, j=1,\ldots,n ]$ jeweils die
Koordinatenringe mit $X_{ij}(v)=v_{ij}, \, v=(v_{ij})\in V$
bzw. $Y_{ij}(w)=w_{ij},\, w=(w_{ij})\in W$. Sei $\varphi^{*}: K[W]\rightarrow
K[V], f\mapsto f\circ \varphi$ der zu $\varphi$ geh"orige Isomorphismus der
Koordinatenringe. Dann gilt
\[
(\varphi^{*}(Y_{ij}))_{i=1,\ldots,m \atop j=1,\ldots,n}=\sigma\cdot (X_{ij})_{i=1,\ldots,m\atop j=1,\ldots,n}\cdot\tau.
\] 
Au"serdem ist 
\begin{eqnarray*}
\varphi^{*}(K[Y_{ij}^{p}: 1\le i\le
m,\,\, 1\le j\le
n])&=&K[\varphi^{*}(Y_{ij})^{p}:1\le i\le
m,\,\, 1\le j\le
n]\\
&=&K[X_{ij}^{p}:1\le i\le
m,\,\, 1\le j\le
n ]
.
\end{eqnarray*}
\end{Lemma}

\Bew Da $\varphi$ Einschr"ankung eines linearen Isomorphismus von $K^{m\times
  n}$ ist, ist klar dass auch $W$ eine Variet"at, $\varphi$ ein
  Isomorphismus von Variet"aten sowie $\varphi^{*}$ ein Isomorphismus der
  Koordinatenringe ist. Sei nun $v=(v_{ij})\in V, \sigma=(\sigma_{ij})\in \GL_{m},
  \tau=(\tau_{ij})\in \GL_{n}$. Dann ist
\begin{eqnarray*}
\varphi^{*}(Y_{ij})(v)&=&Y_{ij}(\varphi(v))=Y_{ij}(\sigma v \tau)\\&=&\sum_{k=1,\ldots,m,
  l=1,\ldots,n} \sigma_{ik}v_{kl}\tau_{lj}\\&=&\sum_{k=1,\ldots,m,
  l=1,\ldots,n} \sigma_{ik}X_{kl}\tau_{lj}(v) \;\quad\myforall v \in V,
\end{eqnarray*}
also gilt die angegebene Transformationsformel. Die letzte Aussage des Satzes
folgt mit dem Frobenius-Homomorphismus. \qed\\

Damit kommen wir zum\\

\noindent \textit{Beweis von Satz \ref{NagataLinRedCrit} (b).} Es sei
$G$ eine zusammenh"angende lineare algebraische Gruppe,
$\chr K = p > 0$, $V=K^{n}=\langle X_{1},\ldots,X_{n} \rangle$ ein $G$-Modul
mit der treuen rationalen Darstellung $\rho: G\rightarrow \GL(V)=\GL_{n}$ und 
\[
W:=F^{p}(V)=\langle
X_{1}^{p},\ldots,X_{n}^{p} \rangle\subseteq S^{p}(V).
\]
Wenn $G$ ein Torus ist, so ist $G$ nach Springer \cite[Theorem 2.5.2
(c)]{SpringerLin} linear reduktiv (der Beweis ist elementar f"uhrbar).
Wenn $G$ linear reduktiv ist, so ist klar, dass $W$ ein Komplement in $S^{p}(V)$
hat. 

Sei nun umgekehrt $U\subseteq S^{p}(V)$  ein $G$-Untermodul mit 
\[
S^{p}(V)=U\oplus
W.
\]
 Wir zeigen, dass dann $G$ ein Torus ist. Angenommen, $G$ sei
\emph{kein} Torus. Nach
Lemma \ref{LemmaEins} ist dann nicht jedes Element von $G$ halbeinfach. Da es
f"ur jedes $x\in G$ nach \cite[Theorem 2.4.8]{SpringerLin} eine eindeutige
Zerlegung $x=x_{h}x_{u}$ mit $x_{h}\in G$ halbeinfach und $x_{u}\in G$
unipotent gibt, enth"alt also $G$ ein vom neutralen Element verschiedenes
unipotentes Element. Dieses liegt dann nach Lemma \ref{LemmaZwei} in einer
nichttrivialen, zusammenh"angenden, abgeschlossenen unipotenten Untergruppe
$H$ von $G$. Dann ist auch $\rho(H)\subseteq \GL_{n}$ eine nichttriviale, zusammenh"angende,
abgeschlossene unipotente Untergruppe von $\GL_{n}$ (Springer
\cite[Proposition 2.2.5, Theorem 2.4.8]{SpringerLin}). Da $S^{p}(V),U,W$ erst recht $H$-
bzw. $\rho(H)$-Moduln sind, 
k"onnen wir ab jetzt
$G=\rho(H)$ annehmen, um die Existenz der Zerlegung $S^{p}(V)=U\oplus W$ f"ur
eine abgeschlossene,
zusammenh"angende, unipotente, nichttriviale Gruppe $G\subseteq \GL_{n}$ zum Widerspruch zu
f"uhren. Nach \cite[2.4.11]{SpringerLin} gibt es ein $\sigma\in \GL_{n}$, so dass alle
Elemente aus $\sigma G \sigma^{-1}$ obere Dreiecksmatrizen sind. Dann sind
$\sigma S^{p}(V)=S^{p}(V)$, $\sigma U$ und $\sigma W = W$ (Frobenius-Homomorphismus!) jeweils
$\sigma G \sigma^{-1}$-Moduln, und damit $S^{p}(V)=\sigma U \oplus W$, d.h. $W$ hat
auch als $\sigma G \sigma^{-1}$-Modul ein Komplement. Wir nehmen
also weiter an, dass alle Elemente aus $G$ bereits obere Dreiecksmatrizen sind. F"ur
$\sigma=(\sigma_{ij})_{i,j=1,\ldots,n}\in G$ gilt also $\sigma_{ii}=1$ f"ur $i=1,\ldots,n$
und $\sigma_{ij}=0$ f"ur $i>j$. Wir k"onnen also $G$ als abgeschlossene
Teilmenge von $K^{\frac{1}{2}n(n-1)}$ auffassen, d.h. wir identifizieren
$\sigma \in G$ mit $(\sigma_{ij})_{1\le i < j \le n}\in K^{\frac{1}{2}n(n-1)}$. Der entsprechende
Koordinatenring von $G$ ist dann
\[
K[G]=K\left[S_{ij}: 1\le i < j \le n  \right] \quad \textrm{mit }
S_{ij}(\sigma)=\sigma_{ij} \textrm{ f"ur } \sigma=(\sigma_{ij})\in G.
\]
(Denn es ist $S_{ii}=1$ und $S_{ij}=0$ f"ur $i> j$.)
Wir setzen
\[
A:=K\left[S_{ij}^{p}: 1\le i < j \le n  \right]=\{f^{p}: f \in K[G]\},
\]
die Unteralgebra von $K[G]$, die aus allen $p$-ten Potenzen besteht.

\noindent {\bf Behauptung:} Es existiert $1\le k < l \le n$ mit $S_{kl}\notin
A$.

{\it Denn} nach Humphreys \cite[Corollary 17.5]{Humphreys} ist $G$ als
unipotente Gruppe nilpotent, also aufl"osbar. Da $G$ auch zusammenh"angend und
kein Torus ist, gibt es nach Springer \cite[Lemma 6.10]{SpringerLin} eine
abgeschlossene Untergruppe $H\le G$ mit $H\cong \Ga$. Daher gibt es Morphismen
$f: \Ga \rightarrow H$ und $g: H\rightarrow \Ga$ mit $g(f(t))=t \;\myforall t
\in \Ga$. Zu $f$ gibt es dann Polynome in einer Variable $f_{ij}\in K[T]$ mit
$f(t)=(f_{ij}(t))_{i,j=1,\ldots,n}\in H \le G$. Da $g((f_{ij}(t))_{i,j=1,\ldots,n})=t$ und
$g$ ebenfalls durch Polynome gegeben ist, gibt es wenigstens ein Indexpaar
$(k,l)$, so dass $T$ als Monom in $f_{kl}$ vorkommt. Dann gilt $k<l$, denn
$f_{ii}=1$ f"ur alle $i$ und $f_{ij}=0$ f"ur alle $i>j$.
 Es folgt $f_{kl}(T)\notin
K[f_{ij}^{p}(T): 1\le i < j \le n]$, denn in keinem Element der rechten Seite
kommt $T$ als Monom vor. Damit ist erst recht $S_{kl}\notin
A$, sonst g"abe es n"amlich $F\in K[G]$ mit $S_{kl}=F^{p}$, und auswerten an
der Stelle $f(t)\in G$ liefert $f_{kl}(t)=F^{p}((f_{ij}(t))_{i,j=1,\ldots,n}) \myforall
t\in \Ga$, also doch $f_{kl}(T)=F^{p}((f_{ij}(T))_{i,j=1,\ldots,n})\in
K[f_{ij}^{p}(T): 1\le i < j \le n]$ - Widerspruch. Dies zeigt die Behauptung.

Aufgrund der Behauptung existiert nun genau ein Indexpaar $(k,l)$, so dass
gilt

\begin{itemize}
\item $S_{kl}\notin A, \quad S_{ij}\in A$  f"ur $i\le j<l$,  $\quad S_{il}\in A$
  f"ur $i<k$  $\quad\quad$(dabei ist $k<l$). \hfill$(*)$  
\end{itemize}
Man erh"alt es, indem man erst ein maximales $l$ w"ahlt so dass $S_{ij}\in A$ f"ur
alle $i\le j<l$ gilt und dann das kleinste $k$ mit $S_{kl}\notin A$. Auf der Menge
$\{(i,j): 1\le i < j \le n\}$ f"uhren wir nun folgende Ordnung ein:
\[
(i_{1},j_{1}) < (i_{2},j_{2}) :\Leftrightarrow j_{1}<j_{2} \textrm{ oder }
j_{1}=j_{2}, \, i_{1}<i_{2}.
\]
F"ur jede unipotente, obere Dreiecksmatrix $\tau\in \GL_{n}$ ist $\tau G\tau^{-1}$
ebenfalls eine zusammenh"angende, unipotente Gruppe oberer Dreiecksmatrizen,
und daher existiert zu jeder solchen Gruppe genau ein Indexpaar $(k,l)$ mit der
Eigenschaft $(*)$. Da die Menge der Indexpaare endlich ist, k"onnen wir ein
$\tau$ mit \emph{maximalem} zugeh"origen Indexpaar gem"a"s obiger Ordnung w"ahlen, und
wir ersetzen $G$ durch $\tau G \tau^{-1}$ (dann ist entsprechend $\tau U$ das neue
Komplement zu $W$). Wir zeigen nun, dass die Maximalit"at von $(k,l)$ folgende
Konsequenz hat:

\begin{equation}\label{depend1}
\lambda_{i}\in K \textrm{ f"ur }i=1,\ldots,n,\,\, \lambda_{k}\ne 0 \Rightarrow
\sum_{i=1}^{n}\lambda_{i}S_{il} \notin A.
\end{equation}
Sei n"amlich stattdessen $\sum_{i=1}^{n}\lambda_{i}S_{il} \in A$. Da
$\lambda_{k}\ne 0$, k"onnen wir O.E. $\lambda_{k}=1$ annehmen. 
Da $S_{il}\in
A$ f"ur $i<k$ nach $(*)$, $S_{ll}=1\in A$ und $S_{il}=0\in A$ f"ur $i>l$, gilt dann auch 
\begin{equation}\label{dependContradiction}
\textrm{\bf Annahme: }\quad\quad
\sum_{i=k}^{l-1}\lambda_{i}S_{il} \in A, \quad\textrm{ mit }\lambda_{k}=1
\quad \textrm{(und }
\lambda_{i}\in K \textrm{ f"ur } i=k,\ldots,l-1\textrm{)}.
\end{equation}
Sei $I_{m}\in K^{m\times m}$ die $m\times m$ Einheitsmatrix. Wir betrachten
nun die unipotente obere Dreiecksmatrix
\[
\tau:=\left(
\begin{array}{ccc}
I_{k-1}\\
&
\begin{array}{ccccc}
1&\lambda_{k+1}&\dots&\dots&\lambda_{l-1}\\
&1&0&\dots&0\\
&&\ddots&\ddots&\vdots\\
&&&1&0\\
&&&&1
\end{array}\\
&&I_{n-l+1}
\end{array}
\right)\in \GL_{n}.
\]
Beim Inversen $\tau^{-1}$ erhalten dann lediglich die $\lambda_{i}$ in obiger
Matrix ein Minuszeichen.  Sei dann $K[Y_{ij}: 1\le i<j \le n]=K[\tau G\tau^{-1}]$ der zur
Variet"at $\tau G\tau^{-1}$ geh"orige Koordinatenring (wobei wieder $Y_{ii}=1$
und $Y_{ij}=0$ f"ur $i>j$). Nach Lemma
\ref{Koordring} gibt es dann einen Isomorphismus $\varphi^{*}:K[Y_{ij}: 1\le
i<j \le n] \rightarrow K[G]$ mit
\[
(\varphi^{*}(Y_{ij}))_{i,j=1,\ldots,n}=\tau\cdot (S_{rs})_{r,s=1,\ldots,n}\cdot \tau^{-1}.
\]
Ist $e_{i}\in K^{n}$ der $i$-te Einheitsvektor, $\delta_{jk}$ das Kronecker-Symbol und setzen wir $\lambda_{i}=0$
f"ur $i<k$ oder $i\ge l$, so ist also f"ur alle $1\le i,j\le n$
\[
\varphi^{*}(Y_{ij})=e_{i}^{T}\tau\cdot (S_{rs})_{r,s=1,\ldots,n}\cdot (e_{j}-\lambda_{j}e_{k}+\delta_{jk}e_{k}).
\]

F"ur $j<l$ bestehen wegen $k<l$ bereits die Komponenten von
\[
(S_{rs})_{r,s=1,\ldots,n}\cdot (e_{j}-\lambda_{j}e_{k}+\delta_{jk}e_{k})
\]
 nur aus
Linearkombinationen von $S_{rs}$ mit $s<l$, liegen also nach $(*)$ in
$A$. Damit gilt auch \[\varphi^{*}(Y_{ij})\in A\quad \textrm{ f"ur } j<l.\]

F"ur $j=l$ und $i<k$ ist

\[
\varphi^{*}(Y_{il})=e_{i}^{T}\cdot (S_{rs})_{r,s=1,\ldots,n}\cdot
e_{l}=S_{il}\in A
\]
wegen $(*)$.

Schliesslich ist f"ur $j=l$ und $i=k$
\[
\varphi^{*}(Y_{kl})=e_{k}^{T}\tau\cdot (S_{rs})_{r,s=1,\ldots,n}\cdot
e_{l}=\sum_{r=k}^{l-1}\lambda_{r}S_{rl}\in A
\]
nach \eqref{dependContradiction}. Nach der letzten Aussage von Lemma
\ref{Koordring} gilt $A=\varphi^{*}(K[Y_{ij}^{p}: 1\le i
< j \le n])$. Da $\varphi^{*}$ ein Isomorphismus ist, folgt nun aus obigen drei Enthaltenseinsrelationen, dass das zu $\tau G \tau^{-1}$ geh"orige,
$(*)$ erf"ullende eindeutige Indexpaar $(k',l')$ gr"o"ser ist als $(k,l)$, was
im Widerspruch zur Maximalit"at von $(k,l)$ steht. Damit ist die Annahme
\eqref{dependContradiction} falsch, d.h. \eqref{depend1} ist richtig.

Sei nun $\sigma=(\sigma_{ij})_{i,j=1,\ldots,n}\in G$ (fest). Wir betrachten
\[
\phi: K[G]\rightarrow K[G],\, f\mapsto \phi(f)
\]
mit
\[
\phi(f): G \rightarrow K,\, x\mapsto f(\sigma x).
\]
Es folgt
\[
\phi(S_{ij})(x)=S_{ij}(\sigma
x)=\sum_{r=1}^{n}(\sigma_{ir}x_{rj})=\sum_{r=1}^{n}\sigma_{ir}S_{rj}(x) \quad
\myforall x\in G,
\]
also
\[
\phi(S_{ij})=\sum_{r=1}^{n}\sigma_{ir}S_{rj}.
\]
Der Frobenius-Homomorphismus liefert sofort $\phi(A)\subseteq A$. Da
$S_{il}\in A$ f"ur $i<k$ nach $(*)$, ist dann also auch $\phi(S_{il})\in A$, d.h.
\[
\phi(S_{il})=\sum_{r=1}^{n}\sigma_{ir}S_{rl} \in A,
\]
also $\sigma_{ik}=0$ nach \eqref{depend1}. Dies gilt f"ur alle $i<k$, und da
$\sigma$ eine unipotente obere Dreiecksmatrix ist, ist die $k$-te Spalte von
$\sigma$ damit gleich dem $k$-ten Einheitsvektor $e_{k}$. Dies gilt f"ur alle
$\sigma \in G$, d.h. $X_{k}$ ist unter $G$ invariant:
\begin{equation}
  \label{XkisInv}
\sigma \cdot X_{k}=X_{k}\quad \myforall \sigma \in G.  
\end{equation}

Wir kehren nun zur"uck zur Zerlegung $S^{p}(V)=U\oplus W$. F"ur $i\ne k$
betrachten wir die eindeutige Zerlegung
\begin{equation}\label{ZerlXiXk}
X_{i}X_{k}^{p-1}=u_{i}+w_{i} \quad \textrm{mit }u_{i}\in U, w_{i}\in W.
\end{equation}
F"ur $\sigma=(\sigma_{ij})_{i,j=1,\ldots,n}\in G$ ist
\begin{equation} \label{GaufXlXk}
\sigma \cdot X_{l}X_{k}^{p-1}\stackrel{~\eqref{XkisInv}}{=}\sum_{i=1}^{l} \sigma_{il}X_{i}X_{k}^{p-1}
\end{equation}
und damit
\begin{eqnarray*}
\sigma \cdot u_{l}&=&\sigma (X_{l}X_{k}^{p-1}-w_{l})=\sum_{i=1}^{l}
\sigma_{il}X_{i}X_{k}^{p-1}-\sigma w_{l}=\sum_{i=1, i \ne
k}^{l}
\sigma_{il}X_{i}X_{k}^{p-1} + \sigma_{kl}X_{k}^{p}-\sigma w_{l}\\
&\stackrel{\eqref{ZerlXiXk}}{=}&\sum_{i=1, i\ne k}^{l}\sigma_{il}u_{i}+\sum_{i=1, i\ne k}^{l}\sigma_{il}w_{i}+\sigma_{kl}X_{k}^{p}-\sigma w_{l}.
\end{eqnarray*}
Die letzten drei Terme liegen dabei alle in $W=\langle
X_{1}^{p},\ldots,X_{n}^{p}\rangle$, und der erste Term in $U$. Da $U$ ein
$G$-Modul ist, also $\sigma u_{l}\in U$, folgt damit aus $U \cap W=\{0\}$:
\begin{equation} \label{GaufU}
\sigma \cdot u_{l}=\sum_{i=1, i\ne k}^{l}\sigma_{il}u_{i}.
\end{equation}
Sei $(X_{k}^{p})^{*}\in S^{p}(V)^{*}$ das Funktional, dass von einem
Polynom in $S^{p}(V)$ gerade den Koeffizienten von $X_{k}^{p}$ liefert. Dann
ist
\[
f: G\rightarrow K,\,  \sigma \mapsto  (X_{k}^{p})^{*}(\sigma\cdot u_{l})
\]
durch Polynome gegeben, also $f\in K[G]$. Wegen \eqref{GaufU} ist
\[
f(\sigma)=\sum_{i=1, i\ne k}^{l}\sigma_{il}(X_{k}^{p})^{*}(u_{i})=\sum_{i=1,
  i\ne k}^{l}(X_{k}^{p})^{*}(u_{i})S_{il}(\sigma) \quad \myforall \sigma \in G,
\]
also mit $\lambda_{i}:=(X_{k}^{p})^{*}(u_{i})$
\begin{equation}\label{falsSil}
f=\sum_{i=1, i\ne k}^{l}\lambda_{i}S_{il}.
\end{equation}
Es ist aber
\begin{eqnarray*}
(f-S_{kl})(\sigma)&=&(X_{k}^{p})^{*}(\sigma\cdot
u_{l})-\sigma_{kl}\stackrel{\eqref{ZerlXiXk}}{=}(X_{k}^{p})^{*}(\sigma\cdot
(X_{l}X_{k}^{p-1}-w_{l}))-\sigma_{kl}\\
&\stackrel{\eqref{GaufXlXk}}{=}& (X_{k}^{p})^{*} \left(\sum_{i=1}^{l}
  \sigma_{il}X_{i}X_{k}^{p-1} -\sigma w_{l}
\right)-\sigma_{kl}=\sigma_{kl}-(X_{k}^{p})^{*}(\sigma w_{l})-\sigma_{kl}\\
&=&-(X_{k}^{p})^{*}(\sigma w_{l}).
\end{eqnarray*}
Da aber $w_{l}\in W=\langle X_{1}^{p},\ldots,X_{n}^{p}\rangle$, sind die
Koeffizienten der Monome von $\sigma \cdot w_{l}$ Polynome in den $p$-ten
Potenzen $\sigma_{ij}^{p}$ der Koeffizienten von $\sigma$, und damit
\[
f-S_{kl}\in A.
\]
Wegen \eqref{falsSil} ist dann aber
\[
\sum_{i=1, i\ne k}^{l}\lambda_{i}S_{il}-S_{kl} \in A,
\]
im Widerspruch zu \eqref{depend1}. Damit war die Annahme, dass $G$ kein
Torus ist falsch, und der Satz ist bewiesen. \qed\\

\subsection {Erste Kohomologie algebraischer Gruppen}\label{firstCohom}\index{Kohomologie!erste}
Sei $V$ ein $G$-Modul. Ein \emph{(1)-Kozyklus} \index{Kozyklus} ist ein Morphismus \[g:
G\rightarrow V\quad \textrm{mit}\quad g_{\sigma \tau}=\sigma g_{\tau}+g_{\sigma}\,\myforall
\sigma,\tau\in G.\] Insbesondere
gilt f"ur das neutrale Element $\iota \in G: g_{\iota}=g_{\iota\iota}=\iota
g_{\iota}+g_{\iota}$, also $g_{\iota}=0$. 
Ist $v\in
V$, so ist durch $\sigma \mapsto (\sigma-1)v:=\sigma v-v$ ebenfalls ein Kozyklus
gegeben, und ein Kozyklus der sich so schreiben l"asst, hei"st ein
\emph{Korand}\index{Korand} oder ein \emph{trivialer Kozyklus}.\index{Kozyklus!trivialer} Die additive Gruppe aller
Kozyklen wird mit $Z^{1}(G,V)$ bezeichnet, die Untergruppe aller Kor"ander mit
$B^{1}(G,V)$, und die Faktorgruppe (1. Kohomologiegruppe)
\index{Kohomologiegruppe} mit
$H^{1}(G,V)=Z^{1}(G,V)/B^{1}(G,V)$. F"ur ein $g \in Z^{1}(G,V)$ definiert man
den \emph{erweiterten} G-Modul\index{erweiterter $G$-Modul} $\tilde{V}:=V\oplus K$ mit der Operation
$\sigma (v,\lambda):=(\sigma v+\lambda g_{\sigma},\lambda)$ f"ur $v\in V,
\lambda \in K$. Man rechnet sofort nach, dass dadurch tats"achlich eine
Operation auf $\tilde{V}$ definiert ist. In $\tilde{V}$ wird $g$ zu einem
Korand, $g\in B^{1}(G,\tilde{V})$, denn $g_{\sigma}=(g_{\sigma},0)=(\sigma-1)(0,1)\,\,\myforall
\sigma\in G$.

Ist $V=K^{n}$ mit Darstellung $\sigma \mapsto A_{\sigma}\in K^{n\times n}$ und
$\sigma\mapsto g_{\sigma}\in K^{n}$ ein Kozyklus, so hat $\tilde{V}$ die
Darstellung 
\begin{equation}\label{darstVSchlange}
\sigma\mapsto \left(\begin{array}{cc}A_{\sigma} &g_{\sigma}\\&1 \end{array}\right).
\end{equation}
(L"ucken in Blockmatrizen wie hier werden wie "ublich mit Nullen aufgef"ullt).

Seien $V,W$ jeweils $G$-Moduln. Dann ist auch $\Hom_{K}(V,W)$ ein $G$-Modul
mit Operation gegeben durch $\sigma \cdot f:=\sigma \circ f \circ \sigma^{-1}$
f"ur $f\in\Hom_{K}(V,W)$, $\sigma\in G$. F"ur den Fixmodul schreiben wir $\Hom_{K}(V,W)^{G}=:\Hom_{G}(V,W)=\Hom_{KG}(V,W)$.\label{homgvwdef}

\begin{SatzDef} \label{HomKVW0}
Sei $V$ ein $G$-Modul und $W$ ein Untermodul. Dann ist
\[
\Hom_{K}(V,W)_{0}:=\left\{f \in \Hom_{K}(V,W): f|_{W}=0 \right\} \cong W
\otimes (V/W)^{*}
\]
ein Untermodul von $\Hom_{K}(V,W)$. 
Damit gilt
\[
\dim_{K} \Hom_{K}(V,W)_{0}=\dim_{K} (W \otimes
(V/W)^{*})=\dim_{K}W\cdot(\dim_{K}V-\dim_{K}W).
\]
\end{SatzDef}

\Bew Nur die angegebene Isomorphie ist beweisbed"urftig. Sei $\rho: V
\rightarrow V/W$ der kanonische Epimorphismus. Dann ist $\Hom_{K}(V/W,W)
\rightarrow\Hom_{K}(V,W)_{0},\, f\mapsto f\circ \rho$ ein Isomorphismus von $G$-Moduln, und
daher
\[
\Hom_{K}(V,W)_{0}\cong \Hom_{K}(V/W,W) \cong W\otimes (V/W)^{*}.
\]
\qed\\

Zur Anwendung der n"achsten Proposition ist die folgende einfache Aussage oft hilfreich:

\begin{Bemerkung} \label{kleinKompl}
Sei $ U \le V \le W$ eine Kette von $G$-Moduln. Wenn $U$ kein Komplement in $V$ hat, so auch nicht in $W$.
\end{Bemerkung}

\Bew Angenommen, $W=U \oplus U'$. Dann gilt auch $V=U \oplus (U' \cap V)$, im Widerspruch zur Voraussetzung: F"ur $v \in V \le W$ gibt es n"amlich $u \in U, u' \in U'$ mit $v=u+u'$, also $u'=v-u \in U' \cap V$. Au"serdem ist $U \cap (U' \cap V) \subseteq U \cap U' =\{0\}$. \qed

\begin{Prop} \label{kompl}
Sei $W$ Untermodul eines $G$-Moduls $V$, sowie $\iota \in \Hom_{K}(V,W)$ mit
$\iota|_{W}=\id_{W}$. Dann ist durch $\sigma \mapsto
g_{\sigma}:=(\sigma-1)\iota$ ein Kozyklus in $Z^{1}(G,\Hom_{K}(V,W)_{0})$
gegeben, welcher genau dann ein Korand ist, wenn $W$ ein ($G$-invariantes)
Komplement in~$V$ hat.
\end{Prop}

\Bew Zun"achst ist $g_{\sigma} \in \Hom_{K}(V,W)_{0}$ f"ur $\sigma \in G$ zu zeigen: F"ur $w \in W$ ist
\begin{equation} \label{g0}
g_{\sigma}(w)=((\sigma-1)\iota)(w)=\sigma\left (\iota(\sigma^{-1}w) \right) -w \stackrel{\iota|_{W}=\textrm{id}_{W}}{=}\sigma\sigma^{-1}w-w=0,
\end{equation}  
also $g_{\sigma}|_{W}=0$ und damit $g_{\sigma} \in \Hom_{K}(V,W)_{0}$. Es ist
klar, dass $g$ ein Kozyklus ist.

Sei nun $g$ ein Korand, d.h. es gibt $f \in \Hom_{K}(V,W)_{0}$ mit
$g_{\sigma}=(\sigma-1)\iota=(\sigma-1)f \, \myforall \sigma\in G$. F"ur $h:=\iota-f$ folgt dann $\sigma h=h$, also $h \in \Hom_{G}(V,W)$, und $\ker h$ ist damit ein Untermodul ($G$-invarianter Untervektorraum) von $V$. F"ur $w \in W$ gilt ferner $h(w)=\iota(w)-f(w)=w-0=w$, also $h|_{W}=\id_{W}$. Da dann $h(V) = W$, folgt also $h(h(v))=h(v) \myforall v \in V$. Also ist $h$ Projektion auf $W$, und in "ublicher Weise folgt nun $V=W \oplus \ker h$:
\[
\myforall v \in V: \quad v=\underbrace{(v-h(v))}_{\in \ker h}+\underbrace{h(v)}_{\in W},
\] 
denn $h(v-h(v))=h(v)-h(h(v))=h(v)-h(v)=0$. Ferner gilt f"ur $v \in W \cap \ker h$, dass $v\stackrel{h|_{W}=\textrm{id}_{W}}{=}h(v)=0$, also ist die Summe direkt.

Gilt umgekehrt $V=W \oplus U$ mit einem $G$-invarianten Teilraum $U$, so w"ahle  
\[
f \in \Hom_{K}(V,W)_{0} \textrm{ mit } f|_{W}=0, \; f|_{U}=\iota|_{U}.
\]
Wir zeigen $g_{\sigma}=(\sigma-1)f$: F"ur $w\in W, u \in U$ ist
\[
\begin{array}{l}
g_{\sigma}(w+u)\stackrel{(\ref{g0})}{=}g_{\sigma}(u)=((\sigma-1)\iota)(u)=\sigma(\iota(\sigma^{-1}u))-\iota(u)\stackrel{f|_{U}=\iota|_{U}}{=}\\
\sigma(f(\sigma^{-1}u))-f(u)=((\sigma-1)f)(u)\stackrel{f|_{W}=0}{=}((\sigma-1)f)(w+u),
\end{array}
\]
also Gleichheit auf $V=W \oplus U$, und $g$ ist ein Korand. \qed\\

Diese Proposition zeigt, dass wenn es einen $G$-Modul gibt, der einen
Untermodul ohne Komplement besitzt (also wenn $G$ nicht linear reduktiv ist),
dann gibt es einen nichttrivialen Kozyklus (der Gruppe $G$). Ist umgekehrt $V$ ein $G$-Modul
und $g\in Z^{1}(G,V)$ ein nichttrivialer Kozyklus, so hat $V$ kein
$G$-invariantes Komplement
in $\tilde{V}=V\oplus K$ - insbesondere ist $G$
nicht linear reduktiv. W"are n"amlich $K(v,\lambda)$ ein solches ($v\in V,
\lambda \in K\setminus\{0\}$), so w"are $\sigma(v,\lambda)-
(v,\lambda)=(\sigma(v)-v+\lambda g_{\sigma},0) \in K(v,\lambda) \cap V=\{0\}
\, \myforall \sigma \in G$, also $\sigma \mapsto g_{\sigma}=(\sigma
-1)(-\frac{1}{\lambda}v)$ doch ein trivialer Kozyklus. Damit ist folgender
Satz gezeigt (Kemper \cite[Proposition 1]{KemperLinRed}, in weniger moderner  Formulierung im
wesentlichen auch bei Nagata \cite[Theorem 4]{NagataComplete}):

\begin{Prop} \label{LinRedKohomologie}
Eine lineare algebraische Gruppe $G$ ist genau dann linear reduktiv, wenn
jeder Kozyklus ein Korand ist, d.h.
$H^{1}(G,V)=0$ f"ur jeden $G$-Modul $V$ gilt.
\end{Prop}

\begin{Korollar}\label{LinearRedFaktorGrp}
Sei $G$ eine lineare algebraische Gruppe und $N\lhd G$ ein
\emph{abgeschlossener} Normalteiler. Ist $G$ linear reduktiv, so auch
$G/N$. Sind umgekehrt $N$ und $G/N$ linear reduktiv, so auch $G$.
\end{Korollar}

Man beachte hierbei, dass man $G/N$ die Struktur einer linearen algebraischen
Gruppe geben kann, siehe Bemerkung \ref{LineareFaktorgruppe}.\\

\Bew Ist $G$ linear reduktiv, so ist jeder $G/N$-Modul via des Homomorphismus
$G\rightarrow G/N$ auch $G$-Modul mit den gleichen Untermoduln und damit
vollst"andig reduzibel, so dass auch
$G/N$ linear reduktiv ist. Seien umgekehrt $N$ und $G/N$ linear reduktiv und
$V$ ein $G$-Modul. Nach der Proposition m"ussen wir $H^{1}(G,V)=0$ zeigen. Sei
also $g\in Z^{1}(G,V)$. Dann ist erst recht $g|_{N}\in Z^{1}(N,V)=B^{1}(N,V)$, denn
$N$ ist linear reduktiv, also $H^{1}(N,V)=0$ nach der Proposition. Damit gibt
es ein $v\in V$ mit $g_{\tau}=(\tau-1)v\,\,\myforall \tau\in N$. F"ur
$h_{\sigma}:=g_{\sigma}-(\sigma-1)v \,\myforall\sigma\in G$ gilt dann 
\begin{equation}\label{hsigmaN0}
h_{\tau}=0 \quad \myforall \tau \in N
\end{equation}
sowie $h\in Z^{1}(G,V)$ mit $h+B^{1}(G,V)=g+B^{1}(G,V)$. Es
gen"ugt also $h\in B^{1}(G,V)$ zu zeigen. F"ur $\sigma\in G,\,\tau\in N$ gilt
$h_{\sigma\tau}=\sigma
h_{\tau}+h_{\sigma}\stackrel{\eqref{hsigmaN0}}{=}h_{\sigma}$. Weiter ist
$\tau\sigma=\sigma\tau'$ mit einem $\tau'\in N$, und damit ist dann auch
$h_{\tau\sigma}=h_{\sigma\tau'}=h_{\sigma}$ nach eben, insgesamt also
\begin{equation}\label{Ninvarianz}
h_{\sigma\tau}=h_{\tau\sigma}=h_{\sigma} \quad \myforall\sigma\in G,\,\tau\in N.
\end{equation}
Mit der Kozyklus Eigenschaft gilt dann
\[
\tau
h_{\sigma}=h_{\tau\sigma}-h_{\tau}\stackrel{~\eqref{Ninvarianz},~\eqref{hsigmaN0}}{=}h_{\sigma}\quad
\myforall\sigma\in G, \tau\in N,
\]
also $h_{\sigma}\in V^{N}\,\,\myforall\sigma\in G$. Wegen \eqref{Ninvarianz} und
Bemerkung \ref{LineareFaktorgruppe} ist also mit $H: G/N\rightarrow
V^{N},\,\,\sigma N\mapsto h_{\sigma}$ dann $H\in
Z^{1}(G/N,V^{N})=B^{1}(G/N,V^{N})$, denn $G/N$ ist linear reduktiv. Also gibt
es $w\in V^{N}$ mit $H(\sigma N)=(\sigma N-1)w$ d.h.
$h_{\sigma}=(\sigma-1)w\,\myforall \sigma\in G$, und damit $h\in
B^{1}(G,V)$. Also ist $G$ linear reduktiv.\qed\\

\begin{Korollar}[Maschke] \label{Maschke}
Eine endliche Gruppe, aufgefasst als lineare algebraische Gruppe "uber einem
K"orper der Charakteristik $p\ge 0$,  ist genau dann linear reduktiv, wenn $p$ kein
Teiler der Gruppenordnung ist.
\end{Korollar}

\Bew "`$\Rightarrow$."' Sei $p=\chr K $ kein Teiler der Gruppenordnung $|G|$, und $g\in Z^{1}(G,V)$ ein
Kozyklus eines $G$-Moduls $V$. Da $|G|$ in $K$ invertierbar ist, ist
\[
v:=\frac{-1}{|G|}\sum_{\tau \in G}g_{\tau}\in V
\]
wohldefiniert. Aus der Kozyklus Eigenschaft $\sigma
g_{\tau}=g_{\sigma\tau}-g_{\sigma}$ folgt
\begin{eqnarray*}
\sigma v - v &=& \frac{-1}{|G|}\sum_{\tau \in G}\left(g_{\sigma\tau}-g_{\sigma}\right)-v\\
&=&v+|G|\cdot \frac{1}{|G|}g_{\sigma}-v = g_{\sigma} \quad \myforall \sigma \in G.
\end{eqnarray*}
Daher ist jeder Kozyklus $g$ ein Korand, also nach Proposition \ref{LinRedKohomologie} $G$ linear
reduktiv.

"`$\Leftarrow.$"' Sei $p$ ein Teiler von $|G|$. Betrachte den
\emph{regul"aren}\index{regul\"arer Modul} $G$-Modul $V$ mit Basis $\{e_{\sigma}\}_{\sigma \in G}$ und
der Operation gegeben durch $\tau e_{\sigma}=e_{\tau\sigma}$. Offenbar ist
dann
\[
e:=\sum_{\sigma \in G}e_{\sigma} \in V^{G}.
\]
W"are $G$ linear reduktiv, so h"atte $Ke$ ein Komplement $U$ in $V$, also
$V=Ke\oplus U$ (mit $U$ Untermodul). Sei $u=\sum_{\sigma \in
  G}\lambda_{\sigma}e_{\sigma} \in U$ mit $\lambda_{\sigma} \in K \, \myforall
\sigma \in G$. Es folgt
\[
\sum_{\tau\in G}\tau u = \sum_{\sigma \in G}\lambda_{\sigma}
\underbrace{\sum_{\tau \in G}e_{\tau \sigma}}_{e}=\sum_{\sigma \in
  G}\lambda_{\sigma}e \in U \quad \textrm{(da $U$ Untermodul)}.
\]
Da $e \notin U$, folgt
\begin{equation}\label{lambdas}
\sum_{\sigma \in G}\lambda_{\sigma}=0 \quad \myforall u=\sum_{\sigma \in
  G}\lambda_{\sigma}e_{\sigma} \in U. 
\end{equation}
Sei
\[
W:=\left\{\sum_{\sigma \in
  G}\lambda_{\sigma}e_{\sigma}\in V: \sum_{\sigma \in G}\lambda_{\sigma}=0 \right\}.
\]
Dann ist $\dim W=\dim V-1$ und $U\subseteq W$ nach \eqref{lambdas}. Da auch
$\dim U=\dim V -1$, folgt also $U=W$. Da aber $|G|\cdot 1 = 0$, ist auch
$e=\sum_{\sigma \in G}1\cdot e_{\sigma}\in W=U$. Dies ist ein Widerspruch zu
$Ke \cap U = \{0\}$. \qed\\

Aus dem Beweis notieren wir

\begin{Korollar}\label{MaschkeNoComplement}
Sei $G$ eine endliche Gruppe, $p = \chr K$ ein Teiler der Gruppenordnung
$|G|$, $V$ die regul"are Darstellung von $G$ mit Basis
$\{e_{\sigma}\}_{\sigma\in G}$ und
$e:=\sum_{\sigma\in G}e_{\sigma}\in V^{G}$. Dann hat der Untermodul $Ke$ kein
Komplement in $V$.
\end{Korollar}

Nun k"onnen wir den Satz von Nagata beweisen:\\

\label{BewNagataTorus}
\noindent{\it Beweis von Satz \ref{NagataTorus}.} Sei $G$ linear
reduktiv. Nach Korollar \ref{NagatasNoComplementCorollar} ist dann $G^{0}$ ein
Torus. Nach Korollar \ref{LinearRedFaktorGrp} ist auch die (endliche, Satz
\ref{zushgkomp}) Faktorgruppe $G/G^{0}$ linear reduktiv, und nach dem Satz von
Maschke gilt damit $p\not| (G:G^{0})$. Ist umgekehrt $G^{0}$ ein Torus und
$p\not| (G:G^{0})$, so sind der Normalteiler $G^{0}$ bzw. die Faktorgruppe
$G/G^{0}$ nach dem bereits zitierten Springer \cite[Theorem 2.5.2(c)]{SpringerLin} bzw. dem
Satz von Maschke linear reduktiv. Nach Korollar \ref{LinearRedFaktorGrp} ist
dann auch $G$ linear reduktiv. \qed\\

Die folgende Proposition (allgemeiner in Kemper \cite[Proposition
2]{KemperLinRed}, allerdings ohne die Formel \eqref{Annulator}) besagt, dass man
jeden nichttrivialen Kozyklus annullieren kann. Sie wird eines der
wichtigsten Hilfsmittel f"ur den Beweis unseres Hauptresultats sein.
Wir verwenden folgende
Notation: F"ur $G$-Moduln $V,W$ und $\pi\in W^{G}$, $g\in Z^{1}(G,V)$
bezeichnen wir mit $\pi \otimes g\in H^{1}(G,W\otimes V)$ die Restklasse des
durch $\sigma\mapsto\pi \otimes g_{\sigma}$ gegebenen Kozyklus aus
$Z^{1}(G,W\otimes V)$.

\begin{Prop} \label{AnnulatorProp}
Sei $V$ ein $G$-Modul, $g \in Z^{1}(G,V)$ und $\tilde{V}$ der entsprechende
erweiterte $G$-Modul. Sei $\{v_{1},\ldots,v_{n+1}\}$ eine Basis von
$\tilde{V}$, derart dass $\{v_{1},\ldots,v_{n}\}$ Basis von $V$ ist und
$\sigma v_{n+1}=v_{n+1}+g_{\sigma} \, \myforall \sigma \in G$ gilt (vgl. die
Darstellung \eqref{darstVSchlange}, S. \pageref{darstVSchlange}). Ist dann
$\{v_{1}^{*},\ldots,v_{n+1}^{*}\}$ die zugeh"orige Dualbasis von $\tilde{V}^{*}$
(also $v_{i}^{*}(v_{j})=\delta_{ij}$), so ist $\pi:=v_{n+1}^{*}$ invariant, und
der Kozyklus $\sigma \mapsto \pi \otimes g_{\sigma}$ aus
$Z^{1}(G,\tilde{V}^{*} \otimes V)$ ist trivial (d.h. $\pi \otimes g=0 \in H^{1}(G,\tilde{V}^{*} \otimes V)$). Genauer gilt
\begin{equation} \label{Annulator}
\pi \otimes g_{\sigma}=-(\sigma-1)(v_{1}^{*}\otimes
v_{1}+\ldots+v_{n}^{*}\otimes v_{n}) \quad \myforall \sigma \in G.
\end{equation}
\end{Prop}

\Bew Bekanntlich ist durch lineare Fortsetzung von $\varphi \otimes v \mapsto
v\varphi(\cdot)$  ein $G$-Modul-Iso\-mor\-phis\-mus $\tilde{V}^{*} \otimes V \rightarrow
\Hom_{K}(\tilde{V},V)$ gegeben, und wir
identifizieren deshalb beide Moduln miteinander. Es gen"ugt dann zu zeigen, dass beide Seiten von
\eqref{Annulator} auf der Basis $\{v_{1},\ldots,v_{n+1}\}$ von $\tilde{V}$
"ubereinstimmen. Da $\{v_{1},\ldots,v_{n}\}$ Basis von $V$ ist, ist $(v_{1}^{*}\otimes
v_{1}+\ldots+v_{n}^{*}\otimes v_{n})|_{V}=\id_{V}$. Da $V$ ein Untermodul von
$\tilde{V}$ ist, folgt $\sigma\cdot\id_{V}=\id_{V}$, so dass die
rechte Seite von  \eqref{Annulator} eingeschr"ankt auf $V$ gleich $0$ ist, und
damit gleich der linken Seite eingeschr"ankt auf $V$ (da $\pi|_{V}=0$). Es
bleibt "Ubereinstimmung auf $v_{n+1}$ zu zeigen. Die linke Seite liefert
$g_{\sigma}$, und die rechte wegen $(v_{1}^{*}\otimes
v_{1}+\ldots+v_{n}^{*}\otimes v_{n})(v_{n+1})=0$ ebenfalls
\[
-\left(\sum_{i=1}^{n}\sigma v_{i}v_{i}^{*}\left(\sigma
 ^{-1}v_{n+1}\right)\right)=-\sigma\left(\sum_{i=1}^{n}v_{i}v_{i}^{*}\left(v_{n+1}+g_{\sigma^{-1}}\right)\right)=-\sigma
 g_{\sigma^{-1}}=g_{\sigma},
\]
wobei im letzten Schritt die Kozyklus Eigenschaft $0=g_{\sigma
  \sigma^{-1}}=\sigma g_{\sigma^{-1}}+g_{\sigma}$ von $g$ ausgenutzt
  wurde. \qed\\

Mit einem $\pi\in\tilde{V}^{*G}$ wie in Proposition \ref{AnnulatorProp}
erhalten wir eine kurze exakte Sequenz von $G$-Moduln
\begin{equation} \label{KozyklusLiefertSequenz}
0\rightarrow V \hookrightarrow
\tilde{V}\stackrel{\pi}{\rightarrow}K\rightarrow 0,
\end{equation}
wir k"onnen also jedem Kozyklus eine solche kurze exakte Sequenz zuordnen. Umgekehrt "`induziert"' eine kurze exakte Sequenz der Form \eqref{KozyklusLiefertSequenz} eindeutig ein Element aus
$H^{1}(G,V)$. Sei n"amlich $v\in \tilde{V}$ mit $\pi(v)=1$. Es ist $\sigma
v-v\in \ker\pi=V$ ($\sigma\in G$), und damit ist mit $h_{\sigma}:=\sigma v -v$ dann $h\in
Z^{1}(G,V)$. (F"ur die Situation aus der Proposition und $v=v_{n+1}$ erhalten wir gerade $h=g$ zur"uck.)
Ist auch $v'\in \tilde{V}$ mit $\pi(v')=1$, so ist wegen
$\pi(v-v')=0$ dann $u:=v-v'\in V$, und f"ur $h'_{\sigma}:=\sigma v' -v'$ gilt
$h_{\sigma}-h'_{\sigma}=(\sigma-1)u$, also $h-h'\in B^{1}(G,V)$. Der exakten
Sequenz kann also wohldefiniert das Element $h+B^{1}(G,V)\in H^{1}(G,V)$
zugeordnet werden.

Wir verfolgen diesen Zusammenhang allgemeiner in Abschnitt
\ref{Kohomologievongruppen}. Dort werden wir als Verallgemeinerung von
Proposition \ref{AnnulatorProp} Annullatoren von Kozyklen in
h"oherer Kohomologie behandeln.\\

Hat $\tilde{V}$ eine Darstellung wie in Gleichung \eqref{darstVSchlange}, so
hat $\tilde{V}^{*}$ die Darstellung
\[
\sigma\mapsto \left(\begin{array}{cc}A_{\sigma^{-1}}^{T} &\\g_{\sigma^{-1}}^{T}&1 \end{array}\right),
\]
woran man nochmal unmittelbar sieht, dass der letzte Basisvektor $\pi$ invariant ist.\\

Wir untersuchen kurz das Zusammenspiel der Propositionen \ref{kompl} und
\ref{AnnulatorProp}: Sei $V$ ein $G$-Modul, $0\ne v\in V^{G}$ so, dass $Kv$ kein
$G$-invariantes Komplement in $V$ hat. Sei weiter $V=Kv\oplus U$ mit einem
Untervektorraum $U$, und $\iota \in \Hom_{K}(V,Kv)$ mit $\iota|_{U}=0$ und
$\iota(v)=v$. Dann ist durch $g_{\sigma}=(\sigma-1)\iota$ ein nichttrivialer
Kozyklus mit Werten in $M:=\Hom_{K}(V,Kv)_{0}\cong Kv\otimes (V/Kv)^{*}\cong
(V/Kv)^{*}$ gegeben. Der zugeh"orige erweiterte Modul ist dann
$\tilde{M}=\Hom_{K}(V,Kv)_{0}+K\iota=\Hom_{K}(V,Kv)\cong V^{*}$. Es folgt
$\tilde{M}^{*}\cong V$, und wegen $\sigma \iota=\iota+g_{\sigma}$ und
$v(\iota):=\iota(v)=1\cdot v$, kann man $\iota$ und $v$ zu einem Paar dualer
Basen in $\tilde{M}$ bzw. $V\cong \tilde{M}^{*}$ wie in Proposition
\ref{AnnulatorProp} erg"anzen, so dass $v=\iota^{*}$. Nach dieser Proposition
wird dann $g\in H^{1}(G,M)$ von $v\in V\cong\tilde{M}^{*}$ annulliert.

{\bf Quintessenz:}\label{Quintessenz} Jede Invariante ohne Komplement tritt als Annullator eines
nichttrivialen Kozyklus auf.

\subsection{Invariantentheorie}
Sei $V$ ein endlich-dimensionaler $K$-Vektorraum. Nach Wahl einer Basis kann man $V$ als $K^{n}$
auffassen.
Dann bezeichnet $K[V]$ die zum Polynomring $K[X_{1},\ldots,X_{n}]$ isomorphe
Algebra der Polynomfunktionen auf $V$. Die Graduierung ist hier durch den
Totalgrad gegeben, und $K[V]_{d}$ bezeichnet dann die Menge der homogenen
Polynome vom Grad $d$ und die $0$.
 Ist $V$ sogar ein $G$-Modul, so operiert
$G$ auf $K[V]$ mittels $\sigma \cdot f:=f \circ \sigma^{-1}$ f"ur $f \in K[V],
\sigma \in G$. Damit ist $K[V]$ isomorph zur \emph{symmetrischen Algebra}\index{Algebra!symmetrische} $S(V^{*})$
des Duals, $K[V]\cong S(V^{*}):=\sum_{d=0}^{\infty}S^{d}(V^{*})$.

\subsubsection{Invariantentheorie reduktiver Gruppen}
Wir kommen zur invariantentheoretischen Charakterisierung der Reduktivit"at.

\begin{Satz}\label{SatzNagataMiyataHaboush}
Sei $G$ eine lineare algebraische Gruppe.

(a) $G$ ist genau dann reduktiv, wenn $G$ \emph{geometrisch
  reduktiv}\index{reduktiv!geometrisch}\index{geometrisch reduktiv} ist,
d.h. wenn f"ur jeden $G$-Modul $V$ und $0 \ne v \in V^{G}$ ein homogenes $f\in K[V]_{+}^{G}$
existiert mit $f(v)\ne 0$.

(b) $G$ ist genau dann linear reduktiv,  \index{linear reduktiv}\index{reduktiv!linear}  wenn f"ur jeden $G$-Modul $V$ und $0
\ne v \in V^{G}$ ein $f\in K[V]_{1}^{G}=V^{*G}$
existiert mit $f(v)\ne 0$.

Insbesondere ist jede linear reduktive Gruppe auch reduktiv.
\end{Satz}

Da wir nun alle Arten von Reduktivit"at beisammen haben, noch eine Bemerkung
zum Begriffsbabylon in der Literatur: In "alteren Arbeiten wie denen
von Nagata wird mit dem Wort "`semi-reduktiv"' \index{semi-reduktiv} unser \emph{geometrisch reduktiv}
bezeichnet (was also  nach Satz \ref{SatzNagataMiyataHaboush} mit unserem reduktiv zusammenf"allt), und mit dem Wort "`reduktiv"' unser \emph{linear
  reduktiv}. (Nagata definiert "`semi-reduktiv"' "uber den Dual und mit
speziellen Kor"andern, man sieht aber leicht, dass seine Definition zu unserer von
"`geometrisch reduktiv"' "aquivalent ist).
Unser "`reduktiv"' wird daher zur besseren Unterscheidung oft als
\emph{gruppentheoretisch reduktiv}\index{reduktiv!gruppentheoretisch} bezeichnet.
Insbesondere bei der Verwendung von "`reduktiv"' ist also
Vorsicht geboten. 
 \\

\Bew (a) Nach Nagata und Miyata \cite[Theorem 2]{NagataMiyata} ist jede
geometrisch reduktive (also "`semi-reduktive"' bei Nagata) Gruppe
reduktiv. Mumford et al. \cite{Mumford} folgend kann man diese Beweisrichtung auch so f"uhren:
Wenn $G$ geometrisch reduktiv ist, so ist nach Nagata \cite{NagataAffine}
$K[X]^{G}$ endlich erzeugt f"ur jede $G$-Variet"at $X$. Nach Popov
\cite{Popov} folgt hieraus jedoch, dass $G$ reduktiv ist.

F"ur die Umkehrung siehe Haboush \cite{Haboush}.

(b) Ist $G$ linear reduktiv und $0 \ne v \in V^{G}$, so existiert ein
Untermodul $U$ von $V$ mit $V=Kv\oplus U$. Dann existiert ein lineares
Funktional $f \in V^{*}$ mit $f(v)=1, f|_{U}=0$, und ein solches ist
$G$-invariant.

Sei umgekehrt die zweite Bedingung erf"ullt. Wir zeigen, dass jeder Kozyklus
$g\in Z^{1}(G,V)$ eines $G$-Moduls $V$
trivial ist, wobei wir die Bezeichnungen aus Proposition \ref{AnnulatorProp}
verwenden. Mit Proposition \ref{LinRedKohomologie} folgt dann die lineare
Reduktivit"at von $G$. Offenbar ist $v_{n+1}^{*} \in \tilde{V}^{*G}$. Daher
existiert nach Voraussetzung
ein $v \in \tilde{V}^{**G}=\tilde{V}^{G}$ mit $v(v_{n+1}^{*})=v_{n+1}^{*}(v) \ne 0$, wobei wir den "ublichen
Isomorphismus $V\rightarrow V^{**}$ verwendet haben. Wir k"onnen
O.E. $v_{n+1}^{*}(v)=1$ annehmen, und dann ist $v=v_{n+1}-u$ mit einem $u\in V$. Da
$v\in \tilde{V}^{G}$, gilt $\sigma v = v \, \myforall \sigma \in G,$ also
$v_{n+1}+g_{\sigma}-\sigma u= v_{n+1}-u$, oder $g_{\sigma}=(\sigma-1)u \,\myforall
\sigma \in G$. Da $u\in V$ ist also $g\in B^{1}(G,V)$ ein Korand. \qed\\

In Charakteristik $0$ fallen die Begriffe reduktiv und linear reduktiv
zusammen:
\begin{Satz}[Nagata und Miyata \cite{NagataMiyata}] \label{NagMiy}
In Charakteristik $0$ ist eine lineare algebraische Gruppe genau dann
reduktiv, wenn sie linear reduktiv ist.
\end{Satz}

Der folgende Beweis ist im wesentlichen der Originalbeweis, wobei wir
jedoch Nagatas Begrifflichkeiten in unsere "aquivalente "ubersetzt
haben. (Insbesondere Nagatas "`semi-reduktiv"' in unser geometrisch reduktiv.)\\

\Bew Wir wissen bereits, dass jede linear reduktive Gruppe auch reduktiv ist. 

Sei
umgekehrt also $G$ reduktiv in Charakteristik $0$. Wir zeigen wieder, dass
jeder Kozyklus $g\in
Z^{1}(G,V)$ eines beliebigen $G$-Moduls $V$ ein Korand ist, was nach
Proposition \ref{LinRedKohomologie} die lineare Reduktivit"at von $G$ zeigt. Wir verwenden wieder die Notation von Proposition
\ref{AnnulatorProp}. Au"serdem identifizieren wir $\tilde{V}=K[\tilde{V}^{*}]_{1}$. Da
$v_{n+1}^{*} \in \tilde{V}^{*G}$, existiert wegen der Reduktivit"at von $G$ ein homogenes $f\in
K[\tilde{V}^{*}]^{G}=S(\tilde{V})^{G}$ vom Grad $d\ge 1$ mit $f(v_{n+1}^{*})=1$ (O.E.). Damit
hat $f$ die folgende Form:
\[
f=v_{n+1}^{d}+a_{d-1}v_{n+1}^{d-1}+\ldots+a_{0} \textrm{ mit } a_{k} \in S^{d-k}(V).
\]
Also ist
\begin{eqnarray*}
\sigma f &=& (v_{n+1}+g_{\sigma})^{d}+\sigma
a_{d-1}(v_{n+1}+g_{\sigma})^{d-1}+\ldots\\
&=&v_{n+1}^{d}+(dg_{\sigma}+\sigma a_{d-1})v_{n+1}^{d-1}+\ldots
\end{eqnarray*}
Da $S(V)$ $G$-invariant ist und $\sigma f=f$, erh"alt man durch
Koeffizientenvergleich bei $v_{n+1}^{d-1}$ also $dg_{\sigma}+\sigma
a_{d-1}=a_{d-1}$ mit $a_{d-1}\in  S^{1}(V)=V$. Da $\chr K = 0$, ist $d$ invertierbar,
also $g_{\sigma}=(\sigma-1)(-\frac{1}{d}a_{d-1}) \, \myforall \sigma \in G$ mit
$-\frac{1}{d}a_{d-1}\in V$. Damit ist
$g\in B^{1}(G,V)$ ein Korand.

\qed\\

Wie wir bereits im Beweis von Satz \ref{SatzNagataMiyataHaboush} bemerkt
haben, ist nach Nagata \cite{NagataAffine} f"ur reduktive Gruppen $G$ und einen
$G$-Modul $V$ der Invariantenring $K[V]^{G}$ endlich erzeugt (das gilt sogar
wenn $V$ blos eine \emph{$G$-Variet"at} ist), und damit also
ein graduierter, affiner Bereich (und nach Popov \cite{Popov} folgt umgekehrt aus
der endlichen Erzeugbarkeit von $K[X]^{G}$ f"ur jede \emph{$G$-Variet"at} $X$
die Reduktivit"at von $G$). Damit k"onnen also Invariantenringe
reduktiver Gruppen mit den in Abschnitt \ref{gradAffineAlgs} entwickelten
Methoden und Begriffen untersucht und beschrieben werden. Dies ist einer der
Gr"unde, warum man heute ausreichend Theorie nur f"ur Invariantenringe
reduktiver Gruppen zur Verf"ugung hat. Das folgende Lemma, das im wesentlichen
von Nagata stammt, ist eines
der
Hauptwerkzeuge bei der Arbeit mit Invariantenringen reduktiver Gruppen:

\begin{Lemma}[Nagata] \label{NagatasRedGrpLemma}
Sei $G$ eine reduktive Gruppe, $V$ ein $G$-Modul und $I\subseteq K[V]^{G}$ ein
Ideal. Dann gilt
\[\sqrt{IK[V]}\cap K[V]^{G}=\sqrt{I}.\]
\end{Lemma}

\Bew (vgl. \cite[Lemma 3]{KemperLinRed}). Nach Nagata \cite[Lemma 5.2.B]{NagataAffine}
(oder Newstead \cite[Lemma 3.4.2]{Newstead}) gilt
\[
IK[V]\cap K[V]^{G} \subseteq \sqrt{I}
\]
f"ur jedes endlich erzeugte Ideal $I$ von $K[V]^{G}$. Da Nagata in derselben
Arbeit zeigt, dass $K[V]^{G}$ eine endlich erzeugte $K$-Algebra, also
noethersch ist, gilt diese Gleichung also f"ur jedes Ideal von $K[V]^{G}$.
Gilt nun  $f\in\sqrt{IK[V]}\cap K[V]^{G}$, so existiert ein $n$ mit $f^{n}\in IK[V]\cap
   K[V]^{G}$, also $f\in\sqrt{I}$ nach obiger Gleichung. Gilt umgekehrt $f\in \sqrt{I}$, so gibt es ein $n$ mit $f^{n}\in I\subseteq IK[V]\cap
   K[V]^{G} $, und dann gilt auch $f\in\sqrt{IK[V]}\cap K[V]^{G}$. 
\qed\\

Damit k"onnen wir die f"ur uns wichtigste 
Eigenschaft reduktiver Gruppen beweisen (Kemper \cite[Lemma 4]{KemperLinRed}):

\begin{Lemma} \label{redphsop}
Ist $G$ reduktiv und bilden $a_{1},\ldots,a_{k}\in K[V]^{G}_{+}$ ein phsop im Polynomring
$K[V]$, so bilden sie auch ein phsop im Invariantenring $K[V]^{G}$.
\end{Lemma}

\Bew Sei $I:=(a_{1},\ldots,a_{k})_{K[V]^{G}}$. Nach Lemma \ref{phsopHeight}
ist $\height (I)=k$ zu zeigen. Sei $\wp \supseteq I$ ein minimales
Primideal. Nach dem Krullschen Hauptidealsatz (\cite[Theorem 13.2]{Eisenbud})
gilt $\height (\wp) \le k$, d.h. wir m"ussen noch $\height (\wp) \ge k$ zeigen.
 Nach Lemma \ref{NagatasRedGrpLemma} gilt
\[
\wp \supseteq \sqrt{I}=\sqrt{IK[V]}\cap K[V]^{G}=\bigcap_{\mathcal{P}\supseteq
  IK[V] \atop \textrm{min. Primideal}} \mathcal{P} \cap K[V]^{G}.
\]
Da $\wp$ prim und rechts ein Schnitt "uber endlich viele Ideale steht
(\cite[Satz III.1.10.c]{Kunz}), umfasst $\wp$ eines dieser Ideale. Daher existiert ein minimales Primideal
$\mathcal{P}\supseteq IK[V]$ mit $\wp \supseteq \mathcal{P} \cap K[V]^{G}
\supseteq I$. Da aber $\wp$ ein minimaler Primteiler von $I$ ist, folgt
$\wp=\mathcal{P} \cap K[V]^{G}$. Da $\mathcal{P}$ ein minimaler Primteiler von
$IK[V]=(a_{1},\ldots,a_{k})_{K[V]}$ ist, gilt nach Krulls Hauptidealsatz
$\height (\mathcal{P})\le k$. Da $a_{1},\ldots,a_{k}$ ein phsop in $K[V]$ ist,
gilt aber auch
$k=\height(a_{1},\ldots,a_{k})_{K[V]}\le \height (\mathcal{P})$ nach Lemma \ref{phsopHeight} und nach Definition der H"ohe, insgesamt
also $\height(\mathcal{P})=k$. Wir erg"anzen nun das phsop zu einem hsop
$a_{1},\ldots,a_{n}$ von $K[V]$ (wobei i.a. $a_{k+i}\not\in K[V]^{G}$), und
setzen $A:=K[a_{1},\ldots,a_{n}]$. Offenbar gilt $\mathcal{P}\cap A \supseteq
(a_{1},\ldots,a_{k})_{A}$, und da $K[V]/A$ eine ganze Erweiterung von
Integrit"atsringen mit $A$ normal ist, gilt "`going-down"' und damit
(\cite[Korollar VI.2.9]{Kunz})
$k=\height(\mathcal{P})=\height(\mathcal{P} \cap A) \ge \height
(a_{1},\ldots,a_{k})_{A}=k$, also $\mathcal{P} \cap
A=(a_{1},\ldots,a_{k})_{A}$, denn auf beiden Seiten stehen Primideale von
$A$. Nun gibt es nach "`going-down"'
(\cite[Theorem 13.9]{Eisenbud}) eine Kette von Primidealen
$\mathcal{P}=\mathcal{P}_{k}\supset \mathcal{P}_{k-1}\supset \ldots\supset
\mathcal{P}_{0}$ von $K[V]$ mit $\mathcal{P}_{i}\cap A=(a_{1},\ldots,a_{i})_{A}$ f"ur
$i=0,\ldots,k$. Dann ist $\wp=\mathcal{P}_{k}\cap K[V]^{G}\supseteq
\mathcal{P}_{k-1}\cap K[V]^{G}\supseteq \ldots\supseteq
\mathcal{P}_{0}\cap K[V]^{G}$ eine Kette von Primidealen in $K[V]^{G}$ mit
echten Inklusionen; Es gilt n"amlich
$a_{i}\in(a_{1},\ldots,a_{i})_{A}\subseteq \mathcal{P}_{i}$ f"ur $1\le i\le k$, also $a_{i}\in\mathcal{P}_{i}\cap K[V]^{G} $, und w"are 
$a_{i}\in\mathcal{P}_{i-1}\cap K[V]^{G}$, so h"atte man auch
$a_{i}\in\mathcal{P}_{i-1}\cap A=(a_{1},\ldots,a_{i-1})_{A}$, ein
Widerspruch (da $A$ Polynomring). Dies zeigt $\height(\wp)\ge k$. \qed\\

Die Invariantentheorie linear reduktiver Gruppen wird entscheidend gepr"agt von

\begin{Satz}[Hochster und Roberts \cite{HochsterRoberts}] \label{HoRo}
Ist $G$ \emph{linear reduktiv}, so ist $K[V]^{G}$ Cohen-Macaulay f"ur jeden $G$-Modul
$V$.
\end{Satz}\index{Satz!Hochster und Roberts}

Zusammen mit dem Satz von Nagata und Miyata \ref{NagMiy} sind also in
Charakteristik $0$ Invariantenringe reduktiver Gruppen stets
Cohen-Macaulay. Bis heute kennt man in Charakteristik $0$ kein Beispiel eines
nicht Cohen-Macaulay Invariantenringes. 

\subsubsection{Roberts' Isomorphismus}\index{Roberts' Isomorphismus}
Wir werden sp"ater massiv von einem Isomorphismus von $\SL_{2}$-Invariantenringen
nach $\Ga$-Invariantenringen gebrauch machen, der (in Charakteristik $0$) auf Roberts
\cite{Roberts} aus dem Jahr 1861 zur"uckgeht. Wir geben hier einen Beweis, der
im wesentlichen eine Spezialisierung des in \cite[Theorem 3.2]{Bryant}
gegebenen Beweises f"ur einen allgemeineren Satz ist. Weitere Referenzen sind
Tyc \cite{Tyc}, Seshadri \cite{Seshadri} und Weitzenb"ock \cite{Weitzenbock}.

Ist $V$ ein $\SL_{2}$-Modul, so ist $V$ auch ein $\Ga$-Modul via 
\begin{equation} \label{sigmat}
t\cdot
v:=\sigma_{t}\cdot v \,\quad 
\textrm{mit}\quad
\sigma_{t}:=
\left(
\begin{array}{cc}
1& t\\
0&1\\ 
\end{array}
\right) \in \SL_{2} \quad
\textrm{f"ur } v\in V,\, t\in \Ga.
\end{equation}

\begin{Satz}[Roberts' Isomorphismus] \label{RobertsIsom}
Sei $V$ ein beliebiger $\SL_{2}$-Modul und $U=K^{2}$ der $\SL_{2}$-Modul mit
der nat"urlichen Darstellung
\[
\left(
\begin{array}{cc}
a& b\\
c&d\\ 
\end{array}
\right)
\left(
\begin{array}{c}
x_{1}\\
x_{2}\\ 
\end{array}
\right)
=\left(
\begin{array}{c}
ax_{1}+bx_{2}\\
cx_{1}+dx_{2}\\ 
\end{array}
\right) \quad \textrm{f"ur }
\left(
\begin{array}{cc}
a&b\\
c&d\\ 
\end{array}
\right) \in \SL_{2}(K), \,\left(
\begin{array}{c}
x_{1}\\
x_{2}\\ 
\end{array}
\right)\in U.
\]
Wir schreiben $f\in K[U\oplus V]^{\SL_{2}}$ mit drei Koordinaten,
 d.h. f"ur die Auswertung von $f$ an einem
$((x_{1},x_{2}),v)\in U\oplus V$ schreiben wir $f(x_{1},x_{2},v)$. Dann ist
durch
\[
\phi: K[U\oplus V]^{\SL_{2}} \rightarrow K[V]^{\Ga}, \quad f
\mapsto \phi(f):=f(1,0,\cdot)
\]
(d.h. f"ur $v\in V$ ist $\phi(f)(v):=f(1,0,v)$) ein Algebren-Isomorphismus gegeben.
\end{Satz}

\Bew Wir pr"ufen zun"achst die Wohldefiniertheit von $\phi$. Sei also $f\in
K[U\oplus V]^{\SL_{2}}$. Zu zeigen ist, dass $\phi(f)\in K[V]$ dann
$\Ga$-invariant ist. F"ur $t\in \Ga$, $\sigma_{t}$ wie oben und $(1,0)\in U$ ist
$\sigma_{t}(1,0)=(1,0)$. Daher ist f"ur $t\in \Ga, v\in V$
\begin{eqnarray*}
(t\cdot \phi(f))(v))&=&\phi(f)(\sigma_{-t}v)=f(1,0,\sigma_{-t}v)=(f\circ
\sigma_{-t})(1,0,v)\\&\stackrel{f\in K[U\oplus V]^{\SL_{2}}}{=}&f(1,0,v)=\phi(f)(v).
\end{eqnarray*}
Da dies f"ur alle $v\in V$ gilt, ist also $t\cdot \phi(f)=\phi(f)$ f"ur $t \in
\Ga$, oder $\phi(f)\in K[V]^{\Ga}$.

Offenbar ist $\phi$ linear und erf"ullt auch $\phi(fg)=\phi(f)\phi(g)$ f"ur
$f,g \in K[U\oplus V]^{\SL_{2}}$,
d.h. $\phi$ ist ein Algebren-Homomorphismus.

Nun zur Injektivit"at von $\phi$. Sei $f\in K[U\oplus V]^{\SL_{2}}$ und
$\phi(f)=0$. Dann ist 
\[
f(1,0,v)=0 \quad \myforall v\in V.
\]
Da f"ur alle $\sigma=\left(
\begin{array}{cc}
a& b\\
c&d\\ 
\end{array}
\right)\in \SL_{2}$ gilt dass $f=f\circ \sigma^{-1}$, ist dann auch 
\[
f(1,0,v)=
\left(f\circ\left(
\begin{array}{cc}
d& -b\\
-c&a\\ 
\end{array}
\right) \right)(1,0,v)=f(d,-c,\sigma^{-1}v)=0 \quad\myforall v\in V.
\]
 Ist $w\in V$
beliebig und setzt man $v:=\sigma w$, so gilt also $f(d,-c,w)=0$ f"ur alle
$w\in V$ und alle $c,d$, die an entsprechender Stelle in Elementen aus
$\SL_{2}$ vorkommen k"onnen, also $(c,d)\ne(0,0)$. Da die beschriebene Menge
aller $(d,-c,w)$ in
$K^{2}\oplus V$ Zariski-dicht liegt, ist also $f=0 \in K[U\oplus V]$. Also ist
$\phi$ injektiv.

Nun zum schwierigsten Teil des Beweises, der Surjektivit"at von
$\phi$. Wir konstruieren ein Urbild zu einem $g\in K[V]^{\Ga}$. Dazu
betrachten wir
\[
G: \SL_{2}\times V \rightarrow K, \quad (\sigma,v)\mapsto
G(\sigma,v):=g(\sigma^{-1}v).
\]
F"ur $\sigma_{t}$ wie in \eqref{sigmat} gilt $g\circ \sigma_{-t}=g$ wegen
$g\in K[V]^{\Ga}$. Also gilt
\begin{equation} \label{eqG}
G(\sigma\sigma_{t},v)=g(\sigma_{-t}\sigma^{-1}v)=g(\sigma^{-1}v)=G(\sigma,v)
\quad \myforall \sigma \in \SL_{2}, t\in \Ga,v \in V.
\end{equation} 
Wir interpretieren diese Gleichung nun f"ur $\sigma=\left(
\begin{array}{cc}
a& b\\
c&d\\ 
\end{array}
\right)\in \SL_{2}$ (also $ad-bc=1$) und unterscheiden die F"alle $a\ne 0$ und $c\ne 0$.

\underline{1.Fall: $a\ne 0$.} Dann setzen wir $t:=-\frac{b}{a}$ und erhalten
aus \eqref{eqG} 
\begin{eqnarray*}
G(\sigma,v)&=&G(\sigma\sigma_{t},v)=G\left(
\left(
\begin{array}{cc}
a& 0\\
c&\frac{1}{a}\\ 
\end{array}
\right)
,v\right)\\&=&\frac{1}{a^{k}}p_{1}(a,c,v) \quad \textrm{mit } k \in \mathbb{N}_{0}
\textrm{ und } p_{1} \in K[U\oplus V]. 
\end{eqnarray*}
Im letzten Schritt wurde verwendet, dass $G(\sigma,v)=g(\sigma^{-1}v)$ ein Polynom in den Koordinaten
von $\sigma$ und $v$ ist, da $g$ ein Polynom ist und die Operation von
$\SL_{2}$ auf $V$ durch Polynome in den Koordinaten $a,b,c,d$ eines $\sigma \in
\SL_{2}$ gegeben ist. Da man f"ur eine Variable $\frac{1}{a}$ einsetzt, kann
man einen Hauptnenner der Form $a^{k}$ ausklammern, und erh"alt dann das
Polynom $p_{1}$. Mit diesem (festen) Polynom gilt obige Gleichung dann f"ur
alle $\sigma,\,v$ mit $a\ne 0$.

\underline{2.Fall: $c\ne 0$.} Dann setzen wir $t:=-\frac{d}{c}$ und erhalten
aus \eqref{eqG} analog wie oben
\begin{eqnarray*}
G(\sigma,v)&=&G(\sigma\sigma_{t},v)=G\left(
\left(
\begin{array}{cc}
a& -\frac{1}{c}\\
c& 0\\ 
\end{array}
\right)
,v\right)\\&=&\frac{1}{c^{l}}p_{2}(a,c,v) \quad \textrm{mit } l \in \mathbb{N}_{0}
\textrm{ und } p_{2} \in K[U\oplus V]. 
\end{eqnarray*}
Sind nun $a\ne 0$ und $c\ne 0$, so sind beide Ausdr"ucke f"ur $G(\sigma,v)$
gleich, und man erh"alt
\[
c^{l}p_{1}(a,c,v)=a^{k}p_{2}(a,c,v) \quad \myforall a\ne 0,c\ne 0,v\in V.
\]
Da die rechts beschriebene Menge jedoch Zariski-dicht liegt und es sich um
eine Polynomgleichung handelt, gilt die Identit"at allgemein, also gilt auch mit
unabh"angigen Variablen $A,C$ und $X=(X_{1},\ldots,X_{n}), n=\dim V$ eines
Polynomrings $K[A,C,X]$, dass
\[
C^{l}p_{1}(A,C,X)=A^{k}p_{2}(A,C,X).
\]
Da $A,C$ teilerfremd sind, folgt $C^{l}|p_{2}(A,C,X)$ und
$A^{k}|p_{1}(A,C,X)$. Daher gilt mit dem Polynom
$f:=p_{1}(A,C,X)/A^{k}=p_{2}(A,C,X)/C^{l}\in K[U\oplus V]$, dass
\begin{equation}  \label{defvonf}
G(\sigma,v)=g(\sigma^{-1}v)
=g\left(
\left(
\begin{array}{cc}
d& -b\\
-c& a\\ 
\end{array}
\right)v\right)
=f(a,c,v) \quad \myforall \sigma\in\SL_{2}, v\in V.  
\end{equation}

Wir behaupten, dass $f$ das gesuchte Urbild zu $g$ ist, also $f\in K[U\oplus
V]^{\SL_{2}}$ und $g=\phi(f)$. Wir zeigen zuerst, dass $f$ invariant ist. Sei
also $(x,y)\in U\setminus(0,0), v\in V, \sigma \in \SL_{2}$. Dann ist
\[
(\sigma^{-1}f)(x,y,v)=f(\sigma(x,y),\sigma v)=f(ax+by,cx+dy,\sigma v),
\]
und wir m"ussen zeigen, dass dies gleich $f(x,y,v)$ ist. Aufgrund der
Zariski-Dichtheit der beschriebenen Koordinatenmenge in $U\oplus V$ folgt dann
allgemein die Invarianz von $f$. Wir unterscheiden wieder $x\ne 0$ und $y\ne
0$.

\underline{1. Fall: $x\ne 0$.} Dann ist $\left(
\begin{array}{cc}
\frac{1}{x}& 0\\
-y& x\\ 
\end{array}
\right)\in \SL_{2}$, und es folgt
\begin{eqnarray*}
f(x,y,v)&\stackrel{~\eqref{defvonf}}{=}&g\left(\left(
\begin{array}{cc}
\frac{1}{x}& 0\\
-y& x\\ 
\end{array}
\right)v\right)=
g\left(\left(
\begin{array}{cc}
\frac{1}{x}& 0\\
-y& x\\ 
\end{array}
\right)
\left(
\begin{array}{cc}
d& -b\\
-c& a\\ 
\end{array}
\right)
\sigma v\right)\\
&=&g\left(
\left(
\begin{array}{cc}
*& *\\
-(cx+dy)& ax+by\\ 
\end{array}
\right)\sigma v\right)\\
&\stackrel{~\eqref{defvonf}}{=}&
f(ax+by,cx+dy,\sigma v).
\end{eqnarray*}

\underline{2. Fall: $y \ne 0$.} Wir erhalten analog
\begin{eqnarray*}
f(x,y,v)&\stackrel{~\eqref{defvonf}}{=}&g\left(\left(
\begin{array}{cc}
0&\frac{1}{y}\\
-y& x\\ 
\end{array}
\right)v\right)=
g\left(\left(
\begin{array}{cc}
0&\frac{1}{y}\\
-y& x\\ 
\end{array}
\right)
\left(
\begin{array}{cc}
d& -b\\
-c& a\\ 
\end{array}
\right)
\sigma v\right)\\
&=&g\left(
\left(
\begin{array}{cc}
*& *\\
-(cx+dy)& ax+by\\ 
\end{array}
\right)\sigma v\right)\\
&\stackrel{~\eqref{defvonf}}{=}&
f(ax+by,cx+dy,\sigma v).
\end{eqnarray*}
Damit ist $f\in K[U\oplus V]^{\SL_{2}}$ gezeigt. Schliesslich gilt mit der
Definition von $\phi$
\[
\phi(f)(v)=f(1,0,v)\stackrel{~\eqref{defvonf}}{=}g\left(\left(
\begin{array}{cc}
1& 0\\
0& 1\\ 
\end{array}
\right)v
\right)=g(v)
\]
f"ur alle $v\in V$, also $\phi(f)=g$, so dass $f$ das gesuchte Urbild ist.\qed\\

Wir haben Roberts' Isomorphismus hier auf der Ebene von Koordinatenringen $K[V]=S(V^{*})$
angegeben. F"ur uns fast wichtiger ist er auf der Ebene der symmetrischen
Algebra $S(V)$. Wir verwenden die Notation im Satz und wie auf
S. \pageref{mynotation} beschrieben. Es hat $U=\langle X,Y \rangle$ die
Darstellung $\sigma \mapsto \left(
\begin{array}{cc}
a& b\\
c&d\\ 
\end{array}
\right)$. Daher hat $U^{*}=\langle X^{*},Y^{*}\rangle$ die Darstellung $\sigma \mapsto \left(
\begin{array}{cc}
d& -c\\
-b&a\\ 
\end{array}
\right)$, also ist $U^{*}\cong \langle Y,-X \rangle$ - man hat also die
Entsprechung $X^{*}\mapsto Y$ und $Y^{*}\mapsto -X$ und wir wollen hier kurz
sogar identifizieren. Roberts' Isomorphismus
\[
S(\langle X^{*},Y^{*} \rangle \oplus V^{*})^{\SL_{2}}=S(\langle Y,-X \rangle
\oplus V^{*})^{\SL_{2}}\rightarrow S(V^{*})^{\Ga}
\]
l"asst sich dann beschreiben durch 
\[
f(X^{*},Y^{*},T)=f(Y,-X,T)\mapsto f(1,0,T),
\]
wobei
$T=(T_{1},\ldots,T_{n})$ eine Basis von $V^{*}$ ist. Ersetzen wir $V^{*}$
durch $V$, so erhalten wir also als Umformulierung:

\begin{Korollar}[Roberts' Isomorphismus f"ur die symmetrische Algebra] \label{RobertsSymmAlg}
Sei $V$ ein $\SL_{2}$-Modul mit Basis $T=(T_{1},\ldots,T_{n})$, $\langle X,Y
\rangle$ der $\SL_{2}$-Modul mit der nat"urlichen Darstellung. Dann ist durch
\[
\phi: S(\langle X,Y \rangle \oplus V)^{\SL_{2}} \rightarrow S(V)^{\Ga}, \quad f(X,Y,T)
\mapsto f(0,1,T)
\]
ein Isomorphismus gegeben.
\end{Korollar}

Wir wollen uns noch die Umkehrabbildung ansehen. Sei wieder $T$ eine Basis von
$V$. Ist dann $g\in K[V^{*}]^{\Ga}$, so hat man f"ur ein Urbild
$f\in K[U\oplus V^{*}]^{\SL_{2}}$ und $(x,y)\in U, x\ne 0, v\in V^{*}$ gem"a"s
\eqref{defvonf}:
\[
f(x,y,v)=g\left(\left(
\begin{array}{cc}
1/x& 0\\
-y&x\\ 
\end{array}
\right)v\right)=\left(\left(
\begin{array}{cc}
x& 0\\
y&1/x\\ 
\end{array}
\right)\cdot g\right)(v).
\]
Die Operation von $\SL_{2}$ ist durch Polynome gegeben. In diese lassen sich
aber nicht nur K"orperelemente, sondern auch neue Variablen einsetzen. Da wir
die Entsprechungen $x\leftrightarrow X^{*}\leftrightarrow Y$ und
$y\leftrightarrow Y^{*}\leftrightarrow -X$ und $T\leftrightarrow v$ haben,
erhalten wir so aus der letzten Gleichung
\[
f(Y,-X,T)=\left(\left(
\begin{array}{cc}
Y& 0\\
-X&1/Y\\ 
\end{array}
\right)\cdot g\right)(T).
\]
Wir d"urfen hier den "Ubergang zu Variablen machen, da links eine Polynomfunktion in
den Koordinaten $x,y,v$ steht und rechts eine rationale Funktion in $x,y,v$
mit h"ochstens einer Potenz von $x$ im Nenner. Multipliziert man mit diesem Nenner, hat man ein
Gleichung von Polynomfunktionen, die f"ur alle $x\ne 0$ gilt, also auf einer
Zariski-dichten Menge. Dann m"ussen
aber auch die entsprechenden Polynome gleich sein, und insbesondere kann der
Nenner gek"urzt werden, so dass auf beiden Seiten Polynome stehen.
Da $f(Y,-X,T)$ das gesuchte Urbild von $g$ ist, erhalten wir also 
\begin{Korollar}[Umkehrung von Roberts' Isomorphismus] \label{UmkehrungRoberts}
Es sei $V$ ein $\SL_{2}$-Modul und $\langle X,Y
\rangle$ der $\SL_{2}$-Modul mit der nat"urlichen Darstellung. Dann ist durch
\[
\phi^{-1}: S(V)^{\Ga}\rightarrow S(\langle X,Y\rangle \oplus V)^{\SL_{2}},
\quad g\mapsto\left(
\begin{array}{cc}
Y& 0\\
-X&1/Y\\ 
\end{array}
\right)\cdot g
\]
die Umkehrung von Roberts' Isomorphismus gegeben. Das angegebene Urbildelement
$\phi^{-1}(g)$ erh"alt man dabei dadurch, indem man in die durch Polynome
gegebene Operation von $\SL_{2}$ auf $g\in S(V)$ statt K"orperelemente
neue unabh"angige Variablen $X,Y$ wie angegeben in die entsprechenden Polynome
einsetzt (Details im Beweis). 
\end{Korollar}

{\bf Beispiel.} Sei $V=\langle X_{1},Y_{1}\rangle \oplus\langle
X_{2},Y_{2}\rangle$ die zweifache Summe der nat"urlichen Darstellung von
$\SL_{2}$. Wir berechnen die Urbilder der $\Ga$-Invarianten $X_{1}$ und
$X_{1}Y_{2}-X_{2}Y_{1}$:
\[
\phi^{-1}(X_{1})=YX_{1}-XY_{1}, \textrm{ und}
\]
\[
\phi^{-1}(X_{1}Y_{2}-X_{2}Y_{1})=(YX_{1}-XY_{1})\cdot
Y_{2}/Y-(YX_{2}-XY_{2})\cdot Y_{1}/Y=X_{1}Y_{2}-X_{2}Y_{1}.
\]
Man verifiziert sofort, dass tats"achlich
\[
\phi(YX_{1}-XY_{1})=X_{1},\quad \phi(X_{1}Y_{2}-X_{2}Y_{1})=X_{1}Y_{2}-X_{2}Y_{1}
\]
gilt.\\

An diesem Beispiel sehen wir auch, dass Roberts' Isomorphismus nicht homogen
ist, zumindest dann, wenn wir beidemale die Standardgraduierung verwenden. So
ist $\deg YX_{1}-XY_{1}=2$, aber $\deg X_{1}=1$.  Zumindest
bildet $\phi$ aber die maximalen homogenen Ideale aufeinander ab, d.h. es gilt
\begin{equation}\label{RobertsErhaltRplus}
\phi(S(\langle X,Y\rangle \oplus V)_{+}^{\SL_{2}})=S(V)_{+}^{\Ga}.
\end{equation}
\Bew Sei $f\in S(\langle X,Y\rangle \oplus V)_{+}^{\SL_{2}}$. Dann gibt es
eine eindeutige Zerlegung $f=g+h$ mit $g\in S(\langle X,Y\rangle)\cdot
S(V)_{+}$ (die Menge der endliche Summen aller Produkte aus beiden Mengen)
und $h\in S(\langle X,Y\rangle)$. Da die letzten beiden Mengen
$\SL_{2}$-invariant sind, sind mit $f$ auch $g$ und $h$ invariant, also $h\in
S(\langle X,Y\rangle)^{\SL_{2}}=K$. Da $f\in S(\langle X,Y\rangle \oplus
V)_{+}$ und $g$ keinen konstanten Anteil hat, ist also $h=0$ und
$f=g\in S(\langle X,Y\rangle)S(V)_{+}$. Da die Anwendung von $\phi$ einfach
ersetzen von $X=0,\,Y=1$ bedeutet, ist also $\phi(f)\in S(V)_{+}$. Dies zeigt
die Inklusion "`$\subseteq$"' in \eqref{RobertsErhaltRplus}. Da weiter
$\phi(K)=K$, also $\phi(S(\langle X,Y\rangle \oplus
V)_{0}^{\SL_{2}})=S(V)_{0}^{\Ga}$ und $\phi$ bijektiv ist, gilt sogar Gleichheit.\qed\\

Wir k"onnen $\phi$ nun
dadurch homogen machen, dass wir auf $S(V)^{\Ga}$ einfach die Graduierung von
$S(\langle X,Y\rangle \oplus V)^{\SL_{2}}$ via $\phi$ vererben.
Nach \eqref{RobertsErhaltRplus} "andert sich bei diesem Wechsel der Graduierung das maximale homogene Ideal von $S(V)^{\Ga}$
nicht, und damit auch nicht $\depth S(V)^{\Ga}=\depth
(S(V)^{\Ga}_{+},S(V)^{\Ga})$ (vgl. Definition \ref{DefOfDepth}). Insbesondere haben wir f"ur jeden $\SL_{2}$-Modul
$V$, dass
\begin{equation}\label{cmdefOfRobertsIsom}
\cmdef S(V)^{\Ga}=\cmdef S(\langle X,Y\rangle \oplus V)^{\SL_{2}}.
\end{equation}

Man beachte, dass Roberts' Isomorphismus nicht besagt, dass jeder
$\Ga$-Invariantenring $S(V)^{\Ga}$ zu einem $\SL_{2}$-Invariantenring isomorph
ist. Als Voraussetzung hierf"ur muss sich die $\Ga$-Darstellung
auf $V$ zu einer $\SL_{2}$-Darstellung erweitern lassen, so dass man die
$\Ga$-Darstellung gem"a"s \eqref{sigmat}, S. \pageref{sigmat} zur"uckerh"alt.
Eine solche fortsetzbare Darstellung von $\Ga$ hei"st dann
\emph{fundamental}. \index{fundamentale Darstellung} Ein Beispiel f"ur eine nicht fundamentale Darstellung in
positiver Charakteristik $p$ ist
\[
\Ga\rightarrow \GL_{3},\quad t\mapsto \left(
\begin{array}{ccc}
1&t&t^{p}\\
&1&\\
&&1
\end{array}
\right),
\]
siehe Fauntleroy \cite[vor Lemma 2]{Fauntleroy}. In Charakteristik $0$ dagegen
ist jede $\Ga$-Darstellung fundamental. Wir geben hier nur einen elementaren
Beweis f"ur $K=\mathbb{C}$, f"ur den allgemeinen Beweis (f"ur den man die
Theorie der Lie-Algebren ben"otigt) siehe \cite[Lemma 10.2]{Grosshans}
oder \cite[III.3.9]{Kraft}. Dann ist also in Charakteristik $0$ jeder
$\Ga$-Invariantenring isomorph zu einem $\SL_{2}$-Invariantenring
(insbesondere also endlich erzeugt - das ist der Satz von Weitzenb"ock~\cite{Weitzenbock}), und damit nach
Hochster und Roberts (Satz \ref{HoRo}) Cohen-Macaulay. In Charakteristik $0$
kann man also nicht nur f"ur reduktive Gruppen, sondern auch f"ur die
einfachste nicht-reduktive Gruppe $\Ga$ kein Beispiel eines nicht Cohen-Macaulay
Invariantenringes finden.

\begin{Satz}
Sei $K=\mathbb{C}$, und $\langle X,Y \rangle$ die Einschr"ankung der
nat"urlichen Darstellung von $\SL_{2}$ auf $\Ga$, d.h. die Darstellung
$\Ga\rightarrow \GL_{2}, \,t\mapsto \left(
\begin{array}{cc}
1& t\\
0&1\\ 
\end{array}
\right)$. Dann ist jede Darstellung von $\Ga$ isomorph zu einer direkten Summe
von symmetrischen Potenzen $S^{m}(\langle X,Y \rangle)$ (f"ur $m=0$ ist dies
gleich $K$,
die triviale Darstellung) und damit
insbesondere fundamental.
\end{Satz}

\Bew Wir berechnen zun"achst die Darstellung von $S^{m}(\langle X,Y
\rangle)=\langle X^{m},X^{m-1}Y,\ldots,Y^{m} \rangle$. Es ist
\begin{eqnarray}
t\cdot X^{m-j}Y^{j}&=&X^{m-j}(tX+Y)^{j}=\sum_{k=0}^{j}X^{m-j} {j \choose k}  (tX)^{j-k}Y^{k}\nonumber\\
&=&\sum_{k=0}^{j}{j\choose k}t^{j-k}X^{m-k}Y^{k} \label{DarstSm}.
\end{eqnarray}

Sei nun $\rho: \mathbb{C}\rightarrow \GL_{n},\, t\mapsto \rho(t)$ eine
beliebige Darstellung von $\Ga$, d.h. $\rho(z+w)=\rho(z)\rho(w)=\rho(w+z)=\rho(w)\rho(z)\, \myforall
z,w\in{\mathbb C}$ und $\rho(0)=I_{n}$, wobei $I_{n}$ die $n\times n$ Einheitsmatrix ist. Da die Eintr"age von $\rho(z)$ Polynome in $z$ sind, ist
die Einschr"ankung $\rho|_{\mathbb R}$ differenzierbar. Mit Hilfe der
Funktionalgleichung f"ur $\rho$ erhalten wir dann
\begin{eqnarray*}
\frac{d\rho}{dt}(t)&=&\lim_{h\rightarrow
  0}\frac{\rho(t+h)-\rho(t)}{h}=\lim_{h\rightarrow
  0}\frac{\rho(h)\rho(t)-\rho(t)}{h}\\&=&\underbrace{\lim_{h\rightarrow
  0}\frac{\rho(h)-I_{n}}{h}}_{=:A}\rho(t)=A\rho(t) \quad\myforall t\in {\mathbb R}.\end{eqnarray*}
Die L"osung der \emph{reellen} Differenzialgleichung $\dot{\rho}=A\rho,\, A\in {\mathbb
  C}^{n\times n}$ mit Anfangswert $\rho(0)=I_{n}$ ist
bekanntlich eindeutig bestimmt und gegeben durch $\rho(t)=\exp(A\cdot t),\,
  t\in {\mathbb R}$. Da aber jede Komponente von $z\mapsto \rho(z)$ und
  $z\mapsto \exp(Az)$ eine holomorphe Funktion ist, und beide Funktionen auf
  ${\mathbb R}$ "ubereinstimmen, stimmen sie auf ganz ${\mathbb C}$ "uberein,
  also gilt $\rho(z)=\exp(Az) \, \myforall z \in{\mathbb C}$.

Wir zeigen als n"achstes, dass alle Eigenwerte von $\rho(z)$ gleich $1$ sind
f"ur jedes $z\in{\mathbb C}$. Sei also $\lambda \in {\mathbb C}$ ein Eigenwert
von $\rho(z)$ und $v\in{\mathbb C}^{n}\setminus\{0\}$ ein zugeh"origer
Eigenvektor, also $\rho(z)v=\lambda v$. Sukzessives multiplizieren mit
$\rho(z)$ ergibt $\rho(nz)v=\lambda^{n}v\, \myforall n \in {\mathbb N}$. Da die
Eintr"age von $\rho(nz)$ Polynome in $nz$ sind und $v$ eine Komponente
$v_{i}$ ungleich $0$ hat, ergibt sich durch Betrachten dieser Komponente von
$\rho(nz)v=\lambda^{n}v$ und Division durch $v_{i}$ die Existenz eines Polynoms $P(n)$ mit
$P(n)=\lambda^{n}\,\, \myforall n\in{\mathbb N}$. Es folgt
\[
P(n+1)=\lambda^{n+1}=\lambda P(n) \quad \myforall n\in{\mathbb N}.
\]
Da also die beiden Polynomfunktionen $x\mapsto P(x+1)$ und $x\mapsto \lambda
P(x)$ auf ${\mathbb N}$ "ubereinstimmen, sind sie auch als Polynome gleich,
also $P(X+1)=\lambda P(X)$. Da aber $P(X+1)$ denselben Grad und Leitkoeffizient
wie $P(X)$ hat, folgt $\lambda=1$.

Da also alle Eigenwerte von $\rho(z)$ gleich $1$ sind, gilt $(\rho(z)-I_{n})^{n}=0\, \myforall z \in{\mathbb C}$, also auch
$\left(\frac{\rho(h)-I_{n}}{h}\right)^{n}=0\,\myforall h \in{\mathbb R}\setminus\{0\}$. Der
Grenz"ubergang $h\rightarrow 0$ liefert $A^{n}=0$, d.h. A ist nilpotent. Die
Jordanbl"ocke von $A$ sind damit von der Form 
\[
J_{m}=\left(
\begin{array}{cccc}
0&1&&\\
&\ddots&\ddots\\
&&\ddots&1\\
&&&0\\
\end{array}
\right)\in{\mathbb C}^{m+1\times m+1}.
\]
Ist $J=SAS^{-1}$ die Jordan-Normalform von A, so ist
$\exp(Jz)=S\exp(Az)S^{-1}=S\rho(z)S^{-1}$, liefert also eine eine zu $\rho$ isomorphe
Darstellung. Die verschiedenen Jordanbl"ocke von $J$ entsprechen dabei
direkten Summanden der Darstellung. Wir beenden den Beweis, indem wir zeigen dass
$z\mapsto \exp(J_{m}z)$ (bei geeigneter Basiswahl) eine Darstellung von
$S^{m}(\langle X,Y\rangle)$ ist. Bekanntlich ist
\begin{equation} \label{DarstJm}
\exp(J_{m}t)=\sum_{k=0}^{\infty}\frac{(J_{m}t)^{k}}{k!}=\left(
  \begin{array}{cccccc}
1&t&\frac{t^{2}}{2!}&\frac{t^{3}}{3!}&\dots&\frac{t^{m}}{m!}\\
&1&t&\frac{t^{2}}{2!}&\ddots&\vdots\\
&&1&t&\ddots&\frac{t^{3}}{3!}\\
&&&\ddots&\ddots&\frac{t^{2}}{2!}\\
&&&&1&t\\
&&&&&1\\
  \end{array}
\right).
\end{equation}
Aus \eqref{DarstSm},
\[
t\cdot \frac{1}{j!}X^{m-j}Y^{j}=\sum_{k=0}^{j} \frac{t^{j-k}}{(j-k)!}\cdot
\frac{1}{k!}X^{m-k}Y^{k},
\]
ersehen wir, dass die Darstellung von $S^{m}(\langle X,Y\rangle)$ bez"uglich
der Basis $\{ \frac{1}{j!}X^{m-j}Y^{j}: j=0,\ldots,m \}$ genau durch \eqref{DarstJm}
gegeben ist. \qed\\

\subsubsection{Die Dimension von Invariantenringen}\label{sectionDimension}
Wir geben hier einige Zusammenh"ange von $\dim K[V]^{G},\, \dim_{K} V$ und $\dim
G$. Sei ferner $K(V):=\Quot(K[V])$ der \emph{K"orper der rationalen Funktionen auf
$V$}. F"ur einen $G$-Modul $V$ operiert $G$ in kanonischer Weise auch auf
$K(V)$. Dann ist
\[
K(V)^{G}:=\{f\in K(V): \sigma \cdot f=f \,\,\myforall \sigma\in G\}
\]
der \emph{Invariantenk"orper}. \index{Invariantenk\"orper} Da er Zwischenk"orper der endlich erzeugten
K"orpererweiterung $K\le K(V)$ ist, ist er nach \cite[Theorem
4.1.5]{NagataFieldTheory} stets endlich erzeugt "uber $K$.

Offenbar gilt $K(V)^{G}\supseteq \Quot (K[V]^{G})$. Das folgende Lemma (siehe
\cite[Exercise 6.10]{Dolgachev}) gibt
ein Kriterium, wann Gleichheit gilt. 

\begin{Lemma} \label{QuotInvRing}
Sei $G$ eine lineare algebraische Gruppe, so dass jeder algebraische
Gruppenhomomorphismus $G\rightarrow \Gm\cong(K\setminus\{0\},\cdot)$ trivial ist. Dann gilt
\[
K(V)^{G}=\Quot K[V]^{G}
\]
f"ur jeden $G$-Modul $V$.
\end{Lemma}

Man beachte dabei, dass ein algebraischer
Homomorphismus $\chi: G\rightarrow \Gm=\{(a,b)\in K^{2}: ab-1=0\}\cong (K\setminus\{0\},\cdot)$ ein
Gruppenhomomorphismus ist, der durch einen Morphismus von Variet"aten gegeben
ist. Schreibt man $\chi=(\chi',\chi'')$, so ist also $\chi'\in K[G]$ und
$\chi'$ induziert einen
Gruppenhomomorphismus $G\rightarrow (K\setminus\{0\},\cdot)$.

Ist umgekehrt durch $\chi'\in K[G]$ ein Gruppenhomomorphismus
$G\rightarrow(K\setminus\{0\},\cdot)$ gegeben, so ist durch $\chi:
G\rightarrow \Gm$, $\sigma\mapsto(\chi'(\sigma),\chi'(\sigma^{-1}))$ ein
algebraischer Gruppenhomomorphismus gegeben, denn
$\chi'(\sigma)\chi'(\sigma^{-1})=\chi'(\sigma\sigma^{-1})=1$ (also
$\chi(\sigma)\in \Gm$), und da $G\rightarrow G, \sigma\mapsto\sigma^{-1}$ ein
Morphismus ist (und damit $(\sigma\mapsto\chi'(\sigma^{-1}))\in K[G]$), ist 
auch $\chi$ ein Morphismus.

Die algebraischen Gruppenhomomorphismen $G\rightarrow \Gm$ entsprechen also
eindeutig den $\chi'\in K[G]$, die einen Homomorphismus $G\rightarrow
(K\setminus\{0\},\cdot)$ induzieren. (Vgl. Bemerkung \ref{WasIstDarstellung};
es ist $\Gm=\GL_{1}$.)
\\

{\noindent\it Beweis des Lemmas.} Sei $0\ne f\in K(V)^{G}$ und $g,h\in K[V]$ teilerfremd mit
$f=\frac{g}{h}$ sowie $\sigma\in G$. Aus $f=\sigma\cdot f$ folgt $g(\sigma \cdot h)=h(\sigma \cdot
g)$, also
$h|g(\sigma \cdot h)$. Da $h,g$ teilerfremd, folgt also $h|\sigma\cdot h$, und da $\deg
h=\deg \sigma\cdot h$ gibt es ein $\chi(\sigma)\in K\setminus\{0\}$ mit $\sigma\cdot 
h=\chi(\sigma) h$. Da die Operation von $\sigma$ auf $h$ durch Polynome in den
Koordinaten von $\sigma$ gegeben ist, gilt $\chi \in K[G]$. Ferner sieht man
aus $\chi(\sigma\tau)h=(\sigma\tau) \cdot h =\sigma(\tau h)= \sigma \cdot
(\chi(\tau)h)=\chi(\sigma)\chi(\tau)h$, dass
$\chi(\sigma\tau)=\chi(\sigma)\chi(\tau)\,\myforall \sigma,\tau\in G$, also dass $\chi$ ein algebraischer
Homomorphismus $G\rightarrow \Gm$ ist. Nach Voraussetzung ist also $\chi=1$
und damit $h\in K[V]^{G}$. Dann ist auch $g=fh\in K[V]^{G}$. Dies zeigt $f\in
\Quot(K[V]^{G})$ und damit die fehlende Inklusion.\qed\\

Dass die Aussage des Lemmas f"ur endliche Gruppen $G$ gilt (wo es durchaus
nichttriviale Homomorphismen $G\rightarrow \Gm$ geben kann, z.B. $G\subseteq
\Gm$ endliche Untergruppe), sieht man an der Gleichung (Bezeichnungen
wie im Beweis)
\[
f=\frac{g}{h}=\frac{g\cdot \prod_{\sigma\in G\setminus\{\iota\}}\sigma
  h}{\prod_{\sigma\in G}\sigma h}.
\]
Dann ist $\prod_{\sigma\in G}\sigma h\in K[V]^{G}$ und mit $f$ ist auch
$f\cdot\prod_{\sigma\in G}\sigma h=g\cdot \prod_{\sigma\in
  G\setminus\{\iota\}}\sigma h$ invariant, also in $K[V]^{G}$. Damit ist $f\in
\Quot K[V]^{G}$.\\

\begin{Satz} \label{DimAbschatzung}
Sei $G$ eine lineare algebraische Gruppe und $V$ ein $G$-Modul. Dann gilt
\[
\trdeg K(V)^{G}/K \ge \dim V -\dim G.
\]
Ist jeder algebraische Homomorphismus $G\rightarrow \Gm\cong(K\setminus\{0\},\cdot)$ trivial,
und ist $K[V]^{G}$ endlich erzeugt,
so gilt insbesondere
\begin{equation}\label{dimVminusdimG}
\dim K[V]^{G}\ge \dim V-\dim G.
\end{equation}
\end{Satz}

Hierbei ist $\dim G=\dim K[G]$ die Krulldimension der affinen Variet"at $G$
und $\dim V$ die Dimension von $V$ als $K$-Vektorraum (was gleich der
Dimension von $V$ als affine Variet"at ist).\\

\Bew Die erste Aussage des Satzes folgt aus Dolgachev \cite[Corollary
  6.2]{Dolgachev}  oder Popov/Vinberg {\it Invariant Theory} in \cite[Corollary zu Lemma 2.4, p. 156 und
  Formel in Section 1.4, p. 151]{Shafarevich}. Unter der Zusatzvoraussetzung gilt dann nach Lemma
\ref{QuotInvRing}, dass $K(V)^{G}=\Quot K[V]^{G}$, und aus $\dim K[V]^{G}=\trdeg
\Quot(K[V]^{G})/K$ f"ur endlich erzeugtes $K[V]^{G}$ folgt dann die Behauptung. \qed\\

\begin{Satz}\label{DimOfGaInvariants}
Es gibt keinen nichttrivialen algebraischen Homomorphismus $\Ga\rightarrow
\Gm$. Insbesondere gilt f"ur jeden $\Ga$-Modul $V$ stets $K(V)^{\Ga}=\Quot
K[V]^{\Ga}$. Ist $K[V]^{\Ga}$ endlich erzeugt, so gilt au"serdem $\dim
K[V]^{\Ga}\ge \dim V-1$.
\end{Satz}

\Bew Sei 
\[
\lambda: \Ga\rightarrow \Gm, \quad a\mapsto \sum_{k=0}^{n}\lambda_{k}a^{k}
\quad \textrm{ mit } \lambda_{k}\in K \textrm{ und } \lambda_{n}\ne 0
\]
ein algebraischer Homomorphismus. Es folgt
\[
1=\lambda(0)=\lambda(a)\lambda(-a)=\left(\sum_{k=0}^{n}\lambda_{k}a^{k}\right)\left(\sum_{k=0}^{n}\lambda_{k}(-a)^{k}
\right) \quad \myforall a\in K. 
\]
Koeffizientenvergleich liefert sofort $n=0$, also
$\lambda(a)=\lambda_{0}=\lambda(0)=1$ f"ur alle $a\in \Ga$.
 
Die restlichen beiden Aussagen folgen sofort aus Lemma \ref{QuotInvRing} und
Satz \ref{DimAbschatzung} mit $\dim \Ga=\dim K^{1}=1$.
\qed\\

Dagegen ist $\exp: ({\mathbb C},+)\rightarrow({\mathbb C}\setminus\{0\},\cdot)$ zwar ein
Gruppenhomomorphismus, aber eben nicht algebraisch.

\begin{Satz}\label{DimOfSL2Invariants}
Es gibt keinen nichttrivialen Gruppenhomomorphismus $\SL_{n}\rightarrow
\Gm$. Insbesondere gilt f"ur jeden $\SL_{n}$-Modul $V$ stets $K(V)^{\SL_{n}}=\Quot
K[V]^{\SL_{n}}$ und \[\dim
K[V]^{\SL_{n}}\ge \dim V-n^{2}+1.\]
\end{Satz}

\Bew Sei $\phi: \SL_{n}\rightarrow \Gm$ ein Gruppenhomomorphismus. Dann ist
$\phi(\SL_{n})\subseteq \Gm$ insbesondere abelsch, und damit liegt die
\emph{Kommutatorgruppe} $\SL_{n}'$ im Kern von $\phi$, also $\SL_{n}'\subseteq \ker
\phi$. Nach Hein \cite[Satz I.6.9]{Hein} gilt aber $\SL_{n}'=\SL_{n}$
(d.h. $\SL_{n}$ ist \emph{perfekt}),\index{Gruppe!perfekt} also
ist $\ker \phi=\SL_{n}$ und damit $\phi$ trivial.

Die restlichen beiden Aussagen folgen nun wieder aus Satz
\ref{DimAbschatzung} mit $\dim \SL_{n}=n^{2}-1$. Da $\SL_{n}$ reduktiv ist, ist
$K[V]^{\SL_{n}}$ stets endlich erzeugt.\qed\\

Dass es keinen \emph{algebraischen} Homomorphismus $\SL_{n}\rightarrow \Gm$
gibt, kann man auch leicht beweisen, ohne dass man die Perfektheit von
$\SL_{n}$ benutzen muss: Sei also $\phi: \SL_{n}\rightarrow \Gm$ ein algebraischer Homomorphismus. Sei
$I_{n}\in K^{n\times n}$ die $n\times n$ Einheitsmatrix und
$E_{ij}=(\delta_{ki}\delta_{lj})_{kl}\in K^{n\times n}$ die Matrix, die genau
in der $i$-ten Zeile und $j$-ten Spalte eine $1$ und sonst nur Nullen
hat. F"ur $i\ne j$ sei $N_{ij}(a):=I_{n}+aE_{ij}$. Bekanntlich wird $\SL_{n}$
erzeugt von der Menge $\{N_{ij}(a): i\ne j, 1\le i,j \le n, a\in K\}$
(Gauss-Algorithmus, siehe auch Hein \cite[Satz I.2.8]{Hein}). Da $N_{ij}(a)N_{ij}(b)=N_{ij}(a+b)$ f"ur $a,b\in K$, ist
die Abbildung \[\phi\circ N_{ij}: \Ga\rightarrow \Gm,\quad a\mapsto
\phi\left(N_{ij}\left(a\right)\right)\] ein algebraischer Homomorphismus, also
nach Satz \ref{DimOfGaInvariants} trivial. Damit gilt
$\phi\left(N_{ij}\left(a\right)\right)=1$ f"ur alle $i\ne j$ und $a\in K$. Da
also $\phi$ auf einem Erzeugendensystem von $\SL_{n}$ konstant gleich $1$ ist,
ist $\phi$ trivial.\\

Dagegen ist etwa $\det: \GL_{n}\rightarrow \Gm$ ein nichttrivialer
algebraischer Homomorphismus. Seien $\langle X_{i},Y_{i} \rangle$ jeweils die
nat"urliche Darstellung der $\GL_{2}$ und $V^{*}:=\bigoplus_{i=1}^{n}\langle
X_{i},Y_{i} \rangle$. Dann ist bekanntlich $K[V]^{\GL_{2}}=K$
(siehe etwa \cite{ConciniProcesi}). Dagegen ist
$\frac{X_{1}Y_{2}-X_{2}Y_{1}}{X_{1}Y_{3}-X_{3}Y_{1}}\in K(V)^{\GL_{2}}$. Es
ist also $K(V)^{\GL_{2}}\ne \Quot K[V]^{\GL_{2}}$. Au"serdem ist $\dim
K[V]^{\GL_{2}}=0 \not\ge \dim V-\dim \GL_{2}=2n-4$ f"ur $n>2$. Gleichung
\eqref{dimVminusdimG} gilt also nicht allgemein.

\begin{Satz}\label{unendlicheDarsteDim}
Sei $G$ eine lineare algebraische Gruppe und $V$ ein $G$-Modul, so dass die
Darstellung $G\rightarrow \GL(V)$ unendliches Bild hat. Ist dann $K[V]^{G}$
endlich erzeugt, so gilt
\[
\dim K[V]^{G}<\dim V.
\]
\end{Satz}

\Bew Sei $n=\dim V$ und $\sigma \mapsto A_{\sigma}\in K^{n\times n}$ die
Darstellung von $G$ auf $V$ bzgl. einer Basis. Nach Voraussetzung ist die Menge $\{A_{\sigma}:\sigma\in
G\}$ unendlich. Dann gibt es auch eine Zeile $i$, so dass $\{e_{i}^{T}A_{\sigma}:\sigma\in
G\}$ (mit $e_{i}\in K^{n}$ $i$-ter Spalteneinheitsvektor) unendlich
ist. Ist $V^{*}=\langle X_{1},\ldots,X_{n}\rangle$ mit Darstellung $\sigma \mapsto
A_{\sigma^{-1}}^{T}$, so ist $A_{\sigma^{-1}}^{T}e_{i}$ der Koordinatenvektor
von $\sigma\cdot X_{i}$.
Also ist $\{\sigma\cdot X_{i}:\sigma\in
G\}$ unendlich. Es ist $K[V]=S(V^{*})=K[X_{1},\ldots,X_{n}]$. Jedenfalls ist
$\dim K[V]^{G}=\trdeg\Quot(K[V]^{G})/K\le\trdeg K(V)/K=n$.
Angenommen, es
w"are $\dim K[V]^{G}=n$. Dann w"are $\trdeg
\Quot(K[V]^{G})/K=n=\trdeg K(V)/K$, also wegen der Additivit"at des
Transzendenzgrades $\trdeg K(V)/\Quot(K[V]^{G})=0$. Insbesondere w"are $X_{i}\in
K(V)$ algebraisch "uber $\Quot(K[V]^{G})$. Daher g"abe es ein $0\ne f\in
\Quot(K[V]^{G})[T]$ mit $f(X_{i})=0$. Da die Koeffizienten von $f$ invariant
unter $G$ sind, folgt dann
\[
0=\sigma\cdot 0 =\sigma \cdot f(X_{i})=f(\sigma\cdot X_{i})\quad \myforall
\sigma\in G.
\]
Also h"atte $f$ die unendlich vielen Nullstellen $\sigma\cdot X_{i},\sigma \in
G$, ein
Widerspruch. \qed\\

Damit erhalten wir die folgende Charakterisierung endlicher Gruppen:

\begin{Korollar}
Sei $G$ eine reduktive Gruppe. $G$ ist genau dann endlich, wenn f"ur jeden
$G$-Modul $V$ $\dim K[V]^{G}=\dim V$ gilt.
\end{Korollar}

\Bew Ist $G$ endlich, so ist $K[V]/K[V]^{G}$ sogar ganz, insbesondere also
$\dim K[V]^{G}=\dim K[V]=\dim V$. Ist $G$ unendlich, so ist f"ur eine nach
Satz \ref{treueDarstellung} existierende treue Darstellung $V$ das Bild von 
$G\rightarrow \GL(V)$ mit $G$ unendlich, und nach obigem Satz gilt dann $\dim
K[V]^{G}<\dim V$. \qed
\begin{Bemerkung}
Jede echte abgeschlossene Teilmenge (insbesondere jede echte abgeschlossene
 Untergruppe) von $\Ga=K^{1}$ ist endlich. Insbesondere hat
 eine nichttriviale rationale Darstellung $\rho: \Ga\rightarrow \GL(V)$ stets
 endlichen Kern und unendliches Bild $\rho(\Ga)$. 
\end{Bemerkung}

\Bew Da $K[\Ga]\cong K[X]$ ein Hauptidealring ist, ist jede echte
abgeschlossene Teilmenge von $\Ga$ Nullstellenmenge eines nichtkonstanten
Polynoms und somit endlich. 
Ist $\rho$ eine nichttriviale rationale Darstellung von $\Ga$, so ist $\ker
\rho=\rho^{-1}(\id_{V})$ also als echte abgeschlossene Untergruppe von $\Ga$ endlich. Dann ist
$\rho(\Ga)\cong \Ga/\ker \rho$ nat"urlich unendlich.\qed\\

Ein Beispiel f"ur eine nichttriviale und nichttreue Darstellung der $\Ga$ ist
zum Beispiel f"ur $\chr K=p,\,\, q=p^{n}>1$ gegeben durch $\rho: \Ga\rightarrow
\GL_{2}$, $a\mapsto \left(\begin{array}{cc}1&a^{q}-a\\&1\end{array}\right)$
mit $\ker \rho={\mathbb F}_{q}$.

\begin{Korollar}\label{exakteDimOfGaInvs}
Sei $V$ ein nichttrivialer $\Ga$-Modul und $K[V]^{\Ga}$ endlich erzeugt ist. Dann gilt
\[
\dim K[V]^{\Ga}=\dim V-1.
\]
\end{Korollar}

\Bew Dies folgt mit der Bemerkung sofort aus den S"atzen \ref{DimOfGaInvariants} und
\ref{unendlicheDarsteDim}. \qed\\

\begin{Korollar}\label{exakteDimOfSL2Invs}
Sei $V'$ ein nichttrivialer $\SL_{2}$-Modul, $\langle X,Y\rangle$ die nat"urliche Darstellung
der $\SL_{2}$ und $V:=V'\oplus\langle X,Y\rangle$. Dann gilt
\[
\dim K[V]^{\SL_{2}}=\dim V-3.
\]
\end{Korollar}

\Bew Nach Roberts' Isomorphismus \ref{RobertsIsom} ist $K[V]^{\SL_{2}}\cong
K[V']^{\Ga}$. Wir zeigen, dass $V'$ auch nichttrivialer $\Ga$-Modul ist, wobei
hier $\Ga=\left\{\left(\begin{array}{cc}1&a\\&1\end{array}\right): a\in K
\right\}\subseteq \SL_{2}$. Sei $\rho: \SL_{2}\rightarrow \GL(V')$ die Darstellung von $V'$. Nach Voraussetzung
ist $\ker \rho\lhd \SL_{2}$ ein echter Normalteiler. Nach
 Hein \cite[Satz 1.2.12, 1.2.10]{Hein} folgt hieraus $\ker \rho\subseteq \{I_{2},-I_{2}\}$. Insbesondere ist $\Ga\cap \ker \rho =\{I_{2}\}$ und damit $V'$ sogar ein
treuer $\Ga$-Modul. Also gilt nach dem letzten Korollar
\[
\dim K[V]^{\SL_{2}}=\dim K[V']^{\Ga}=\dim V'-1=\dim V-3.
\]
\qed\\
An dem Beweis sehen wir auch, dass jede nichttriviale, fundamentale (also auf $\SL_{2}$ fortsetzbare) Darstellung der $\Ga$ treu ist; Insbesondere ist
die oben angegebene nicht treue Darstellung der $\Ga$ nicht fundamental.
\begin{BemRoman}\label{exakteDimOfSL2FrobInvs}
Das Korollar gilt auch, wenn $V=\langle X^{p},Y^{p}\rangle \oplus V'$ ist. Wir
werden n"amlich in Satz \ref{IsomOfFrobInvs} sehen, dass
\[
S(\langle X^{p},Y^{p}\rangle \oplus V')^{\SL_{2}}\cong S(\langle X,Y\rangle
\oplus V')^{\SL_{2}}\cap K[X^{p},Y^{p}]\otimes S(V')\]
gilt. Au"serdem ist 
$R_{1}:=S(\langle X,Y\rangle
\oplus V')^{\SL_{2}}$ ganz "uber $R_{2}:=S(\langle X,Y\rangle
\oplus V')^{\SL_{2}}\cap K[X^{p},Y^{p}]\otimes S(V')$, denn f"ur $f\in R_{1}$
ist $f^{p}\in R_{2}$. Dann ist also $\dim R_{2}=\dim R_{1}=\dim V'-1$, und
damit auch $\dim S(V)^{\SL_{2}}=\dim R_{2}=\dim V-3$.
\end{BemRoman}

\newpage
\section{Kohomologie von Gruppen}\label{Kohomologievongruppen}
Dieser Abschnitt ist im wesentlichen als Exkurs zu verstehen, da die Resultate
im weiteren nicht verwendet werden. Ziel ist es, eine Verallgemeinerung von
Proposition \ref{AnnulatorProp} f"ur h"ohere Kohomologie zu geben. Das Ergebnis
k"onnte zwar im Zusammenhang etwa mit \cite[Theorem 1.4 oder Corollary
1.6]{KemperOnCM} von Bedeutung sein, aber da sich die gemeinsamen
Voraussetzungen schwer unter einen Hut bringen lassen, ist das Ergebnis mehr
von "`akademischem Interesse"'. Wer an dieser
Verallgemeinerung nicht interessiert ist, kann diesen Abschnitt daher
(evtl. nach Einf"uhrung der Grundbegriffe in Abschnitt
\ref{GrundbegriffeKohom}) "uberspringen. Im Gegensatz zu der in der Literatur
"ublichen abstrakten Einf"uhrung der $n$-ten Kohomologiegruppen als $\Ext_{KG}^{n}(K,V)$
verwenden wir den expliziten und elementaren Zugang "uber $n$-Kozyklen.

Eine weitere Besonderheit dieses Abschnitts ist eine andere Bedeutung bereits
verwendeter Begriffe. Zun"achst ist $G$ hier eine \emph{beliebige} (also nicht
notwendig lineare algebraische) Gruppe. Tr"agt $G$ zus"atzlich eine
algebraische Struktur, wird diese ignoriert. Einen \emph{$G$-Modul} (oder auch
$KG$-Modul) nennen wir
in diesem Abschnitt einen (nicht notwendig endlich-dimensionalen)
$K$-Vektorraum $V$, auf dem $G$ linear operiert (auch f"ur eine lineare
algebraische Gruppe muss die Operation nicht durch einen Morphismus gegeben sein).

\subsection{Koketten, Kozyklen und Kor"ander} \label{GrundbegriffeKohom}\index{Kokette}\index{Kozyklus}\index{Korand}
Sei $G$ eine Gruppe und $V$ ein $G$-Modul. Es sei
\[
C^{n}(G,V):=\{g: G^{n}\rightarrow V\}, \quad n\ge 1
\]
die Menge aller Abbildungen von $G^{n}$ nach $V$
und
\[
C^{0}(G,V):=V.
\]
Ist $v\in C^{0}(G,V)$ und $f: V\rightarrow W$ eine Abbildung in einen
$G$-Modul $W$, so schreiben wir
auch $f\circ v$ f"ur $f(v)$. Dann ist $f\circ v\in C^{0}(G,W)$.
Ein Element $g\in C^{n}(G,V)$ mit $n\ge 0$ hei"st \emph{$n$-Kokette (mit
  Koeffizienten in $V$)}. Offenbar
ist $C^{n}(G,V)$ in kanonischer Weise $KG$-Modul (f"ur $g\in C^{n}(G,V)$,
$\sigma\in G$ ist $\sigma g$ definiert durch $(\sigma g)(x):=\sigma(g(x))$
f"ur alle $x\in G^{n}$), inbesondere also auch
$K$-Vektorraum und additive Gruppe. Wir betrachten f"ur $n\ge 0$ die
$K$-lineare Abbildung\label{DefVonPartnV}
\[
\partial_{n}^{V}: C^{n}(G,V)\rightarrow C^{n+1}(G,V),\quad g\mapsto \partial_{n}^{V}g
\]
mit
\begin{eqnarray}
\partial_{n}^{V}g(\sigma_{1},\ldots,\sigma_{n+1})&:=&\sigma_{1}g(\sigma_{2},\ldots,\sigma_{n+1})+\sum_{i=1}^{n}(-1)^{i}g(\sigma_{1},\ldots,\sigma_{i-1},\sigma_{i}\sigma_{i+1},\sigma_{i+2},\ldots,\sigma_{n+1})\nonumber\\&&+(-1)^{n+1}g(\sigma_{1},\ldots,\sigma_{n})\quad
\textrm{f"ur }(\sigma_{1},\ldots,\sigma_{n+1})\in G^{n+1}.\label{defpartdn}
\end{eqnarray}
Falls keine
Verwechslungsgefahr besteht, schreiben wir auch $\partial_{n}$ statt $\partial_{n}^{V}$.
Wir schreiben die Formel
f"ur $n=0,1,2$ und jeweils $g\in C^{n}(G,V)$ explizit auf:
\begin{eqnarray*}
  \partial_{0}g(\sigma)&=&\sigma g -g \quad \myforall \sigma\in G\quad
  (\textrm{hier ist } g\in C^{0}(G,V)=V)\\
  \partial_{1}g(\sigma,\tau)&=&\sigma
  g(\tau)-g(\sigma\tau)+g(\sigma)\quad \myforall \sigma,\tau\in G\\
  \partial_{2}g(\sigma_{1},\sigma_{2},\sigma_{3})&=&\sigma_{1}g(\sigma_{2},\sigma_{3})-g(\sigma_{1}\sigma_{2},\sigma_{3})+g(\sigma_{1},\sigma_{2}\sigma_{3})-g(\sigma_{1},\sigma_{2})
  \quad \forall\sigma_{1},\sigma_{2},\sigma_{3}\in G.
\end{eqnarray*}
Wir definieren
\[
Z^{n}(G,V):=\ker \partial_{n}^{V}\subseteq C^{n}(G,V), \quad n\ge 0,
\]
die additive Gruppe der \emph{$n$-Kozyklen (mit Koeffizienten in $V$)}. Im
n"achsten Unterabschnitt zeigen wir
$\partial_{n+1}^{V}\circ \partial_{n}^{V}=0 \,\,\myforall n\ge 0$, daher ist die
Gruppe der \emph{$n$-Kor"ander (mit Koeffizienten in $V$)},
\[
B^{n}(G,V):=\im \partial_{n-1}^{V}\subseteq Z^{n}(G,V), \quad n\ge 1,
\]
eine Untergruppe von $Z^{n}(G,V)$. Die Faktorgruppe
\[
H^{n}(G,V):=Z^{n}(G,V)/B^{n}(G,V), \quad n\ge 1
\]
hei"st \emph{$n$-te Kohomologiegruppe von $G$ (mit Koeffizienten in V)}.\index{Kohomologiegruppe} Wir
setzen auch $H^{0}(G,V):=Z^{0}(G,V)=V^{G}$ und $B^{0}(G,V):=0$. Ein
$n$-Kozyklus $g\in Z^{n}(G,V)$ hei"st \emph{trivial}, falls $g\in B^{n}(G,V)$,
sonst \emph{nichttrivial}. \index{Kozyklus!trivialer}

Ein $0$-Kozyklus ist also eine Invariante $g\in V^{G}$. Ein $1$-Korand ist
eine Abbildung $g: G\rightarrow V$, f"ur die es ein $v\in V=C^{0}(G,V)$ gibt
mir $g(\sigma)=\sigma v -v \,\,\myforall \sigma\in G$. Ein $1$-Kozyklus ist eine
Abbildung $g: G\rightarrow V$, die die Funktionalgleichung
$g(\sigma\tau)=\sigma g(\tau)+g(\sigma)\,\myforall\sigma,\tau\in G$ erf"ullt. In
Abschnitt \ref{firstCohom}, wo wir die erste Kohomologie \emph{algebraischer}
Gruppen definiert haben, haben wir zus"atzlich gefordert, dass $g$ ein
Morphismus ist. Um die $n$-te Kohomologie \emph{algebraischer} Gruppen zu
definieren, w"urde man f"ur eine lineare algebraische Gruppe $G$ und eine
rationale Darstellung $V$ "uberall $C^{n}(G,V)$ durch
\[
C_{mor}^{n}(G,V):=\{f: G^{n}\rightarrow V,\quad f \textrm{ Morphismus}\}
\]
ersetzen (f"ur endliche Gruppen ist dies gleich $C^{n}(G,V)$). In diesem Exkurs werden wir dies jedoch nicht weiter verfolgen.

\subsection{Der Kokomplex}
Wir m"ussen noch zeigen, dass
\[
\ldots\stackrel{\partial_{n-1}^{V}}{\rightarrow}C^{n}(G,V)\stackrel{\partial_{n}^{V}}{\rightarrow}C^{n+1}(G,V)\stackrel{\partial_{n+1}^{V}}{\rightarrow}C^{n+2}(G,V)\stackrel{\partial_{n+2}^{V}}{\rightarrow}\ldots
\]
ein \emph{Kokomplex}\index{Kokomplex} ist, d.h. dass 
\[
\partial_{n+1}^{V}\circ\partial_{n}^{V}=0 \in \Hom_{K}(C^{n}(G,V),C^{n+2}(G,V))
\]
 ist. (Bei
einem Kokomplex ist die Indizierung nach rechts aufsteigend, bei einem
Komplex\index{Komplex} 
absteigend). Dies l"asst sich direkt durch Anwenden der Definitionsgleichung
zeigen, doch ist die entstehende Rechung etwas un"ubersichtlich. Stattdessen
gehen wir "ahnlich wie in Benson \cite[section 3.4]{BensonCohomology1} in
mehreren Schritten vor.\\

\Bew (i) Wir betrachten f"ur $n\ge 1$ folgenden Untervektorraum von $C^{n}(G,V)$,
\[
L^{n}(G,V):=\{h\in C^{n}(G,V):
h(\sigma\sigma_{1},\ldots,\sigma\sigma_{n})=\sigma h(\sigma_{1},\ldots,\sigma_{n})
\quad\textrm{f"ur } \sigma,\sigma_{1},\ldots,\sigma_{n}\in G\}.
\]
Dann ist durch
\[
\Psi_{n}: C^{n}(G,V)\rightarrow L^{n+1}(G,V), \quad g\mapsto\Psi_{n}g \quad
(n\ge 0)
\]
mit
\[
\Psi_{n}g(\sigma_{0},\ldots,\sigma_{n}):=\sigma_{0}g(\sigma_{0}^{-1}\sigma_{1},\sigma_{1}^{-1}\sigma_{2},\ldots,\sigma_{n-1}^{-1}\sigma_{n})
\quad\myforall \sigma_{0},\ldots,\sigma_{n}\in G
\]
eine $K$-lineare Abbildung gegeben; F"ur die Wohldefiniertheit "uberpr"ufen wir
\begin{eqnarray*}
  \Psi_{n}g(\sigma\sigma_{0},\ldots,\sigma\sigma_{n})&=&\sigma\sigma_{0}g(\sigma_{0}^{-1}\sigma_{1},\sigma_{1}^{-1}\sigma_{2},\ldots,\sigma_{n-1}^{-1}\sigma_{n})\\
&=&\sigma \Psi_{n}g(\sigma_{0},\ldots,\sigma_{n}) \quad \myforall \sigma,\sigma_{0},\ldots,\sigma_{n}\in G,
\end{eqnarray*}
also tats"achlich $\Psi_{n}g\in L^{n+1}(G,V)$. Weiter betrachten wir die
lineare Abbildung
\[
\Phi_{n}: L^{n+1}(G,V)\rightarrow C^{n}(G,V), \quad h\mapsto \Phi_{n}h \quad
(n\ge 0)
\]
mit
\[
\Phi_{n}h(\sigma_{1},\ldots,\sigma_{n}):=h(1,\sigma_{1},\sigma_{1}\sigma_{2},\ldots,\sigma_{1}\sigma_{2}\cdot\ldots\cdot\sigma_{n})
\quad \myforall \sigma_{1},\ldots,\sigma_{n}\in G.
\]
Wir bezeichnen hier mit $1\in G$ das neutrale Element von $G$.

(ii) Wir zeigen
\begin{equation}\label{isomCLn}
\Psi_{n}\circ \Phi_{n}=\id_{L^{n+1}(G,V)} \quad \textrm{und}\quad
\Phi_{n}\circ \Psi_{n}=\id_{C^{n}(G,V)} \quad (n\ge 0).
\end{equation}
(Die zweite Gleichung ben"otigen wir nicht und dient nur der Vollst"andigkeit.)\\
F"ur $\sigma_{0},\ldots,\sigma_{n}\in G$ und $h\in L^{n+1}(G,V), g\in
C^{n}(G,V)$ gilt n"amlich
\begin{eqnarray*}
\Psi_{n}\Phi_{n}h(\sigma_{0},\ldots,\sigma_{n})&=&\sigma_{0}\Phi_{n}h(\sigma_{0}^{-1}\sigma_{1},\sigma_{1}^{-1}\sigma_{2},\ldots,\sigma_{n-1}^{-1}\sigma_{n})\\
&=&\sigma_{0}h(1,\sigma_{0}^{-1}\sigma_{1},\sigma_{0}^{-1}\sigma_{2},\ldots,\sigma_{0}^{-1}\sigma_{n})\\
&=&h(\sigma_{0},\ldots,\sigma_{n})
\end{eqnarray*}
(im letzten Schritt haben wir $h\in L^{n+1}(G,V)$ verwendet), also
$\Psi_{n}\Phi_{n}h=h$, und analog
\begin{eqnarray*}
  \Phi_{n}\Psi_{n}g(\sigma_{1},\ldots,\sigma_{n})&=&\Psi_{n}g(1,\sigma_{1},\sigma_{1}\sigma_{2},\ldots,\sigma_{1}\cdot\ldots\cdot\sigma_{n})\\
&=&g(\sigma_{1},\ldots,\sigma_{n}),
\end{eqnarray*}
also $\Phi_{n}\Psi_{n}g=g$.

(iii) Wir betrachten nun die lineare Abbildung
\[
\delta_{n}: L^{n}(G,V)\rightarrow L^{n+1}(G,V),\quad h\mapsto \delta_{n}h
\quad (n\ge 1)
\]
mit
\[
\delta_{n}h(\sigma_{0},\ldots,\sigma_{n}):=\sum_{i=0}^{n}(-1)^{i}h(\sigma_{0}\ldots,\hat{\sigma_{i}},\ldots,\sigma_{n})
\quad \myforall \sigma_{0},\ldots,\sigma_{n}\in G,
\]
wobei das Dach $\hat{}$ bedeutet, dass der Eintrag gestrichen wird. Man pr"uft
leicht nach, dass $\delta_{n}h\in L^{n+1}(G,V)$ gilt. Als n"achstes zeigen wir
\begin{equation}
  \label{deltandeltanplus1gleich0}
\delta_{n+1}\circ\delta_{n}=0\in \Hom_{K}(L^{n}(G,V),L^{n+2}(G,V)) \quad (n\ge
1),  
\end{equation}
dass also ein Kokomplex vorliegt: F"ur $\sigma_{0},\ldots,\sigma_{n+1}\in G$ gilt
\begin{eqnarray*}
\delta_{n+1}\delta_{n}h(\sigma_{0},\ldots,\sigma_{n+1})&=&\sum_{i=0}^{n+1}(-1)^{i}\delta_{n}h
(\sigma_{0}\ldots,\hat{\sigma_{i}},\ldots,\sigma_{n+1})\\
&=&\sum_{i=0}^{n+1}(-1)^{i}\Big(\sum_{j=0}^{i-1} (-1)^{j}h
(\sigma_{0},\ldots,\hat{\sigma_{j}}\ldots,\hat{\sigma_{i}},\ldots,\sigma_{n+1})\\
&&+  \sum_{j=i+1}^{n+1}(-1)^{j+1}h
(\sigma_{0},\ldots,\hat{\sigma_{i}}\ldots,\hat{\sigma_{j}},\ldots,\sigma_{n+1})\Big)=0.
\end{eqnarray*}

(iv) Wir zeigen nun
\begin{equation}
  \label{deltandn}
\partial_{n}^{V}=\Phi_{n+1}\circ\delta_{n+1}\circ\Psi_{n} \quad (n\ge 0).  
\end{equation}
 F"ur $g\in
C^{n}(G,V)$ und $\sigma_{1},\ldots,\sigma_{n+1}\in G$ gilt
\begin{eqnarray*}
  \Phi_{n+1}\delta_{n+1}\Psi_{n}g(\sigma_{1},\ldots,\sigma_{n+1})&=&\delta_{n+1}\Psi_{n}g(1,\sigma_{1},\ldots,\sigma_{1}\sigma_{2}\cdot\ldots\cdot\sigma_{n+1})\\
&=&\Psi_{n}g(\sigma_{1},\sigma_{1}\sigma_{2},\ldots,\sigma_{1}\sigma_{2}\cdot\ldots\cdot\sigma_{n+1})\\
&&+\sum_{i=1}^{n+1}(-1)^{i}\Psi_{n}g(1,\sigma_{1},\ldots,\widehat{\sigma_{1}\cdot\ldots\cdot\sigma_{i}},\ldots,\sigma_{1}\sigma_{2}\cdot\ldots\cdot\sigma_{n+1})\\&=&\sigma_{1}g(\sigma_{2},\ldots,\sigma_{n+1})+\sum_{i=1}^{n}(-1)^{i}g(\sigma_{1},\ldots,\sigma_{i}\sigma_{i+1},\ldots,\sigma_{n+1})\\
&&+(-1)^{n+1}g(\sigma_{1},\sigma_{2},\ldots,\sigma_{n})\\
&=&\partial_{n}^{V}g(\sigma_{1},\sigma_{2},\ldots,\sigma_{n+1}).
\end{eqnarray*}

(v) Damit folgt dann f"ur $n\ge 0$
\begin{eqnarray*}
\partial_{n+1}^{V}\circ\partial_{n}^{V}&\stackrel{~\eqref{deltandn}}{=}&\Phi_{n+2}\circ\delta_{n+2}\circ\underbrace{\Psi_{n+1}\circ\Phi_{n+1}}_{=\id_{L^{n+2}(G,V)}~\eqref{isomCLn}}\circ\delta_{n+1}\circ\Psi_{n}\\
&=&\Phi_{n+2}\circ\underbrace{\delta_{n+2}\circ\delta_{n+1}}_{=0~\eqref{deltandeltanplus1gleich0}}\circ\Psi_{n}\\&=&0\,\,\,\in \Hom_{K}(C^{n}(G,V),C^{n+2}(G,V)).  
\end{eqnarray*}
Dies wollten wir zeigen. \qed\\

\subsection{Die bar resolution}\label{sectiondiebarresolutin}
In diesem Abschnitt geben wir eine freie Aufl"osung von $K$ als $KG$-Modul an,
die sogenannte \emph{bar resolution}.\index{bar resolution} Wie "ublich operiert hier $G$ trivial auf $K$.
Sei
\[
P_{n}:=\left\{f: G^{n}\rightarrow KG: f(x)\ne 0 \textrm{ nur f"ur endlich
  viele } x\in
  G^{n}\right\} \quad (n\ge 1)
\]
der freie $KG$-Modul mit Basis $\{[\sigma_{1},\ldots,\sigma_{n}]:
\sigma_{1},\ldots,\sigma_{n}\in G\}$, wobei
\[
[\sigma_{1},\ldots,\sigma_{n}]: G^{n}\rightarrow KG,\quad
(\tau_{1},\ldots,\tau_{n})\mapsto \delta_{\sigma_{1},\tau_{1}}\cdot\ldots\cdot\delta_{\sigma_{n},\tau_{n}}.
\]
(Hier ist $\delta$ das Kronecker-Symbol.)
Weiter sei \[P_{0}:=KG \textrm{  mit einzigem Basiselement } []=1\in KG.\]

Wir betrachten die Folge von Abbildungen $d_{n}\in \Hom_{KG}(P_{n},P_{n-1})$
($n\ge 1$)
gegeben durch $KG$-lineare Fortsetzung von
\begin{eqnarray}
d_{n}([\sigma_{1},\ldots,\sigma_{n}])&:=&\sigma_{1}[\sigma_{2},\ldots,\sigma_{n}]+\sum_{i=1}^{n-1}(-1)^{i}[\sigma_{1},\ldots,\sigma_{i}\sigma_{i+1},\ldots,\sigma_{n}]\nonumber\\&&+(-1)^{n}[\sigma_{1},\ldots,\sigma_{n-1}], \label{defdn}
\end{eqnarray}
sowie $d_{0}\in\Hom_{KG}(P_{0},K)$ gegeben durch $KG$-lineare Fortsetzung von
\begin{equation}\label{defd0}
d_{0}([]):=1.
\end{equation}
 Ziel dieses Abschnitts ist es zu zeigen, dass
\begin{equation}\label{thebarresolution}
\ldots\stackrel{d_{n+1}}{\rightarrow}
P_{n}\stackrel{d_{n}}{\rightarrow}P_{n-1}\stackrel{d_{n-1}}{\rightarrow}P_{n-2}\stackrel{d_{n-2}}{\rightarrow} \ldots\stackrel{d_{2}}{\rightarrow}P_{1}\stackrel{d_{1}}{\rightarrow}P_{0}\stackrel{d_{0}}{\rightarrow}K\rightarrow 0
\end{equation}
eine exakte Sequenz ist (d.h. es gilt $\ker d_{n}=\im d_{n+1} \,\myforall n\ge 0$), also eine freie Aufl"osung des $KG$-Moduls $K$. Diese
Aufl"osung hei"st \emph{bar resolution}.\\

\Bew F"ur jeden $KG$-Modul $V$ und $n\ge 0$ definieren wir eine lineare Abbildung
\begin{equation}
  \label{deffvonomega}
\omega_{n}^{V}: C^{n}(G,V)\rightarrow \Hom_{KG}(P_{n},V), \quad
g\mapsto\omega_{n}^{V}(g)
\end{equation}
mit
\[
\omega_{n}^{V}(g): P_{n}\rightarrow V, \quad
[\sigma_{1},\ldots,\sigma_{n}]\mapsto g(\sigma_{1},\ldots,\sigma_{n})\quad
(KG-\textrm{linear fortgesetzt}).
\]
Dann ist $\omega_{n}^{V}$ bijektiv, denn f"ur
\begin{equation}
  \label{defvonen}
e_{n}\in C^{n}(G,P_{n}) \textrm{ mit }
e_{n}(\sigma_{1},\ldots,\sigma_{n}):=[\sigma_{1},\ldots,\sigma_{n}]  
\end{equation}
gilt 
\begin{equation}
  \label{omegaTog}
  \omega_{n}^{V}(g)\circ e_{n}=g\quad \myforall g\in C^{n}(G,V)
\end{equation}
und
\begin{equation}
  \label{omegatophi}
  \omega_{n}^{V}(\varphi\circ e_{n})=\varphi \quad\myforall \varphi\in \Hom_{KG}(P_{n},V).
\end{equation}
Insbesondere f"ur $\varphi=\id_{P_{n}}$ folgt
\begin{equation}
  \label{omegaenidn}
\omega_{n}^{P_{n}}(e_{n})=\id_{P_{n}}.  
\end{equation}
Anhand der Definitionen \eqref{defpartdn} und \eqref{defdn} pr"uft man leicht,
dass f"ur $g\in C^{n}(G,V)$ ($n\ge 0$) die Gleichung
\begin{equation}
  \label{omegaequation}
\omega_{n+1}^{V}(\partial_{n}^{V}g)=\omega_{n}^{V}(g)\circ d_{n+1}  
\end{equation}
gilt. Insbesondere haben wir f"ur $e_{n}\in C^{n}(G,P_{n})$
\begin{eqnarray*}
\omega_{n+2}^{P_{n}}(\underbrace{\partial_{n+1}^{P_{n}}\partial_{n}^{P_{n}}e_{n}}_{=0})&\stackrel{~\eqref{omegaequation}}{=}&\omega_{n+1}^{P_{n}}(\partial_{n}^{P_{n}}e_{n})\circ
d_{n+2}\\
&\stackrel{~\eqref{omegaequation}}{=}&\omega_{n}^{P_{n}}(e_{n})\circ
d_{n+1}\circ d_{n+2}\\
&\stackrel{~\eqref{omegaenidn}}{=}&\id_{P_{n}}\circ
d_{n+1}\circ d_{n+2}=d_{n+1}\circ d_{n+2},
\end{eqnarray*}
und damit $d_{n+1}\circ d_{n+2}=0$ $(n\ge 0)$. Da auch $d_{0}\circ d_{1}=0$
(es ist $d_{0}(d_{1}([\sigma]))=d_{0}(\sigma[]-[])=1-1=0 \,\myforall \sigma\in
G$), ist also \eqref{thebarresolution} ein Komplex, also $\im d_{n+1}\subseteq
\ker d_{n} \,\myforall n\ge 0$. 

Um die umgekehrte Inklusion zu zeigen, betrachten wir die $K$-linearen
Abbildungen
\[
t_{n}:P_{n}\rightarrow P_{n+1},\quad
\sigma_{0}[\sigma_{1},\ldots,\sigma_{n}]\mapsto
[\sigma_{0},\sigma_{1},\ldots,\sigma_{n}] \quad (K\textrm{-linear fortgesetzt}).
\]
Da $\{\sigma_{0}[\sigma_{1},\ldots,\sigma_{n}]:\sigma_{0},\ldots,\sigma_{n}\in
G\}$ eine Basis von $P_{n}$ als $K$-Vektorraum ist, ist $t_{n}$
wohldefiniert. Dann gilt
\[
\id_{P_{n}}=d_{n+1}\circ t_{n}+t_{n-1}\circ d_{n}\quad (n\ge 1),
\]
denn unter der rechten Abbildung erhalten wir
\begin{eqnarray*}
\sigma_{0}[\sigma_{1},\ldots,\sigma_{n}]&\mapsto&d_{n+1}([\sigma_{0},\sigma_{1},\ldots,\sigma_{n}])+t_{n-1}\Big(\sigma_{0}\big(\sigma_{1}[\sigma_{2},\ldots,\sigma_{n}]\\&&+\sum_{i=1}^{n-1}(-1)^{i}[\sigma_{1},\ldots,\sigma_{i}\sigma_{i+1},\ldots,\sigma_{n}]+(-1)^{n}[\sigma_{1},\ldots,\sigma_{n-1}]
  \big)\Big)\\
&=&\sigma_{0}[\sigma_{1},\ldots,\sigma_{n}]+\sum_{i=0}^{n-1}(-1)^{i+1}[\sigma_{0},\ldots,\sigma_{i}\sigma_{i+1},\ldots,\sigma_{n}]\\
&&+(-1)^{n+1}[\sigma_{0},\ldots,\sigma_{n-1}]+[\sigma_{0}\sigma_{1},\sigma_{2},\ldots,\sigma_{n}]\\
&&+\sum_{i=1}^{n-1}(-1)^{i}[\sigma_{0},\sigma_{1},\ldots,\sigma_{i}\sigma_{i+1},\ldots,\sigma_{n}]+(-1)^{n}[\sigma_{0},\ldots,\sigma_{n-1}]\\
&=&\sigma_{0}[\sigma_{1},\ldots,\sigma_{n}],
\end{eqnarray*}
d.h. die rechte Abbildung ist gleich der Identit"at. Damit k"onnen wir nun
$\ker d_{n}\subseteq \im d_{n+1}$ f"ur $n\ge 1$ zeigen: F"ur $x\in \ker d_{n}$
ist
\[
x=\id_{P_{n}}(x)=d_{n+1}\circ t_{n}(x)+\underbrace{t_{n-1}\circ
  d_{n}(x)}_{=0}=d_{n+1}\circ t_{n}(x)\in \im d_{n+1}.
\]
Bleibt noch der Fall $n=0$. F"ur $x=\sum_{i=1}^{k}\lambda_{i}\sigma_{i}[]\in
\ker d_{0}$ mit $\lambda_{i}\in K, \sigma_{i}\in G$ folgt
$d_{0}(x)=\sum_{i=1}^{k}\lambda_{i}=0$. Mit $y:=\sum_{i=1}^{k}\lambda_{i}[\sigma_{i}]$
folgt
\[
d_{1}(y)=\sum_{i=1}^{k}\lambda_{i}(\sigma_{i}[]-[])=\sum_{i=1}^{k}\lambda_{i}\sigma_{i}[]=x,
\]
also auch $\ker d_{0}\subseteq \im d_{1}$. Damit ist nun bewiesen, dass
\eqref{thebarresolution} eine exakte Sequenz ist.\qed\\

\subsection{Beschreibung der Kohomologie durch $\Ext$}
Ist $V$ ein $G$-Modul, so erh"alt man aus dem Komplex \eqref{thebarresolution}
den Kokomplex
\begin{equation}
  \label{KokomplexOfBarResol}
  0\stackrel{d_{0}^{*}}{\rightarrow}\Hom_{KG}(P_{0},V)\stackrel{d_{1}^{*}}{\rightarrow}\ldots\stackrel{d_{n-1}^{*}}{\rightarrow} \Hom_{KG}(P_{n-1},V)\stackrel{d_{n}^{*}}{\rightarrow} \Hom_{KG}(P_{n},V)\stackrel{d_{n+1}^{*}}{\rightarrow}\ldots
\end{equation}
mit Differentialen
\[
d_{n}^{*}: \Hom_{KG}(P_{n-1},V)\rightarrow  \Hom_{KG}(P_{n},V), \quad f\mapsto
d_{n}^{*}(f)=f\circ d_{n} \quad (n\ge 1),
\]
und $d_{0}^{*}:=0$. Man bezeichnet mit
\[
\Ext_{KG}^{n}(K,V):=\ker d_{n+1}^{*}/\im d_{n}^{*} \quad (n\ge 0)
\]
die $n$-te Kohomologiegruppe dieses Kokomplexes. Wir zeigen die Isomorphie
der additiven Gruppen
\[
H^{n}(G,V)\cong \Ext_{KG}^{n}(K,V),
\]
womit (im wesentlichen) der Bogen zu der in der Literatur "ublichen Beschreibung gespannt ist.

Aus Gleichung \eqref{omegaequation} folgt sofort
\[
\omega_{n+1}^{V}\circ \partial_{n}^{V}=d_{n+1}^{*}\circ \omega_{n}^{V}, \quad (n\ge 0)
\]
und aus der Bijektivit"at von $\omega_{n}^{V} \,\,\myforall n\ge 0$ folgt hieraus
\[
\omega_{n}^{V}(\ker \partial_{n}^{V})=\ker d_{n+1}^{*}\quad \textrm{ und }\quad
\omega_{n+1}^{V}(\im \partial_{n}^{V})=\im d_{n+1}^{*} \quad \myforall n\ge 0.
\]
Daher induziert $\omega_{n}^{V}$ einen Gruppenisomorphismus von $\ker \partial_{n}^{V}/\im
\partial_{n-1}^{V}=H^{n}(G,V)$ auf $\ker d_{n+1}^{*}/\im
d_{n}^{*}=\Ext_{KG}^{n}(K,V)$ 
$(n \ge 1)$, und
$\omega_{0}^{V}$ einen Gruppenisomorphismus von $\ker \partial_{0}^{V}=H^{0}(G,V)$ auf $\ker d_{1}^{*}=\Ext_{KG}^{0}(K,V)$.

\subsection{Von kurzen zu langen exakten Sequenzen}
Wir wollen der Vollst"andigkeit halber noch eine h"aufig benutzte Eigenschaft der
Kohomologiegruppen beschreiben (wir verwenden diese nicht weiter, so dass dieser
Abschnitt "ubersprungen werden kann):

\begin{Satz}
Eine kurze exakte Sequenz von $G$-Moduln
\[
0\rightarrow U \stackrel{\varepsilon}{\hookrightarrow}V
\stackrel{\pi}{\rightarrow}W\rightarrow 0
\]
induziert eine lange exakte Sequenz von Kohomologiegruppen
\begin{eqnarray*}
&&0\rightarrow U^{G} \stackrel{\varepsilon_{0}}{\rightarrow}V^{G}
\stackrel{\pi_{0}}{\rightarrow}W^{G}
\stackrel{\gamma_{0}}{\rightarrow}H^{1}(G,U)
\stackrel{\varepsilon_{1}}{\rightarrow}
H^{1}(G,V)\stackrel{\pi_{1}}{\rightarrow}H^{1}(G,W)\stackrel{\gamma_{1}}{\rightarrow}H^{2}(G,U)\stackrel{\varepsilon_{2}}{\rightarrow}\ldots\\
&&\ldots\stackrel{\pi_{n-1}}{\rightarrow} H^{n-1}(G,W)  \stackrel{\gamma_{n-1}}{\rightarrow} H^{n}(G,U)\stackrel{\varepsilon_{n}}{\rightarrow}H^{n}(G,V)\stackrel{\pi_{n}}{\rightarrow}H^{n}(G,W)\stackrel{\gamma_{n}}{\rightarrow} H^{n+1}(G,U)\stackrel{\varepsilon_{n+1}}{\rightarrow}\ldots
\end{eqnarray*}
mit folgenden "Ubergangshomomorphismen ($n\ge 0$):
\begin{eqnarray*}
\varepsilon_{n}:& H^{n}(G,U)\rightarrow H^{n}(G,V),& g+B^{n}(G,U)\mapsto
\varepsilon \circ g + B^{n}(G,V) \quad (g\in Z^{n}(G,U))\\
\pi_{n}:& H^{n}(G,V)\rightarrow H^{n}(G,W),& g+B^{n}(G,V)\mapsto
\pi \circ g + B^{n}(G,W)\quad (g\in Z^{n}(G,V))\\
\gamma_{n}:& H^{n}(G,W)\rightarrow H^{n+1}(G,U),& g+B^{n}(G,W)\mapsto
\partial_{n}^{V}h+B^{n+1}(G,U)\quad (g\in Z^{n}(G,W)),\\&&\textrm{wenn } h\in C^{n}(G,V) \textrm{ mit } g=\pi \circ h.
\end{eqnarray*}
\end{Satz}

\Bew (i) Wir zeigen zun"achst die Wohldefiniertheit von $\pi_{n}$. Ist $g\in
Z^{n}(G,V)$, also $\partial_{n}^{V}g=0$, so ist $\partial_{n}^{W}(\pi\circ
g)=\pi\circ\partial_{n}^{V}g=0$, also $\pi \circ g \in Z^{n}(G,W)$. Ist $g\in
B^{n}(G,V)$, also $g=\partial_{n-1}^{V}h$ mit $h\in C^{n-1}(G,V)$, so ist $\pi
\circ g=\pi \circ \partial_{n-1}^{V}h=\partial_{n-1}^{W}(\pi \circ h)$ mit $\pi
\circ h\in C^{n-1}(G,W)$. Dies zeigt $\pi \circ g \in B^{n}(G,W)$. Es folgt
die Wohldefiniertheit von $\pi_{n}$. Die Wohldefiniertheit von
$\varepsilon_{n}$ zeigt man genauso. Es ist klar, dass $\pi_{n}$ und
$\varepsilon_{n}$ Gruppenhomomorphismen sind.

(ii) Wir zeigen die Exaktheit an der Stelle $H^{n}(G,V)$. Aus $\pi \circ
\varepsilon=0$ folgt sofort $\pi_{n}\circ \varepsilon_{n}=0$, also $\im
\varepsilon_{n}\subseteq \ker \pi_{n}$. F"ur die umgekehrte Inklusion sei
$g\in Z^{n}(G,V)$ mit $\pi_{n}(g+B^{n}(G,V))=0$, d.h. $\pi \circ g \in
B^{n}(G,W)$. Dann existiert ein $h'\in C^{n-1}(G,W)$ mit $\pi \circ g =
\partial_{n-1}^{W}h'$. Da $\pi$ surjektiv, existiert zu $h'$ ein $h\in
C^{n-1}(G,V)$ mit $h'=\pi \circ h$. Es folgt $\pi \circ g=\partial_{n-1}^{W}(\pi
\circ h)=\pi \circ \partial_{n-1}^{V}h$, also $\pi \circ
(g-\partial_{n-1}^{V}h)=0$. Aus $\ker \pi = U$ folgt also
$f:=g-\partial_{n-1}^{V}h\in C^{n}(G,U)$. Wegen $U\subseteq V$ ist dann
\[
\partial_{n}^{U}f=\partial_{n}^{V}g-\partial_{n}^{V}\partial_{n-1}^{V}h=0
\]
wegen $g\in Z^{n}(G,V)=\ker \partial_{n}^{V}$ und
$\partial_{n}\circ\partial_{n-1}=0$. Also ist $f\in Z^{n}(G,U)$ und 
\[
\varepsilon_{n}(f+B^{n}(G,U))=f+B^{n}(G,V)=g-\partial_{n-1}^{V}h+B^{n}(G,V)=g+B^{n}(G,V).
\]
Es folgt $\ker \pi_{n}\subseteq \im \varepsilon_{n}$.

(iii) Wir zeigen die Wohldefiniertheit von $\gamma_{n}$. Ist $g\in
Z^{n}(G,W)$, so existiert wegen der Surjektivit"at von $\pi$ jedenfalls ein
$h\in C^{n}(G,V)$ mit $g=\pi\circ h$. Es folgt $0=\partial_{n}^{W}g=\pi \circ
\partial_{n}^{V}h$, also $\partial_{n}^{V}h\in C^{n+1}(G,U)$. Es folgt
$\partial_{n+1}^{U}(\partial_{n}^{V}h)=\partial_{n+1}^{V}\partial_{n}^{V}h=0$,
also $\partial_{n}^{V}h\in Z^{n+1}(G,U)$. Ist auch $h'\in C^{n}(G,V)$ mit
$g=\pi \circ h'$, so ist $0=\pi \circ (h-h')$, also $h-h'\in C^{n}(G,U)$. Dann
ist $\partial_{n}^{V}h-\partial_{n}^{V}h'=\partial_{n}^{U}(h-h')\in
B^{n+1}(G,U)$ und damit
$\partial_{n}^{V}h+B^{n+1}(G,U)=\partial_{n}^{V}h'+B^{n+1}(G,U)$. 
Daher ist die Abbildung
\begin{eqnarray*}
\Gamma_{n}: Z^{n}(G,W)\rightarrow H^{n+1}(G,U),&&g\mapsto
\partial_{n}^{V}h+B^{n+1}(G,U)\quad (g\in Z^{n}(G,W))\\
&&\textrm{ wenn }h\in C^{n}(G,V) \textrm{ mit } g=\pi \circ h
\end{eqnarray*}
wohldefiniert. Man sieht auch leicht, dass $\Gamma_{n}$
ein Gruppenhomomorphismus ist. Wir zeigen, dass $B^{n}(G,W)$ im Kern von $\Gamma_{n}$
liegt, was dann zeigt, dass $\gamma_{n}$
wohldefiniert und ein Gruppenhomomorphismus ist.
Ist $g\in
B^{n}(G,W)$, so gibt es $f\in C^{n-1}(G,W)$ mit $g=\partial_{n-1}^{W}f$. Zu
$f$ gibt es dann $t\in C^{n-1}(G,V)$ mit $f=\pi\circ t$. Dann ist $g=\partial_{n-1}^{W}(\pi\circ t)=\pi\circ
\partial_{n-1}^{V}t$, und f"ur $h:=\partial_{n-1}^{V}t\in C^{n}(G,V)$ gilt
dann $g=\pi \circ h$ und
$\partial_{n}^{V}h=\partial_{n}^{V}\partial_{n-1}^{V}t=0$, also $\Gamma_{n}(g)=0$. Dies zeigt
$B^{n}(G,W)\subseteq \ker \Gamma_{n}$.

(iv) Exaktheit an der Stelle $H^{n}(G,W)$. Sei $g\in Z^{n}(G,V)$. Dann ist
$\pi\circ g\in Z^{n}(G,W)$. Es folgt
\[
\gamma_{n}(\pi\circ g + B^{n}(G,W))=\partial_{n}^{V}g+B^{n+1}(G,U)=0+B^{n+1}(G,U),
\]
da $\partial_{n}^{V}g=0$. Also ist $\gamma_{n}\circ\pi_{n}=0$ und $\im
\pi_{n}\subseteq \ker \gamma_{n}$. F"ur die umgekehrte Inklusion sei $g\in
Z^{n}(G,W)$ mit $\gamma_{n}(g+B^{n}(G,W))=0$, d.h. f"ur $h\in C^{n}(G,V)$ mit
$g=\pi\circ h$ gilt
$\partial_{n}^{V}h\in B^{n+1}(G,U)$. Dann gibt es $f\in C^{n}(G,U)$ mit
$\partial_{n}^{V}h=\partial_{n}^{U}f$, also $\partial_{n}^{V}(h-f)=0$ oder
$h-f\in Z^{n}(G,V)$. Da $\pi \circ f=0$ (wegen $f\in C^{n}(G,U)$), folgt also $g=\pi \circ h=\pi\circ
(h-f)$ und damit $g+B^{n}(G,W)=\pi_{n}(h-f+B^{n}(G,V))$. Dies zeigt $\ker
\gamma_{n}\subseteq \im \pi_{n}$.

(v) Exaktheit an der Stelle $H^{n}(G,U)$. Sei $f\in Z^{n-1}(G,W)$ und $h\in
C^{n-1}(G,V)$ mit $f=\pi \circ h$. Dann ist
$\gamma_{n-1}(f+B^{n-1}(G,W))=\partial_{n-1}^{V}h+B^{n}(G,U)$ und damit
$\varepsilon_{n}(\gamma_{n-1}(f+B^{n-1}(G,W)))=\partial_{n-1}^{V}h+B^{n}(G,V)=0+B^{n}(G,V)$,
da $h\in C^{n-1}(G,V)$. Also ist $\varepsilon_{n}\circ\gamma_{n-1}=0$ und $\im
\gamma_{n-1}\subseteq\ker\varepsilon_{n}$. Sei umgekehrt $g\in Z^{n}(G,U)$ mit
$g+B^{n}(G,U)\in \ker \varepsilon_{n}$, also $g\in B^{n}(G,V)$. Dann gibt es
$f\in C^{n-1}(G,V)$ mit $g=\partial_{n-1}f$. F"ur $\pi\circ f\in C^{n-1}(G,W)$
gilt dann $\partial_{n-1}^{W}(\pi\circ f)=\pi \circ \partial_{n-1}^{V}f=\pi
\circ g=0$, da $g\in Z^{n}(G,U)$. Also ist $\pi \circ f\in Z^{n-1}(G,W)$ und
$\gamma_{n-1}(\pi\circ f +
B^{n-1}(G,W))=\partial_{n-1}^{V}f+B^{n}(G,U)=g+B^{n}(G,U)$. Dies zeigt $\ker
\varepsilon_{n}\subseteq \im \gamma_{n-1}$. 

Insgesamt haben wir die Exaktheit der langen Sequenz gezeigt. \qed\\

\subsection{Von exakten Sequenzen zu Kozyklen}
In diesem Abschnitt wollen wir eine Konstruktion angeben, bei der gewissen exakten
Sequenzen ein bis auf einen Korand eindeutiger Kozyklus zugewiesen wird.

\begin{SatzDef} \label{SatzVomIndKoz}
Sei $n\ge 1$,
\begin{equation} \label{allgemeineKozyklusSeq}
0\rightarrow
U_{n}\stackrel{\pi_{n}}{\rightarrow}U_{n-1}\stackrel{\pi_{n-1}}{\rightarrow}U_{n-2}\stackrel{\pi_{n-2}}{\rightarrow}\ldots\stackrel{\pi_{1}}{\rightarrow}U_{0}\stackrel{\pi_{0}}{\rightarrow}U_{-1}
\end{equation}
eine exakte Sequenz von $G$-Moduln und $w\in U_{-1}^{G} \cap \im\pi_{0}$. F"ur $k=0,\ldots,n$ sei $g_{k}\in C^{k}(G,U_{k})$ mit
\[
\pi_{0}(g_{0})=w\quad \textrm{ und }\quad \pi_{k}\circ
g_{k}=\partial_{k-1}^{U_{k-1}}g_{k-1} \quad \textrm{ f"ur }
k=1,\ldots,n. \quad (*)
\]
Dann ist $g_{n}\in Z^{n}(G,U_{n})$. Eine Folge $(g_{k})$ mit den Eigenschaften $(*)$
existiert immer, und die Restklasse $g_{n}+B^{n}(G,U_{n})$ ist unabh"angig von
der speziellen Wahl einer Folge $(g_{k})$ mit den Eigenschaften $(*)$.
Ist $g_{n}'\in g_{n}+B^{n}(G,U_{n})$ ein weiterer Repr"asentant der Restklasse, so existieren
dazu $g_{k}'\in C^{k}(G,U_{k}),\, k=0,\ldots,n-1$ mit der $(*)$ entsprechenden Eigenschaft,
d.h. jeder Kozyklus der induzierten Restklasse kann durch die Konstruktion
auch "`realisiert"' werden.
 
Existiert ein $v\in U_{0}^{G}$ mit $\pi_{0}(v)=w$, oder f"ur ein $0\le k_{0}\le n$
ein $G$-Homomorphismus $\mu_{k_{0}}: U_{k_{0}-1}\rightarrow U_{k_{0}}$ mit
$\pi_{k_{0}}\circ\mu_{k_{0}}=\id_{U_{k_{0}-1}}$ (d.h. ein \emph{rechts-Splitting}
an der Stelle $k_{0}$) oder $\mu_{k_{0}}\circ \pi_{k_{0}}=\id_{U_{k_{0}}}$
(d.h. ein \emph{links-Splitting} \index{Splitting}
an der Stelle $k_{0}$), so ist $g_{n}\in B^{n}(G,U_{n})$.

Im
Falle einer \emph{exakten Sequenz der L"ange $n$ in
  Standardform},\index{exakte Sequenz!Standardform}
\index{exakte Sequenz!L\"ange $n$} 
\begin{equation}
  \label{exakteSequenzdesKozyklus}
0\rightarrow
U_{n}\stackrel{\pi_{n}}{\rightarrow}U_{n-1}\stackrel{\pi_{n-1}}{\rightarrow}U_{n-2}\stackrel{\pi_{n-2}}{\rightarrow}\ldots\stackrel{\pi_{1}}{\rightarrow}U_{0}\stackrel{\pi_{0}}{\rightarrow}K\rightarrow
0
\end{equation}
und $w=1\in K$
nennen wir $g_{n}$ einen \emph{von der exakten Sequenz
  \eqref{exakteSequenzdesKozyklus} induzierten Kozyklus}.
\end{SatzDef}

In Beispiel \ref{Bsp2Kozyklus} werden wir sehen, dass die angegebenen
Kriterien f"ur das Vorliegen eines Korands nur hinreichend, aber nicht
notwendig sind. Falls man ein links/rechts Splitting hat, dass nicht am
linken/rechten Rand liegt, so ist das Kriterium relativ trivial, da dann die
links/rechts benachbarte Abbildung die Nullabbildung ist.\\

Interpretiert man $\pi_{n}$ als Inklusion $U_{n}\hookrightarrow U_{n-1}$, so
folgt aus $(*)$,  das $g_{n}=\partial_{n-1}^{U_{n-1}}g_{n-1}$, also $g_{n}\in
B^{n}(G,U_{n-1})$. In $U_{n-1}$ wird $g_{n}$ also zu einem Korand.\\

\Bew (i) Wir zeigen zun"achst die Existenz einer Folge $(g_{k})$ mit der
Eigenschaft $(*)$. Da $w\in \im \pi_{0}$ (dies gilt insbesondere
f"ur den Fall der Sequenz \eqref{exakteSequenzdesKozyklus} und $w=1\in K=K^{G}$), existiert ein $g_{0}\in
C^{0}(G,U_{0})=U_{0}$ mit $\pi_{0}(g_{0})=w$. Dann gilt f"ur $\sigma\in G$,
dass $\partial_{0}^{U_{0}}g_{0}(\sigma)=\sigma g_{0}-g_{0}$, also
$\pi_{0}(\partial_{0}^{U_{0}}g_{0}(\sigma))=\sigma w-w=0$ ($w\in U_{-1}^{G}$). Da $\ker \pi_{0}=\im
\pi_{1}$, existiert also $g_{1}\in C^{1}(G,U_{1})$ mit
$\partial_{0}^{U_{0}}g_{0}=\pi_{1}\circ g_{1}$. Wir haben also $(*)$ f"ur $k=1$.

Es sei nun $g_{0},\ldots,g_{k}$ mit $(*)$ bereits konstruiert. Es folgt
\begin{equation}\label{keinSpezialfall}
\pi_{k}\circ\partial_{k}^{U_{k}}g_{k}=\partial_{k}^{U_{k-1}}(\pi_{k}\circ
g_{k})\stackrel{(*)}{=}\partial_{k}^{U_{k-1}}\partial_{k-1}^{U_{k-1}}g_{k-1}=0,
\end{equation}
(da $\partial_{k}\circ\partial_{k-1}=0$) und wegen $\ker \pi_{k}=\im \pi_{k+1}$ gibt es also $g_{k+1}\in
C^{k+1}(G,U_{k+1})$ mit $\pi_{k+1}\circ
g_{k+1}=\partial_{k}^{U_{k}}g_{k}$. Dies ist $(*)$ f"ur $k+1$.

(ii) Ist nun $(g_{k})$ ein Folge mit $(*)$, so gilt insbesondere f"ur $k=n$, dass
\[
\pi_{n}\circ \partial_{n}^{U_{n}}g_{n}=\partial_{n}^{U_{n-1}}(\pi_{n}\circ
g_{n})\stackrel{(*)}{=}\partial_{n}^{U_{n-1}}\partial_{n-1}^{U_{n-1}}g_{n-1}=0.
\]
(Der
Leser mag einwenden, dass diese Gleichung nur ein Spezialfall von
\eqref{keinSpezialfall} f"ur $k=n$ ist. Auf \eqref{keinSpezialfall} werden
wir aber durch \emph{Konstruktion} einer \emph{speziellen} Folge gef"uhrt, die
$(*)$ erf"ullen \emph{soll}, w"ahrend wir vorige Gleichung f"ur \emph{jede} Folge zeigen,
die $(*)$ erf"ullt. Die Rechnung ist nat"urlich die gleiche.)
Da $\pi_{n}$ injektiv, folgt dann aber $\partial_{n}^{U_{n}}g_{n}=0$, also
$g_{n}\in Z^{n}(G,U_{n})$.

(iii) Sei nun $h_{k}\in C^{k}(G,U_{k})$ f"ur $k=0,\ldots,n$ eine weitere Folge mit
der Eigenschaft $(*)$, also
\[
\pi_{0}(h_{0})=w\quad \textrm{ und }\quad \pi_{k}\circ
h_{k}=\partial_{k-1}^{U_{k-1}}h_{k-1} \quad \textrm{ f"ur }
k=1,\ldots,n.
\]
Wir m"ussen $g_{n}+B^{n}(G,U_{n})=h_{n}+B^{n}(G,U_{n})$ zeigen.
 
Wir zeigen zun"achst: F"ur $k=1,\ldots,n$ existiert ein $f_{k}\in
C^{k-1}(G,U_{k})$ mit
\[
\pi_{k}\circ(g_{k}-h_{k})=\pi_{k}\circ \partial_{k-1}^{U_{k}}f_{k}.\quad (**)
\]
Da $\pi_{0}(g_{0}-h_{0})=w-w=0$ und $\ker \pi_{0}=\im \pi_{1}$, gibt es
$f_{1}\in C^{0}(G,U_{1})$ mit $\pi_{1}\circ f_{1} = g_{0}-h_{0}$. Dann ist
\[
\pi_{1}\circ(g_{1}-h_{1})\stackrel{(*)}{=}\partial_{0}^{U_{0}}(g_{0}-h_{0})=\partial_{0}^{U_{0}}(\pi_{1}\circ
f_{1})=\pi_{1}\circ \partial_{0}^{U_{1}} f_{1}.
\]
Dies ist $(**)$ f"ur $k=1$. Sei nun $f_{1},\ldots,f_{k}$ mit $(**)$ bereits
konstruiert. Dann ist
\[
\pi_{k}\circ (g_{k}-h_{k}-\partial_{k-1}^{U_{k}}f_{k})=0,
\]
und wegen $\ker \pi_{k}=\im \pi_{k+1}$ existiert ein $f_{k+1}\in
C^{k}(G,U_{k+1})$ mit
\[
g_{k}-h_{k}-\partial_{k-1}^{U_{k}}f_{k}=\pi_{k+1}\circ f_{k+1}.
\]
Anwenden von $\partial_{k}^{U_{k}}$ liefert
\[
\partial_{k}^{U_{k}}g_{k}-\partial_{k}^{U_{k}}h_{k}-\partial_{k}^{U_{k}}\partial_{k-1}^{U_{k}}f_{k}=\partial_{k}^{U_{k}}(\pi_{k+1}\circ f_{k+1})=
\pi_{k+1}\circ\partial_{k}^{U_{k+1}} f_{k+1},
\]
und wegen $\partial_{k}^{U_{k}}\circ\partial_{k-1}^{U_{k}}=0$ und $(*)$ folgt
\[
\pi_{k+1}\circ(g_{k+1}-h_{k+1})=\pi_{k+1}\circ\partial_{k}^{U_{k+1}} f_{k+1}.
\]
Dies ist $(**)$ f"ur $k+1$. Insbesondere f"ur $k=n$ folgt aus $(**)$ und der
Injektivit"at von $\pi_{n}$, dass
\[
g_{n}-h_{n}=\partial_{n-1}^{U_{n}}f_{n}\in B^{n}(G,U_{n}).
\]
Damit ist also die Restklasse $g_{n}+B^{n}(G,U_{n})=h_{n}+B^{n}(G,U_{n})$ von der speziellen Wahl der
Folge~$(g_{k})$ mit $(*)$ unabh"angig.

(iv) Sei nun $g_{n}'\in g_{n}+B^{n}(G,U_{n})$ ein weiterer Repr"asentant der
Restklasse, d.h. $g_{n}'=g_{n}+\partial_{n-1}^{U_{n}}f$ mit $f\in
C^{n-1}(G,U_{n})$. Wir setzen $g_{n-1}':=g_{n-1}+\pi_{n}\circ f\in
C^{n-1}(G,U_{n-1})$ und $g_{k}':=g_{k}$ f"ur $k=0,\ldots,n-2$. Dann hat die
Folge der $g_{k}'$ mit $k=0,\ldots,n$ die $(*)$ entsprechende Eigenschaft und
induziert damit $g_{n}'$. Nach Konstruktion m"ussen wir dabei $(*)$ nur noch
f"ur $k=n$ und $k=n-1$ testen. Da
\begin{eqnarray*}
\pi_{n}\circ g_{n}'&=&\pi_{n}\circ
(g_{n}+\partial_{n-1}^{U_{n}}f)\stackrel{{(*)}}{=}\partial_{n-1}^{U_{n-1}}g_{n-1}+\partial_{n-1}^{U_{n-1}}(\pi_{n}\circ
f)\\& =& \partial_{n-1}^{U_{n-1}}(g_{n-1}+\pi_{n}\circ
f)=\partial_{n-1}^{U_{n-1}}g_{n-1}',
\end{eqnarray*}
also $(*)$ f"ur die Folge $(g_{k}')$ und $k=n$ gilt, sowie
\begin{eqnarray*}
\pi_{n-1}\circ g_{n-1}'&=&\pi_{n-1}\circ(g_{n-1}+\pi_{n}\circ
f)\stackrel{\pi_{n-1}\circ \pi_{n}=0}{=}\pi_{n-1}\circ
g_{n-1}\\&\stackrel{(*)}{=}&\partial_{n-2}^{U_{n-2}}g_{n-2}=\partial_{n-2}^{U_{n-2}}g_{n-2}'
\end{eqnarray*}
gilt, also $(*)$ f"ur die Folge $(g_{k}')$ und $k=n-1$ gilt, erf"ullt also auch die Folge
$(g_{k}')$ f"ur $k=0,\ldots,n$ die Bedingung $(*)$ und induziert den
Repr"asentant der Restklasse $g_{n}'$. (Im Fall $n=1$ und $k=n-1=0$ erh"alt man
analog $\pi_{0}(g_{0}')=\pi_{0}(g_{0}+\pi_{1}\circ f)=w$).

(v) Wir nehmen nun die Existenz eines $v\in U_{0}^{G}$ mit $\pi_{0}(v)=w$ an. Wir
definieren $g_{k}'\in C^{k}(G,U_{k})$ f"ur $k=0,\ldots,n$ durch $g_{0}':=v$ und
$g_{k}':=0$ f"ur $k=1,\ldots,n$. Dann gilt $\pi_{0}(g_{0}')=w$. Weiter gilt
\[
\partial_{0}^{U_{0}}g_{0}'(\sigma)=\sigma v - v\stackrel{v\in
  U_{0}^{G}}{=}0=\pi_{1}\circ g_{1}'(\sigma),
\]
und offenbar auch $\pi_{k}\circ g_{k}'=0=\partial_{k-1}^{U_{k-1}}g_{k-1}'$ f"ur
$k=2,\ldots,n$. Dies ist $(*)$ f"ur die Folge $(g_{k}')$, und nach der bereits 
bewiesenen Unabh"angigkeit der induzierten Restklasse folgt
$g_{n}+B^{n}(G,U_{n})=g_{n}'+B^{n}(G,U_{n})=0+B^{n}(G,U_{n})$, 
d.h. $g_{n}$ ist ein $n$-Korand.

(vi) Sei nun f"ur ein $0\le k_{0}\le n$
ein $G$-Homomorphismus $\mu_{k_{0}}: U_{k_{0}-1}\rightarrow U_{k_{0}}$ mit
$\pi_{k_{0}}\circ\mu_{k_{0}}=\id_{U_{k_{0}-1}}$ gegeben, also ein rechts-Splitting. Im Fall $k_{0}=0$
setzen wir $v:=\mu_{0}(w)\in U_{0}^{G}$ (dann gilt $\pi_{0}(v)=w$) und kommen
zum vorherigen Fall (v). Sei
daher nun $1\le k_{0}\le n$. Sei weiter
$(g_{k})_{k=0,\ldots,n}$ eine Folge mit $(*)$. Aus $\ker \pi_{k_{0}-1}=\im
\pi_{k_{0}}=U_{k_{0}-1}$ (da $\pi_{k_{0}}\circ\mu_{k_{0}}=\id_{U_{k_{0}-1}}$) folgt $\pi_{k_{0}-1}=0$ und daraus sofort, dass man
$g_{k}=0\in B^{k}(G,U_{k})$ f"ur $k\ge k_{0}-1$ w"ahlen kann.

Wir geben noch
ein alternatives Argument. (Der eilige Leser kann zu (vii) springen).  Dieses funktioniert auch dann, wenn $\mu_{k_{0}}$ nur
auf einem "`hinreichend gro"sen"' Untermodul $U'\subseteq U_{k_{0}-1}$
definiert ist und $\pi_{k_{0}}\circ\mu_{k_{0}}=\id_{U'}$ erf"ullt; Hinreichend
gro"s hei"st hier, dass $\mu_{k_{0}}\circ
\partial_{k_{0}-1}^{U_{k_{0}-1}}g_{k_{0}-1}$ definiert ist, die Bedingung
h"angt also von der Wahl der Folge $(g_{k})$ ab. Wir setzen
$g_{k_{0}}':=\mu_{k_{0}}\circ \partial_{k_{0}-1}^{U_{k_{0}-1}}g_{k_{0}-1}\in
C^{k_{0}}(G,U_{k_{0}})$ und $g_{k}':=0\in C^{k}(G,U_{k})$ f"ur
$k=k_{0}+1,\ldots,n$, und behaupten, dass die Folge
$g_{0},\ldots,g_{k_{0}-1},g_{k_{0}}',g_{k_{0}+1}',\ldots,g_{n}'$ die $(*)$
entsprechende Eigenschaft mit $g_{n}'\in B^{n}(G,U_{n})$ erf"ullt. Dann gilt $g_{n}-g_{n}'\in
B^{n}(G,U_{n})$ nach dem bereits bewiesenen, und damit auch $g_{n}\in
B^{n}(G,U_{n})$. Offenbar ist $(*)$ nur noch f"ur $k=k_{0}$ und $k=k_{0}+1$
(wenn $k_{0}<n$) zu
pr"ufen. F"ur $k=k_{0}$ erhalten wir
\[
\pi_{k_{0}}\circ g_{k_{0}}'=\pi_{k_{0}}\circ \mu_{k_{0}}\circ
\partial_{k_{0}-1}^{U_{k_{0}-1}}g_{k_{0}-1}=\id_{U_{k_{0}-1}}\circ\partial_{k_{0}-1}^{U_{k_{0}-1}}g_{k_{0}-1}=
\partial_{k_{0}-1}^{U_{k_{0}-1}}g_{k_{0}-1},
\]
also $(*)$ f"ur $k=k_{0}$, und f"ur $k=k_{0}+1$ erhalten wir
\[
\partial_{k_{0}}^{U_{k_{0}}}g_{k_{0}}'=\partial_{k_{0}}^{U_{k_{0}}}\left(\mu_{k_{0}}\circ
\partial_{k_{0}-1}^{U_{k_{0}-1}}g_{k_{0}-1}\right)=\mu_{k_{0}}\circ
\partial_{k_{0}}^{U_{k_{0}-1}}\partial_{k_{0}-1}^{U_{k_{0}-1}}g_{k_{0}-1}=0=\pi_{k_{0}+1}\circ
g_{k_{0}+1}',
\]
da $\partial_{k_{0}}^{U_{k_{0}-1}}\circ\partial_{k_{0}-1}^{U_{k_{0}-1}}=0$ und
$g_{k_{0}+1}'=0$,
also $(*)$ f"ur $k=k_{0}+1$. Ist $k_{0}<n$, so ist $g_{n}'=0\in
B^{n}(G,U_{n})$. F"ur $k_{0}=n$ erhalten wir nach Definition ebenfalls
$g_{n}'=g_{k_{0}}'=\mu_{k_{0}}\circ
\partial_{k_{0}-1}^{U_{k_{0}-1}}g_{k_{0}-1}=\partial_{k_{0}-1}^{U_{k_{0}}}(\mu_{k_{0}}\circ
g_{k_{0}-1}
)\in B^{n}(G,U_{n})$.

(vii) Wir kommen zum Fall, dass f"ur ein $0\le k_{0}\le n$
ein $G$-Homomorphismus $\mu_{k_{0}}: U_{k_{0}-1}\rightarrow U_{k_{0}}$ mit
$\mu_{k_{0}}\circ\pi_{k_{0}}=\id_{U_{k_{0}}}$ gegeben ist, also ein
links-Splitting. (Auch hier gen"ugt es, dass $\mu_{k_{0}}$ auf einem
"`hinreichend gro"sen"' Untermodul definiert ist, vgl. die Bemerkung oben).
Sei wieder
$(g_{k})_{k=0,\ldots,n}$ eine Folge mit $(*)$. 
Im
Fall $k_{0}=0$ folgt aus $\pi_{0}(g_{0})=w$ nach Anwendung von $\mu_{0}$, dass
$g_{0}=\mu_{0}(w)\in U_{0}^{G}$ (da $w\in U_{-1}^{G}$); Nach dem
bereits betrachteten Kriterium mit $v:=g_{0}$ folgt $g_{n}\in B^{n}(G,U_{n})$.
Sei also $k_{0}\ge 1$. Wir setzen $g_{k}':=0\in
C^{k}(G,U_{k})$ f"ur $k=k_{0}+1,\ldots,n$ und zeigen, dass die Folge
$g_{0},\ldots,g_{k_{0}},g_{k_{0}+1}',\ldots,g_{n}'$ die $(*)$ entsprechende
Eigenschaft erf"ullt. Dabei ist $(*)$ diesmal offenbar nur f"ur $k=k_{0}+1$
(falls $k_{0}<n$) zu pr"ufen. 
Aus $\pi_{k_{0}}\circ
g_{k_{0}}\stackrel{(*)}{=}\partial_{k_{0}-1}^{U_{k_{0}-1}}g_{k_{0}-1}$ folgt 
nach Anwenden von $\mu_{k_{0}}$ unter Beachtung von
$\mu_{k_{0}}\circ \pi_{k_{0}}=\id_{U_{k_{0}}}$ dann
$g_{k_{0}}=\partial_{k_{0}-1}^{U_{k_{0}}}(\mu_{k_{0}}\circ g_{k_{0}-1})$. Also gilt
$g_{k_{0}}\in B^{k_{0}}(G,U_{k_{0}})$. Ist $k_{0}=n$, so haben wir damit
bereits die Behauptung. Wenn $k_{0}<n$, so gilt $\pi_{k_{0}+1}\circ
g_{k_{0}+1}'=0=\partial_{k_{0}}^{U_{k_{0}}}g_{k_{0}}$ (da $g_{k_{0}+1}'=0$ und
 $g_{k_{0}}\in B^{k_{0}}(G,U_{k_{0}})$), also $(*)$. Wir
erhalten $g_{n}=g_{n}-g_{n}'\in B^{n}(G,U_{n})$, also die Behauptung. \qed

\begin{BemRoman}\label{gnAlsnKorand}
W"ahlt man jeweils die konstruierte Folge (mit der Mischung aus $g_{k}'$ und $g_{k}$) in einem der
Kriterien f"ur das Vorliegen eines Korands als alternative Folge $(h_{k})$ zur
Konstruktion der Folge $(f_{k})$ gem"a"s $(**)$, so erh"alt man wegen
 $\pi_{n}$ injektiv im Fall $h_{n}=0$ aus $(**)$ also insbesondere
$g_{n}=\partial_{n-1}f_{n}$ (wenn $h_{n}\ne 0$, so konnten wir $g_{n}$ im
Beweis sowieso schon so schreiben). Da die Konstruktion der Folge $(f_{k})$ explizit
ist (sofern man Urbilder unter den $\pi_{k}$ explizit angeben kann), zeigt der
Beweis also auch, wie man dann $g_{n}$ explizit als $n-$Korand schreiben kann.
\end{BemRoman}

\begin{BspRoman}\label{Bsp2Kozyklus}
Wir betrachten die exakte Sequenz von $\Ga$-Moduln
\[
0\rightarrow K \stackrel{\pi_{2}}{\rightarrow}K^{2}
\stackrel{\pi_{1}}{\rightarrow}K^{2}\stackrel{\pi_{0}}{\rightarrow}K\rightarrow
0
\]
mit $\pi_{2}(x):=(x,0)$, $\pi_{1}(x,y)=(y,0)$ und
$\pi_{0}(x,y)=y$ f"ur $x,y\in K$.
Die Operation auf $K$ ist dabei trivial, und die Operation auf $K^{2}$ gegeben
durch $a\cdot (x,y):=(x+ay,y)\, \myforall a\in \Ga, (x,y)\in K^{2}$. Offenbar ist $\pi_{0}(0,1)=1$, also w"ahlen wir
$g_{0}=(0,1)\in C^{0}(\Ga,K^{2})$. Sei ab jetzt $x,y\in K,\,\, a,b\in \Ga$.
Es ist $\partial_{0}g_{0}(a)=(a,1)-(0,1)=(a,0)=\pi_{1}(0,a)$, und wir w"ahlen daher
$g_{1}\in C^{1}(\Ga,K^{2})$ mit $g_{1}(a):=(0,a)$. Es ist
\[
\partial_{1}g_{1}(a,b)=a\cdot
g_{1}(b)-g_{1}(a+b)+g_{1}(a)=(ab,b)-(0,a+b)+(0,a)=(ab,0)=\pi_{2}(ab).
\]
Daher w"ahlen wir $g_{2}\in C^{2}(\Ga,K)$ mit $g_{2}(a,b)=ab$. Dann ist
$g_{2}\in Z^{2}(\Ga,K)$. 

Ist $\chr K\ne 2$, so ist mit $h\in C^{1}(\Ga,K)$ gegeben durch $h(a):=-\frac{1}{2} a^{2}$
\[
g_{2}(a,b)=ab=-\frac{1}{2}b^{2} +\frac{1}{2}(a+b)^{2}-\frac{1}{2}a^{2}=a\cdot
h(b)-h(a+b)+h(a)=\partial_{1}h(a,b),
\]
also $g_{2}\in B^{2}(\Ga,K)$. Dennoch gibt es kein $(x,y)\in (K^{2})^{\Ga}$ mit
$\pi_{0}(x,y)=1$. Dann gibt es erst recht kein $\mu_{0}: K\rightarrow K^{2}$
mit $\pi_{0}\circ \mu_{0}=\id$. Aus $\mu_{0}\circ \pi_{0}=\id$ folgte der
Widerspruch $\pi_{1}=0$. G"abe es $\mu_{1}: K^{2}\rightarrow K^{2}$ mit $\pi_{1}\circ
\mu_{1}=\id$ bzw. $\mu_{1}\circ \pi_{1}=\id$, so w"are $\pi_{0}=0$ bzw. $\pi_{2}=0$. G"abe es
$\mu_{2}: K^{2}\rightarrow K$ mit $\pi_{2}\circ
\mu_{2}=\id$, so w"are $\pi_{1}=0$. Es bleibt der Fall $\mu_{2}\circ
\pi_{2}=\id$. Dann w"are aber $\ker \mu_{2}$ ein $\Ga$-stabiles Komplement zu
$\im \pi_{2}=K\cdot
(1,0)$ in $K^{2}$ - ein Widerspruch.
Unsere Kriterien f"ur das Vorliegen eines Korands sind also
nur hinreichend, aber nicht notwendig.

Sei nun $\chr K=2$. Angenommen, es g"abe $h\in
C^{1}(\Ga,K)$ mit $g_{2}=\partial_{1}h.$ Dann ist
\[
g_{2}(a,b)=ab=a\cdot h(b)-h(a+b)+h(a)=h(a)+h(b)+h(a+b)\quad \myforall a,b\in K.
\]
F"ur $a=b$ folgt hieraus $a^{2}=h(0) \,\,\myforall a\in K$, ein Widerspruch. Dies
zeigt $g_{2} \not\in B^{2}(\Ga,K)$ f"ur $\chr K=2$.
\end{BspRoman}

\subsection{Pushout und Pullback}
F"ur die Konstruktion von exakten Sequenzen aus Kozyklen ben"otigen wir die
Konstruktionen "`Pushout"' und "`Pullback"'. Diese sind Spezialf"alle von
direktem und inversem Limes. F"ur diese allgemeineren Konstrukte siehe etwa
\cite[S. 39ff]{Rotman}. Wir beschr"anken uns hier auf $KG$-Moduln. Die
Formulierung f"ur $R$-Moduln mit einem beliebigen Ring $R$ ist w"ortlich
dieselbe. 

\subsubsection{Pushout}
\begin{SatzDef}\label{DefPushout}
Seien $X,Y_{1},Y_{2}$ jeweils $G$-Moduln und $\varphi_{1}: X\rightarrow Y_{1}$,
$\varphi_{2}: X\rightarrow Y_{2}$ jeweils $G$-Homomorphismen. Der
\emph{Pushout}\index{Pushout} von $\varphi_{1}$ und $\varphi_{2}$ ist der $G$-Modul
\[
Z:=(Y_{1}\oplus Y_{2})/W  \quad \textrm{ mit } \quad W:=\left\{(\varphi_{1}(x),-\varphi_{2}(x))\in Y_{1}\oplus
  Y_{2}: x\in X\right\}
\]
zusammen mit den Homomorphismen
\[
\pi_{1}: Y_{1}\rightarrow Z, \quad y_{1}\mapsto (y_{1},0)+W
\]
und
\[
\pi_{2}: Y_{2}\rightarrow Z, \quad y_{2}\mapsto (0,y_{2})+W.
\]
Dann gilt $\pi_{1}\circ \varphi_{1}=\pi_{2}\circ \varphi_{2}$, 
\begin{center}
~\begin{xy}
\xymatrix{
X \ar[r]^{\varphi_{2}} \ar[d]_{\varphi_{1}} & Y_{2} \ar@{.>}[d]^{\pi_{2}} \\
Y_{1} \ar@{.>}[r]_{\pi_{1}} & Z,\\
}
\end{xy}
\end{center}
und es gilt
folgende universelle Eigenschaft: F"ur jeden weiteren $G$-Modul $Z'$ und
Homomorphismen $\pi_{1}': Y_{1}\rightarrow Z'$ und $\pi_{2}': Y_{2}\rightarrow
 Z'$ mit $\pi_{1}'\circ\varphi_{1}=\pi_{2}'\circ\varphi_{2}$ existiert genau
ein Homomorphismus $\pi: Z\rightarrow Z'$, der das folgende
\emph{Pushout-Diagramm}\index{Pushout!Diagramm} kommutativ macht:
\begin{center}
~\begin{xy}
\xymatrix{
X \ar[r]^{\varphi_{2}} \ar[d]_{\varphi_{1}} & Y_{2} \ar[d]^{\pi_{2}} \ar@/^/[ddr]^{\pi_{2}'}&\\
Y_{1} \ar[r]_{\pi_{1}}\ar@/_/[drr]_{\pi_{1}'} & Z \ar@{.>}[dr]|-{\pi}&\\
&&Z'
}
\end{xy}
\end{center}
Weiter gilt: Ist $\varphi_{1}$ injektiv/surjektiv, so ist auch  $\pi_{2}$ injektiv/surjektiv.
\end{SatzDef}

\Bew Offenbar ist $W$ ein $G$-Untermodul von $Y_{1}\oplus Y_{2}$ und daher $Z$
ein $G$-Modul. 
F"ur $x\in X$ ist
\begin{eqnarray*}
\pi_{1}\circ\varphi_{1}(x)&=&(\varphi_{1}(x),0)+W=(\varphi_{1}(x),0)-(\varphi_{1}(x),-\varphi_{2}(x))+W\\&=&(0,\varphi_{2}(x))+W=\pi_{2}\circ\varphi_{2}(x),
\end{eqnarray*}
also $\pi_{1}\circ\varphi_{1}=\pi_{2}\circ\varphi_{2}$. Sind weiter
$Z',\pi_{1}',\pi_{2}'$ wie beschrieben, so gilt f"ur
\[
\pi':  Y_{1}\oplus Y_{2}\rightarrow Z', \quad (y_{1},y_{2})\mapsto
\pi_{1}'(y_{1})+\pi_{2}'(y_{2})
\]
wegen $\pi_{1}'\circ \varphi_{1}=\pi_{2}'\circ \varphi_{2}$ offenbar $W\subseteq \ker \pi'$. Daher induziert $\pi'$ eine Abbildung
\[
\pi: Z=(Y_{1}\oplus Y_{2})/W \rightarrow Z', \quad \textrm{ mit }\quad
(y_{1},y_{2})+W\mapsto \pi_{1}'(y_{1})+\pi_{2}'(y_{2}),
\]
die offenbar $\pi\circ \pi_{1}=\pi_{1}'$ und $\pi\circ \pi_{2}=\pi_{2}'$
erf"ullt. Ist $\bar{\pi}$ eine weitere solche Abbildung, so gilt offenbar
\[
\bar{\pi}((y_{1},y_{2})+W)=\bar{\pi}(\pi_{1}(y_{1})+\pi_{2}(y_{2}))=\pi_{1}'(y_{1})+\pi_{2}'(y_{2})=\pi((y_{1},y_{2})+W),
\]
also $\pi=\bar{\pi}$.

Sei nun $\varphi_{1}$ injektiv, und $y_{2}\in Y_{2}$ mit $\pi_{2}(y_{2})=0$,
d.h. $(0,y_{2})\in W$. Dann gibt es $x\in X$ mit $\varphi_{1}(x)=0$ und
$\varphi_{2}(x)=y_{2}$. Aus der Injektivit"at von $\varphi_{1}$ folgt $x=0$
und damit dann auch $y_{2}=\varphi_{2}(x)=0$, d.h. $\pi_{2}$ ist injektiv.

Sei nun $\varphi_{1}$ surjektiv. Zu $(y_{1},y_{2})+W\in Z$ (mit $y_{i}\in
Y_{i},\, i=1,2$) existiert dann ein $x\in X$ mit
$y_{1}=\varphi_{1}(x)$. Dann gilt
\begin{eqnarray*}
\pi_{2}(y_{2}+\varphi_{2}(x))&=&(0,y_{2}+\varphi_{2}(x))+W=(0,y_{2}+\varphi_{2}(x))+(\varphi_{1}(x),-\varphi_{2}(x))+W\\&=&(\varphi_{1}(x),y_{2})+W
=(y_{1},y_{2})+W,  
\end{eqnarray*}
also ist auch $\pi_{2}$ surjektiv. \qed\\

{\it Der Rest dieses Unterabschnitts wird nur f"ur eine alternative Konstruktion
  gebraucht und kann "ubersprungen werden.}\\

\begin{samepage}
\begin{Korollar}\label{PushoutExakteSeq}
Im folgenden Diagramm von $G$-Moduln
\begin{center}
~\begin{xy}
\xymatrix{
Y_{1} \ar[r]^{\varepsilon_{1}} \ar@{->>}[d]_{\varphi_{1}} &
Y_{2}\ar[r]^{\varepsilon_{2}} \ar@{.>>}[d]_{\varphi_{2}}&Y_{3}
\ar@{.>>}[d]_{\varphi_{3}}\\
X_{1} \ar@{.>}[r]_{\pi_{1}}& X_{2}\ar@{.>}[r]_{\pi_{2}}& X_{3}
}
\end{xy}
\end{center}
sei die obere Zeile exakt und $\varphi_{1}$ surjektiv. Weiter sei $X_{2}$ der
Pushout von $\varphi_{1}$ und $\varepsilon_{1}$ mit Homomorphismen $\pi_{1}$
und $\varphi_{2}$, sowie $X_{3}$ der
Pushout von $\varphi_{2}$ und $\varepsilon_{2}$ mit Homomorphismen $\pi_{2}$
und $\varphi_{3}$. Dann sind auch $\varphi_{2}, \varphi_{3}$ surjektiv, und
auch die untere Zeile ist exakt.
\end{Korollar}
\end{samepage}

\Bew Die Surjektivit"at von $\varphi_{2}$ und $\varphi_{3}$ folgt direkt aus
dem Satz. Aufgrund der Konstruktion ist das angegebene Diagramm jedenfalls
kommutativ. Daher gilt
\[
\pi_{2}\circ\pi_{1}\circ\varphi_{1}=\varphi_{3}\circ\underbrace{\varepsilon_{2}\circ\varepsilon_{1}}_{=0}=0,
\]
und wegen der Surjektivit"at von $\varphi_{1}$ folgt auch $\pi_{2}\circ
\pi_{1}=0$, also $\im \pi_{1}\subseteq \ker \pi_{2}$.

Sei umgekehrt $x_{2}\in X_{2}$ mit $\pi_{2}(x_{2})=0$. Nach Definition des
Pushouts ist $X_{3}=(X_{2}\oplus Y_{3})/W$ mit
$W=\left\{(\varphi_{2}(y_{2}),-\varepsilon_{2}(y_{2}))\in X_{2}\oplus
  Y_{3}: y_{2}\in Y_{2} \right\}$, und $\pi_{2}(x_{2})=(x_{2},0)+W\in
X_{3}$. Aus $\pi_{2}(x_{2})=0$ folgt daher $(x_{2},0)\in W$, d.h. es gibt
$y_{2}\in Y_{2}$ mit $x_{2}=\varphi_{2}(y_{2})$, $\varepsilon_{2}(y_{2})=0$. Mit
der Exaktheit der oberen Sequenz gibt es also wegen $y_{2}\in \ker
\varepsilon_{2}=\im \varepsilon_{1}$ ein $y_{1}\in Y_{1}$ mit
$y_{2}=\varepsilon_{1}(y_{1})$. Damit ist
$x_{2}=\varphi_{2}(y_{2})=\varphi_{2}(\varepsilon_{1}(y_{1}))=\pi_{1}(\varphi_{1}(y_{1}))\in
\im \pi_{1}$, also $\ker \pi_{2}\subseteq\im \pi_{1}.$

\qed 
\begin{BspRoman} \label{nichtexaktPushPull}
In diesem (Gegen-)Beispiel wollen wir zeigen, dass die Surjektivit"at von
$\varphi_{1}$ in Korollar \ref{PushoutExakteSeq} notwendig ist, damit die
induzierte Sequenz "uberhaupt ein Komplex ist. Sei $K^{n}$ jeweils ein
trivialer $KG$-Modul.
Wir betrachten mit den
Bezeichnungen des Korollars das folgende kommutative Diagramm
\begin{center}
~\begin{xy}
\xymatrix{
K \ar[r]^{\varepsilon_{1}} \ar@{->}[d]_{\varphi_{1}} &
K^{2}\ar[r]^{\varepsilon_{2}} \ar@{.>}[d]_{\varphi_{2}}&K
\ar@{.>}[d]_{\varphi_{3}}\\
K^{2} \ar@{.>}[r]_{\pi_{1}}& K^{3}\ar@{.>}[r]_{\pi_{2}}& K^{2}
}
\end{xy}
\end{center}
mit 
\[
\varepsilon_{1}(x_{1})=(x_{1},0),\quad \varepsilon_{2}(x_{1},x_{2})=x_{2}
\]
\[
\varphi_{1}(x_{1})=(x_{1},0),\quad
\varphi_{2}(x_{1},x_{2})=(x_{1},x_{2},0),\quad \varphi_{3}(x_{2})=(x_{2},0)
\]
\[
\pi_{1}(x_{1},x_{3})=(x_{1},0,x_{3}),\quad
\pi_{2}(x_{1},x_{2},x_{3})=(x_{2},x_{3})
\]
jeweils f"ur $x_{1},x_{2},x_{3}\in K$. Man beachte, dass hier
$\varphi_{1}$ nicht surjektiv ist. Dabei sind in der unteren Zeile die rechten
beiden Moduln Pushouts wie im Korollar beschrieben.  Jedoch ist $\pi_{2}\circ
\pi_{1}(x_{1},x_{3})=(0,x_{3})$, d.h. die untere Zeile ist nicht einmal ein Komplex.
\end{BspRoman}
\subsubsection{Pullback}
Wir kommen zu der zum Pushout "`dualen"' Konstruktion.
\begin{SatzDef}\label{DefPullback}
Seien $X_{1},X_{2},Y$ jeweils $G$-Moduln und $\varphi_{1}: X_{1}\rightarrow Y$ sowie
$\varphi_{2}: X_{2}\rightarrow Y$ jeweils $G$-Homomorphismen. Der
\emph{Pullback}\index{Pullback} von $\varphi_{1}$ und $\varphi_{2}$ ist der $G$-Modul
\[
Z:=\left\{(x_{1},x_{2})\in X_{1}\oplus X_{2}: \quad\varphi_{1}(x_{1})=\varphi_{2}(x_{2})\right\}
\]
zusammen mit den Homomorphismen
\[
\pi_{1}: Z\rightarrow X_{1}, \quad (x_{1},x_{2})\mapsto x_{1}
\]
und
\[
\pi_{2}: Z\rightarrow X_{2}, \quad (x_{1},x_{2})\mapsto x_{2}.
\]
Dann gilt $\varphi_{1}\circ\pi_{1}=\varphi_{2}\circ\pi_{2}$, 
\begin{center}
~\begin{xy}
\xymatrix{
Z \ar@{.>}[r]^{\pi_{2}} \ar@{.>}[d]_{\pi_{1}}&X_{2} \ar[d]^{\varphi_{2}}\\
X_{1} \ar[r]_{\varphi_{1}}&Y,
}
\end{xy}
\end{center}
und es gilt
folgende universelle Eigenschaft: F"ur jeden weiteren $G$-Modul $Z'$ und
Homomorphismen $\pi_{1}': Z'\rightarrow X_{1}$ und $\pi_{2}': Z'\rightarrow
X_{2}$ mit $\varphi_{1}\circ\pi_{1}'=\varphi_{2}\circ\pi_{2}'$ existiert genau
ein Homomorphismus $\pi: Z'\rightarrow Z$, der das folgende
\emph{Pullback-Diagramm}\index{Pullback!Diagramm} kommutativ macht:
\begin{center}
~\begin{xy}
\xymatrix{
Z' \ar@{.>}[rd]|-{\pi} \ar@/^/[rrd]^{\pi_{2}'} \ar@/_/[rdd]_{\pi_{1}'}&&\\
& Z \ar[r]^{\pi_{2}} \ar[d]_{\pi_{1}}&X_{2} \ar[d]^{\varphi_{2}}\\
&X_{1} \ar[r]_{\varphi_{1}}&Y
}
\end{xy}
\end{center}
Weiter gilt: Ist $\varphi_{2}$ injektiv/surjektiv, so ist auch  $\pi_{1}$ injektiv/surjektiv.
\end{SatzDef}

\Bew Offenbar ist $Z$ ein $G$-Untermodul von $X_{1}\oplus X_{2}$. 
Aufgrund der Definition von $Z$ gilt $\varphi_{1}\circ
\pi_{1}=\varphi_{2}\circ\pi_{2}$. Sind weiter $Z',\pi_{1}',\pi_{2}'$ wie
beschrieben, so ist
\[
\pi: Z'\rightarrow Z, \quad z' \mapsto (\pi_{1}'(z'),\pi_{2}'(z'))
\]
wohldefiniert, denn
$\varphi_{1}(\pi_{1}'(z'))=\varphi_{2}(\pi_{2}'(z'))$, also
$(\pi_{1}'(z'),\pi_{2}'(z'))\in Z$. Ferner gilt f"ur $z'\in Z'$ offenbar
$\pi_{1}(\pi(z'))=\pi_{1}((\pi_{1}'(z'),\pi_{2}'(z')))=\pi_{1}'(z')$, also
$\pi_{1}\circ \pi=\pi_{1}'$ und analog $\pi_{2}\circ \pi=\pi_{2}'$. Ist nun
$\bar{\pi}: Z'\rightarrow Z$ eine weitere Abbildung mit diesen Eigenschaften,
so gilt f"ur $z'\in Z'$ und $\bar{\pi}(z')=(x_{1},x_{2})$ mit $x_{i}\in X_{i}, \,
i=1,2$ dann $x_{1}=\pi_{1}((x_{1},x_{2}))=\pi_{1}(\bar{\pi}(z'))=\pi_{1}'(z')$ und
analog $x_{2}=\pi_{2}'(z')$. Also ist
$\bar{\pi}(z')=(\pi_{1}'(z'),\pi_{2}'(z'))=\pi(z')$, also $\bar{\pi}=\pi$.

Sei nun $\varphi_{2}$ injektiv. F"ur $z=(x_{1},x_{2})\in \ker \pi_{1}$ folgt
dann $0=\pi_{1}(z)=x_{1}$. Da $z\in Z$ gilt
$\varphi_{2}(x_{2})=\varphi_{1}(x_{1})=0$. Aus der Injektivit"at von
$\varphi_{2}$ folgt dann $x_{2}=0$, und damit $z=(x_{1},x_{2})=(0,0)=0\in
Z$. Also ist $\ker \pi_{1}=0$ und $\pi_{1}$ injektiv. 

Sei nun $\varphi_{2}$ surjektiv. Zu $x_{1}\in X_{1}$ gibt es dann $x_{2}\in
X_{2}$ mit $\varphi_{1}(x_{1})=\varphi_{2}(x_{2})$. Also ist $(x_{1},x_{2})\in
Z$ und $\pi_{1}((x_{1},x_{2}))=x_{1}$, d.h. $\pi_{1}$ ist surjektiv.\qed\\

{\it Der Rest dieses Unterabschnitts wird nur f"ur eine alternative Konstruktion
  gebraucht und kann "ubersprungen werden.}\\
\begin{Korollar}\label{PullbackExakteSeq}
In dem folgenden kommutativen Diagramm von $G$-Moduln
\begin{center}
 ~\begin{xy}
\xymatrix{
X_{1} \ar@{.>}[r]^{\pi_{1}} \ar@{^{(}.>}[d]_{\varphi_{1}} &
X_{2}\ar@{.>}[r]^{\pi_{2}} \ar@{^{(}.>}[d]_{\varphi_{2}} &
X_{3}\ar@{^{(}->}[d]^{\varphi_{3}}\\
Y_{1} \ar[r]_{\varepsilon_{1}} & Y_{2} \ar[r]_{\varepsilon_{2}} & Y_{3}
}
\end{xy}
\end{center}
sei die untere Zeile exakt und $\varphi_{3}$ injektiv. Weiter sei $X_{2}$ der
Pullback von $\varepsilon_{2}$ und $\varphi_{3}$ mit Homomorphismen
$\varphi_{2}$ und $\pi_{2}$, sowie $X_{1}$ der Pullback von $\varepsilon_{1}$
und $\varphi_{2}$ mit Homomorphismen $\varphi_{1}$ und $\pi_{1}$. Dann sind
auch $\varphi_{1}$ und $\varphi_{2}$ injektiv und die obere Zeile ist
ebenfalls exakt.
\end{Korollar}

\Bew Die Injektivit"at von $\varphi_{2}$ und $\varphi_{1}$ folgt direkt aus
dem Satz. Aufgrund der Konstruktion ist das Diagramm kommutativ. Damit ist
\[
\varphi_{3}\circ \pi_{2}\circ\pi_{1}=\underbrace{\varepsilon_{2}\circ \varepsilon_{1}}_{=0}
\circ \varphi_{1}=0,
\]
und wegen der Injektivit"at von $\varphi_{3}$ folgt auch
$\pi_{2}\circ\pi_{1}=0$, also $\im \pi_{1}\subseteq \ker \pi_{2}$. 

Sei nun
umgekehrt $x_{2}\in \ker \pi_{2}$. Da $x_{2}\in X_{2}=\left\{(y_{2},x_{3})\in
    Y_{2}\oplus X_{3}: \quad \varepsilon_{2}(y_{2})=\varphi_{3}(x_{3})
  \right\}$, gibt es also $y_{2}\in Y_{2},\, x_{3}\in X_{3}$ mit
  $x_{2}=(y_{2},x_{3})$ und $\varepsilon_{2}(y_{2})=\varphi_{3}(x_{3})$. Aus
  $0=\pi_{2}(x_{2})=\pi_{2}((y_{2},x_{3}))=x_{3}$ (Definition von $\pi_{2}$)
  folgt  $0=\varphi_{3}(x_{3})=\varepsilon_{2}(y_{2})$, also $y_{2}\in \ker
  \varepsilon_{2}=\im \varepsilon_{1}$ (Exaktheit der unteren Sequenz). Daher
  gibt es $y_{1}\in Y_{1}$ mit $y_{2}=\varepsilon_{1}(y_{1})$. Wir setzen
  $x_{1}:=(y_{1},x_{2})\in Y_{1}\oplus X_{2}$. Da
  $\varepsilon_{1}(y_{1})=y_{2}=\varphi_{2}((y_{2},x_{3}))=\varphi_{2}(x_{2})$,
  gilt sogar $x_{1}\in X_{1}$. Dann ist
  $\pi_{1}(x_{1})=\pi_{1}((y_{1},x_{2}))=x_{2}$, also $x_{2}\in \im \pi_{1}$
  und damit $\ker \pi_{2}\subseteq \im
  \pi_{1}$. \qed\\

\begin{BspRoman} \label{nichtexaktPullback}
Auch f"ur den Pullback wollen wir anhand eines Gegenbeispiels zeigen, dass die
Injektivit"at von $\varphi_{3}$ in obigem Korollar eine notwendige Voraussetzung ist. Es seien
jeweils wieder $K^{n}$ triviale $KG$-Moduln. Wir betrachten
\begin{center}
 ~\begin{xy}
\xymatrix{
K^{2} \ar@{.>}[r]^{\pi_{1}} \ar@{.>}[d]_{\varphi_{1}} &
K^{3}\ar@{.>}[r]^{\pi_{2}} \ar@{.>}[d]_{\varphi_{2}} &
K^{2}\ar@{->}[d]^{\varphi_{3}}\\
K \ar[r]_{\varepsilon_{1}} & K^{2} \ar[r]_{\varepsilon_{2}} & K
}
\end{xy}
\end{center}
mit 
\[
\varepsilon_{1}(x_{1})=(x_{1},0),\quad \varepsilon_{2}(x_{1},x_{2})=x_{2}
\]
\[
\varphi_{1}(x_{1},x_{3})=x_{1},\quad
\varphi_{2}(x_{1},x_{2},x_{3})=(x_{1},x_{2}),\quad \varphi_{3}(x_{2},x_{3})=x_{2}
\]
\[
\pi_{1}(x_{1},x_{3})=(x_{1},0,x_{3}),\quad
\pi_{2}(x_{1},x_{2},x_{3})=(x_{2},x_{3})
\]
jeweils f"ur $x_{1},x_{2},x_{3}\in K$. Man beachte, dass hier
$\varphi_{3}$ nicht injektiv ist. Dabei sind in der oberen Zeile die linken
beiden Moduln Pullbacks wie im Korollar beschrieben. Jedoch ist $\pi_{2}\circ
\pi_{1}(x_{1},x_{3})=(0,x_{3})$, d.h. die obere Zeile ist nicht einmal ein Komplex.
\end{BspRoman}

\subsection{Homomorphe Bilder von Kozyklen}
Seien $U,V$ jeweils $G$-Moduln,  $\varepsilon: U\rightarrow V$ ein
Homomorphismus von $G$-Moduln und $g\in Z^{n}(G,U)$ ein Kozyklus,
d.h. $\partial_{n}^{U}g=0$. Dann ist offenbar auch $\varepsilon\circ g\in
Z^{n}(G,V)$ ein Kozyklus, denn es gilt auch $\partial_{n}^{V}(\varepsilon\circ
g)=\varepsilon\circ\partial_{n}^{U}g=0$. Ein h"aufiger Fall ist hierbei, dass
$U\subseteq V$ ein Untermodul und $\varepsilon: U\hookrightarrow V$ die
Inklusion ist.

Wir wollen in diesem Abschnitt zeigen, wie
man aus einer exakten Sequenz, die $g\in Z^{n}(G,U)$ induziert, eine exakte
Sequenz machen kann, die $\varepsilon\circ g\in Z^{n}(G,V)$ induziert.

\begin{Satz}\label{HomomorpheBilderSatz}
Sei  $n\ge 1$,
\[ 
0\rightarrow
U_{n}\stackrel{\pi_{n}}{\rightarrow}U_{n-1}\stackrel{\pi_{n-1}}{\rightarrow}U_{n-2}\stackrel{\pi_{n-2}}{\rightarrow}\ldots\stackrel{\pi_{1}}{\rightarrow}U_{0}\stackrel{\pi_{0}}{\rightarrow}U_{-1}
\]
eine exakte Sequenz von $G$-Moduln und $w\in U_{-1}^{G} \cap \im\pi_{0}$. F"ur $k=0,\ldots,n$ sei $g_{k}\in C^{k}(G,U_{k})$ mit
\[
\pi_{0}(g_{0})=w\quad \textrm{ und }\quad \pi_{k}\circ
g_{k}=\partial_{k-1}^{U_{k-1}}g_{k-1} \quad \textrm{ f"ur }
k=1,\ldots,n \quad\quad (*)
\]
eine Folge, die $g_{n}\in Z^{n}(G,U_{n})$ gem"a"s Definition \ref{SatzVomIndKoz} induziert. 

Sei $U_{n}'$ ein weiterer $G$-Modul und $\varepsilon_{n}: U_{n}\rightarrow
U_{n}'$ ein Homomorphismus von $G$-Moduln.


Sei weiter $U_{n-1}'$ der Pushout von $\varepsilon_{n}$ und $\pi_{n}$ mit
Homomorphismen $\pi_{n}'$ und $\varepsilon_{n-1}$. Gem"a"s der universellen
Eigenschaft des Pushouts bilde man weiter die Abbildung $\pi_{n-1}':
U_{n-1}'\rightarrow U_{n-2}$, die folgendes Diagramm kommutativ macht (Details
im Beweis):


\begin{center}
~\begin{xy}
\xymatrix{
0 \ar[r]&U_{n}\ar@{^{(}->}[r]^{\pi_{n}}\ar@{->}[d]_{\varepsilon_{n}}
&U_{n-1}\ar[r]^{\pi_{n-1}}\ar@{.>}[d]_{\varepsilon_{n-1}}&U_{n-2}\ar[r]^{\pi_{n-2}}&\ldots\ar[r]^{\pi_{2}}&U_{1}\ar[r]^{\pi_{1}}&U_{0}\ar[r]^{\pi_{0}}&U_{-1}.\\
0 \ar@{.>}[r]&U_{n}'\ar@{^{(}.>}[r]^{\pi_{n}'} \ar@/_4pc/[rru]^{0} &U_{n-1}'\ar@{.>}[ru]_{\pi_{n-1}'}
}
\end{xy}
\end{center}
Dann ist die nach unten abknickende, zu $U_{n}'$ f"uhrende Sequenz exakt und
induziert (mit gleichem $w\in U_{-1}^{G}\cap \im\pi_{0}$) gem"a"s Definition \ref{SatzVomIndKoz} den Kozyklus $\varepsilon_{n}\circ
g_{n}\in Z^{n}(G,U_{n}')$.
\end{Satz}

\Bew
 Da 
$\pi_{n}$ injektiv ist, ist nach Satz \ref{DefPushout} auch $\pi_{n}'$
injektiv und damit die untere Sequenz an der Stelle $U_{n}'$ exakt. Ferner gilt mit der Nullabbildung $0:
U_{n}'\rightarrow \im(\pi_{n-1})$, dass $0\circ
\varepsilon_{n}=\pi_{n-1}\circ\pi_{n}$, denn die obere Sequenz ist exakt. Nach
der universellen Eigenschaft des Pushouts, angewendet auf die Abbildungen $0:
U_{n}'\rightarrow \im(\pi_{n-1})$ und $\pi_{n-1}: U_{n-1}\rightarrow \im (\pi_{n-1})$
(beachte die Einschr"ankung des Bildbereichs!) existiert also wie behauptet eine Abbildung
$\pi_{n-1}': U_{n-1}'\rightarrow \im(\pi_{n-1})\subseteq U_{n-2}$, die das Diagramm kommutativ
macht. Nach Definition gilt dann $\im(\pi_{n-1}')\subseteq \im
(\pi_{n-1})$. Aufgrund der Kommutativit"at des Diagramms gilt $\pi_{n-1}'\circ
\varepsilon_{n-1}=\pi_{n-1}$, und damit umgekehrt auch $\im(\pi_{n-1})\subseteq
\im(\pi_{n-1}')$, folglich $\im(\pi_{n-1}')=\im(\pi_{n-1})=\ker(\pi_{n-2})$.
Damit ist die nach unten abgebogene Sequenz an der Stelle $U_{n-2}$ exakt. 

Weiter folgt aus der Kommutativit"at des Diagramms $0=\pi_{n-1}'\circ \pi_{n}'$, also $\im
\pi_{n}'\subseteq \ker \pi_{n-1}'$. Nach Konstruktion des Pushouts als
Faktormodul $U_{n-1}'=(U_{n}'\oplus U_{n-1})/W$, wobei dann $\pi_{n}'$ und
$\varepsilon_{n-1}$ die entsprechenden Projektionen sind (siehe Definition
\ref{DefPushout}), gilt
$U_{n-1}'=\pi_{n}'(U_{n}')+\varepsilon_{n-1}(U_{n-1})$. Insbesondere gibt es
zu $x\in \ker \pi_{n-1}'\subseteq U_{n-1}'$ dann $x_{1}\in U_{n}'$ und $x_{2}\in
U_{n-1}$ mit $x=\pi_{n}'(x_{1})+\varepsilon_{n-1}(x_{2})$. Anwenden von $\pi_{n-1}'$ auf
diese Gleichung unter Beachtung der Kommutativit"at des Diagramms liefert
$0=0+\pi_{n-1}(x_{2})$, also $x_{2}\in \ker \pi_{n-1}=\im \pi_{n}$. Also gibt es $x_{3}\in
U_{n}$ mit $x_{2}=\pi_{n}(x_{3})$, und wir erhalten
$x=\pi_{n}'(x_{1})+\varepsilon_{n-1}(x_{2})=\pi_{n}'(x_{1})+\varepsilon_{n-1}(\pi_{n}(x_{3}))=\pi_{n}'(x_{1})+\pi_{n}'(\varepsilon_{n}(x_{3}))\in
\im \pi_{n}'$. Es gilt also auch $\ker \pi_{n-1}'\subseteq \im \pi_{n}'$ und damit
insgesamt $\ker \pi_{n-1}'=\im \pi_{n}'$. Also ist die nach unten
abgebogene Sequenz auch an der Stelle $U_{n-1}'$ exakt. Damit ist sie dann "uberall exakt.

Wir setzen nun $g_{n}':=\varepsilon_{n}\circ g_{n}\in C^{n}(G,U_{n}')$,
$g_{n-1}':=\varepsilon_{n-1}\circ g_{n-1}\in C^{n-1}(G,U_{n-1}')$ und
$g_{k}':=g_{k}\in C^{k}(G,U_{k})$ f"ur $k=0,\ldots,n-2$, und zeigen, dass die
Folge $(g_{k}')$ die entsprechende Eigenschaft $(*)$ aus Satz
\ref{SatzVomIndKoz} besitzt. Gem"a"s Definition induziert die Sequenz dann den
Kozyklus $g_{n}'=\varepsilon_{n}\circ g_{n}\in
Z^{n}(G,U_{n}')$ wie behauptet. Dabei ist $(*)$ offenbar nur noch f"ur $k=n$ und
$k=n-1$ zu pr"ufen. Wir verfizieren
\begin{eqnarray*}
\pi_{n}'\circ g_{n}'&=&\pi_{n}'\circ \varepsilon_{n}\circ
g_{n}=\varepsilon_{n-1}\circ\pi_{n}\circ g_{n}\stackrel{(*)}{=}\varepsilon_{n-1}\circ
\partial_{n-1}^{U_{n-1}}g_{n-1}\\&=&\partial_{n-1}^{U_{n-1}'}(\varepsilon_{n-1}\circ
g_{n-1})=\partial_{n-1}^{U_{n-1}'}g_{n-1}',
\end{eqnarray*}
also die $(*)$ entsprechende Eigenschaft f"ur die Folge $(g_{k}')$ und $k=n$, und ebenso
\[
\pi_{n-1}'\circ g_{n-1}'=\pi_{n-1}'\circ \varepsilon_{n-1}\circ
g_{n-1}=\pi_{n-1}\circ g_{n-1}\stackrel{(*)}{=}
\partial_{n-2}^{U_{n-2}}g_{n-2}=\partial_{n-2}^{U_{n-2}}g_{n-2}',
\]
also die $(*)$ entsprechende Eigenschaft f"ur $k=n-1$. F"ur $n=1$ und $k=0$
erh"alt man entsprechend $\pi_{0}'(g_{0}')=\pi_{0}(g_{0})=w$.  \qed\\

\subsection{Der generische $n$-Kozyklus}\label{AbschnGenKoz}
Wir betrachten nochmals die bar resolution
\[
\ldots\stackrel{d_{n+1}}{\rightarrow}
P_{n}\stackrel{d_{n}}{\rightarrow}P_{n-1}\stackrel{d_{n-1}}{\rightarrow}P_{n-2}\stackrel{d_{n-2}}{\rightarrow} \ldots\stackrel{d_{2}}{\rightarrow}P_{1}\stackrel{d_{1}}{\rightarrow}P_{0}\stackrel{d_{0}}{\rightarrow}K\rightarrow 0
\]
und verwenden die Bezeichnungen aus Abschnitt
\ref{sectiondiebarresolutin}. Aus den Definitionen von $d_{n}$
(S. \pageref{defdn}, ~\eqref{defdn}), $\partial_{n}^{P_{n}}$
(S. \pageref{defpartdn}, ~\eqref{defpartdn}) und $e_{n}$
(S. \pageref{defvonen}, \eqref{defvonen}) erhalten wir sofort
\begin{equation} \label{enKozyklus}
d_{0}(e_{0})=1\quad\textrm{ und }\quad d_{n}\circ e_{n}=\partial_{n-1}^{P_{n-1}}e_{n-1}\quad \myforall n\ge 1.  
\end{equation}
Wir betrachten nun 
\begin{equation}\label{defbarPn}
\overline{P_{n}}:=P_{n}/\im{d_{n+1}}\end{equation} mit der kanonischen
Projektion 
\begin{equation}\label{projpn}
p_{n}: P_{n}\rightarrow \overline{P_{n}}
\end{equation}
 sowie
\begin{equation}\label{defbaren}
\overline{e_{n}}:=p_{n}\circ e_{n}\in C^{n}(G,\overline{P_{n}}).
\end{equation}
Dann gilt
\[
\partial_{n}^{\overline{P_{n}}}\overline{e_{n}}=\partial_{n}^{\overline{P_{n}}}(p_{n}\circ
e_{n})=p_{n}\circ
\partial_{n}^{P_{n}}e_{n}\stackrel{~\eqref{enKozyklus}}{=}p_{n}\circ
d_{n+1}\circ e_{n}=0,
\]
(da $p_{n}\circ d_{n+1}=0$), also $\overline{e_{n}}\in Z^{n}(G,\overline{P_{n}})$. Wir nennen diesen
"`allgemeinsten"' Kozyklus den \emph{generischen $n$-Kozyklus}\index{Kozyklus!generischer}.

Da $\im d_{n+1}=\ker d_{n}$ existiert eine \emph{injektive} Abbildung 
\begin{equation}\label{defbardn}
\overline{d_{n}}:
\overline{P_{n}}\rightarrow P_{n-1}\quad \textrm{ mit }\quad \overline{d_{n}}\circ
p_{n}=d_{n},
\end{equation}
und wir erhalten so die \emph{generische exakte Sequenz des generischen
  $n$-Kozyklus}\index{exakte Sequenz!generische}
\begin{equation}\label{generischeexakteSeq}
0\rightarrow
\overline{P_{n}}\stackrel{\overline{d_{n}}}{\rightarrow}P_{n-1}\stackrel{d_{n-1}}{\rightarrow}P_{n-2}\stackrel{d_{n-2}}{\rightarrow} \ldots\stackrel{d_{2}}{\rightarrow}P_{1}\stackrel{d_{1}}{\rightarrow}P_{0}\stackrel{d_{0}}{\rightarrow}K\rightarrow 0.
\end{equation}
Wegen \eqref{enKozyklus} und 
\begin{equation}\label{dnenbardnbaren}
\overline{d_{n}}\circ \overline{e_{n}}\stackrel{~\eqref{defbaren}}{=}\overline{d_{n}}\circ
p_{n}\circ e_{n}\stackrel{~\eqref{defbardn}}{=}d_{n}\circ e_{n}\stackrel{~\eqref{enKozyklus}}{=}\partial_{n-1}^{P_{n-1}}e_{n-1}\end{equation} erf"ullt
die Folge $e_{0},\ldots,e_{n-1},\overline{e_{n}}$ die Voraussetzung $(*)$ von
Definition \ref{SatzVomIndKoz}, so dass $\overline{e_{n}}$ von der generischen
Sequenz induziert wird.

\begin{Satz}\label{GenericUniverselleEig}
Der generische $n$-Kozyklus $\overline{e_{n}}\in Z^{n}(G,\overline{P_{n}})$ besitzt
folgende universelle Eigenschaft: Zu jedem $G$-Modul $V$ und jedem $n$-Kozyklus
$g\in Z^{n}(G,V)$ existiert genau eine $KG$-lineare Abbildung $\phi_{g}^{n}:
\overline{P_{n}}\rightarrow V$ mit $g=\phi_{g}^{n}\circ \overline{e_{n}}$.
\end{Satz}

\Bew Wir betrachten gem"a"s der Konstruktion \eqref{deffvonomega},
S. \pageref{deffvonomega} die $KG$-lineare Abbildung
\[
\omega_{n}^{V}(g): P_{n}\rightarrow V,\quad\textrm{ mit }\quad
[\sigma_{1},\ldots,\sigma_{n}]\mapsto g(\sigma_{1},\ldots,\sigma_{n}).
\]
Nach Gleichung \eqref{omegaequation} gilt 
\[
\omega_{n}^{V}(g)\circ
d_{n+1}=\omega_{n+1}^{V}(\partial_{n}^{V}g)=\omega_{n+1}^{V}(0)=0,
\]
da $g\in Z^{n}(G,V)$. Also ist $\omega_{n}^{V}(g)|_{\im d_{n+1}}=0$, und damit
existiert eine Faktorisierung 
\begin{equation}\label{Faktorisierunphign}
\phi_{g}^{n}: \overline{P_{n}}=P_{n}/\im
d_{n+1}\rightarrow V\quad \textrm{mit}\quad
\omega_{n}^{V}(g)=\phi_{g}^{n}\circ p_{n}\end{equation} (siehe \eqref{projpn}).
Wir erhalten 
\[
\phi_{g}^{n}\circ \overline{e_{n}}\stackrel{~\eqref{defbaren}}{=}\phi_{g}^{n}\circ
p_{n}\circ
e_{n}\stackrel{~\eqref{Faktorisierunphign}}{=}\omega_{n}^{V}(g)\circ
e_{n}\stackrel{~\eqref{omegaTog}}{=}g
\]
wie gew"unscht. Die Eindeutigkeit folgt, weil $\overline{P_{n}}$ als $KG$-Modul von
$\overline{e_{n}}(G^{n})$ erzeugt wird.\qed\\

\subsection{Von Kozyklen zu exakten Sequenzen}\label{Kozykluszusequenz}
In Satz \ref{SatzVomIndKoz} haben wir gesehen, dass jeder exakten Sequenz der
Form \eqref{exakteSequenzdesKozyklus} ein bis auf einen Korand eindeutiger
Kozyklus zugeordnet werden kann. Wir wollen nun zeigen, dass jedem Kozyklus
auch eine ihn induzierende Sequenz der Form \eqref{exakteSequenzdesKozyklus}
zugeordnet werden kann.

\begin{Satz}\label{neuerKozyklusIndizierer}
Sei $n\ge 1$, $U_{n}$ ein $G$-Modul und $g_{n}\in Z^{n}(G,U_{n})$ ein $n$-Kozyklus. Dann
existiert eine exakte Sequenz der Form \eqref{exakteSequenzdesKozyklus}, die
$g_{n}$ induziert. Genauer kann man hierf"ur die nach unten abknickende Sequenz
des folgenden Diagramms w"ahlen,
\begin{center}
~\begin{xy}
\xymatrix{
0 \ar[r]&\overline{P_{n}}\ar@{^{(}->}[r]^{\overline{d_{n}}}\ar@{->}[d]_{\phi_{g_{n}}^{n}}
&P_{n-1}\ar[r]^{d_{n-1}}\ar@{.>}[d]_{\varepsilon_{n-1}}&P_{n-2}\ar[r]^{d_{n-2}}&\ldots\ar[r]^{d_{2}}&P_{1}\ar[r]^{d_{1}}&P_{0}\ar[r]^{d_{0}}&K\ar[r]&0.\\
0 \ar@{.>}[r]&U_{n}\ar@{^{(}.>}[r]^{d_{n}'} \ar@/_4pc/[rru]^{0} &U_{n-1}'\ar@{.>}[ru]_{d_{n-1}'}
}
\end{xy}
\end{center}
Hierbei ist die obere Zeile die generische exakte Sequenz des generischen
 $n$-Kozyklus \eqref{generischeexakteSeq},
 $\phi_{g_{n}}^{n}: \overline{P_{n}}\rightarrow U_{n}$ die $KG$-lineare Abbildung aus
 Satz \ref{GenericUniverselleEig}, und die untere Zeile wird mit Hilfe des
 Pushouts wie in Satz \ref{HomomorpheBilderSatz} gebildet.
\end{Satz}

\Bew Die obere Zeile induziert nach Abschnitt \ref{AbschnGenKoz} den
generischen $n$-Kozyklus $\overline{e_{n}}\in Z^{n}(G,\overline{P_{n}})$. Nach Satz
\ref{GenericUniverselleEig} gilt f"ur den $KG$-Homomorphismus $\phi_{g_{n}}^{n}:
\overline{P_{n}}\rightarrow U_{n}$ dann $g_{n}=\phi_{g_{n}}^{n}\circ \overline{e_{n}}$. Nach
Satz \ref{HomomorpheBilderSatz} induziert dann die nach unten abknickende
Sequenz den Kozyklus $g_{n}$. \qed\\

Dieser Satz liefert eine eher unhandliche exakte Sequenz zur"uck. Insbesondere
sind f"ur unendliches $G$ die $P_{k}$ unendlichdimensional, w"ahrend man den
Kozyklus dennoch oft durch eine Sequenz aus endlichdimensionalen Moduln
induzieren kann.


\begin{BspRoman}\label{KorandInduzierer}
Um eine induzierende Sequenz f"ur einen Korand $g_{n}\in
B^{n}(G,U_{n})$ anzugeben, k"onnen wir nach Satz \ref{SatzVomIndKoz} ohne Einschr"ankung
$g_{n}=0\in B^{n}(G,U_{n})$ annehmen (denn eine Sequenz induziert jeden
Repr"asentanten der Restklasse). Dann ist $\phi_{g_{n}}^{n}=0$, und nach
Konstruktion des Pushouts gilt dann $U_{n-1}'=(U_{n}\oplus
P_{n-1})/\{(0,-\overline{d_{n}}(x)): x\in \overline{P_{n}}\}\cong U_{n}\oplus
(P_{n-1}/\im d_{n})\stackrel{~\eqref{defbarPn}}{=} U_{n}\oplus
\overline{P_{n-1}}$. Es ist dann $d_{n}'$ durch die Einbettung
$U_{n}\hookrightarrow U_{n}\oplus
\overline{P_{n-1}}$ gegeben, und es
existiert ein "`links-Splitting"' $\mu_{n}: U_{n-1}'\rightarrow U_{n}$ mit
$\mu_{n}\circ d_{n}'=\id_{U_{n}}$, wie in einem der Kriterien von Satz
\ref{SatzVomIndKoz} f"ur das Vorliegen eines Korands gefordert. Auch wenn dieses
Kriterium nicht notwendig f"ur das Vorliegen eines
Korands ist (wie wir in Beispiel~\ref{Bsp2Kozyklus} gesehen haben), gibt es also
wenigstens eine induzierende Sequenz zu einem gegebenen Korand, die dieses Kriterium erf"ullt.
\end{BspRoman}


\begin{BspRoman}\label{einsKozyklenallgemein}
Wir untersuchen die Konstruktion im Fall $n=1$ noch etwas expliziter. Sei also
$g=g_{1}\in Z^{1}(G,U_{1})$. Wir haben also das Diagramm
\begin{center}
~\begin{xy}
\xymatrix{
0 \ar[r]&\overline{P_{1}}\ar@{^{(}->}[r]^{\overline{d_{1}}}\ar@{->}[d]_{\phi_{g}^{1}}
&KG\ar[r]^{d_{0}}\ar@{.>}[d]_{\varepsilon_{0}}&K\ar[r]&0.\\
0 \ar@{.>}[r]&U_{1}\ar@{^{(}.>}[r]^{d_{1}'} \ar@/_4pc/[rru]^{0} &U_{0}\ar@{.>}[ru]_{d_{0}'}
}
\end{xy}
\end{center}

Dann ist $U_{0}:=(U_{1}\oplus KG)/W$ mit 
\[
W=\left\langle\left\{(g(\sigma),-(\sigma-\iota))\in U_{1}\oplus KG:\quad
    \sigma\in G\right\} \right\rangle_{KG}.
\]
($\iota\in G$ das neutrale Element.)
Da 
\[
(\tau
g(\sigma),-\tau(\sigma-\iota))=(g(\tau\sigma),-(\tau\sigma-\iota))-(g(\tau),-(\tau-\iota))
\]
(aufgrund der Kozyklus-Eigenschaft), gilt sogar
\begin{equation}
  \label{WistKVR}
W=\left\langle\left\{(g(\sigma),-(\sigma-\iota))\in U_{1}\oplus KG:\quad
    \sigma\in G\right\} \right\rangle_{K}. 
\end{equation}

 Wir definieren
$\tilde{U_{1}}:=U_{1}\oplus K$ mit $G$-Operation 
\[
\sigma \cdot
(u,\lambda):=(\sigma u+\lambda g(\sigma), \lambda) \quad \textrm{f"ur } u\in
U_{1},\,\lambda\in K,\, \sigma\in G
\]
 (wie in Abschnitt
\ref{firstCohom}) und behaupten, dass durch
\[
F: \tilde{U_{1}}\rightarrow U_{0},\quad (u,\lambda)\mapsto (u,\lambda\cdot\iota)+W
\]
ein Isomorphismus von $KG$-Moduln gegeben ist. Die $K$-Linearit"at ist
klar. Ferner ist
\begin{eqnarray*}
F(\sigma \cdot
(u,\lambda))&=&F((\sigma u+\lambda g(\sigma), \lambda))=(\sigma u+\lambda
g(\sigma), \lambda\cdot\iota) +W\\&=&(\sigma u+\lambda
g(\sigma), \lambda\cdot\iota)-\lambda(g(\sigma),-(\sigma-\iota))+W=\\
&=&(\sigma u,\lambda\sigma)+W=\sigma F(u,\lambda) \quad\myforall u\in
U_{1},\lambda\in K,
\end{eqnarray*}
d.h. $F$ ist $KG$-linear. $F$ ist injektiv, denn f"ur $u\in
U_{1},\lambda\in K$ mit $F((u,\lambda))=0$ folgt $(u,\lambda\cdot\iota)\in W$. Nach \eqref{WistKVR}
gibt es dann $\sigma_{i}\in G,\lambda_{i}\in K,\, i=1,\ldots,n$ mit
\begin{equation}
  \label{darstulambda}
(u,\lambda\cdot\iota)=\sum_{i=1}^{n}(\lambda_{i}g(\sigma_{i}),-\lambda_{i}(\sigma_{i}-\iota)).  
\end{equation}

Dabei k"onnen wir O.E. $\sigma_{i}\ne \sigma_{j}$ f"ur $i\ne j$ annehmen. Da
$g(\iota)=g(\iota\iota)=\iota g(\iota)+g(\iota)$, also $g(\iota)=0$, k"onnen
wir au"serdem $\sigma_{i}\ne \iota\,\myforall i$ annehmen. Da in der zweiten
Komponente der linken Seite von \eqref{darstulambda} kein Term mit
$\sigma_{i}\ne \iota$ vorkommt, rechts aber $\lambda_{i}\sigma_{i}$, folgt
$\lambda_{i}=0\,\myforall i$. Mit \eqref{darstulambda} folgt $(u,\lambda)=0$,
also ist $F$ injektiv.

Zum Beweis der Surjektivit"at gen"ugt es wegen der $KG$-Linearit"at von $F$ ein
Urbild f"ur $(u,\lambda\sigma)+W$ ($u\in U_{1},\lambda\in K, \sigma\in G$)
anzugeben. Da
\begin{eqnarray*}
  F((u+\lambda g(\sigma),\lambda))&=&(u+\lambda
  g(\sigma),\lambda\iota)+W\\&=&(u+\lambda
  g(\sigma),\lambda\iota)-\lambda(g(\sigma),-(\sigma-\iota))+W\\
&=&(u,\lambda\sigma)+W,
\end{eqnarray*}
 ist $F$ also auch surjektiv.

Der Abbildung $d_{1}': U_{1}\rightarrow U_{0}, \, u\mapsto (u,0)+W$
entspricht via $F$ dann die Abbildung $\pi_{1}:
U_{1}\rightarrow \tilde{U_{1}},\, u\mapsto (u,0)$.

Nach dem Beweis zur universellen Eigenschaft des Pushouts (Satz
\ref{DefPushout}) gilt f"ur $d_{0}'$ ferner
\[
d_{0}'((u,\lambda\cdot\iota)+W)=0(u)+d_{0}(\lambda\cdot\iota)=\lambda \quad
\myforall u\in U_{1},\lambda\in K.
\]
Also entspricht $d_{0}'$ via $F$ der Abbildung $\pi_{0}:
  \tilde{U_{1}}\rightarrow K, \, (u,\lambda)\mapsto \lambda$, und wir erhalten
  die exakte Sequenz
  \[
 0\rightarrow
 U_{1}\stackrel{\pi_{1}}{\rightarrow}\tilde{U_{1}}\stackrel{\pi_{0}}{\rightarrow}K\rightarrow
 0.
 \]
 Damit entspricht unsere Konstruktion einer exakten Sequenz aus einem $1$-Kozyklus genau
 der im Text nach Proposition \ref{AnnulatorProp} angegebenen.
\end{BspRoman}

\subsubsection*{Konstruktion durch sukzessive Pushout-Bildung}
{\it In diesem erg"anzenden Abschnitt wollen wir eine weitere
  Konstruktion einer exakten Sequenz angeben, die einen
  vorgegebenen Kozyklus induziert. Diese funktioniert zwar nur unter
  Zusatzvoraussetzungen, f"uhrt aber im Allgemeinen zu "`kleineren"' Moduln in
  der Sequenz. Dieser Abschnitt kann "ubersprungen werden.\\}

Ist $U_{n}'\subseteq U_{n}$ ein Untermodul und $g_{n}\in Z^{n}(G,U_{n})$, so
dass $g_{n}(G^{n})\subseteq U_{n}'$, so ist offenbar auch $g_{n}\in
Z^{n}(G,U_{n}')$. Mit Hilfe von Satz \ref{neuerKozyklusIndizierer} kann man dann also auch eine
exakte Sequenz konstruieren, die $g_{n}\in Z^{n}(G,U_{n}')$ induziert.
Ist $g_{n}\in Z^{n}(G,U_{n})$ nichttrivial, so gilt dies
erst recht f"ur $g_{n}\in Z^{n}(G,U_{n}')$. Der kleinste Untermodul
$U_{n}'\subseteq U_{n}$ mit $g_{n}(G^{n})\subseteq U_{n}'$ ist offenbar
gegeben durch $\langle g_{n}(G^{n})\rangle_{KG}$, und wegen $\partial_{n}g_{n}=0$ und der Definition von
$\partial_{n}$ ist
\[
\langle g_{n}(G^{n})\rangle_{K}=\langle g_{n}(G^{n})\rangle_{KG},
\]
also der von den $g_{n}(\sigma_{1},\ldots,\sigma_{n}),\,\, \sigma_{1},\ldots,\sigma_{n}\in
G$ erzeugte $K$-Unterraum sogar ein $KG$-Untermodul. Ist also $g_{n}\in Z^{n}(G,U_{n})$ nichttrivial, so erst recht
$g_{n}\in Z^{n}(G,\langle g_{n}(G^{n})\rangle_{K})$. Es ist daher keine gro"se
Einschr"ankung, wenn wir $g_{n}\in Z^{n}(G,U_{n})$ \emph{nichttrivial} und
\begin{equation}
  \label{gnsurjektiv}
U_{n}=\langle g_{n}(G^{n})\rangle_{K}
\end{equation}
fordern. (Gegebenenfalls kann man mit Satz \ref{HomomorpheBilderSatz} dann
   nachtr"aglich $U_{n}$ in einen gr"o"seren Modul einbetten und erh"alt so auch
   f"ur den gr"o"seren Modul eine exakte Sequenz).
   Unter diesen Voraussetzungen wollen wir in diesem Abschnitt eine
alternative Konstruktion einer
exakten Sequenz \eqref{exakteSequenzdesKozyklus}
(S. \pageref{exakteSequenzdesKozyklus}) angeben, die $g_{n}$ induziert.

Wir verwenden wieder die Notation von Abschnitt \ref{sectiondiebarresolutin},
S. \pageref{sectiondiebarresolutin}, insbesondere ben"otigen wir die
numerierten Gleichungen.

Offenbar ist \eqref{gnsurjektiv} genau dann erf"ullt, wenn
\[
\varphi_{n}:=\omega_{n}^{U_{n}}(g_{n})\in \Hom_{KG}(P_{n},U_{n})
\]
(siehe \eqref{deffvonomega}, S. \pageref{deffvonomega}) surjektiv ist. Durch sukzessive Pushout-Bildung
erhalten wir aus der bar resolution \eqref{thebarresolution} (S. \pageref{thebarresolution}) und $\varphi_{n}$ das folgende
kommutative Diagramm:
\begin{center}
  ~\begin{xy}
\xymatrix{
P_{n+1}\ar[r]^{d_{n+1}}&P_{n}\ar[r]^{d_{n}}\ar@{->>}[d]_{\varphi_{n}}&P_{n-1}\ar[r]^{d_{n-1}}\ar@{.>>}[d]_{\varphi_{n-1}}&P_{n-2}\ar[r]^{d_{n-2}}\ar@{.>>}[d]_{\varphi_{n-2}}&\ldots\ar[r]^{d_{3}}&P_{2}\ar[r]^{d_{2}}\ar@{.>>}[d]_{\varphi_{2}}&P_{1}\ar[r]^{d_{1}}\ar@{.>>}[d]_{\varphi_{1}}&P_{0}\ar[r]^{d_{0}}\ar@{.>>}[d]_{\varphi_{0}}&K\ar[r]\ar@{.>>}[d]_{\varphi_{-1}}&0\ar@{.>>}[d]_{\varphi_{-2}}\\
0\ar[r]&U_{n}\ar@{.>}[r]_{\pi_{n}}&U_{n-1}\ar@{.>}[r]_{\pi_{n-1}}&U_{n-2}\ar@{.>}[r]_{\pi_{n-2}}&\ldots\ar@{.>}[r]_{\pi_{3}}&U_{2}\ar@{.>}[r]_{\pi_{2}}&U_{1}\ar@{.>}[r]_{\pi_{1}}&U_{0}\ar@{.>}[r]_{\pi_{0}}&U_{-1}\ar@{.>}[r]_{\pi_{-1}}&U_{-2}
}
\end{xy}
\end{center}
Hierbei ist f"ur $k=n-1,\ldots,-2$ jeweils $U_{k}$ der Pushout von
$\varphi_{k+1}$ und $d_{k+1}$ mit Homomorphismen $\pi_{k+1}$ und
$\varphi_{k}$.
\begin{Satz} \label{KozyklusLiefertSeq}
Sei $n\ge 1$ und $g_{n}\in Z^{n}(G,U_{n})$ ein \emph{nichttrivialer} $n$-Kozyklus, so dass
\[
U_{n}=\left\langle \left\{g_{n}(\sigma_{1},\ldots,\sigma_{n}):\quad
  \sigma_{1},\ldots,\sigma_{n}\in G\right\}\right\rangle_{K}.
\]
Dann ist die untere Zeile des obigen kommutativen Diagramms eine exakte Sequenz der Form
\begin{center}
  ~\begin{xy}
\xymatrix{
0\ar[r]&U_{n}\ar[r]^{\pi_{n}}&U_{n-1}\ar[r]^{\pi_{n-1}}&U_{n-2}\ar[r]^{\pi_{n-2}}&\ldots\ar[r]^{\pi_{3}}&U_{2}\ar[r]^{\pi_{2}}&U_{1}\ar[r]^{\pi_{1}}&U_{0}\ar[r]^{\pi_{0}}&K\ar[r]&0
}
\end{xy}
\end{center}
und induziert (gem"a"s Satz/Definition \ref{SatzVomIndKoz}) den Kozyklus $g_{n}$.
\end{Satz}

Da die $\varphi_{k}$ zwar surjektiv sind, in der Regel aber nicht injektiv
sein werden, sind die $U_{k}$ also als Faktormoduln der $P_{k}$ "`kleiner"'
als die Moduln der nach Satz \ref{neuerKozyklusIndizierer} konstruierten
Sequenz.\\

\Bew Da die obere Sequenz exakt und $\varphi_{n}$ nach Voraussetzung
surjektiv ist, ist nach Korollar \ref{PushoutExakteSeq} auch die untere Sequenz
(bis evtl. an der Stelle $U_{n}$) exakt und alle $\varphi_{k}$ ebenfalls
surjektiv. F"ur die Exaktheit an der Stelle $U_{n}$ zeigen wir, dass $\pi_{n}$
injektiv ist. Sei also $u\in U_{n}$ mit $\pi_{n}(u)=0\in U_{n-1}$. Nach
Konstruktion ist $U_{n-1}=(U_{n}\oplus P_{n-1})/W$ mit
$W=\left\{(\varphi_{n}(x),-d_{n}(x))\in U_{n}\oplus P_{n-1}:\quad x\in P_{n}
\right\}$, und es ist $\pi_{n}(u)=(u,0)+W$. Aus $\pi_{n}(u)=0$ folgt also
$(u,0)\in W$, d.h. es gibt ein $x\in P_{n}$ mit $u=\varphi_{n}(x)$ und
$0=d_{n}(x)$. Wegen der Exaktheit der oberen Sequenz und $x\in \ker d_{n}=\im
d_{n+1}$ gibt es $y\in P_{n+1}$ mit $x=d_{n+1}(y)$. Dann ist also
\[
u=\varphi_{n}(x)=\varphi_{n}\circ d_{n+1}(y)=\omega_{n}^{U_{n}}(g_{n})\circ
d_{n+1}(y)\stackrel{~\eqref{omegaequation}}{=}\omega_{n+1}^{U_{n}}(\partial_{n}^{U_{n}}g_{n})(y)=0,
\]
da $\partial_{n}^{U_{n}}g_{n}=0$ wegen $g_{n}\in Z^{n}(G,U_{n})$. Also ist
$\ker \pi_{n}=0$ und die untere Sequenz exakt.

Aus der Surjektivit"at von $\varphi_{-2}$ folgt au"serdem $U_{-2}=0$.

Mit $e_{k}$ wie in \eqref{defvonen} gilt 
\[
\varphi_{n}\circ
e_{n}=\omega_{n}^{U_{n}}(g_{n})\circ e_{n}\stackrel{~\eqref{omegaTog}}{=}g_{n},
\]
 und wir setzen daher auch f"ur $k=n-1,\ldots,0$
 \begin{equation}
   \label{gkdefphik}
g_{k}:=\varphi_{k}\circ e_{k} \in C^{k}(G,U_{k}).   
 \end{equation}
Dann gilt
\begin{equation}
  \label{omegagkphik}
\omega_{k}^{U_{k}}(g_{k})=\omega_{k}^{U_{k}}(\varphi_{k}\circ
e_{k})\stackrel{~\eqref{omegatophi}}{=}\varphi_{k} \quad \myforall k=0,\ldots,n.  \end{equation}
Mit der Kommutativit"at des Diagramms erhalten wir
\begin{eqnarray*}
\pi_{k}\circ g_{k}&\stackrel{~\eqref{gkdefphik}}{=}&\pi_{k}\circ\varphi_{k}\circ
e_{k}=\varphi_{k-1}\circ d_{k}\circ
e_{k}\stackrel{~\eqref{omegagkphik}}{=}\omega_{k-1}^{U_{k-1}}(g_{k-1})\circ
d_{k}\circ
e_{k}\\&\stackrel{~\eqref{omegaequation}}{=}&\omega_{k}^{U_{k-1}}(\partial_{k-1}^{U_{k-1}}g_{k-1})\circ
e_{k}\stackrel{~\eqref{omegaTog}}{=} \partial_{k-1}^{U_{k-1}}g_{k-1} \quad
\myforall k=1,\ldots,n,
\end{eqnarray*}
also
\begin{equation}
  \label{zweiterteilvonstern}
  \pi_{k}\circ g_{k}=\partial_{k-1}^{U_{k-1}}g_{k-1} \quad
\myforall k=1,\ldots,n,
\end{equation}
d.h. die $g_{k}$ erf"ullen den zweiten Teil der Bedingung $(*)$ von Satz
\ref{SatzVomIndKoz} (S. \pageref{SatzVomIndKoz}). Wenn wir nun noch
$U_{-1}\cong K$ und $\pi_{0}(g_{0})=1$ zeigen, so ist die untere Sequenz von der
behaupteten Form und induziert nach Satz \ref{SatzVomIndKoz} den Kozyklus
$g_{n}$. Wir zeigen $\pi_{0}(g_{0})\ne 0$. Dann ist insbesondere $U_{-1}\ne 0$
(da $\pi_{0}(g_{0})\in U_{-1}$), und aus der Surjektivit"at der $KG$-linearen Abbildung $\varphi_{-1}:
K\rightarrow U_{-1}$ folgt dann $U_{-1}\cong K$. Da man den Isomorphismus so
w"ahlen kann, dass $\pi_{0}(g_{0})\ne 0$ gerade $1\in K$ entspricht, sind wir
dann fertig.

Wir nehmen also stattdessen $\pi_{0}(g_{0})=0$ an, und zeigen f"ur
$k=1,\ldots,n$ die Existenz eines
\begin{equation}
  \label{hkGleichung}
h_{k}\in C^{k-1}(G,U_{k})\quad \textrm{ mit } \quad\pi_{k}\circ g_{k}=\pi_{k}\circ
\partial_{k-1}^{U_{k}}h_{k}.  
\end{equation}
Da $g_{0}\in \ker \pi_{0}=\im \pi_{1}$, gibt es ein $h_{1}\in U_{1}=C^{0}(G,U_{1})$
mit $g_{0}=\pi_{1}\circ h_{1}$.
Es folgt
\[
\pi_{1}\circ
g_{1}\stackrel{~\eqref{zweiterteilvonstern}}=\partial_{0}^{U_{0}}g_{0}=\partial_{0}^{U_{0}}(\pi_{1}\circ
h_{1})=\pi_{1}\circ \partial_{0}^{U_{1}}h_{1},
\]
also \eqref{hkGleichung} f"ur $k=1$.

Sei nun $h_{1},\ldots,h_{k}$ bereits konstruiert. Aus \eqref{hkGleichung} folgt
\[
\pi_{k}\circ(g_{k}-\partial_{k-1}^{U_{k}}h_{k})=0.
\]
Da $\ker \pi_{k}=\im \pi_{k+1}$, gibt es also $h_{k+1}\in C^{k}(G,U_{k+1})$
mit
\begin{equation}
  \label{nachsthk}
g_{k}-\partial_{k-1}^{U_{k}}h_{k}=\pi_{k+1}\circ h_{k+1}.
\end{equation}
Mit Hilfe von 
$\partial_{k}^{U_{k}}\circ \partial_{k-1}^{U_{k}}=0$ erhalten wir hieraus
\[
\pi_{k+1}\circ g_{k+1}\stackrel{~\eqref{zweiterteilvonstern}}{=}\partial_{k}^{U_{k}}g_{k}\stackrel{~\eqref{nachsthk}}{=}\partial_{k}^{U_{k}}(\pi_{k+1}\circ
h_{k+1})=\pi_{k+1}\circ \partial_{k}^{U_{k+1}}h_{k+1},
\]
also \eqref{hkGleichung} f"ur $k+1$.

Mit der Injektivit"at von $\pi_{n}$ erhalten wir aus
\eqref{hkGleichung} insbesondere f"ur $k=n$, dass
\[g_{n}=\partial_{n-1}^{U_{n}}h_{n}\in B^{n}(G,U_{n}),\] im Widerspruch zur
vorausgesetzten Nichttrivialit"at von $g_{n}$. Also ist $\pi_{0}(g_{0})\ne
0$. \qed\\

\begin{BspRoman}
Wir wollen noch kurz zeigen, dass auch diese Konstruktion im Fall $n=1$ zum
selben Ergebnis f"uhrt wie in Abschnitt \ref{firstCohom}. Dabei sehen wir auch
gut, wo die Surjektivit"at von $\varphi_{n}$ ins Spiel kommt. Wir haben dann also das
Diagramm
\begin{center}
  ~\begin{xy}
\xymatrix{
P_{2}\ar[r]^{d_{2}}&P_{1}\ar[r]^{d_{1}}\ar@{->>}[d]_{\varphi_{1}}&KG\ar[r]^{d_{0}}\ar@{.>>}[d]_{\varphi_{0}}&K\ar[r]\ar@{.>>}[d]_{\varphi_{-1}}&0\\
0\ar[r]&U_{1}\ar@{.>}[r]_{\pi_{1}}&U_{0}\ar@{.>}[r]_{\pi_{0}}&U_{-1}\ar@{.>}[r]&0.
}
\end{xy}
\end{center}
Wir k"onnen hier an Beispiel \ref{einsKozyklenallgemein} ankn"upfen (und
verwenden die dortigen Bezeichnungen): Die
Konstruktion von $U_{0}=(U_{1}\oplus KG)/W$ ist identisch, und es ist $d_{1}'=\pi_{1}$. Es bleibt
also noch die Konstruktion von $U_{-1}$ und $\pi_{0}$.
   Um $\pi_{0}$ explizit anzugeben, schreiben wir $(u,x)_{W}:=(u,x)+W$ f"ur $u\in
   U_{1},\,x\in KG$, und es ist $U_{-1}=(U_{0}\oplus K)/W_{0}$. In  $U_{-1}$ gilt 
  \begin{equation}
    \label{GleichunginUminus1}
  ((0,\sigma)_{W},0)\equiv
  ((0,0)_{W},1)\quad \textrm{f"ur }\sigma\in G.  
  \end{equation}
  Es ist  $\pi_{0}: U_{0}\rightarrow U_{-1},\, (u,x)_{W}\mapsto ((u,x)_{W},0)+W_{0}$. Nach Voraussetzung ist $U_{1}=\langle g(\sigma):\sigma\in
  G\rangle_{K}$. F"ur ein $u\in U_{1}$ gibt es also eine Darstellung
  $u=\sum_{i=1}^{n}\lambda_{i}g(\sigma_{i})$ mit $\lambda_{i}\in K,\sigma_{i}\in
  G\, \myforall i$ (hier geht also die Surjektivit"at von $\varphi_{1}$ ein!). Da $(g_{\sigma},0)_{W}=(0,\sigma-\iota)_{W}\myforall
  \sigma\in G$, ist also f"ur $\lambda\in K$
  \[
  (u,\lambda\iota)_{W}=\left(\sum_{i=1}^{n}\lambda_{i}g(\sigma_{i}),\lambda\iota\right)_{W}=\left(0,\lambda\iota+\sum_{i=1}^{n}\lambda_{i}(\sigma_{i}-\iota)\right)_{W}.
  \] 
  Daher ist
  \begin{eqnarray*}
  \pi_{0}((u,\lambda\iota)_{W})&=&\left(\left(0,\lambda\iota+\sum_{i=1}^{n}\lambda_{i}(\sigma_{i}-\iota)\right)_{W},0\right)+W_{0}\\&\stackrel{~\eqref{GleichunginUminus1}}{=}&((0,0)_{W},\lambda\cdot
  1+\sum_{i=1}^{n}\lambda_{i}\cdot(1-1)))+W_{0}\\&=&((0,0)_{W},\lambda_{})+W_{0}.
  \end{eqnarray*}
Also entspricht $\pi_{0}$ via $F$ der Abbildung $\tilde{U_{1}}\rightarrow K,
\, (u,\lambda)\mapsto \lambda$ wie erwartet.
\end{BspRoman}

\subsection{Induktion von Kozyklen durch Standardsequenzen}
Gem"a"s Satz \ref{neuerKozyklusIndizierer} k"onnen wir jeden Kozyklus durch
eine "`Standardsequenz"', also eine exakte Sequenz der Form \eqref{exakteSequenzdesKozyklus} (S. \pageref{exakteSequenzdesKozyklus})
 induzieren.
Dabei sind die $P_{k}$ jedoch oft "`unn"otig gro"s"'. Ist der Kozyklus bereits
durch
eine exakte Sequenz der Form \eqref{allgemeineKozyklusSeq}
(S. \pageref{allgemeineKozyklusSeq}) gegeben, so kann man
mit Hilfe des Pullbacks diese Sequenz auf eine Standardform bringen, die
denselben Kozyklus induziert. Wir verwenden dazu jeweils ein zu den Konstruktionen aus
Abschnitt \ref{Kozykluszusequenz} "`duales"' Verfahren. Beim ersten Verfahren
bilden wir einmal einen Pullback und h"angen diesen mittels dessen
universeller Eigenschaft an die alte Sequenz an. Bei der zweiten Methode (die
wieder nur unter einer Zusatzvoraussetzung funktioniert) bilden wir sukzessive
Pullbacks. Der Vorteil dieser Verfahren gegen"uber einer Anwendung von Satz
\ref{neuerKozyklusIndizierer} besteht darin, dass die Moduln in der gegebenen
Sequenz vermutlich einfachere Struktur haben als in der bar resolution, und
bei den in diesem Abschnitt vorgestellten Verfahren bleibt diese Struktur eher erhalten.
Besteht insbesondere eine gegebene Sequenz \eqref{allgemeineKozyklusSeq} nur aus
endlich-dimensionalen Moduln, so gilt dies auch f"ur die daraus in diesem Abschnitt konstruierten Sequenzen.\\  

\begin{Satz}
Ein $n$-Kozyklus ($n\ge 1$) werde durch eine exakte Sequenz der Form
\eqref{allgemeineKozyklusSeq} (S. \pageref{allgemeineKozyklusSeq})  gegeben
(mit vorgegebenem $w\in \im\pi_{0}\cap U_{-1}^{G}$). Sei $U_{0}'$ der Pullback
von
$\pi_{0}$ und $\varphi_{-1}: K\rightarrow U_{-1},\,\, \lambda\mapsto
\lambda\cdot w$ mit Homomorphismen $\varphi_{0}$ und $\pi_{0}'$, und $\pi_{1}'$
die nach der universellen Eigenschaft des Pullbacks existierende Abbildung,
die folgendes Diagramm kommutativ macht (Details im Beweis):
\begin{center}
~\begin{xy}
\xymatrix{
 &&&&&U_{0}'\ar@{.>}[r]^{\pi_{0}'}
\ar@{.>}[d]_{\varphi_{0}}&K\ar@{.>}[r]
\ar[d]_{\varphi_{-1}}&0\\
0 \ar[r]&U_{n}\ar[r]^<<<<<{\pi_{n}}&U_{n-1}\ar[r]^{\pi_{n-1}}&\ldots\ar[r]^{\pi_{2}}&U_{1}\ar@{.>}[ru]^{\pi_{1}'}\ar[r]^{\pi_{1}}\ar@/^4pc/[rru]^{0}&U_{0}\ar[r]^{\pi_{0}}&U_{-1}.
}
\end{xy}
\end{center}
Dann ist die nach oben abknickende "`Standardsequenz"' exakt und induziert denselben
Kozyklus wie die untere exakte Sequenz.
\end{Satz}

\Bew Mit der Nullabbildung $0: U_{1}\rightarrow K$ gilt  $\varphi_{-1}\circ
0=\pi_{0}\circ \pi_{1}$, und gem"a"s der universellen Eigenschaft des Pullback
gibt es eine Abbildung $\pi_{1}': U_{1}\rightarrow U_{0}'$, die das
Diagramm kommutativ macht. Wir zeigen zun"achst, dass die nach oben
abknickende Sequenz exakt ist. Nach Definition des Pullbacks ist
\[
U_{0}'=\{(u_{0},\lambda)\in U_{0}\oplus K:\quad
\pi_{0}(u_{0})=\varphi_{-1}(\lambda)=\lambda\cdot w\},
\]
und $\varphi_{0}$ bzw. $\pi_{0}'$ sind dann die Projektionen auf die erste
bzw. zweite Koordinate. Nach Voraussetzung gibt es ein $u_{0}\in U_{0}$ mit
$\pi_{0}(u_{0})=w$. Dann ist $u_{0}':=(u_{0},1)\in U_{0}'$ und
$\pi_{0}'(u_{0}')=1$, also $\pi_{0}'$ surjektiv und die abknickende Sequenz an
der Stelle $K$ exakt. Aus der Kommutativit"at des Diagramms folgt
$0=\pi_{0}'\circ \pi_{1}'$, also $\im\pi_{1}'\subseteq \ker \pi_{0}'$.
Ist umgekehrt $(u,\lambda)\in \ker \pi_{0}'\subseteq
U_{0}'$, also $\lambda=0$, so folgt aus $\pi_{0}(u)=\lambda\cdot w=0$ also
$u\in \ker\pi_{0}=\im\pi_{1}$. Also gibt es $u_{1}\in U_{1}$ mit
$\pi_{1}(u_{1})=u$. Nach Konstruktion der Abbildung $\pi_{1}'$ (siehe der
Beweis zu Satz \ref{DefPullback}) gilt dann
$\pi_{1}'(u_{1})=(\pi_{1}(u_{1}),0(u_{1}))=(u,0)=(u,\lambda)$, also auch
$\ker\pi_{0}'\subseteq \im\pi_{1}'$, und wir haben Exaktheit an der Stelle
$U_{0}'$. Weiter gilt f"ur $u_{2}\in U_{2}$ stets $\pi_{1}'\circ
\pi_{2}(u_{2})=(\pi_{1}(\pi_{2}(u_{2})),0(\pi_{2}(u_{2})))=(0,0)$ (da
$\pi_{1}\circ\pi_{2}=0$) und damit $\im\pi_{2}\subseteq\ker\pi_{1}'$. Ist
umgekehrt $u_{1}\in \ker\pi_{1}'$, also $(\pi_{1}(u_{1}),0)=(0,0)$, so folgt
$u_{1}\in\ker\pi_{1}=\im\pi_{2}$ und somit auch
$\ker\pi_{1}'\subseteq\im\pi_{2}$. Dies zeigt Exaktheit an der Stelle $U_{1}$,
und damit ist die Sequenz insgesamt exakt. 

Wir m"ussen noch zeigen, dass beide Sequenzen den gleichen Kozyklus
induzieren. Sei dazu $g_{0}'\in U_{0}'$ mit $\pi_{0}'(g_{0}')=1$, $g_{1}\in
C^{1}(G,U_{1})$ mit $\pi_{1}'\circ g_{1}=\partial_{0}^{U_{0}'}g_{0}'$ und f"ur
$k=2,\ldots,n$ jeweils $g_{k}\in C^{k}(G,U_{k})$ mit $\pi_{k}\circ
g_{k}=\partial_{k-1}^{U_{k-1}}g_{k-1}$, d.h. die angegebene Folge erf"ullt die
Eigenschaft $(*)$ aus Satz \ref{SatzVomIndKoz}. Dann induziert die nach oben abknickende Sequenz also
gem"a"s Satz \ref{SatzVomIndKoz} den Kozyklus $g_{n}\in Z^{n}(G,U_{n})$. Wir
setzen nun $g_{0}:=\varphi_{0}(g_{0}')\in U_{0}$ und zeigen, dass die Folge der
$g_{0},\ldots,g_{n}$ ebenfalls die Eigenschaft $(*)$ aus Satz
\ref{SatzVomIndKoz} erf"ullt. Es gilt
$\pi_{0}(g_{0})=\pi_{0}(\varphi_{0}(g_{0}'))=\varphi_{-1}(\pi_{0}'(g_{0}'))=\varphi_{-1}(1)=w$
und $\pi_{1}\circ g_{1}=\varphi_{0}\circ \pi_{1}'\circ g_{1}=\varphi_{0}\circ
\partial_{0}^{U_{0}'}g_{0}'=\partial_{0}^{U_{0}}(\varphi_{0}(g_{0}'))=\partial_{0}^{U_{0}}g_{0}$.
F"ur die anderen Werte von $k$ gilt $(*)$ sowieso. Also induziert auch die
urspr"ungliche Sequenz den Kozyklus $g_{n}\in Z^{n}(G,U_{n})$, und damit beide
Sequenzen den gleichen Kozyklus.\qed\\

\subsubsection*{Konstruktion durch sukzessive Pullback-Bildung}
{\it Diese alternative Konstruktion kann wieder "ubersprungen werden.}

Wir wollen hier die exakte Standardsequenz mit Hilfe von Korollar
 \ref{PullbackExakteSeq} konstruieren. 
 In der Situation von Satz
\ref{SatzVomIndKoz} ist dabei $\varphi_{-1}: K\rightarrow U_{-1},\quad
\lambda\mapsto \lambda\cdot w$ f"ur $w\ne 0$ ein injektiver
$KG$-Homomorphismus, und dann k"onnen wir das Korollar anwenden. Durch sukzessive Pullbackbildung erh"alt man so das
kommutative Diagramm

\begin{center}
~\begin{xy}
\xymatrix{
0 \ar@{.>}[r]^{\varepsilon_{n+1}} \ar@{.>}[d]_{\varphi_{n+1}} &
V_{n}\ar@{.>}[r]^{\varepsilon_{n}} \ar@{.>}[d]_{\varphi_{n}}&V_{n-1}\ar@{.>}[d]_{\varphi_{n-1}}\ar@{.>}[r]^{\varepsilon_{n-1}}
&\ldots\ar@{.>}[r]^{\varepsilon_{2}}&V_{1}\ar@{.>}[r]^{\varepsilon_{1}}
\ar@{.>}[d]_{\varphi_{1}}&V_{0}\ar@{.>}[r]^{\varepsilon_{0}}
\ar@{.>}[d]_{\varphi_{0}}&K\ar@{.>}[r]
\ar[d]_{\varphi_{-1}}&0\\
0 \ar[r]_{\pi_{n+1}}&U_{n}\ar[r]_{\pi_{n}}&U_{n-1}\ar[r]_{\pi_{n-1}}&\ldots\ar[r]_{\pi_{2}}&U_{1}\ar[r]_{\pi_{1}}&U_{0}\ar[r]_{\pi_{0}}&U_{-1}.
}
\end{xy}
\end{center}
Dabei ist $V_{i}$ f"ur $i=0,\ldots,n+1$ jeweils Pullback von $\pi_{i}$ und $\varphi_{i-1}$ mit
Homomorphismen $\varphi_{i}$ und $\varepsilon_{i}$. Nach dem Korollar sind
dann mit $\varphi_{-1}$ auch die $\varphi_{i}$ mit $i=0,\ldots,n+1$
injektiv und die obere Sequenz exakt (Surjektivit"at von $\varepsilon_{0}$
zeigen wir gleich). Aus der Injektivit"at von
$\varphi_{n+1}: V_{n+1}\rightarrow 0$ folgt $V_{n+1}=0$. Sei $g_{k}\in
C^{k}(G,U_{k})$ mit $(*)$ wie in Satz \ref{SatzVomIndKoz}, also
\[
\pi_{0}(g_{0})=w\quad \textrm{ und }\quad \pi_{k}\circ
g_{k}=\partial_{k-1}^{U_{k-1}}g_{k-1} \quad \textrm{ f"ur }
k=1,\ldots,n. \quad (*)
\]
Wir konstruieren $h_{k}\in C^{k}(G,V_{k})$, $\,k=0,\ldots,n$ mit der $(*)$
entsprechenden Eigenschaft 
\[
\varepsilon_{0}(h_{0})=1\quad \textrm{ und }\quad \varepsilon_{k}\circ
h_{k}=\partial_{k-1}^{V_{k-1}}h_{k-1} \quad \textrm{ f"ur }
k=1,\ldots,n, \quad (**)
\]
sowie \[\varphi_{k}\circ h_{k}=g_{k}. \quad (***)\] Wir setzen
zun"achst $h_{0}:=(g_{0},1)\in U_{0}\oplus K$. Dann ist
$\pi_{0}(g_{0})=w=\varphi_{-1}(1)$, nach Definition des Pullbacks also sogar $h_{0}\in V_{0}=C^{0}(G,V_{0})$ sowie
$\varepsilon_{0}(h_{0})=1$ und $\varphi_{0}\circ h_{0}=g_{0}$. Wegen
$\varepsilon_{0}(h_{0})=1$ ist dann insbesondere auch $\varepsilon_{0}$
surjektiv, d.h. die obere exakte Sequenz ist tats"achlich von der angegebenen
Form. Sei nun $h_{0},\ldots,h_{k-1}$ mit $(**)$ und $(***)$ bereits konstruiert. Wir setzen dann
\[
h_{k}:=(g_{k},\partial_{k-1}^{V_{k-1}}h_{k-1})\in C^{k}(G,U_{k}\oplus
V_{k-1}).
\]
Da \[
\varphi_{k-1}\circ\partial_{k-1}^{V_{k-1}}h_{k-1}
=\partial_{k-1}^{U_{k-1}}(\varphi_{k-1}\circ
h_{k-1})\stackrel{(***)}{=}\partial_{k-1}^{U_{k-1}}g_{k-1}\stackrel{(*)}{=}\pi_{k}\circ g_{k},
\]
ist sogar $h_{k}\in C^{k}(G,V_{k})$. Weiter gilt $\varepsilon_{k}\circ
h_{k}=\partial_{k-1}^{V_{k-1}}h_{k-1}$ und $\varphi_{k}\circ h_{k}=g_{k}$ nach
Definition der Homomorphismen des Pullbacks. Dies sind $(**)$ und $(***)$ f"ur
$k$.

Insbesondere gilt dann nach Satz \ref{SatzVomIndKoz} $h_{n}\in Z^{n}(G,V_{n})$,
und es ist $g_{n}=\varphi_{n}\circ h_{n}$ mit einer injektiven Abbildung
$\varphi_{n}$. Falls $\langle g_{n}(G^{n})\rangle_{KG}=U_{n}$, so ist
$\varphi_{n}$ sogar surjektiv, also $U_{n}\cong V_{n}$. 

In diesem Fall
induziert die obere Sequenz also (bis auf Isomorphie) den Kozyklus $g_{n}$.

Ist $\varphi_{n}$ nicht surjektiv, so kann man noch Satz
\ref{HomomorpheBilderSatz} anwenden, um eine Sequenz zu erhalten, die mit
$U_{n}$ endet und $g_{n}=\varphi_{n}\circ h_{n}\in Z^{n}(G,U_{n})$ induziert.

\subsection{"Aquivalenz von Sequenzen}
{\it Dieser Abschnitt dient nur der Vollst"andigkeit und kann "ubersprungen werden.}\\

Ein Element $g_{n}\in H^{n}(G,U_{n})$ kann offenbar durch viele verschiedene,
auf $U_{n}$ endende exakte Sequenzen der L"ange $n$ in Standardform
\eqref{exakteSequenzdesKozyklus} induziert werden. 

\begin{Def}
Wir nennen zwei exakte
Sequenzen  der L"ange $n\ge 1$ in Standardform
\begin{equation}\label{dieersteSeq}
E_{1}:\quad 0\rightarrow
U_{n}\stackrel{\pi_{n}}{\rightarrow}U_{n-1}\stackrel{\pi_{n-1}}{\rightarrow}U_{n-2}\stackrel{\pi_{n-2}}{\rightarrow}\ldots\stackrel{\pi_{1}}{\rightarrow}U_{0}\stackrel{\pi_{0}}{\rightarrow}K\rightarrow
0
\end{equation}
und
\begin{equation}\label{diezweiteSeq}
E_{2}:\quad 0\rightarrow
U_{n}\stackrel{\pi_{n}'}{\rightarrow}U_{n-1}'\stackrel{\pi_{n-1}'}{\rightarrow}U_{n-2}'\stackrel{\pi_{n-2}'}{\rightarrow}\ldots\stackrel{\pi_{1}'}{\rightarrow}U_{0}'\stackrel{\pi_{0}'}{\rightarrow}K\rightarrow
0
\end{equation}
 mit gleichem Ende $U_{n}$ \emph{"aquivalent}, in Zeichen $E_{1}\sim E_{2}$, \index{exakte Sequenz!\"aquivalente}
 wenn sie die gleiche Restklasse $g_{n}+B^{n}(G,U_{n})\in H^{n}(G,U_{n})$
 induzieren. 
\end{Def}
Die so definierte "Aquivalenz ist offenbar symmetrisch, reflexiv und
transitiv. Obwohl die Gesamtheit der zugrundeliegenden Objekte keine Menge
sondern eine "`Klasse"' bildet, wollen wir bei einer Relation mit diesen drei
Eigenschaften von einer \emph{"Aquivalenzrelation} sprechen. 

Ziel dieses Abschnitts ist es, die "Aquivalenz von Sequenzen in Standardform
mit gleichem Ende $U_{n}$ auf eine weitere Art zu charakterisieren.

\begin{Def}\label{YonedaPraaq}
Zwei exakte Sequenzen \eqref{dieersteSeq}, \eqref{diezweiteSeq} der L"ange
$n\ge 1$ in Standardform mit gleichem Ende $U_{n}$
hei"sen \emph{Yoneda-pr"a"aquivalent}, in Zeichen $E_{1}\sim_{pY}E_{2}$, \index{exakte Sequenz!Yoneda-pr\"a\"aquivalent}, wenn es f"ur $k=0,\ldots,n-1$ $KG$-Homomorphismen $\varphi_{k}:
U_{k}\rightarrow U_{k}'$ gibt, so dass folgendes Diagramm kommutativ wird:

\begin{center}
~\begin{xy}
\xymatrix{
E_{1}:&0 \ar[r]  &
U_{n}\ar[r]^{\pi_{n}} \ar[d]_{\id_{U_{n}}}&U_{n-1}\ar[d]_{\varphi_{n-1}}\ar[r]^{\pi_{n-1}}
&\ldots\ar[r]^{\pi_{2}}&U_{1}\ar[r]^{\pi_{1}}
\ar[d]_{\varphi_{1}}&U_{0}\ar[r]^{\pi_{0}}
\ar[d]_{\varphi_{0}}&K\ar[r]
\ar[d]_{\id_{K}}&0\\
E_{2}:&0 \ar[r]&U_{n}\ar[r]_{\pi_{n}'}&U_{n-1}'\ar[r]_{\pi_{n-1}'}&\ldots\ar[r]_{\pi_{2}'}&U_{1}'\ar[r]_{\pi_{1}'}&U_{0}'\ar[r]_{\pi_{0}'}&K\ar[r]&0.
}
\end{xy}
\end{center}
\end{Def}

Die Yoneda-Pr"a"aquivalenz ist offenbar reflexiv und transitiv, im Allgemeinen
aber nicht symmetrisch. 

\begin{Lemma}\label{LemmaYonedaPraeq}
Sind zwei exakte Sequenzen \eqref{dieersteSeq}, \eqref{diezweiteSeq} der
L"ange $n\ge 1$ in Standardform mit gleichem Ende $U_{n}$ Yoneda-pr"a"aquivalent,
also $E_{1}\sim_{pY}E_{2}$, so sind sie auch "aquivalent, also $E_{1}\sim E_{2}$.
\end{Lemma}

\Bew Wir verwenden die Bezeichnungen von Definition \ref{YonedaPraaq}. Sei
$g_{k}\in C^{k}(G,U_{k})$ f"ur $k=0,\ldots,n$ eine Folge mit $(*)$ aus
Definition \ref{SatzVomIndKoz}, also 
\[
\pi_{0}(g_{0})=1\quad\textrm{und}\quad\pi_{k}\circ
g_{k}=\partial_{k-1}^{U_{k-1}}g_{k-1} \quad \textrm{ f"ur } k=1,\ldots, n,
\quad (*)\] d.h. $E_{1}$
induziert $g_{n}\in Z^{n}(G,U_{n})$. Wir setzen
zus"atzlich $\varphi_{n}:=\id_{U_{n}}$ und
$g_{k}':=\varphi_{k}\circ g_{k}\in C^{k}(G,U_{k}')$ f"ur $k=0,\ldots,n$ und
zeigen, dass die Folge der $(g_{k}')$ die $(*)$ entsprechende Eigenschaft
erf"ullt. Dann induziert $E_{2}$ ebenfalls den Kozyklus $g_{n}'=g_{n}\in
Z^{n}(G,U_{n})$, und damit gilt dann $E_{1}\sim E_{2}$.  Aus der
Kommutativit"at des Diagramms folgt
\[\pi_{0}'(g_{0}')=\pi_{0}'(\varphi_{0}(g_{0}))=\id_{K}(\pi_{0}(g_{0}))\stackrel{(*)}{=}1\]
sowie 
\begin{eqnarray*}
\pi_{k}'\circ g_{k}'&=&\pi_{k}'\circ \varphi_{k}\circ
g_{k}=\varphi_{k-1}\circ \pi_{k}\circ
g_{k}\stackrel{(*)}{=}\varphi_{k-1}\circ\partial_{k-1}^{U_{k-1}}g_{k-1}\\&=&\partial_{k-1}^{U_{k-1}'}(\varphi_{k-1}\circ
g_{k-1})=\partial_{k-1}^{U_{k-1}'}g_{k-1}' \quad \myforall k=1,\ldots,n,
\end{eqnarray*}
und dies ist $(*)$ f"ur die Folge $(g_{k}')$. \qed\\

\begin{Def}\label{DefYonedeAq}
Zwei exakte Sequenzen $E, E'$ in Standardform mit gleichem Ende $U_{n}$ ($n\ge
1$)
hei"sen \emph{Yoneda-"aquivalent}, in Zeichen $E\sim_{Y} E'$,
\index{exakte Sequenz!Yoneda-\"aquivalent} wenn es eine Zahl $m\ge 0$ und exakte
Sequenzen $E_{1},\ldots, E_{m}$ der L"ange $n$ in Standardform mit gleichem
Ende $U_{n}$ gibt, so dass f"ur die Folge $E_{0}:=E, E_{1},E_{2},\ldots, E_{m},
E_{m+1}:=E'$ gilt: F"ur jedes $k=0,\ldots,m$ gilt $E_{k}\sim_{pY} E_{k+1}$
oder $E_{k+1}\sim_{pY} E_{k}$, d.h. je zwei benachbarte Glieder sind evtl. nach
Vertauschen Yoneda-pr"a"aquivalent. 
\end{Def}

\begin{Satz}\label{SatzVonderYoneadAquiv}
Zwei exakte Sequenzen $E, E'$ in Standardform mit gleichem Ende $U_{n}$ ($n\ge
1$) sind genau dann "aquivalent, wenn
sie Yoneda-"aquivalent sind, also
\[
E\sim E' \quad \Leftrightarrow \quad E\sim_{Y}E'.
\]
Insbesondere ist die Yoneda-"Aquivalenz eine "Aquivalenzrelation.
\end{Satz}

\Bew Sei zun"achst $E\sim_{Y} E'$, d.h. es gibt exakte Sequenzen
$E_{1},\ldots, E_{m}$ wie in De\-finition~\ref{DefYonedeAq}. Wir setzen wieder
$E_{0}:=E,\,\, E_{m+1}:=E'$. Nach Lemma
\ref{LemmaYonedaPraeq} und weil $\sim$ symmetrisch ist, gilt dann $E_{k}\sim E_{k+1}$ f"ur $k=0,\ldots,m$, und
da $\sim$ eine "Aquivalenzrelation ist gilt dann auch $E\sim E'$.

Sei umgekehrt $E\sim E'$, d.h. $E$ und $E'$ induzieren dieselbe Restklasse
$g_{n}+B^{n}(G,U_{n})\in H^{n}(G,U_{n})$. Nach Satz \ref{SatzVomIndKoz} induzieren beide
Sequenzen dann auch jeden Repr"asentanten dieser Restklasse, d.h. beide
Sequenzen induzieren insbesondere auch $g_{n}\in Z^{n}(G,U_{n})$.
 Im
n"achsten Lemma zeigen wir, dass es eine nur von $U_{n}$ und $g_{n}\in Z^{n}(G,U_{n})$ abh"angige Sequenz $E_{1}$ gibt mit
$E_{1}\sim_{pY}E$ und $E_{1}\sim_{pY}E'$. Nach Definition~\ref{DefYonedeAq}
(mit $m=1$) gilt dann $E\sim_{Y} E'$. \qed\\

\begin{Lemma}
Sei $n\ge 1$ und $E$ eine exakte Sequenz der L"ange $n$ in Standardform mit
Ende $U_{n}$, die den Kozyklus $g_{n}\in
Z^{n}(G,U_{n})$ induziert. Sei weiter $E'$ die in Satz
\ref{neuerKozyklusIndizierer} aus $g_{n}$ und $U_{n}$ konstruierte nach unten
abknickende exakte Sequenz, die ebenfalls $g_{n}$ induziert. Dann gilt
$E'\sim_{pY}E$.
\end{Lemma}

\Bew Sei $E$ gegeben wie in Gleichung \eqref{dieersteSeq}. Wir betrachten das
folgende Diagramm, dessen oberen beiden Zeilen sich gem"a"s der Konstruktion
aus Satz \ref{neuerKozyklusIndizierer} aus $g_{n}$ und $U_{n}$ ergeben. Die
dritte Zeile ist die gegebene Sequenz $E$.

\begin{center}
~\begin{xy}
\xymatrix{
&0 \ar[r]&\overline{P_{n}}\ar@{^{(}->}[r]^{\overline{d_{n}}}\ar@{->}[d]_{\phi_{g_{n}}^{n}}
&P_{n-1} \ar@/^2pc/[dd]^>>>>>>{\varphi_{n-1}}\ar[r]^{d_{n-1}}\ar@{.>}[d]_{\varepsilon_{n-1}}&P_{n-2}\ar[dd]^{\varphi_{n-2}}\ar[r]^{d_{n-2}}&\ldots\ar[r]^{d_{2}}&P_{1}\ar[dd]_{\varphi_{1}}\ar[r]^{d_{1}}&P_{0}\ar[dd]_{\varphi_{0}}\ar[r]^{d_{0}}&K\ar[r]\ar[dd]_{\id_{K}}&0\\
E':&0 \ar@{.>}[r]&U_{n}\ar@{^{(}.>}[r]^{d_{n}'}\ar[d]_{\id_{U_{n}}}
&U_{n-1}'\ar@{.>}[ru]_{d_{n-1}'}\ar@{.>}[d]_{\varphi_{n-1}'}\\
E: &0 \ar[r]&U_{n}\ar[r]^{\pi_{n}}
&U_{n-1}\ar[r]^{\pi_{n-1}}&U_{n-2}\ar[r]^{\pi_{n-2}}&\ldots\ar[r]^{\pi_{2}}&U_{1}\ar[r]^{\pi_{1}}&U_{0}\ar[r]^{\pi_{0}}&K\ar[r]&0.
}
\end{xy}
\end{center}
Sei weiter $g_{k}\in C^{k}(G,U_{k})$ f"ur $k=0,\ldots,n$ mit der $(*)$
entsprechenden Eigenschaft aus Satz~\ref{SatzVomIndKoz}, also
\[
\pi_{0}(g_{0})=1\quad\textrm{ und } \quad \pi_{k}\circ
g_{k}=\partial_{k-1}^{U_{k-1}}g_{k-1}\,\,\textrm{ f"ur } k=1,\ldots,n.\quad
(*)
\]
Wir verwenden wieder die Notation aus Abschnitt
\ref{sectiondiebarresolutin}, und gem"a"s Gleichung \eqref{deffvonomega}
(S. \pageref{deffvonomega}) definieren wir
\[\varphi_{k}:=\omega_{k}^{U_{k}}(g_{k})\quad \textrm{ f"ur }  k=0,\ldots,n.\] Mit $e_{k}\in
C^{k}(G,P_{k})$ wie in \eqref{defvonen} gilt dann
\begin{eqnarray*}
\pi_{0}(\varphi_{0}(e_{0}))&=&\pi_{0}(
\omega_{0}^{U_{0}}(g_{0})\circ
e_{0})\stackrel{~\eqref{omegaTog}}{=}\pi_{0}(
g_{0})\stackrel{(*)}{=}1\\&\stackrel{~\eqref{defd0}, S. \pageref{defd0}}{=}&\id_{K}(d_{0}(e_{0}))
\end{eqnarray*}
und
\begin{eqnarray}\label{FetteRechnungZumKommutieren1}
\pi_{k}\circ \varphi_{k}\circ e_{k}&=&\pi_{k}\circ
\omega_{k}^{U_{k}}(g_{k})\circ
e_{k}\stackrel{~\eqref{omegaTog}}{=}\pi_{k}\circ
g_{k}\stackrel{(*)}{=}\partial_{k-1}^{U_{k-1}}g_{k-1}\\&\stackrel{~\eqref{omegaTog}}{=}&\partial_{k-1}^{U_{k-1}}(\omega_{k-1}^{U_{k-1}}(g_{k-1})\circ
e_{k-1})=\partial_{k-1}^{U_{k-1}}(\varphi_{k-1}\circ
e_{k-1})\label{FetteRechnungZumKommutieren2}\\&=&\varphi_{k-1}\circ
\partial_{k-1}^{P_{k-1}}e_{k-1}\stackrel{~\eqref{enKozyklus},\, S. \pageref{enKozyklus}}{=}\varphi_{k-1}\circ
d_{k}\circ e_{k} \quad \textrm{ f"ur } k=1,\ldots,n.\label{FetteRechnungZumKommutieren3}
\end{eqnarray}
Da die Bilder $e_{k}(G^{k})$ jeweils $P_{k}$ als $KG$-Modul erzeugen und alle
betrachteten Abbildungen $KG$-linear sind, gilt also
\begin{equation}\label{rechtskommuts}
\pi_{0}\circ \varphi_{0}=\id_{K}\circ d_{0}\quad\textrm{und}\quad
 \pi_{k}\circ\varphi_{k}=\varphi_{k-1}\circ d_{k}\,\,\textrm{ f"ur
 }k=1,\ldots, n,
\end{equation}
insbesondere kommutiert das Diagramm
"`rechtsseitig"' der Abbildung $\varphi_{n-1}$. Weiter gilt mit $\overline{e_{n}}$
wie in \eqref{defbaren}, S. \pageref{defbaren}
\begin{eqnarray*}
\pi_{n}\circ \phi_{g_{n}}^{n}\circ\overline{e_{n}}&\stackrel{\textrm{Satz }
  \ref{GenericUniverselleEig}}{=}&\pi_{n}\circ
  g_{n}\stackrel{~\eqref{FetteRechnungZumKommutieren1},~\eqref{FetteRechnungZumKommutieren3}}=\varphi_{n-1}\circ d_{n}\circ
  e_{n}\\&\stackrel{~\eqref{dnenbardnbaren}, \, S. \pageref{dnenbardnbaren}}{=}&\varphi_{n-1}\circ \overline{d_{n}}\circ \overline{e_{n}},
\end{eqnarray*}
und da $\overline{e_{n}}(G^{n})$ den $KG$-Modul $\overline{P_{n}}$ erzeugt, also
\[
\pi_{n}\circ \phi_{g_{n}}^{n}=\varphi_{n-1}\circ \overline{d_{n}}.
\]
Da $U_{n-1}'$ als Pushout von $\phi_{g_{n}}^{n}$ und $\overline{d_{n}}$ definiert
ist, liefert dessen universelle Eigenschaft also eine Abbildung
$\varphi_{n-1}': U_{n-1}'\rightarrow U_{n-1}$ mit
\begin{equation}\label{kommutganzlinks}
\varphi_{n-1}'\circ d_{n}'=\pi_{n}=\pi_{n}\circ \id_{U_{n}},
\end{equation}
also Kommutativit"at des Diagramms "`links unten"',  und 
\begin{equation}\label{senkrechtkomm}
\varphi_{n-1}'\circ \varepsilon_{n-1}=\varphi_{n-1}.
\end{equation}
Da wir in Satz \ref{neuerKozyklusIndizierer} $d_{n-1}'$ ebenfalls aus der universellen
Eigenschaft des Pushouts erhalten haben, gilt ausserdem
\begin{equation}\label{dn1kommutiert}
d_{n-1}'\circ \varepsilon_{n-1}=d_{n-1}.
\end{equation}
Wir m"ussen als letztes noch die Kommutatvit"at des "`Trapezes"' zeigen,
n"amlich $\pi_{n-1}\circ \varphi_{n-1}'=\varphi_{n-2}\circ d_{n-1}'$. Sei dazu
$y\in U_{n-1}'$ beliebig. Nach Konstruktion des Pushouts (siehe Definition
\ref{DefPushout}) gilt $U_{n-1}'=d_{n}'(U_{n})+\varepsilon_{n-1}(P_{n-1})$,
d.h. zu $y$ gibt es $u\in U_{n},\,x\in P_{n-1}$ mit 
\begin{equation}\label{Defvonydneps}
y=d_{n}'(u)+\varepsilon_{n-1}(x).
\end{equation}
Dann gilt
\begin{eqnarray*}
\pi_{n-1}\circ \varphi_{n-1}'(y)&\stackrel{~\eqref{Defvonydneps}}{=}&\pi_{n-1}\circ
\varphi_{n-1}'(d_{n}'(u)+\varepsilon_{n-1}(x))\\&\stackrel{~\eqref{kommutganzlinks},~\eqref{senkrechtkomm}}{=}&\pi_{n-1}\circ\pi_{n}(u)+\pi_{n-1}\circ\varphi_{n-1}(x)\\&\stackrel{\pi_{n-1}\circ\pi_{n}=0,\,\,
\eqref{rechtskommuts}}{=}&\varphi_{n-2}\circ
d_{n-1}(x)\\&\stackrel{~\eqref{dn1kommutiert}}{=}&\varphi_{n-2}\circ
d_{n-1}'\circ \varepsilon_{n-1}(x)\\&\stackrel{d_{n-1}'\circ
  d_{n}'=0}{=}&\varphi_{n-2}\circ d_{n-1}'(d_{n}'(u)+\varepsilon_{n-1}(x))\\&\stackrel{~\eqref{Defvonydneps}}{=}&\varphi_{n-2}\circ d_{n-1}'(y)
\end{eqnarray*}
f"ur alle $y\in U_{n-1}'$, d.h. $\pi_{n-1}\circ \varphi_{n-1}'=\varphi_{n-2}\circ d_{n-1}'$. Also kommutiert das gesamte Diagramm,  und damit
gilt $E'\sim_{pY}E$ nach Definition \ref{YonedaPraaq}. \qed\\

\begin{BspRoman}\label{Yonedafuer1}
Wir wollen die (Yoneda-)"Aquivalenz f"ur den Fall $n=1$ untersuchen und
zeigen, dass sie hier mit der Yoneda-Pr"a"aquivalenz "ubereinstimmt, dass hier
also insbesondere bereits die Yoneda-Pr"a"aquivalenz eine "Aquivalenzrelation
ist. Wir betrachten dazu zwei Yoneda-pr"a"aquivalente exakte Sequenzen $E\sim_{pY}
E'$, haben also ein kommutatives Diagramm
\begin{center}
~\begin{xy}
\xymatrix{
E:&0 \ar[r]  & U_{1}\ar[r]^{\pi_{1}}
\ar[d]_{\id_{U_{1}}}&U_{0}\ar[r]^{\pi_{0}}
\ar[d]_{\varphi_{0}}&K\ar[r]
\ar[d]_{\id_{K}}&0\\
E':&0 \ar[r]&U_{1}\ar[r]_{\pi_{1}'}&U_{0}'\ar[r]_{\pi_{0}'}&K\ar[r]&0.
}
\end{xy}
\end{center}
Wir zeigen, dass dann $\varphi_{0}$ ein Isomorphismus ist, woraus dann auch
$E'\sim_{pY}E$ folgt. Sei $u_{0}\in \ker \varphi_{0}$. Dann ist
$0=\pi_{0}'(\varphi_{0}(u_{0}))=\id_{K}(\pi_{0}(u_{0}))$, also
$u_{0}\in\ker\pi_{0}=\im \pi_{1}$. Dann gibt es $u_{1}\in U_{1}$ mit
$u_{0}=\pi_{1}(u_{1})$. Es folgt
$0=\varphi_{0}(u_{0})=\varphi_{0}(\pi_{1}(u_{1}))=\pi_{1}'(\id_{U_{1}}(u_{1}))$. Da $\pi_{1}'$
injektiv ist, folgt $u_{1}=0$ und damit $u_{0}=\pi_{1}(u_{1})=0$, also $\ker
\varphi_{0}=0$ und $\varphi_{0}$ ist injektiv. Zum Beweis der Surjektivit"at
von $\varphi_{0}$ sei $u_{0}'\in U_{0}'$. Da $\pi_{0}$ surjektiv ist, gibt es dann
$u_{0}\in U_{0}$ mit $\pi_{0}(u_{0})=\pi_{0}'(u_{0}')$. Dann ist
$\pi_{0}'(u_{0}'-\varphi_{0}(u_{0}))=\pi_{0}'(u_{0}')-\pi_{0}(u_{0})=0$, also
$u_{0}'-\varphi_{0}(u_{0})\in \ker \pi_{0}'=\im \pi_{1}'$. Es gibt also
$u_{1}'\in U_{1}$ mit $\pi_{1}'(u_{1}')=u_{0}'-\varphi_{0}(u_{0})$. Es folgt
$\varphi_{0}(u_{0}+\pi_{1}(u_{1}'))=\varphi_{0}(u_{0})+\pi_{1}'(u_{1}')=u_{0}'$,
also die Surjektivit"at von $\varphi_{0}$.

Im Fall $n=1$ gibt es also zu jedem Kozyklus im Wesentlichen nur eine exakte
Sequenz, die diesen induziert.
\end{BspRoman}

\begin{BemRoman}
Aus Satz \ref{SatzVomIndKoz}, Beispiel \ref{KorandInduzierer} und Satz
\ref{SatzVonderYoneadAquiv} folgt: Ein von einer exakten Sequenz mit Ende
$U_{n}$ induzierter
Kozyklus ist genau dann trivial, wenn die exakte Sequenz (Yoneda-)"aquivalent
ist zu einer exakten Sequenz $0\rightarrow
U_{n}\stackrel{\pi_{n}}{\rightarrow}U_{n-1}\rightarrow \ldots$, f"ur die ein
links-Splitting $\mu_{n}: U_{n-1}\rightarrow U_{n}$ mit $\mu_{n}\circ
\pi_{n}=\id_{U_{n}}$ existiert. Im Fall $n=1$ folgt aus Beispiel
\ref{Yonedafuer1}, dass dies genau dann der Fall ist, wenn f"ur eine beliebige
(und dann jede) den Kozyklus induzierende Sequenz ein solches Splitting existiert.
\end{BemRoman}

\subsection{Annullation von $n$-Kozyklen}
Sind $V,W$ jeweils $KG$-Moduln und $g\in C^{n}(G,V)$ sowie $\varphi\in (W^{*})^{G}$, so
schreiben wir $\varphi\cdot g\in C^{n}(G,\Hom_{K}(W,V))$ f"ur die Abbildung $G^{n}\rightarrow \Hom_{K}(W,V),\,\,x\mapsto \varphi\cdot g(x)$ mit
$\varphi\cdot g(x): W\rightarrow V,\,\, w\mapsto \varphi(w)g(x)$. Ist $g\in
Z^{n}(G,V)$, also $\partial_{n}^{V}g=0$,  so folgt wegen $\varphi\in (W^{*})^{G}$ dann
auch $\partial_{n}^{\Hom_{K}(W,V)}(\varphi\cdot g)=0$, also $\varphi\cdot
g\in Z^{n}(G,\Hom_{K}(W,V))$. Ist $g\in B^{n}(G,V)$, also
$g=\partial_{n-1}^{V}f$ mit $f\in C^{n-1}(G,V)$, so gilt entsprechend auch
$\varphi\cdot g=\varphi\cdot
\partial_{n-1}^{V}f=\partial_{n-1}^{\Hom_{K}(W,V)}(\varphi \cdot f)\in B^{n}(G,\Hom_{K}(W,V))$. Die
Multiplikation mit $\varphi$ induziert also eine Abbildung
$H^{n}(G,V)\rightarrow H^{n}(G,\Hom_{K}(W,V))$. Unter anderem werden wir in
diesem Abschnitt zeigen, dass die \emph{Augmentationsabbildung}\index{Augmentationsabbildung} $d_{0}\in
(KG)^{*}, \,\, \sum_{\sigma\in G}\lambda_{\sigma}\cdot \sigma\mapsto
\sum_{\sigma\in G}\lambda_{\sigma}$ (wobei nur endlich viele
$\lambda_{\sigma}\in K$ ungleich $0$ sind)
stets die Nullabbildung induziert.

Ist $W$ \emph{endlich-dimensional}, so gilt bekanntlich die Isomorphie
$W^{*}\otimes_{K} V\cong \Hom_{K}(W,V)$, gegeben durch lineare Fortsetzung von
$\varphi\otimes v\rightarrow \varphi\cdot v$. Entsprechend ist dann $\varphi\otimes
g\in Z^{n}(G,W^{*}\otimes_{K} V)$ definiert.
\\

Sind $U,V,W$ jeweils $KG$-Moduln und $\psi\in \Hom_{G}(U,V)$ (also $\psi\circ
\sigma=\sigma\circ \psi\,\,\myforall \sigma\in G$), so ist die
Abbildung $\psi_{*}: \Hom_{K}(W,U)\rightarrow \Hom_{K}(W,V),\,\, f\mapsto
\psi\circ f$ ein Element aus $\Hom_{G}(\Hom_{K}(W,U),\Hom_{K}(W,V))$. Denn
f"ur $\sigma\in G,\,\, f\in\Hom_{K}(W,U) $ gilt $\psi_{*}(\sigma\cdot
f)=\psi\circ\sigma\circ f\circ \sigma^{-1}=\sigma \circ \psi\circ f\circ
\sigma^{-1}=\sigma\cdot \psi_{*}(f)$. Das folgende Lemma ist dann
wohlbekannt. 

\begin{Lemma}
Sind $A,B,C,W$ jeweils $KG$-Moduln und ist
\[A\stackrel{\varepsilon}{\rightarrow}B\stackrel{\pi}{\rightarrow}C\]
eine exakte Sequenz von $KG$-Moduln, so ist auch
\[\Hom_{K}(W,A)\stackrel{\varepsilon_{*}}{\longrightarrow}\Hom_{K}(W,B)\stackrel{\pi_{*}}{\longrightarrow}\Hom_{K}(W,C)\]
eine exakte Sequenz von $KG$-Moduln.
\end{Lemma}

\Bew (Der Vollst"andigkeit halber.) F"ur $f\in \Hom_{K}(W,A)$ gilt
$\pi_{*}\circ \varepsilon_{*}(f)=\pi\circ\varepsilon \circ f=0$, da
$\pi\circ\varepsilon=0$, also $\im \varepsilon_{*}\subseteq \ker \pi_{*}$. Sei
umgekehrt $g\in \ker \pi_{*}$, also $\pi\circ g=0$. Dann ist $\im g\subseteq
\ker \pi=\im \varepsilon$. Da wir von Vektorr"aumen sprechen, gibt es dann ein
$f\in \Hom_{K}(W,A)$ mit $g=\varepsilon\circ f=\varepsilon_{*}(f)\in \im
\varepsilon_{*}$, also auch $\ker \pi_{*}\subseteq \im \varepsilon_{*}$. \qed\\


Wir kommen zu der angek"undigten Verallgemeinerung von Proposition
\ref{AnnulatorProp}, dem eigentlichen Ziel und Hauptresultat dieses Abschnitts
"uber Kohomologie von Gruppen.

\begin{Satz}\label{AnnulatorDurchSeq}
Sei $n\ge 1$ und
\begin{center}
  ~\begin{xy}
\xymatrix{
0\ar[r]&U_{n}\ar[r]^{\pi_{n}}&U_{n-1}\ar[r]^{\pi_{n-1}}&U_{n-2}\ar[r]^{\pi_{n-2}}&\ldots\ar[r]^{\pi_{3}}&U_{2}\ar[r]^{\pi_{2}}&U_{1}\ar[r]^{\pi_{1}}&U_{0}\ar[r]^{\pi_{0}}&K\ar[r]&0
}
\end{xy}
\end{center}
eine exakte Sequenz von $G$-Moduln, die einen Kozyklus $g\in Z^{n}(G,U_{n})$
induziert. Dann ist $\pi_{0}\in (U_{0}^{*})^{G}$, und es gilt $\pi_{0}\cdot
g\in B^{n}(G,\Hom_{K}(U_{0},U_{n}))$. Anders fomuliert: $\pi_{0}\cdot
g=0\in H^{n}(G,\Hom_{K}(U_{0},U_{n}))$, d.h. die Restklasse von $g$ in $H^{n}(G,U_{n})$ wird von
$\pi_{0}$ annulliert.

Falls $\dim_{K}U_{0}<\infty$, so gilt entsprechend $\pi_{0}\otimes g\in
B^{n}(G,U_{0}^{*}\otimes_{K} U_{n})$. 
\end{Satz}

\Bew Da $\pi_{0}$ ein $G$-Homomorphimsmus ist, ist $\pi_{0}\in
(U_{0}^{*})^{G}$. 
Sei $g_{n}:=g$ und $g_{k}\in C^{k}(G,U_{k})$ f"ur $k=0,\ldots,n$ die
zugeh"orige $g$ induzierende Folge mit
$(*)$ wie in Satz/Definition \ref{SatzVomIndKoz}, d.h. 
\[
\pi_{0}(g_{0})=1\quad \textrm{ und }\quad \pi_{k}\circ
g_{k}=\partial_{k-1}^{U_{k-1}}g_{k-1} \quad \textrm{ f"ur }
k=1,\ldots,n. \quad (*)
\]
Wir wenden nun auf die exakte Sequenz den "`Funktor"' $\Hom_{K}(U_{0},\cdot)$ an und erhalten so die nach
dem Lemma ebenfalls exakte Sequenz
\begin{center}
  ~\begin{xy}
\xymatrix{
0\ar[r]&\Hom_{K}(U_{0},U_{n})\ar[r]^{\pi_{n*}}&\Hom_{K}(U_{0},
U_{n-1})\ar[r]^{\pi_{n-1*}}&\Hom_{K}(U_{0},U_{n-2})\ar[r]^<<<<<<<{\pi_{n-2*}}&\ldots\ldots\ldots\quad\quad
&&\\\ldots\ar[r]^<<<<<{\pi_{3*}} &\Hom_{K}(U_{0},U_{2})\ar[r]^{\pi_{2*}}&\Hom_{K}(U_{0},U_{1})\ar[r]^{\pi_{1*}}&\Hom_{K}(U_{0},U_{0})\ar[r]^<<<<<<<{\pi_{0*}}&\Hom_{K}(U_{0},K)\quad(**)
}
\end{xy}
\end{center}
mit $\pi_{k*}$ definiert wie oben durch Verkettung mit $\pi_{k}$.
F"ur die Folge der \[\pi_{0}\cdot g_{k}\in C^{k}(G,\Hom_{K}(U_{0},
U_{k}))\quad\textrm{ mit }\quad k=0,\ldots,n\] gilt dann offenbar
\[
\pi_{0*}(\pi_{0}\cdot g_{0})=\pi_{0}\cdot \pi_{0}(g_{0})\stackrel{(*)}{=}\pi_{0}\cdot 1=\pi_{0}
\]
und
\begin{eqnarray*}
\pi_{k*}\circ
(\pi_{0}\cdot g_{k})&=&\pi_{0}\cdot (\pi_{k}\circ
g_{k})\stackrel{(*)}{=}\pi_{0}\cdot \partial_{k-1}^{U_{k-1}}g_{k-1}\\&=& \partial_{k-1}^{\Hom_{K}(U_{0},
  U_{k-1})}(\pi_{0}\cdot g_{k-1}) \quad 
\myforall k=1,\ldots,n. 
\end{eqnarray*}
Dies ist die Eigenschaft $(*)$ von Satz/Definition \ref{SatzVomIndKoz} f"ur
die exakte Sequenz $(**)$ (mit $w:=\pi_{0}\cdot 1\in \Hom_{K}(U_{0},
K)^{G}$), und damit induziert die exakte Sequenz $(**)$ den $n$-Kozyklus $\pi_{0}\cdot g_{n}\in
Z^{n}(G,\Hom_{K}(U_{0}, U_{n}))$. 

Wir betrachten nun $\id_{U_{0}}\in\Hom_{K}(U_{0},U_{0})^{G}$. Es gilt
$\pi_{0*}(\id_{U_{0}})=\pi_{0}\circ \id_{U_{0}}=\pi_{0}=\pi_{0}\cdot 1$.
Nach Satz/Definition \ref{SatzVomIndKoz} (mit $v=\id_{U_{0}}$) ist also
$\pi_{0}\cdot g_{n}\in B^{n}(G,\Hom_{K}(U_{0}, U_{n}))$. 

Der Zusatz f"ur $\dim_{K}U_{0}<\infty$ folgt nun aus der Isomorphie
$\Hom_{K}(U_{0},U_{n})\cong U_{0}^{*}\otimes_{K} U_{n}$. \qed\\

Im Falle $n=1$ und eines Kozyklus $g\in Z^{1}(G,U)$ mit $\dim_{K}U<\infty$, der von einer kurzen
exakten Sequenz
\[
0\rightarrow U \hookrightarrow \tilde{U} \stackrel{\pi}{\rightarrow} K
\rightarrow 0
\]
induziert wird (siehe Beispiel \ref{einsKozyklenallgemein}), liefert der Satz
$\pi \otimes g\in B^{1}(G,\tilde{U}^{*}\otimes_{K} U)$. Mit der in Bemerkung
\ref{gnAlsnKorand} beschriebenen Konstruktion ergibt sich auch die Formel
\eqref{Annulator} (S. \pageref{Annulator}).
 Wir erhalten also Proposition
\ref{AnnulatorProp} zur"uck.

\begin{Korollar}\label{augmentationannulliert}
Sei $G$ eine beliebige Gruppe und 
\[
d_{0}: KG\rightarrow K, \quad \sum_{\sigma\in G}\lambda_{\sigma}\cdot\sigma\mapsto
\sum_{\sigma\in G}\lambda_{\sigma}
\]
(mit nur endlich  vielen $\lambda_{\sigma}\in K$ ungleich $0$)
die \emph{Augmentationsabbildung}.\index{Augmentationsabbildung}
Dann annulliert $d_{0}$ jeden Kozyklus eines jeden
$G$-Moduls $U$, d.h. f"ur alle $g\in Z^{n}(G,U)$ ($n\ge 1$) gilt $d_{0}\cdot g\in B^{n}(G,\Hom_{K}(KG,U))$.
\end{Korollar}

\Bew Sei zun"achst $n\ge 2$ und $U_{n}:=U$. Nach Satz
\ref{neuerKozyklusIndizierer} wird $g\in Z^{n}(G,U_{n})$ von einer exakten
Sequenz der Form
\begin{center}
  ~\begin{xy}
\xymatrix{
0\ar[r]&U_{n}\ar[r]^{d_{n}'}&U_{n-1}\ar[r]^{d_{n-1}'}&P_{n-2}\ar[r]^{d_{n-2}}&\ldots\ar[r]^{d_{1}}&KG\ar[r]^{d_{0}}&K\ar[r]&0
}
\end{xy}
\end{center}
induziert, die also "`auf die bar resolution endet"' (da $n\ge 2$, kommt
$P_{0}=KG$ tats"achlich vor). Nach Satz \ref{AnnulatorDurchSeq} gilt also
$d_{0}\cdot g\in B^{n}(G,\Hom_{K}(KG,U))$.

Es bleibt der Fall $n=1$. Wir erledigen ihn durch explizite Rechnung. Sei also
$g\in Z^{1}(G,U)$. Wir definieren dazu $f\in \Hom_{K}(KG,U)$ durch $K$-lineare Fortsetzung von
\[
G\rightarrow U,\quad \tau\mapsto-g_{\tau}
\]
($G$ ist $K$-Basis von $KG$).
Wir berechnen nun
\[
\partial_{0}^{\Hom_{K}(KG,U)}f(\sigma)=(\sigma-1)f=\sigma\circ
f\circ\sigma^{-1}-f\in \Hom_{K}(KG,U), 
\]
indem wir die Bilder dieser Abbildung auf der Basis $G$ von $KG$ angeben. F"ur
$\tau\in G$ gilt
\begin{eqnarray*}
((\partial_{0}f)(\sigma))(\tau)&=&(\sigma\circ
f\circ\sigma^{-1})(\tau)-f(\tau)=\sigma(-g_{\sigma^{-1}\tau})+g_{\tau}\\&=&-(\sigma
g_{(\sigma^{-1}\tau)}+g_{\sigma}-g_{\sigma})+g_{\tau}\\
&=&-g_{\sigma(\sigma^{-1}\tau)}+g_{\sigma}+g_{\tau}=g_{\sigma}=(d_{0}\cdot g_{\sigma})(\tau),
\end{eqnarray*}
wobei wir die Kozyklus Eigenschaft von $g$ verwendet haben. Also gilt $\partial_{0}f(\sigma)=d_{0}\cdot g_{\sigma}$
oder $d_{0}\cdot g=\partial_{0}f\in B^{1}(G,\Hom_{K}(KG,U))$. \qed\\

\noindent Wir geben noch einen\\

{\it\noindent 2. Beweis.} Wir betrachten den generischen
$n$-Kozyklus $\overline{e_{n}}\in Z^{n}(G,\overline{P_{n}})$ \eqref{defbaren},
S. \pageref{defbaren}. Er wird induziert von der exakten Sequenz
\eqref{generischeexakteSeq}, und nach Satz \ref{AnnulatorDurchSeq} gilt also
$d_{0}\cdot \overline{e_{n}}\in B^{n}(G,\Hom_{K}(KG,\overline{P_{n}}))$,
d.h. es gibt $f_{n-1}\in C^{n-1}(G,\Hom_{K}(KG,\overline{P_{n}}))$ mit

\begin{equation}\label{annulgenericnkoz}
d_{0}\cdot
 \overline{e_{n}}=\partial_{n-1}^{\Hom_{K}(KG,\overline{P_{n}})}f_{n-1}.
\end{equation}
 Sei
nun $g_{n}\in Z^{n}(G,U_{n})$ mit einem $G$-Modul $U_{n}$.
Wir
betrachten die Abbildung $\phi_{g_{n}}^{n}\in
\Hom_{KG}(\overline{P_{n}},U_{n})$ aus Satz \ref{GenericUniverselleEig} mit
\begin{equation}\label{SchreibeNochmalgnphignen}
g_{n}=\phi_{g_{n}}^{n}\circ \overline{e_{n}}.
\end{equation}
 Wie wir in der Diskussion am
Anfang dieses Abschnitts gesehen haben, ist dann die Abbildung
\[(\phi_{g_{n}}^{n})_{*}: \Hom_{K}(KG,\overline{P_{n}})\rightarrow
\Hom_{K}(KG,U_{n}),\quad \psi\mapsto
(\phi_{g_{n}}^{n})_{*}(\psi)=\phi_{g_{n}}^{n}\circ \psi\] ebenfalls ein
$KG$-Homomorphismus, also 
\begin{equation}\label{phignnistkghom}
(\phi_{g_{n}}^{n})_{*}\in
\Hom_{KG}(\Hom_{K}(KG,\overline{P_{n}}),\Hom_{K}(KG,U_{n})).
\end{equation}
 Damit erhalten
wir
\begin{eqnarray*}
d_{0}\cdot g_{n}&\stackrel{~\eqref{SchreibeNochmalgnphignen}}{=}&d_{0}\cdot (\phi_{g_{n}}^{n}\circ
\overline{e_{n}})=(\phi_{g_{n}}^{n})_{*}\circ (d_{0}\cdot
\overline{e_{n}})\\
&\stackrel{~\eqref{annulgenericnkoz}}{=}&(\phi_{g_{n}}^{n})_{*}\circ\partial_{n-1}^{\Hom_{K}(KG,\overline{P_{n}})}f_{n-1}\\
&\stackrel{~\eqref{phignnistkghom}}{=}&\partial_{n-1}^{\Hom_{K}(KG,U_{n})}((\phi_{g_{n}}^{n})_{*}\circ f_{n-1})
\end{eqnarray*}
mit $(\phi_{g_{n}}^{n})_{*}\circ f_{n-1}\in C^{n-1}(G,\Hom_{K}(KG,U_{n}))$,   also $d_{0}\cdot g_{n}\in B^{n}(G,\Hom_{K}(KG,U_{n}))$. \qed\\

\begin{BemRoman}
Ist $|G|<\infty$, so ist $KG$ selbstdual als $G$-Modul. Ist n"amlich
$\delta_{\tau}\in \Hom_{K}(KG,K)$ f"ur $\tau\in G$ gegeben durch $K$-lineare Fortsetzung von 
\[
\delta_{\tau}(\sigma)=\left\{\begin{array}{cl}
1 & \textrm{ f"ur }\sigma=\tau\\
0&\textrm{ sonst,}\end{array}\right.
\]
 so bilden die $\delta_{\tau}$ eine Basis von
$(KG)^{*}=\Hom_{K}(KG,K)$. Ferner ist durch $K$-lineare Fortsetzung von $\tau\mapsto
\delta_{\tau}$ ein $G$-Isomorphismus $KG\rightarrow (KG)^{*}$ gegeben, denn f"ur $\sigma\in G$ gilt
$\sigma\cdot \delta_{\tau}=\delta_{\tau}\circ
\sigma^{-1}=\delta_{\sigma\tau}$. Die Augmentationsabbildung
$d_{0}=\sum_{\sigma\in G}\delta_{\sigma}\in\Hom_{K}(KG,K)$ entspricht also der "`Spur"' $\pi:=\sum_{\sigma\in
  G}\sigma\in KG$. Da $\dim_{K}KG<\infty$, also $\Hom_{K}(KG,U)\cong
(KG)^{*}\otimes_{K} U\cong KG\otimes_{K} U$ f"ur jeden $G$-Modul $U$ gilt, haben wir also mit dem Korollar f"ur jedes
$g\in Z^{n}(G,U)$, dass $\pi\otimes g\in B^{n}(G,KG\otimes_{K} U)$. Dieses
Resultat "uber die Annullation von Kozyklen endlicher Gruppen ist bekannt,
siehe etwa Kemper \cite[Lemma 1.7]{KemperOnCM} oder \cite[Proposition 2.3]{KemperDepthCohom}.
\end{BemRoman}

Wir geben noch ein letztes Beispiel.

\begin{BspRoman}
Wir greifen Beispiel \ref{Bsp2Kozyklus} wieder auf, also den von der exakten
Sequenz
\[
0\rightarrow K \stackrel{\pi_{2}}{\rightarrow}K^{2}
\stackrel{\pi_{1}}{\rightarrow}K^{2}\stackrel{\pi_{0}}{\rightarrow}K\rightarrow
0
\]
induzierten $2$-Kozyklus $g\in Z^{2}(\Ga,K)$
 mit $g(a,b)=ab\,\,\myforall a,b\in \Ga$. F"ur $\chr K=2$ ist $g$ nichttrivial,
 aber nach dem Satz gilt $\pi_{0}\otimes g\in B^{2}(\Ga,(K^{2})^{*}\otimes_{K}
 K)$. Tats"achlich kann $\pi_{0}\otimes g$ mit $h\in Z^{2}(\Ga,K^{2})$ mit
 $h(a,b)=(ab,0)\,\,\myforall a,b\in \Ga$ identifiziert werden, und dann ist
 $h=\partial_{1}f$, mit $f\in C^{1}(\Ga,K^{2})$ gegeben durch
 $f(a)=(0,a)\,\,\myforall a\in \Ga$ - also wie vom Satz behauptet $h\in B^{2}(\Ga,K^{2})$.
\end{BspRoman}

\newpage
\section{Folgen von Invariantenringen mit unbeschr"ankt wachsendem Cohen-Macaulay-Defekt}\label{Hauptabschnitt}

In diesem Abschnitt behandeln wir eines der Hauptresultate dieser Arbeit -
insbesondere werden wir f"ur jede reduktive, nicht linear reduktive Gruppe $G$ einen $G$-Modul
$V$ konstruieren mit $\lim_{k \to \infty}\cmdef K[\bigoplus_{i=1}^{k}V]^{G}=\infty$.\\ 

Im folgenden Lemma ist $G$ eine beliebige (nicht notwendig
reduktive) lineare algebraische Gruppe.

\begin{Lemma} \label{ZweiReg}
Homogene Elemente $a_{1},a_{2}\in K[V]_{+}$ bilden genau dann ein phsop in
$K[V]$, wenn sie teilerfremd sind, was genau dann der Fall ist, wenn sie eine
regul"are Sequenz in $K[V]$ bilden.

Gilt dann zus"atzlich $a_{1},a_{2}\in K[V]^{G}$, so bilden sie auch eine
regul"are Sequenz in $K[V]^{G}$ (und dann dort auch ein phsop, falls
$K[V]^{G}$ endlich erzeugt, z.B. $G$ reduktiv).
\end{Lemma}

\Bew Ist $a_{1},a_{2}$ ein phsop
in $K[V]$, so ist es dort auch eine regul"are Sequenz, denn der Polynomring
$K[V]$ ist Cohen-Macaulay. Ist nun $d\in K[V]$ ein gemeinsamer Teiler von
$a_{1},a_{2}$, also $a_{1}=dt_{1}, a_{2}=dt_{2}$ mit $t_{1},t_{2}\in K[V]$, so folgt
$t_{1}a_{2}=dt_{1}t_{2}=t_{2}a_{1}\in (a_{1})$. Aufgrund der Regularit"at hat
man also $t_{1}\in (a_{1})$, oder $t_{1}=a_{1}t'$ mit $t'\in K[V]$. Es folgt
$t_{1}=dt_{1}t'$ oder $1=dt'$. Also ist der gemeinsame Teiler $d$ eine Einheit
und damit $a_{1},a_{2}$ teilerfremd.

Seien nun $a_{1},a_{2}$ teilerfremd. Wir zeigen, dass dann $a_{1},a_{2}$ eine
regul"are Sequenz in $K[V]$ bzw. unter der Zusatzvoraussetzung auch in
$K[V]^{G}$ ist. Dann bilden die beiden Elemente dort auch jeweils
ein  phsop (Satz~\ref{regPhsop} (a)).

 Sei je nachdem $R=K[V]$ oder $R=K[V]^{G}$, und $h_{2}\in R$ mit $h_{2}a_{2}\in
(a_{1})_{R}$. Dann existiert $h_{1}\in R$ mit
$h_{2}a_{2}=h_{1}a_{1}$. Da $a_{1},a_{2}$ in $K[V]$ teilerfremd
sind, folgt also
aus $a_{1}|h_{2}a_{2}$, dass $a_{1}|h_{2}$ in $K[V]$; D.h. es gibt ein $t\in
K[V]$ mit $h_{2}\stackrel{(*)}{=}a_{1}t$, d.h. $h_{2}\in (a_{1})_{K[V]}$. Dies zeigt die
$K[V]$-Regularit"at. Im Fall $R=K[V]^{G}$ sind $h_{2}\in R$ und $a_{1}$
invariant, und Anwendung eines $\sigma \in G$ auf Gleichung $(*)$ liefert
$h_{2}=a_{1}(\sigma t)$. 
Da $K[V]$ nullteilerfrei ist,  ist auch $t=h_{2}/a_{1}=\sigma t$ invariant,
also $t \in K[V]^{G}$. Damit ist $h_{2}=a_{1}t\in (a_{1})_{K[V]^{G}}$, also
$a_{1},a_{2}$ auch $K[V]^{G}$-regul"ar. \qed\\

\subsection{Die Charakterisierung linear reduktiver Gruppen nach Kemper}

Das Haupthilfsmittel f"ur die
 Konstruktion von Invariantenringen mit beliebig gro"sem
Cohen-Ma\-cau\-lay-Defekt ist das folgende Lemma, das im Beweis von Kemper \cite[Proposition
6]{KemperLinRed} steckt, und zur Konstruktion von
nicht-Cohen-Macaulay-Invariantenringen (also Invariantenringen mit positivem
Cohen-Macaulay-Defekt) f"uhrt. Aus diesem folgt dann auch die Charakterisierung
linear reduktiver Gruppen nach Kemper als diejenigen reduktiven Gruppen, deren
Invariantenringe stets Cohen-Macaulay sind.


\begin{Lemma} \label{Hauplemma}
Sei $G$ eine lineare algebraische Gruppe, $V$ ein $G$-Modul, und es existiere ein $0 \ne g \in
H^{1}(G,K[V])$. Seien weiter  $a_{1},a_{2}\in K[V]^{G}_{+}$ homogen und teilerfremd in
$K[V]$ (d.h. ein phsop in
$K[V]$), die beide  $g$ annullieren, also $a_{i}g=0 \in H^{1}(G,K[V])$ f"ur
$i=1,2$. 

Dann gibt es ein $m\in K[V]^{G}$ mit $m \notin
(a_{1},a_{2})_{K[V]^{G}}$,  so dass f"ur jedes weitere $a_{3}\in K[V]^{G}$ mit
$a_{3}g=0 \in H^{1}(G,K[V])$ gilt, dass $ma_{3}\in(a_{1},a_{2})_{K[V]^{G}}.$

Ist $G$ reduktiv und bilden $a_{1},a_{2},a_{3}$ mit obigen Eigenschaften ein phsop in
$K[V]$, so bilden sie aufgrund der Reduktivit"at von $G$ auch eines in
$K[V]^{G}$, aber dort keine regul"are Sequenz, insbesondere ist also
$K[V]^{G}$ nicht Cohen-Macaulay.
\end{Lemma}

\Bew Wir kl"aren zun"achst die einfache Folgerung. Ist $G$ reduktiv und bilden $a_{1},a_{2},a_{3}$ ein phsop in $K[V]$, so wegen Lemma \ref{redphsop}
auch eines in $K[V]^{G}$. Da aber nach dem ersten Teil des Satzes $ma_{3}\in
(a_{1},a_{2})_{K[V]^{G}}$, aber $m\notin (a_{1},a_{2})_{K[V]^{G}}$, bildet das phsop
dort keine regul"are Sequenz. Also ist $K[V]^{G}$ nach Satz \ref{regPhsop}
(b) nicht Cohen-Macaulay.

Nun zum eigentlichen Beweis des Lemmas. Sei $(\sigma\mapsto g_{\sigma})\in
Z^{1}(G,K[V])$ der zu $g$ geh"orige Kozyklus.
Nach Voraussetzung sind die Ko\-zy\-klen $(\sigma \mapsto a_{i} g_{\sigma})\in
Z^{1}(G, K[V]), i=1,2,3$ trivial, also gibt es $b_{i} \in K[V]$ mit 
\begin{equation*} 
a_{i} g_{\sigma}=(\sigma -1)b_{i} \quad \myforall \sigma \in G,\quad i=1,2,3.
\end{equation*}
Sei 
\begin{equation*} 
u_{ij}=a_{i}b_{j}-a_{j}b_{i} \textrm{ f"ur } 1 \le i < j \le 3.
\end{equation*}
Offenbar ist $u_{ij} \in K[V]^{G}$ ($\sigma u_{ij}=a_{i}(b_{j}+a_{j}g_{\sigma})-a_{j}(b_{i}+a_{i}g_{\sigma})=u_{ij}$), und es gilt
\[
u_{23}a_{1}-u_{13}a_{2}+u_{12}a_{3}=
\left|
\begin{array}{ccc}
a_{1} & a_{2} & a_{3}\\
a_{1} & a_{2} & a_{3}\\
b_{1} & b_{2} & b_{3}\\
\end{array}
\right|=0.
\]
Wir setzen nun $m:=u_{12}=a_{1}b_{2}-a_{2}b_{1}$. Dabei h"angt $m$ nur
von $a_{1}$ und $a_{2}$ ab, und nach obiger Gleichung gilt
$ma_{3}\in(a_{1},a_{2})_{K[V]^{G}}$. Wir m"ussen zeigen, dass $m\notin
(a_{1},a_{2})_{K[V]^{G}}$, und nehmen dazu das Gegenteil an, also $m=u_{12}\in
(a_{1},a_{2})_{K[V]^{G}}$. Dann gibt es $f_{1},f_{2} \in K[V]^{G}$ mit 
\begin{equation*} 
u_{12}=a_{1}b_{2}-a_{2}b_{1}=f_{1}a_{1}+f_{2}a_{2}.
\end{equation*}
Aus $a_{1}(b_{2}-f_{1})=a_{2}(f_{2}+b_{1})$ und der Teilerfremdheit von
$a_{1},a_{2}$ folgt dann, dass $a_{1}$ Teiler von $f_{2}+b_{1}$ ist, also $f_{2}+b_{1}=a_{1} \cdot h$ mit $h \in K[V]$. Nun ist
\[
a_{1} \cdot (\sigma -1)h=(\sigma -1)(a_{1}h)=(\sigma -1)(f_{2}+b_{1})=(\sigma -1)b_{1}=a_{1}g_{\sigma} \qquad \myforall \sigma \in G,
\]
also $g_{\sigma}=(\sigma -1)h$. Damit ist die zugeh"orige Restklasse $g\in H^{1}(G,K[V])$ gleich 0, was im Widerspruch zur Voraussetzung steht. Also war
die Annahme falsch, und es gilt $m\notin
(a_{1},a_{2})_{K[V]^{G}}$.\qed\\

Aus $a_{1}g_{\sigma}=(\sigma-1)b_{1}\myforall \sigma\in G$ mit $b_{1}\in K[V]$
folgt "ubrigens, dass $(g_{\sigma})_{\sigma\in G}$
"`gradbeschr"ankt"' 
ist, also dass es Zahlen $0\le m\le n<\infty$ gibt mit
$(g_{\sigma})_{\sigma\in G}\in
Z^{1}(G,\oplus_{k=m}^{n}S^{k}(V))$. Au"serdem folgt aus $g_{\sigma}=\frac{1}{a_{1}}(\sigma-1)b_{1}$, dass $g$ in jedem Fall
durch einen Morphismus gegeben ist (Polynomdivison, etwa bzgl. graduierter
lexikographischer Ordnung, gibt die Koeffizienten von $g_{\sigma}$ als
Linearkombination der Koeffizienten von $(\sigma-1)b_{1}$), auch falls man diese Forderung an die
Elemente aus $Z^{1}(G,K[V])$ (wie in Abschnitt \ref{Kohomologievongruppen})
zun"achst nicht stellen will.

\begin{Korollar} \label{NCMRing}
Sei $G$ eine reduktive Gruppe, $g\in Z^{1}(G,U)$ ein nichttrivialer Kozyklus,
und $\tilde{U}$ der zugeh"orige erweiterte $G$-Modul (siehe Abschnitt
\ref{firstCohom}). Ist dann
\[
V^{*}:=U\oplus\tilde{U}^{*}\oplus\tilde{U}^{*}\oplus\tilde{U}^{*} \quad \textrm{
  d.h. } V:=U^{*}\oplus\tilde{U}\oplus\tilde{U}\oplus\tilde{U},
\]
so ist $K[V]^{G}$ nicht Cohen-Macaulay.
\end{Korollar}

\Bew Es ist
$K[V]=S(V^{*})=S(U\oplus\tilde{U}^{*}\oplus\tilde{U}^{*}\oplus\tilde{U}^{*})$.
Dann ist auch $g\in Z^{1}(G,S(V^{*}))$ ein nichttrivialer Kozyklus, denn $U$
ist direkter Summand von $S(V^{*})$. Nach
Proposition \ref{AnnulatorProp} gibt es $\pi \in \tilde{U}^{*G}\setminus\{0\}$, so dass $\pi
\otimes g=0\in H^{1}(G,\tilde{U}^{*}\otimes U)$. Sind $a_{1},a_{2},a_{3}$ die
drei Kopien von $\pi$ in den entsprechenden direkten Summanden von $S(V^{*})$,
so bilden diese also ein annullierendes phsop des Kozyklus $g$ in $K[V]$, und nach
vorhergehendem Lemma ist $K[V]^{G}$ nicht Cohen-Macaulay. \qed\\ 

Man beachte, dass bei dieser Konstruktion $V$ bzw. $V^{*}$ niemals
vollst"andig reduzibel sein kann: Die Nichttrivialit"at des Kozyklus $g$ ist
n"amlich "aquivalent dazu, dass $U\subseteq \tilde{U}$
bzw. $K\pi\subseteq\tilde{U}^{*}$ kein Komplement hat. Mit Hilfe verfeinerter
Methoden werden wir sp"ater dennoch Beispiele mit $V$ vollst"andig reduzibel
und nicht Cohen-Macaulay Invariantenring angeben.\\

Es folgt die Charakterisierung linear reduktiver Gruppen nach Kemper:

\begin{Satz}[Kemper \cite{KemperLinRed}] \label{KempersBigTheorem}
Eine reduktive Gruppe $G$ ist genau dann linear reduktiv, wenn $K[V]^{G}$ f"ur
jeden $G$-Modul $V$ Cohen-Macaulay ist.
\end{Satz}

\Bew Ist $G$ linear reduktiv, so ist $K[V]^{G}$ nach Hochster und Roberts
stets Cohen-Macaulay. Ist $G$ reduktiv, aber nicht linear reduktiv, so gibt es
nach Proposition \ref{LinRedKohomologie} einen $G$-Modul $U$ mit
$H^{1}(G,U)\ne 0$, d.h. es existiert ein nichttrivialer Kozyklus $g\in
Z^{1}(G,U)$. Nach dem Korollar existiert dann aber ein $G$-Modul $V$ mit nicht
Cohen-Macaulay Invariantenring $K[V]^{G}$. \qed\\

Wir wollen abschliessend noch eine Begr"undung angeben, warum man mit Lemma
\ref{Hauplemma} die meisten bekannten Beispiele von nicht Cohen-Macaulay
Invariantenringen verstehen kann. Die nicht Cohen-Macaulay Eigenschaft ergibt sich n"amlich in
vielen F"allen aus einem phsop der L"ange $3$, welches keine regul"are
Sequenz ist. (Das folgende Lemma (allgemeiner in \cite[Theorem 1.4]{KemperOnCM}) wird nicht weiter verwendet und kann
ggf. "ubersprungen werden.)

\begin{Lemma}
Sei $K[V]^{G}$ endlich erzeugt und $a_{1},a_{2},a_{3}\in K[V]^{G}$ ein phsop in
$K[V]$. Ist $a_{1},a_{2},a_{3}$ keine regul"are Sequenz in $K[V]^{G}$, so gibt
es ein $0\ne g\in H^{1}(G,K[V])$ mit $a_{i}g=0\in H^{1}(G,K[V]), \, i=1,2,3$.
\end{Lemma}

\Bew Nach Lemma \ref{ZweiReg} sind $a_{1},a_{2}$ teilerfremd in $K[V]$ und
bilden eine regul"are Sequenz in $K[V]^{G}$. Da aber $a_{1},a_{2},a_{3}$
nicht $K[V]^{G}$-regul"ar sind, gibt es
dann also $r_{1},r_{2},r_{3}\in K[V]^{G}$ mit
$r_{1}a_{1}+r_{2}a_{2}+r_{3}a_{3}=0$ und
$r_{3}\not\in(a_{1},a_{2})_{K[V]^{G}}$. Da $K[V]$ Cohen-Macaulay und
$a_{1},a_{2},a_{3}$ dort ein phsop ist, gibt es aber nach Satz \ref{regPhsop} $s_{1},s_{2}\in K[V]$ mit
\[
r_{3}=s_{1}a_{1}+s_{2}a_{2}.
\]
Wir wenden hierauf $(\sigma-1)$ mit $\sigma\in G$ an. Wegen
$r_{3},a_{1},a_{2}\in K[V]^{G}$ erhalten wir
\[
0=a_{1}(\sigma-1)s_{1}+a_{2}(\sigma-1)s_{2}.
\]
Da $a_{1},a_{2}$ in $K[V]$ teilerfremd sind, folgt $a_{2}|(\sigma-1)s_{1}$. Daher
ist
\[
g_{\sigma}:=\frac{(\sigma-1)s_{1}}{a_{2}}=-\frac{(\sigma-1)s_{2}}{a_{1}}\in
K[V] \quad \myforall \sigma\in G.
\]
Offenbar ist $g\in Z^{1}(G,K[V])$ und $a_{1}g,a_{2}g\in B^{1}(G,K[V])$. 
Weiter ist
\begin{eqnarray*}
0&=&r_{1}a_{1}+r_{2}a_{2}+r_{3}a_{3}=r_{1}a_{1}+r_{2}a_{2}+(s_{1}a_{1}+s_{2}a_{2})a_{3}\\&=&(r_{1}+s_{1}a_{3})a_{1}+(r_{2}+s_{2}a_{3})a_{2}.
\end{eqnarray*}
Da $a_{1},a_{2}$ teilerfremd sind, folgt $a_{2}|(r_{1}+s_{1}a_{3})$, also gibt es $h\in K[V]$ mit
\[
r_{1}+s_{1}a_{3}=a_{2}h.
\]
Es folgt
\[
a_{3}g_{\sigma}=\frac{(\sigma-1)(a_{3}s_{1})}{a_{2}}=\frac{(\sigma-1)(r_{1}+a_{3}s_{1})}{a_{2}}=\frac{(\sigma-1)(a_{2}h)}{a_{2}}=(\sigma-1)h
\,\myforall\sigma\in G,
\]
also auch $a_{3}g\in B^{1}(G,K[V])$.
W"are
auch $g\in B^{1}(G,K[V])$, so g"abe es $v\in K[V]$ mit $g_{\sigma}=(\sigma-1)v
\,\myforall \sigma\in G$. Es folgt $(\sigma-1)s_{1}=(\sigma-1)(a_{2}v)$ und
$(\sigma-1)s_{2}=(\sigma-1)(-a_{1}v)$ f"ur alle $\sigma\in G$, also
$s_{1}-a_{2}v,s_{2}+a_{1}v\in K[V]^{G}$. Damit w"are dann doch
\[
r_{3}=s_{1}a_{1}+s_{2}a_{2}=(s_{1}-a_{2}v)a_{1}+(s_{2}+a_{1}v)a_{2}\in (a_{1},a_{2})_{K[V]^{G}},
\]
Widerspruch! \qed\\

\subsection{Der Hauptsatz - eine untere Schranke f"ur den Cohen-Macaulay-Defekt}
Wir leiten in diesem Abschnitt eine untere Schranke f"ur den
Cohen-Macaulay-Defekt her, die es uns erlauben wird, f"ur jede reduktive,
nicht linear reduktive Gruppe einen Invariantenring mit beliebig gro"sem
Cohen-Macaulay-Defekt zu konstruieren, und dies l"asst sich dann sogar mit
Vektorinvarianten erreichen. Der folgende (technische) Satz, zusammen mit
seinen unmittelbaren, weniger technischen Folgerungen und den Anwendungen im
n"achsten Abschnitt mit konkreten Beispielen stellt das Hauptresultat dieser Arbeit dar.

\begin{Hauptsatz} \label{BigMainTheorem}
Sei $G$ eine reduktive Gruppe, $V$ ein $G$-Modul, und es existiere ein $0 \ne g \in
H^{1}(G,K[V])$. Seien weiter  $a_{1},\ldots,a_{k}\in K[V]^{G}$ mit $k\ge 2$ ein phsop in
$K[V]$ mit $a_{i}g=0\in H^{1}(G,K[V])$ f"ur $i=1,\ldots,k$. Dann gilt
\[
\depth(a_{1},\ldots,a_{k})_{K[V]^{G}}=2, \quad \textrm{ also } \cmdef (a_{1},\ldots,a_{k})_{K[V]^{G}}=k-2,
\]
und damit
\[
\cmdef K[V]^{G} \ge k-2.
\]
\end{Hauptsatz}

Wir werden auch die folgende {\bf allgemeinere Formulierung} brauchen: {\it Sei
$G$ eine beliebige lineare algebraische Gruppe, $V$ ein $G$-Modul so, dass
$K[V]^{G}$ endlich erzeugt ist, und  $0 \ne g \in
H^{1}(G,K[V])$. Seien weiter $a_{1},\ldots,a_{k}\in K[V]_{+}^{G}$ homogen mit
$a_{i}g=0\in H^{1}(G,K[V])$. Ist $a_{1},a_{2}$ teilerfremd in $K[V]$, so gilt
$\depth(a_{1},\ldots,a_{k})_{K[V]^{G}}=2$. Ist zus"atzlich
$a_{1},\ldots,a_{k}$ ein phsop in $K[V]^{G}$, so gilt $\cmdef
(a_{1},\ldots,a_{k})_{K[V]^{G}}=k-2$ und damit $\cmdef K[V]^{G} \ge k-2$.}\\

\Bew Nach Lemma \ref{ZweiReg} bildet $a_{1},a_{2}$ eine regul"are Sequenz in $K[V]^{G}$.  
Nach Lemma \ref{Hauplemma} existiert ein $m \in K[V]^{G}$ mit $m\notin
(a_{1},a_{2})_{K[V]^{G}}$, aber $ma_{i}\in (a_{1},a_{2})_{K[V]^{G}}\, \myforall
i=1,\ldots,k$. Satz \ref{SatzShankWehlau} mit $M=R=K[V]^{G}$ liefert dann
$\depth(a_{1},\ldots,a_{k})_{K[V]^{G}}=2$. Da $a_{1},\ldots,a_{k}$ ein phsop in
$K[V]$ ist, ist es wegen Lemma \ref{redphsop} also auch eines in $K[V]^{G}$
(in der allgemeineren Formulierung ist dies eine Voraussetzung),
und es folgt
$\height(a_{1},\ldots,a_{k})_{K[V]^{G}}=k$ (Lemma \ref{phsopHeight}). Damit
gilt nach Definition \ref{DefCMDef} $\cmdef(a_{1},\ldots,a_{k})_{K[V]^{G}}=k-2$, und mit Satz
\ref{cmdefAbsch} also $\cmdef K[V]^{G}\ge k-2$. \qed\\

Der Beweis ist aufgrund der in Abschnitt \ref{Grundlagen} gemachten Vorbereitungen
ziemlich kurz. Ich m"ochte noch einen Beweis angeben, der zwar etwas l"anger ist, weil er im Prinzip den Beweis von Satz
\ref{cmdefAbsch} in einer etwas "ubersichtlicheren Situation mit eingebaut hat, mir aber ebenfalls lehrreich erscheint und
n"aher an meinem urspr"unglichen Beweis f"ur diesen Satz liegt.\\

\noindent\textit{2. Beweis.} Genau wie beim ersten Beweis folgert man zun"achst
$
\depth(a_{1},\ldots,a_{k})_{K[V]^{G}}=2.
$

Sei nun $n=\dim K[V]^{G}$, und man erg"anze $a_{1},\ldots,a_{k}$ zu einem
hsop $a_{1},\ldots,a_{k},a_{k+1},\ldots,a_{n}$ von $K[V]^{G}$. Es gilt dann
\begin{eqnarray*}
\depth K[V]^{G} &=& \depth (a_{1},\ldots,a_{k},a_{k+1},\ldots,a_{n})_{K[V]^{G}} \quad
\textrm{ (Korollar \ref{KorDepthR})}\\
&\le&\depth (a_{1},\ldots,a_{k})_{K[V]^{G}}+(n-k) \quad \textrm{(Korollar
  \ref{depthkKorollar})}\\
&=&2+(n-k),
\end{eqnarray*}
also
\begin{samepage}
\[
k-2\le n-\depth K[V]^{G}=\dim K[V]^{G}-\depth K[V]^{G} = \cmdef K[V]^{G}.
\]
\qed\\
\end{samepage}

Wir ziehen die Korollar \ref{NCMRing} entsprechende Folgerung; F"ur $k=3$
erhalten wir sogar genau das Korollar \ref{NCMRing} zur"uck.

\begin{Korollar} \label{depthToInfty}
Sei $G$ eine reduktive Gruppe, $g\in Z^{1}(G,U)$ ein nichttrivialer Kozyklus,
und $\tilde{U}$ der zugeh"orige erweiterte $G$-Modul (siehe Abschnitt
\ref{firstCohom}). Ist dann
\[
 V^{*}:=U\oplus\bigoplus_{i=1}^{k}\tilde{U}^{*} \quad \textrm{ d.h. } V:=U^{*}\oplus\bigoplus_{i=1}^{k}\tilde{U}
\]
so gilt
\[
\cmdef K[V]^{G}\ge k-2.
\]
\end{Korollar}

\Bew Analog wie der Beweis von Korollar \ref{NCMRing}.
Die $k$ Kopien $a_{1},\ldots,a_{k}\in K[V]^{G}$ des nach Proposition \ref{AnnulatorProp}
existierenden Annullators $\pi \in \tilde{U}^{*G}\setminus\{0\}$ des (eingebetteten)
Kozyklus $0\ne g\in H^{1}(G,K[V])$  erf"ullen die
Voraussetzungen des Hauptsatzes. \qed\\

Damit erhalten wir auch eine neue Charakterisierung linear reduktiver
Gruppen. Zum Vergleich beachte man, dass man Satz \ref{KempersBigTheorem} auch
so formulieren kann:

\textit{Eine reduktive Gruppe $G$ ist genau dann linear reduktiv, wenn $\cmdef
  K[V]^{G}=0$ f"ur jeden $G$-Modul $V$ ist.}

\begin{Korollar}
Eine reduktive Gruppe $G$ ist genau dann linear reduktiv, wenn es eine "`globale
  Cohen-Macaulay-Defekt Schranke"' gibt, d.h. eine Zahl
  $B\in {\mathbb N}$ mit $\cmdef
  K[V]^{G} \le B$ f"ur jeden $G$-Modul $V$.
\end{Korollar}

\Bew Ist $G$ linear reduktiv, so ist der Cohen-Macaulay-Defekt stets $0$. Ist
$G$ reduktiv, aber nicht linear reduktiv, so gibt es einen $G$-Modul $U$ mit
$H^{1}(G,U)\ne 0$ (Proposition \ref{LinRedKohomologie}), so dass die
Voraussetzungen von Korollar \ref{depthToInfty} erf"ullt sind. Setzt man dort
$k>B+2$, so erh"alt man mit dem dortigen Modul $V$ also $\cmdef K[V]^{G}>B$,
also existiert keine globale Cohen-Macaulay-Defekt Schranke $B$. \qed\\

Man kann den Cohen-Macaulay-Defekt sogar mit Vektorinvarianten gegen
unendlich treiben:

\begin{Korollar} \label{depthInfVects}
Ist $G$ reduktiv, aber nicht linear reduktiv, so existiert ein $G$-Modul $V$
mit
\[
\cmdef K\left[\oplus_{i=1}^{k}V\right]^{G}\ge k-2 \quad \myforall k\in{\mathbb N},
\]
insbesondere also
\[
\lim_{k\rightarrow \infty}  \cmdef K\left[\oplus_{i=1}^{k}V\right]^{G} =\infty.
\]
\end{Korollar}

\Bew Da $G$ nicht linear reduktiv ist, gibt es nach Proposition \ref{LinRedKohomologie}  einen $G$-Modul $U$ mit
nichttrivialem Kozyklus wie in Korollar
\ref{depthToInfty} gefordert. Der dortige Beweis bleibt richtig, wenn man
anstatt $U^{*}\oplus\bigoplus_{i=1}^{k}\tilde{U}$ den Modul
$W:=\bigoplus_{i=1}^{k}\left( U^{*}\oplus\tilde{U}\right)$ betrachtet - weitere
nichttriviale Kozyklen in $K[W]$ st"oren ja nicht. Man w"ahlt einen einzigen
der $k$ nichttrivialen Kozyklen aus, und dieser wird dann von den $k$ Kopien der
Invarianten $\pi \in \tilde{U}^{*G}\subseteq S(W^{*})^{G}$  annulliert. Mit
$V:=U^{*}\oplus\tilde{U}$ gilt dann also die Behauptung. \qed\\

\begin{BemRoman}\label{BemerkungNachDepthInvVects} 
Die Frage, ob in den Korollaren \ref{depthToInfty} oder
\ref{depthInfVects} jeweils auch ein \emph{treuer} $G$-Modul $V$ mit den
genannten Eigenschaften existiert, l"asst sich sehr leicht mit "`Ja"'
beantworten: Man nehme einfach einen beliebigen, treuen $G$-Modul $M$ (der ja
stets existiert, Satz \ref{treueDarstellung})  als
Summand hinzu, setze also etwa $V':=V\oplus M$. Nichttriviale Kozyklen bleiben
auch in $K[V']$ nichttrivial, und entsprechendes gilt f"ur annullierende
phsops im Polynomring; und mit $M$ ist nat"urlich auch $V'$ treu. Damit gelten
die Aussagen der beiden Korollare dann auch f"ur $V'$ statt $V$.
\end{BemRoman}

Zum Vergleich geben wir noch das entsprechende Resultat f"ur endliche Gruppen
an:

\begin{Satz}[Gordeev, Kemper] \label{FiniteCMdef}
Sei $\chr K=p>0$, $G$ eine endliche Gruppe, $V$ ein $G$-Modul, $N$ der Kern der Darstellung
$G\rightarrow \GL(V)$ und $p$ ein Teiler des Index $(G:N)$.  Dann gilt
\[
\lim_{k\rightarrow \infty}  \cmdef K\left[\oplus_{i=1}^{k}V\right]^{G} =\infty.
\]
Genauer: Ist $r_{0}>0$ mit $H^{r_{0}}(G/N,K)\ne 0$ (ein solches $r_{0}$ gibt
 es stets) oder $H^{r_{0}}(G/N,V^{*})\ne 0$ oder $H^{r_{0}}(G/N,K[V])\ne 0$, so gilt
\[
\cmdef K\left[\oplus_{i=1}^{k}V\right]^{G} \ge k-r_{0}-1\quad \myforall k \in
{\mathbb N}.
\]
\end{Satz}

Es gilt sogar folgende {\bf Verallgemeinerung:} {\it Es erf"ulle $V_{1}:=V$
  obige Voraussetzungen, wobei zus"atzlich $V$ \emph{treu} sei. Ist $(V_{i})_{i\ge
  2}$ eine beliebige Folge \emph{treuer} $G$-Moduln, so gilt
\[
\cmdef K\left[\oplus_{i=1}^{k}V_{i}\right]^{G} \ge k-r_{0}-1\quad \myforall k \in
{\mathbb N}.
\]}
Man beachte, dass $G/N$ treu auf $V$ operiert und $K[V]^{G/N}=K[V]^{G}$
gilt. Der Satz gibt damit das bestm"ogliche Resultat, denn wenn $p$ kein
Teiler von $(G:N)$ ist,
so kommt die $G$-Operation letztlich von der linear reduktiven Gruppe $G/N$,
und dann gilt nat"urlich nach Hochster und Roberts $\cmdef K\left[\oplus_{i=1}^{k}V\right]^{G}=0$ f"ur
alle $k$.

In Gordeev, Kemper  \cite[Corollary 5.15]{Gordeev} wird der Satz aus
einem allgemeineren Resultat gefolgert; Wir geben einen etwas elementareren
Beweis, der auf den Resultaten der Arbeit \cite{KemperOnCM} basiert.\\

\Bew  Nach dem
Vorspann k"onnen wir annehmen, dass $G$ treu auf $V$ operiert und $p$ ein
Teiler der Gruppenordnung $|G|$ ist. Es gen"ugt dann offenbar, die
Verallgemeinerung zu beweisen (man kann ja $V_{i}:=V\,\, \myforall i$ w"ahlen). Setze $R:=K\left[\oplus_{i=1}^{k}V_{i}\right]$.
Nach \cite[Theorem
4.1.3]{BensonCohomology2} existiert wegen $p\mid |G|$  ein $r_{0}>0$, so dass die $r_{0}$-te
Kohomologiegruppe $H^{r_{0}}(G,K)\ne 0$ ist. Ist diese oder eine der anderen
beiden Bedingungen erf"ullt, so ist auch
$H^{r_{0}}(G,R)\ne 0$ (denn es sind $K$ bzw. $V^{*}$ bzw. $K[V]$ jeweils
direkte, $G$-stabile Summanden von $R$ mit $G$-stabilem Komplement). Insbesondere gibt es dann ein minimales $1\le r \le
r_{0}$ mit der Eigenschaft $H^{r}(G,R)\ne 0$. Man kann dann $0\ne g\in
H^{r}(G,R)$ w"ahlen, so dass der zugeh"orige nichttriviale $r$-Kozyklus nur
Werte in einer homogenen Komponente annimmt (indem man etwa eine geeignete homogene
Komponente eines gegebenen nichttrivialen Kozyklus $g$ w"ahlt - da $G$ endlich, sind nur endlich
viele Projektionen von $g$ auf eine homogene Komponente ungleich Null, und
w"aren alle solchen Projektionen trivial, so auch $g$). Insbesondere ist dann das
"`Annullationsideal"'
\[
I_{g}:=\left\{a\in R^{G}: ag=0 \in H^{r}(G,R)\right\} \lhd R^{G}
\]
homogen. Nach \cite[Corollary 1.6]{KemperOnCM} gilt dann
\begin{equation} \label{depthIgr}
\depth I_{g}\le r+1 \le r_{0}+1.
\end{equation}
F"ur $\sigma \in G$ sei $A_{i,\sigma}$ die zugeh"orige lineare Abbildung auf $V_{i}$, und $B_{\sigma}$
die zugeh"orige lineare Abbildung auf $\oplus_{i=1}^{k}V_{i}$. Da $G$ treu auf
$V_{i}$ operiert, gilt f"ur jedes
$\sigma \ne \iota$ ($\iota\in G$ neutrales Element), dass $\rang
\left(A_{i,\sigma}-\id_{V_{i}}\right)\ge 1$, und damit $\rang
\left(B_{\sigma}-\id_{\oplus_{i=1}^{k}V_{i}}\right) \ge k$. Nach \cite[Lemma 2.1,
Lemma 1.7]{KemperOnCM} gibt es dann $a_{1},\ldots,a_{k}\in I_{g}$ mit
$\height(a_{1},\ldots,a_{k})_{R^{G}}=k$, und es folgt $\height(I_{g})\ge
k$. Mit Gleichung \eqref{depthIgr} folgt also
\[
\cmdef I_{g}=\height(I_{g})-\depth(I_{g})\ge k -r_{0}-1,
\]
und aus Satz \ref{cmdefAbsch} dann $\cmdef R^{G}\ge k-r_{0}-1$.  Hieraus
folgen die Behauptungen. \qed\\

Wir geben eine h"ubsche Anwendung. 

\begin{Korollar}[Campbell, Geramita, Hughes, Kemper, Shank,
  Wehlau]\label{depthZpn} $\textrm{ }$\\ Sei $K$ ein algebraisch abgeschlossener
 K"orper der Charakteristik $p>0$ und $G$ eine endliche Gruppe, die einen
 Normalteiler vom Index $p$ enth"alt (z.B. $G$ eine $p$-Gruppe, oder $G$
 nilpotent und $p||G|$). Dann gilt f"ur jeden treuen $G$-Modul
 $V$
\[
\cmdef K\left[\oplus_{i=1}^{k}V\right]^{G} \ge k-2,
\]
insbesondere ist der Invariantenring f"ur $k\ge 3$ nicht Cohen-Macaulay.
\end{Korollar}

Ist $G$ eine $p$-Gruppe, so existiert eine Unterguppe $N\subseteq G$ mit
$(G:N)=p$, und eine solche ist automatisch auch Normalteiler vom Index
$p$. Ist $G$ nilpotent, also direktes Produkt ihrer Sylow-Untergruppen, so
ist das direkte Produkt eines Normalteilers der $p$-Sylowgruppe vom Index $p$
(existiert nach eben)
und der anderen Sylowgruppen ein Normalteiler von $G$ vom Index $p$.\\

\Bew (vgl. \cite[Corollary 21]{CampEtAl} sowie \cite[Theorem
1.2]{CampNonCMVectInvs} und \cite[Theorem 2.3]{KemperOnCM}).
Sei ${\mathbb F}_{p}\subseteq K$ der Primk"orper von $K$. 
Sei $N\lhd G$ ein Normalteiler mit $(G:N)=p$. Dann ist $G/N\cong Z_{p}\cong
({\mathbb F_{p}},+)\subseteq (K,+)$. ($Z_{p}$ ist die zyklische Gruppe der
Ordnung $p$). Wir schreiben $G$ multiplikativ und bezeichnen f"ur $a\in G$ mit $\bar{a}\in K$ das $aN\in G/N$ entsprechende
Element gem"a"s obiger Einbettung. Dann ist $\overline{ab}=\bar{a}+\bar{b}$
($a,b\in G$).

Offenbar ist durch $g: G\rightarrow K, a\mapsto g_{a}=\bar{a}$ ein
nichttrivialer Kozyklus gegeben ($K$ tr"agt triviale $G$-Operation, hier
bezeichnet mit $\circ$): Zum einen
ist
\[
g_{ab}=\overline{ab}=\bar{a} + a \circ \bar{b} = a\circ g_{b} + g_{a}\quad
\myforall a,b\in G,
\] 
also $g$ ein Kozyklus, und f"ur alle $v\in K$ und $a\in G$ ist $(a-1)\circ
v=v-v=0$, also $g$ nichttrivial. Dies zeigt $H^{1}(G,K)\ne 0$, und mit
$r_{0}=1$ in Satz \ref{FiniteCMdef} folgt die Behauptung. \qed\\

Im Gegensatz zu Satz \ref{FiniteCMdef} besagt Korollar \ref{depthInfVects}  nur, dass "uberhaupt eine
Darstellung existiert, deren Cohen-Macaulay-Defekt der Vektorinvarianten
gegen unendlich geht. Man k"onnte vermuten, dass dennoch die entsprechende
Aussage f"ur treue rationale Darstellungen reduktiver, nicht linear reduktiver Gruppen
gilt, aber dem ist nicht so. Es gilt n"amlich:

\begin{Satz} \label{VectInvsAreCM}
Sei $\chr K$ beliebig und $V=K^{n}$. Ist dann $G$ eine der Gruppen
$\SL_{n}(K), \GL_{n}(K)$ oder $\SO_{n}(K)$ (dann $\chr K\ne 2$), so ist $K\left[\oplus_{i=1}^{k}V\right]^{G}$
Cohen-Macaulay f"ur alle $k\in {\mathbb N}$.
\end{Satz} 

\Bew F"ur die Gruppe $\SO_{n}$ siehe die erst in K"urze erscheinende
Arbeit \cite{Raghavan}. F"ur $G=\SL_{n}$ wird
$K\left[\oplus_{i=1}^{k}V\right]^{G}$ nach de Concini und Procesi
\cite{ConciniProcesi} auch in positiver Charakteristik von den sogenannten
"`Pl"ucker-Invarianten"' erzeugt, und bildet daher den Koordinatenring einer
Grassman-Variet"at. Ein solcher ist nach Hochster \cite[Corollary
3.2]{Hochster} Cohen-Macaulay. (Vgl. auch Hochster und Eagon \cite[p.1029,
mitte]{HochsterEagon}). F"ur $G=\GL_{n}$ gilt dagegen nach
\cite{ConciniProcesi} stets $K\left[\oplus_{i=1}^{k}V\right]^{G}=K$, so dass
diese Invariantenringe trivialerweise ebenfalls stets Cohen-Macaulay sind. \qed\\ 

Die Gruppen $\SL_{n}, \GL_{n}$ mit $n\ge 2$ und $\SO_{n}$ mit $n\ge 3$ und
$p\ne 2$ sind in positiver Charakteristik zwar
reduktiv, aber nicht linear reduktiv (siehe Bemerkung \ref{ExpliziteBspFuerKlassGruppen}), und operieren treu auf $V=K^{n}$. Dennoch
sind die zugeh"origen Vektorinvarianten nach diesem Satz stets
Cohen-Macaulay. Daher l"asst sich Satz \ref{FiniteCMdef} nicht auf unendliche
Gruppen verallgemeinern.\\

Wir haben im Beweis von Satz \ref{FiniteCMdef} ein Resultat "uber die Tiefe
des sogenannten Annullationsideals eines $r$-Kozyklus (f"ur endliche Gruppen) verwendet. Mit unseren
Methoden k"onnen wir dieses Resultat f"ur $r=1$ auch auf unendliche Gruppen
verallgemeinern. Nochmal die Definition f"ur den Fall $r=1$: 
 Ist $g\in H^{1}(G,K[V])$, so wird
das \emph{Annullationsideal}\index{Annullationsideal} definiert als
\[
I_{g}:=\left\{a\in K[V]^{G}: ag=0 \in H^{1}(G,K[V])\right\},
\]
und man sieht sofort, dass dies tats"achlich ein Ideal in $K[V]^{G}$ ist. Wir nennen einen Kozyklus $g\in Z^{1}(G,K[V])$ (bzw. seine Restklasse in
$H^{1}(G,K[V])$ \emph{homogen},\index{Kozyklus!homogen} wenn es ein $d\in
{\mathbb N}_{0}$ gibt mit $g_{\sigma}\in K[V]_{d}\, \myforall \sigma \in G$. F"ur
einen homogenen Kozyklus ist offenbar $I_{g}$ ein homogenes Ideal. Der
folgende Satz ist eine Version des Hauptsatzes (mit praktisch demselben Beweis) mit schw"acheren Forderungen,
die entsprechend auch eine schw"achere Aussage liefert.

\begin{Satz} \label{depthKillerideal}
Sei $G$ eine (nicht notwendig reduktive!) lineare algebraische
  Gruppe und  $V$ ein $G$-Modul, so dass $K[V]^{G}$ endlich erzeugt ist. Sei
  ferner $0\ne g\in H^{1}(G,K[V])$ \emph{homogen}. Falls $I_{g}$
  zwei homogene, in $K[V]$ teilerfremde Elemente positiven Grades enth"alt, so
  gilt
\[
\depth(I_{g})=2.
\]
\end{Satz}

\Bew Aufgrund der Voraussetzungen ist $I_{g}\ne K[V]^{G}$ ein homogenes Ideal, denn
$g\ne 0$ (also $1\notin I_{g}$) und $g$ ist homogen. Seien $a_{1},a_{2}\in I_{g}$
homogen, positiven Grades und teilerfremd. Nach Lemma \ref{ZweiReg} bilden
$a_{1},a_{2}$ eine regul"are Sequenz in $I_{g}$, und nach Lemma
\ref{Hauplemma} gibt es ein $m\in K[V]^{G}$ mit $m \notin
(a_{1},a_{2})_{K[V]^{G}}$ aber $am \in(a_{1},a_{2})_{K[V]^{G}}$ f"ur alle
$a\in I_{g}$. Nach Satz \ref{SatzShankWehlau} (der "`"aquivalenten Formulierung"') folgt $\depth I_{g}=2$. \qed\\

Kennt man zus"atzlich eine untere Schranke f"ur $\height I_{g}$,  $\height
I_{g}\ge k$, etwa weil man ein in $I_{g}$ liegendes phsop von
$K[V]^{G}$ der L"ange $k$  kennt, so erh"alt man mit Satz \ref{cmdefAbsch} die Absch"atzung
\[
\cmdef R \ge \height(I_{g})-\depth(I_{g})\ge k-2.
\]
Insbesondere erh"alt man so auch nochmal die Aussage des Hauptsatzes, wenn $G$
reduktiv ist und man ein annullierendes phsop in $K[V]$ hat.\\

Man vergleiche diesen Satz mit dem entsprechenden Resultat f"ur endliche
  Gruppen (nur den Spezialfall f"ur
  die erste Kohomologie), Kemper \cite[Corollary 1.6]{KemperOnCM}, welches wir
  (f"ur allgemeines $r$) im Beweis von Satz \ref{FiniteCMdef} verwendet haben: 

\begin{Korollar}
Ist $G$ eine endliche Gruppe, $0 \ne g \in H^{1}(G,K[V])$ homogen, so ist 
\[
\depth(I_{g})=\min\{2,\height(I_{g})\}.
\]
\end{Korollar}

\Bew Im Fall $\height(I_{g})=0$ ist $I_{g}=\{0\}$, also auch $\depth(I_{g})=0=\min\{2,0\}$. Im
Fall $\height(I_{g})=1$ ist $I_{g}\ne \{0\}$. Da $K[V]^{G}$ nullteilerfrei, ist dann
auch $1\le \depth(I_{g})\le \height(I_{g})=1$, also $\depth(I_{g})=1=\min\{2,1\}$. Ist nun
$\height(I_{g})\ge 2$, so enth"alt $I_{g}$ ein phsop $a_{1},a_{2}$ von $K[V]^{G}$ der L"ange $2$ (Lemma
\ref{phsopHeight} (b)). Da $K[V]$ eine ganze Erweiterung des normalen Rings
$K[V]^{G}$ (\cite[Proposition 2.3.11]{DerksenKemper}) ist, ist $a_{1},a_{2}$ nach Korollar \ref{ganzPhsop} auch ein phsop
in $K[V]$, also dort teilerfremd (Lemma ~\ref{ZweiReg}). Mit Satz
\ref{depthKillerideal} folgt also
$\depth(I_{g})=2=\min\{2,\height(I_{g})\}$. \qed\\

Die Verallgemeinerung des Korollars auf unendliche Gruppen scheitert (zumindest
bei diesem Beweis) offenbar
daran, dass es dann phsops in $K[V]^{G}$ gibt, die kein phsop in $K[V]$
bilden. Siehe etwa \cite[Remark 5]{KemperLinRed} f"ur ein Beispiel.

\newpage
\section{Anwendungen des Hauptsatzes}\label{Anwendungen}
Im Folgenden wollen wir konkret Darstellungen $V$ f"ur eine reduktive, nicht
linear reduktive Gruppe $G$ angeben, die die
Voraussetzungen des Hauptsatzes \ref{BigMainTheorem} mit beliebig vorgegebenem $k$
erf"ullen, f"ur die also $\cmdef K[V]^{G}\ge k-2$ gilt. Mit den Ergebnissen
aus Abschnitt \ref{sectionDimension} l"asst sich oft auch $\dim K[V]^{G}$ exakt
angeben, so dass wir eine Absch"atzung $\depth K[V]^{G}\le \dim K[V]^{G}-k+2$ erhalten.
F"ur die allgemeine
Konstruktion m"ussen wir dabei nur die Ergebnisse der letzten Abschnitte
zusammensetzen, was allerdings zu einer gro"sen Vektorraumdimension des
Darstellungsmoduls $V$ f"uhrt. Wir f"uhren daher auch verfeinerte
Betrachtungen durch, die zu teils erheblich niedrigeren Dimensionen
f"uhren. Dabei verwenden wir Roberts' Isomorphismus, der au"serdem auch Beispiele f"ur eine Klasse von nicht reduktiven Gruppen liefert.

\subsection{Motivation: Additive Gruppen endlicher K"orper}
Sei $K$ ein algebraisch abgeschlossener K"orper der Charakteristik $p$. F"ur
jede Potenz $q$ von $p$ ist
\[
G:=\{a \in K: a^{q}-a =0\}\subseteq K^{1}
\]
zusammen mit der Addition "`$+$"' von $K$ eine lineare algebraische Gruppe,
die isomorph ist zur additiven Gruppe $({\mathbb F}_{q},+)$ des endlichen
K"orpers mit $q$ Elementen. Sei $\langle X,Y \rangle$ der $G$-Modul mit der
nat"urlichen Darstellung
\[
a\mapsto \left(
\begin{array}{cc}
1&a\\
&1
\end{array}
\right).
\]
Dabei ist $\langle X,Y \rangle$ selbstdual, genauer ist $\langle X,Y
\rangle^{*}=\langle X^{*},Y^{*} \rangle \cong \langle Y,-X \rangle$.
Wir betrachten nun die $n$-fache direkte Summe der nat"urlichen Darstellung,
\[
V^{*}:=\bigoplus_{i=1}^{n}\langle X_{i},Y_{i} \rangle. 
\]
Als n"achstes ben"otigen wir einen nichttrivialen Kozyklus in $K[V]\cong
S(V^{*})$. Ein solcher wird gegeben durch
\[
g: G \rightarrow K[V]=S(V^{*}),\quad a \mapsto a\cdot 1_{K[V]}\in K[V]_{0}.
\]
Dabei ist $g$ ein Kozyklus in $Z^{1}(G,K[V])$, denn es gilt $g_{a\circ b}=g_{a+b}=(a+b)\cdot 1_{K[V]}=g_{b}+g_{a}=a\cdot
g_{b}+g_{a}$ (die Operation von $G$ auf $K[V]_{0}=K$ ist trivial). Au"serdem
ist $g$ auch nichttrivial, denn wenn es ein $v\in K[V]$ g"abe mit
$g_{a}=(a-1)v\,\, \myforall a\in G$, dann k"onnte man $v$
 auch homogen vom Grad $0$ w"ahlen, also $v\in K=K[V]_{0}$. Dann ist aber
$a\cdot v -v =v-v=0$, im Widerspruch etwa zu $g_{1}=1$. 

Weiter sind die $X_{i}, i=1,\ldots,n$ Annullatoren des Kozyklus $g$, denn es gilt
\[
X_{i}g_{a}=aX_{i}=(a-1)\cdot Y_{i}\quad \myforall a \in G,
\]
also $X_{i}g\in B^{1}(G,K[V])$.
Da $X_{1},\ldots,X_{n}$ als Variablen au"serdem ein phsop bilden, liefert der
Hauptsatz \ref{BigMainTheorem} also sofort
\begin{equation} \label{FiniteAdditiveGroups}
\cmdef K[V]^{G}=\cmdef K\left[\bigoplus_{i=1}^{n}\langle X_{i},Y_{i}
\rangle\right]^{({\mathbb F}_{q},+)} \ge n-2.
\end{equation}
(Wir h"atten dies nat"urlich auch aus Korollar \ref{depthZpn} folgern k"onnen).
Da $G$ eine endliche Gruppe ist, gilt $\dim K[V]^{G}=\dim_{K} V = 2n$, und
damit
\[
\depth K[V]^{G}=\dim K[V]^{G}-\cmdef K[V]^{G}\le 2n - (n-2)=n+2.
\]
Die rechte Seite ist offenbar gleich $\dim_{K} V^{G}+2$.
Nach dem folgenden Resultat, das auch noch etwas allgemeiner gilt, gilt hier sogar
Gleichheit (dabei ist $V^{\sigma}=\{v\in V: \sigma v = v\}$):

\begin{Satz}[Campbell, Ellingsrud, Hughes, Kemper, Shank, Skjelbred, Wehlau]
  \label{EllingsrudSkjelbred} $\textrm{ }$\\
Sei $\chr K=p>0$, $G$ eine \emph{$p$-Gruppe} und $V$ ein $G$-Modul. Falls $G$
\emph{nicht} von den Elementen $\sigma \in G$ mit $\dim V^{\sigma} > \dim
V^{G}$ erzeugt wird (z.B. wenn $G$ eine zyklische $p$-Gruppe ist), dann gilt
\[
\depth K[V]^{G}=\min\{\dim_{K}V^{G}+2, \dim_{K}V\}.
\]
\end{Satz}\index{Satz!Ellingsrud, Skjelbred}

\Bew Siehe Campbell, Hughes, Kemper, Shank, Wehlau \cite{CampEtAl}, Theorem 1, Theorem
2 und Remark 4 (a). Im zyklischen Fall folgt aus $\dim V^{\sigma}>\dim V^{G}$
jedenfalls, dass $\sigma$ nicht $G$ erzeugt, und daher in der eindeutig
bestimmten Untergruppe vom Index $p$ liegt, d.h. alle solchen $\sigma$ k"onnen
h"ochstens diese Untergruppe erzeugen. Der zyklische Fall wurde bereits in der
bahnbrechenden Arbeit von Ellingsrud und Skjelbred
\cite{Ellingsrud} behandelt. \qed\\

In unserem Fall $V\cong V^{*}=\bigoplus_{i=1}^{n}\langle X_{i},Y_{i}
\rangle$ ist offenbar $\dim V^{G}=n$, und das einzige $\sigma \in G$ mit $\dim
V^{\sigma}>\dim V^{G}$ ist das neutrale Element $\sigma=\iota$, von welchem $G$
nicht erzeugt wird. Damit ist der Satz anwendbar und liefert $\depth
K[V]^{G}=\min\{n+2,2n\}$. Insbesondere gilt dann in
\eqref{FiniteAdditiveGroups} Gleichheit, wenn $n\ge 2$ ist.\\

\label{MotivationsBeispiel}
Was passiert, wenn man $q$ "`gegen unendlich"' gehen l"asst, also zu
der gesamten additiven Gruppe $\Ga=(K,+)$ des K"orpers "ubergeht? $\Ga$ ist
nicht reduktiv, und der Hauptsatz nicht anwendbar. Obwohl die Gleichung
\eqref{FiniteAdditiveGroups}  unabh"angig von $q$ ist, gilt sie nun nicht
mehr. Stattdessen gilt stets (in beliebiger Charakteristik)
\[
\cmdef K\left[\bigoplus_{i=1}^{n}\langle X_{i},Y_{i} \rangle\right]^{\Ga}=0,
\]
denn nach Roberts' Isomorphismus \ref{RobertsIsom}  gilt ja die Isomorphie
\[
K\left[\bigoplus_{i=1}^{n}\langle X_{i},Y_{i} \rangle\right]^{\Ga}\cong K\left[\bigoplus_{i=1}^{n+1}\langle X_{i},Y_{i}
  \rangle\right]^{\SL_{2}}
\]
 (wobei auf der rechten Seite $\langle X_{i},Y_{i} \rangle$ die
nat"urliche Darstellung der $\SL_{2}$ ist), und der rechte Invariantenring ist
nach Satz \ref{VectInvsAreCM} Cohen-Macaulay\footnote{Der Invariantenring ist
  also nicht nach Hochster und Roberts Cohen-Macaulay, sondern nach Hochster
  und Roberts' Isomorphismus!}. Das Argument des Hauptsatzes geht
dabei deshalb schief, weil dann $X_{1},\ldots,X_{i}$ f"ur $i\ge 3$ kein phsop
mehr im Invariantenring bildet. 

Wir werden in Satz \ref{ZweitesHauptResultat} sehen, wie man durch Hinzunahme des Summanden $\langle
X_{0}^{p},Y_{0}^{p} \rangle$ zu $V$ die Ungleichung
\eqref{FiniteAdditiveGroups} auch f"ur die volle additive Gruppe $\Ga$ (in
positiver Charakteristik $p$) retten kann.

\subsection{Konstruktion von Darstellungen mit gro"sem
  Cohen-Macaulay-Defekt des Invariantenrings}\label{ExpliziterKonstruktion}
Wir geben nun f"ur jede reduktive Gruppe, f"ur die das "uberhaupt m"oglich ist,
\emph{explizit} eine Folge von Darstellungen an, f"ur die der
Cohen-Macaulay-Defekt des zugeh"origen Invariantenrings beliebig gro"s wird.

\begin{Satz} \label{ExpliziteDarstellungen}
Sei $G$ eine reduktive, aber nicht linear
reduktive Gruppe in positiver Charakteristik $p=\chr K$. Wir unterscheiden
zwei F"alle:
\begin{enumerate}
\item Falls die Zusammenhangskomponente $G^{0}$ kein Torus ist, so sei 
  $V$ ein treuer $G$-Modul (existiert nach Satz \ref{treueDarstellung}),
  \[U:=\Hom_{K}(S^{p}(V),F^{p}(V))_{0}\] (siehe Definition \ref{HomKVW0}) und $\iota \in
  \Hom_{K}(S^{p}(V),F^{p}(V))$ mit $\iota|_{F^{p}(V)}=\id_{F^{p}(V)}$. 

\item Falls die Zusammenhangskomponente $G^{0}$ ein Torus ist, so teilt $p$ die
  Ordnung der endlichen Gruppe $H:=G/G^{0}$. (Dies
  ist insbesondere dann der Fall, wenn $G$ endlich ist. Dann ist $|G^{0}|=1$ und H=G). Jeder
  $H$-Modul ist auch $G$-Modul mittels des kanonischen Homomorphismus
  $G\rightarrow H=G/G^{0}$. Sei
  dann $V:=\oplus_{\sigma \in H} Ke_{\sigma}$ die regul"are Darstellung von
  $H$, also $\sigma e_{\tau}=e_{\sigma\tau} \,\, \myforall \sigma,\tau \in H$. Mit
  $e:=\sum_{\sigma\in H}e_{\sigma}$ betrachte dann \[U:=\Hom_{K}(V,Ke)_{0}\] und
  $\iota \in \Hom_{K}(V,Ke)$ mit $\iota|_{Ke}=\id_{Ke}$.
\end{enumerate}
Es sei dann jeweils 
  $g: G\rightarrow U, \, \sigma\mapsto (\sigma-1)\iota$ der zugeh"orige
  Kozyklus $g\in Z^{1}(G,U)$. Weiter sei dann $\tilde{U}$ der zu $g$ und $U$ geh"orige
  erweiterte $G$-Modul (siehe Abschnitt \ref{firstCohom}), d.h. es ist
  \[\tilde{U}=U\oplus K\cdot \iota.\]

Dann gilt
\[
\cmdef K\left[U^{*}\oplus \bigoplus_{i=1}^{k}\tilde{U}\right]^{G} \ge k-2.
\]
Insbesondere ist der Invariantenring f"ur $k\ge 3$ nicht Cohen-Macaulay.
\end{Satz}

\Bew Da $G$ nicht linear reduktiv ist, ist nach Satz \ref{NagataTorus}
entweder $G^{0}$ kein Torus oder falls doch, $p$ ein Teiler von
$|G/G^{0}|$. Daher tritt einer der beiden F"alle auf. Wenn gezeigt ist, dass der
Kozyklus $g$ nichttrivial ist, folgt die Behauptung des Satzes dann sofort
aus Korollar \ref{depthToInfty}.

Im 1. Fall hat $F^{p}(V)$ kein Komplement in $S^{p}(V)$ nach Korollar
\ref{NagatasNoComplementCorollar}.

Im 2. Fall hat $Ke$ kein Komplement in $V$ nach Korollar
\ref{MaschkeNoComplement}.

\noindent Daher ist in beiden F"allen nach Proposition \ref{kompl} der Kozyklus $g$ nichttrivial. Dies war zu zeigen.\qed\\

Im 2. Fall kann man alternativ auch Satz \ref{FiniteCMdef} zur Konstruktion
verwenden. Es ist z.B. $U\cong (V/Ke)^{*}$ (siehe Satz \ref{HomKVW0}) f"ur
$|H|\ne 2$ ein treuer $H$-Modul mit
$H^{1}(H,U)\ne 0$, und mit $r_{0}=1$ in Satz \ref{FiniteCMdef} folgt $\cmdef K\left[\bigoplus_{i=1}^{k}U^{*}\right]^{H}\ge k-2$.\\

\noindent Wir berechnen noch die Dimension der Darstellung $M=U^{*}\oplus
\bigoplus_{i=1}^{k}\tilde{U}$.

{\it 1. Fall:} Es ist $\dim_{K}S^{p}(V)={n+p-1 \choose p}$,
wobei $n=\dim_{K}V=\dim_{K}F^{p}(V)$, also
\[\dim_{K}U=\dim_{K}F^{p}(V)(\dim_{K}S^{p}(V)-\dim_{K}F^{p}(V))=n\left({n+p-1 \choose p}-n \right).\]
Da $\dim_{K}\tilde{U}=\dim_{K}(U)+1$, ist $\dim_{K}M=(k+1)\dim_{K}U +k$, also
\begin{equation}\label{DimOfM1}
\dim_{K} M=n(k+1)\left({n+p-1 \choose p}-n
\right)+k.
\end{equation}
Insbesondere f"ur $n=2$ ist dies
\[
\dim_{K} M=2(k+1)(p-1)+k.
\]

{\it 2. Fall:} Sei $m:=|G/G^{0}|$. Es ist $\dim_{K}V=m$, $\dim_{K}Ke=1$ und
damit $\dim_{K}U=1\cdot(m-1)$, also
\begin{equation}\label{DimOfM2}
\dim_{K} M=(k+1)(m-1)+k.
\end{equation}\\

\begin{BemRoman}\label{smallCocycles} Man kann die Dimension von $M$ etwas dr"ucken, denn es ist
$U$ von der Form $\Hom_{K}(V,W)_{0}=W\otimes (V/W)^{*}$ (siehe Satz
\ref{HomKVW0}). Daher kann man den Summand $U^{*}$ von $M$ durch $W^{*}\oplus
(V/W)$ ersetzen. Denn dann ist $K[M]=S(M^{*})=S(W\oplus (V/W)^{*}\oplus\ldots)$,
und wegen \[
U\cong W\otimes (V/W)^{*} \le S^{2}(W\oplus (V/W)^{*})\cong S^{2}(W)\oplus
S^{2}((V/W)^{*}) \oplus W\otimes (V/W)^{*}
\]
 hat
auch hier $S(M^{*})$ einen nichttrivialen Kozyklus, der nun Werte in der
zweiten symmetrischen Potenz annimmt. Der Beweis von Korollar
\ref{depthToInfty} gilt dann weiterhin mit dieser Ersetzung. Im ersten Fall
von Satz \ref{ExpliziteDarstellungen}
ist au"serdem $W^{*}=F^{p}(V)^{*}$ mit $V$ ein treuer $G$-Modul, der damit als
direkter Summand in $M$ auftaucht, so dass dann $M$ auch treu ist; Der Modul
$U$ ist dagegen nicht immer treu. Im 2. Fall allein deshalb nicht, weil es
sich im Allgemeinen um die Darstellung einer Faktorgruppe handelt. Ein Beispiel f"ur
den 1. Fall mit $U$ nicht treu ist gegeben durch $G=\SL_{2}$ mit
nat"urlicher Darstellung $V=\langle X,Y \rangle$ in Charakteristik~$3$. Wir werden in Abschnitt
\ref{SL2char3} sehen, dass dann $U\cong \langle X,Y \rangle \otimes\langle
X^{3},Y^{3} \rangle$ ist, d.h. die Darstellung von $U$ ist gegeben durch
\[
\sigma=\left(\begin{array}{cc}a&b\\c&d\end{array}\right) \mapsto
\left(\begin{array}{cc}a&b\\c&d\end{array}\right)\otimes
\left(\begin{array}{cc}a^{3}&b^{3}\\c^{3}&d^{d}\end{array}\right)\in
K^{4\times 4}
\]
(Kronecker Produkt von Matrizen), insbesondere wird das negative der
Einheitsmatrix $\sigma=-I_{2}\in \SL_{2}$ auf die Identit"at $I_{4}\in
K^{4\times 4}$ abgebildet, d.h. $U$ ist nicht treu. Auch $\tilde{U}$ ist hier
nicht treu, wie die Darstellung \eqref{US3} (S. \pageref{US3}) (wieder f"ur $\sigma=-I_{2}\in \SL_{2}$)
zeigt; Siehe auch Bemerkung~\ref{BemerkungNachDepthInvVects} f"ur eine andere L"osung
dieses Problems.
\end{BemRoman}

Wir formulieren noch kurz das entsprechende explizite Resultat f"ur
Vektorinvarianten. Dieses ergibt sich sofort aus Korollar \ref{depthInfVects}
(bzw. dem Beweis)
und Satz \ref{ExpliziteDarstellungen} zusammen mit der Bemerkung. Wir
beschr"anken uns auf den 1. Fall, der 2. Fall wird genauso behandelt.

\begin{Korollar}
Seien die Voraussetzungen des 1. Falls von Satz \ref{ExpliziteDarstellungen}
erf"ullt. Dann ist
\[
W:=F^{p}(V)^{*}\oplus \left(S^{p}(V)/F^{p}(V)\right) \oplus \tilde{U}
\]
ein \emph{treuer} G-Modul mit
\[
\cmdef K\left[\bigoplus_{i=1}^{k}W\right]^{G}\ge k-2.
\]
\end{Korollar}

Der zweite Fall von Satz \ref{ExpliziteDarstellungen} handelt im wesentlichen
von endlichen Gruppen.  F"ur diese gibt es in der Literatur bessere Ergebnisse
als unser Satz liefert, siehe etwa Satz \ref{EllingsrudSkjelbred} f"ur $p$-Gruppen. Interessant sind daher f"ur uns die
Gruppen, f"ur die der erste Fall zutrifft, f"ur die also die
Zusammenhangskomponente kein Torus ist. Insbesondere trifft dies f"ur
praktisch alle Klassischen Gruppen zu. Hier die vollst"andige "Ubersicht der f"ur
die Klassischen Gruppen auftretenden F"alle:

\begin{Bemerkung} \label{ExpliziteBspFuerKlassGruppen}
Sei $\chr K=p,\, n\ge 2$. 
\begin{enumerate}
\item Die Gruppen $\SL_{n}, \, \GL_{n},\, \Sp_{n}$ sind
zusammenh"angend  ($n$ gerade f"ur $\Sp_{n}$) und kein
Torus. 
\item F"ur $p\ne 2, n\ge 3$ ist $\SO_{n}$ zusammenh"angend, kein Torus und die
Zusammenhangskomponente von $\On_{n}$. 
\item F"ur $p\ne 2,\, n=2$ ist $\SO_{2}$ ein Torus, die
  Zusammenhangskomponente von $\On_{2}$ und es ist $\On_{2}/\SO_{2}=Z_{2}$;
  insbesondere sind beide Gruppen linear reduktiv.
\item F"ur $p=2$ sind die orthogonalen Gruppen speziell definiert (siehe der
  folgende Beweis), und es gilt $\SO_{n}=\On_{n}$. Es ist dann f"ur $n\ge 3$ die Zusammenhangskomponente von
  $\SO_{n}$ kein Torus. Ferner ist $\SO_{2}$ semidirektes Produkt eines
  Torus mit der $Z_{2}$, und damit nicht linear reduktiv. 
\end{enumerate}
Insbesondere also liefert Satz
\ref{ExpliziteDarstellungen} f"ur jede der in $1,2,4$ genannten Gruppen $G$
im Fall $p>0$ zu gegebenem $k$
\emph{explizit} eine Darstellung $V$ mit $\cmdef K[V]^{G}\ge k-2.$

F"ur die in $3.$ genannten Gruppen dagegen ist jeder Invariantenring Cohen-Macaulay.
\end{Bemerkung}

\Bew Um zu zeigen dass die Zusammenhangskomponente $G^{0}$ einer linearen
algebraischen Gruppe $G$ kein Torus ist, gen"ugt es, eine nichttriviale,
abgeschlossene, zusammenh"angende, unipotente Untergruppe $U\subseteq G$
anzugeben. Dann ist n"amlich $U\subseteq G^{0}$, und damit ist $G^{0}$ kein
Torus. Denn Elemente eines Torus sind halbeinfach, und nur das neutrale Element
ist zugleich halbeinfach und unipotent, d.h. man h"atte sonst $U=\{\iota\}$. Hierf"ur gen"ugt es wiederum, einen
nichttrivialen algebraischen Homomorphismus $\rho: \Ga \rightarrow G$ anzugeben. Denn
dann ist $\rho(\Ga)\subseteq G$ eine nichttriviale, zusammenh"angende
abgeschlossene unipotente Untergruppe (\cite[2.2.5, 2.4.8]{SpringerLin}). 

{\it 1.} Es sind $\SL_{n}\subseteq K^{n^{2}}$ bzw. $\GL_{n}\subseteq
K^{n^{2}+1}$ Nullstellenmengen der irreduziblen Polynome $\det-1$ bzw. $E\cdot
\det-1$ (die zu $E$ geh"orige Koordinate liefert hier das Inverse der
Determinante der zugeh"origen Matrix) und damit als Variet"aten irreduzibel, also
zusammenh"angend. F"ur den Zusammenhang von $\Sp_{n}$ siehe \cite[Exercise
2.2.9]{SpringerNew}. Mit $n=2m$,
$J_{2}:=\left(\begin{array}{cc}0&1\\-1&0\end{array}\right)$ und
$J_{2m}:=\diag(J_{2},\ldots,J_{2})\in K^{2m\times 2m}$ (Blockdiagonalmatrix)
ist 
\[
\Sp_{n}=\left\{A\in K^{n\times n}: A^{T}J_{n}A=J_{n}\right\}. 
\]
Man verifiziert sofort, dass $\Sp_{2}=\SL_{2}$, und mit $A\mapsto
\diag(A,I_{2},\ldots,I_{2})\in K^{n\times n}$ hat man eine Einbettung
$\SL_{2}\subseteq \Sp_{n}$ ($I_{2}\in K^{2\times 2}$ Einheitsmatrix). Analog hat man Einbettungen $\SL_{2}\subseteq
\SL_{n},\GL_{n}$. Da 
\begin{equation} \label{HomGa}
\Ga\rightarrow \SL_{2},\quad a\mapsto \left(\begin{array}{cc}
1&a\\
0&1\end{array}\right)
\end{equation}
ein injektiver algebraischer Homomorphismus ist, ist keine der drei
zusammenh"angenden Gruppen
ein Torus.

{\it 2.} F"ur den Zusammenhang von $\SO_{n}$ siehe \cite[Exercise
2.2.2]{SpringerNew}. Da $\On_{n}/\SO_{n}\cong Z_{2}$, ist $\SO_{n}=\On_{n}^{0}$. Offenbar
hat man eine Einbettung $\SO_{3}\subseteq \SO_{n}$. Mit $i$ einer festen L"osung von
$X^{2}+1=0$ und $\sqrt{2}$ einer festen L"osung von $X^{2}-2=0$ in $K$  pr"uft man
leicht (aber etwas m"uhsam) nach, dass durch
\begin{equation} \label{HomSOn}
\rho:\,\Ga\rightarrow \SO_{3},\quad a\mapsto \left(\begin{array}{ccc}
1&i\sqrt{2}a&\sqrt{2}a\\
-i\sqrt{2}a&1+a^{2}&-ia^{2}\\
-\sqrt{2}a&-ia^{2}&1-a^{2}
\end{array}\right)
\end{equation}
ein Homomorphismus gegeben ist, so dass $\SO_{n}$ kein Torus ist.

{\it 3.} Wir m"ussen nur noch zeigen, dass $\SO_{2}$ ein Torus ist. Sei
$A=\left(\begin{array}{cc}a&b\\c&d\end{array}\right)\in \SO_{2}$. Aus
$a^{2}+b^{2}=1=a^{2}+c^{2}$ folgt $c=\pm b$ und analog $d=\pm a$. Eine kurze
Fallunterscheidung zeigt, dass
\[
\SO_{2}=\left\{\left(\begin{array}{cc}a&b\\-b&a\end{array}\right): a,b\in K,\,
a^{2}+b^{2}=1\right\}.
\]
Ist $i=\sqrt{-1}$ ein feste Wurzel von $X^{2}+1=0$ in $K$, so zeigt
\[
\left(\begin{array}{cc}1&1\\i&-i\end{array}\right)^{-1}
\left(\begin{array}{cc}a&b\\-b&a\end{array}\right)
\left(\begin{array}{cc}1&1\\i&-i\end{array}\right)
=
\frac{i}{2}\left(\begin{array}{cc}-i&-1\\-i&1\end{array}\right)
\left(\begin{array}{cc}a+ib&a-ib\\ia-b&-ia-b\end{array}\right)=
\]
\[
=
\frac{i}{2}\left(\begin{array}{cc}-ia+b-ia+b&-ia-b+ia+b\\-ia+b+ia-b&-ia-b-ia-b\end{array}\right)=
\left(\begin{array}{cc}a+ib&0\\0&a-ib\end{array}\right),
\] 
dass die $\SO_{2}$ diagonalisierbar ist. Da die $a+ib$ mit $a^{2}+b^{2}=1$ auch ganz
$K\setminus\{0\}$ aussch"opfen (f"ur $t\in K\setminus\{0\}$ setze
$a:=\frac{1+t^{2}}{2t}$. Mit $b=\sqrt{1-a^{2}}$ und geeignetem Vorzeichen ist
dann $a+ib=t$), ist $\SO_{2}\cong \Gm$ ein Torus. Da $\On_{2}/\SO_{2}=Z_{2}$
 sind also $\SO_{2}$ und $\On_{2}$ f"ur $p\ne 2$ nach Satz
\ref{NagataTorus} linear reduktiv.

{\it 4.} Wir erinnern zun"achst an die Definition der orthogonalen Gruppen in
Charakteristik $2$ (vgl. Carter \cite[1.6]{Carter}).
F"ur $n=2m$ betrachte
\[
f=X_{1}X_{2}+X_{3}X_{4}+\ldots+X_{n-1}X_{n}.
\]
F"ur $n=2m+1$ betrachte
\[
f=X_{1}X_{2}+X_{3}X_{4}+\ldots+X_{n-2}X_{n-1}+X_{n}^{2}.
\]
Dann ist
\[
\On_{n}=\SO_{n}=\left\{\sigma \in \GL_{n}: f\circ \sigma = f \right\}.
\]
Offenbar hat man Einbettungen $\On_{3}\subseteq \On_{2m+1}$ bzw.
$\On_{4}\subseteq \On_{2m}$ f"ur $m\ge 1$ bzw. $m\ge 2$.

F"ur $\On_{3}$ betrachte 
\begin{equation} \label{HomOn3}
\rho: \Ga\rightarrow \GL_{3},\quad a\mapsto \left(\begin{array}{ccc}
1&&\\
a^{2}&1&\\
a&&1
\end{array}\right).
\end{equation}
Offenbar ist $\rho$ ein Homomorphismus, und wegen
\[
f(\rho(a)(x_{1},x_{2},x_{3}))=f(x_{1},a^{2}x_{1}+x_{2},ax_{1}+x_{3})=x_{1}x_{2}+a^{2}x_{1}^{2}+a^{2}x_{1}^{2}+x_{3}^{2}=x_{1}x_{2}+x_{3}^{2}
\]
ist $\rho(a)\in \On_{3}$. Damit enth"alt $\On_{2m+1}^{0}$ eine
zusammenh"angende unipotente Untergruppe, ist also kein Torus.

F"ur $\On_{4}$ betrachte den Homomorphismus
\begin{equation} \label{HomOn4}
\rho: \Ga\rightarrow \GL_{4},\quad a\mapsto \left(\begin{array}{cccc}
1&&&\\
&1&a&\\
&&1\\
a&&&1
\end{array}\right).
\end{equation}
Es folgt
\begin{eqnarray*}
f(\rho(a)(x_{1},x_{2},x_{3},x_{4}))&=&f(x_{1},x_{2}+ax_{3},x_{3},ax_{1}+x_{4})\\&=&x_{1}x_{2}+ax_{1}x_{3}+ax_{1}x_{3}+x_{3}x_{4}=x_{1}x_{2}+x_{3}x_{4},
\end{eqnarray*}
also auch hier $\rho(a)\in \On_{4}$, und $\On_{2m}^{0}$ ist kein Torus f"ur
$m \ge 2$.

$\On_{2}$ schliesslich besteht aus allen Elementen der Form
$\left(\begin{array}{cc}a&0\\0&a^{-1}\end{array}\right)$ und
$\left(\begin{array}{cc}0&b\\-b^{-1}&0\end{array}\right)$, und ist damit
isomorph zu $\Gm \rtimes Z_{2}$. 

Jede der in $1,2,4$ genannten Gruppen $G$ ist also im Fall $p>0$ nicht linear
reduktiv (Satz \ref{NagataTorus}). Au"ser im Fall $\SO_{2}$ f"ur $p=2$ ist
n"amlich $G^{0}$
kein Torus, und  da dann jeweils $K^{n}$ ein treuer $G$-Modul ist, liefert
Satz \ref{ExpliziteDarstellungen} nach Fall $1.$ zu gegebenem $k$ explizit einen
$G$-Modul $V$ mit $\cmdef K[V]^{G}\ge k-2$. Die Gruppe $G=\SO_{2},\, p=2$,
f"ur die $G^{0}$ ein Torus ist, aber $p=2$ Teiler von $|G/G^{0}|$ ist, wird
entsprechend von Fall $2.$ in Satz \ref{ExpliziteDarstellungen} abgedeckt.  

Da die Gruppen in $3.$ linear reduktiv sind, folgt die letzte Aussage aus dem
Satz von Hochster und Roberts \ref{HoRo}.
\qed\\

\subsection{Nichttriviale Kozyklen auf elementarem Weg und Beispiele f"ur
  endliche Gruppen}
Die f"ur unser Konstruktionsverfahren ben"otigten nichttrivialen Kozyklen
h"angen im wesentlichen an dem relativ komplizierten Satz von Nagata, Korollar 
\ref{NagatasNoComplementCorollar}. Hat man aber eine Gruppe konkret
vorgegeben, z.B. $\SL_{n}$, und dazu einen $G$-Modul $V$ sowie einen Kozyklus
$g\in Z^{1}(G,V)$, so kann man relativ leicht feststellen, ob ein Kozyklus
nichttrivial ist. Man macht dazu f"ur einige $\sigma_{i}\in G, \, i=1,\ldots,m$ den
Ansatz
\[
(\sigma_{i} -1)v=g_{\sigma},\quad i=1,\ldots,m,
\]
und erh"alt so ein inhomogenes lineares Gleichungssystem f"ur $v$. Bei
geschickter (meist naheliegender) Wahl der $\sigma_{i}$ erh"alt man dann
entweder einen Widerspruch, was dann zeigt, dass $g$ nichttrivial ist, oder
eine L"osung $v_{0}$ f"ur $v$ (evtl. gibt es mehrere L"osungen). Dann muss man noch verifizieren, ob allgemein
$g_{\sigma}=(\sigma-1)v_{0}\,\myforall \sigma \in G$ gilt, was dann zeigt dass $g$
trivial ist. Ist dies nicht der Fall, so kann man eine andere L"osung
probieren oder mittels Erh"ohung von $m$ die L"osungen weiter einschr"anken und
von vorne beginnen. F"ur die Wahl der $\sigma_{i}$ ist (falls existent) etwa
ein Homomorphismus $\rho: \Ga \rightarrow G$ n"utzlich, und man w"ahlt dann
einige Bilder von $\rho$ als die $\sigma_{i}$. Falls $N$ der gr"osste
 Grad eines der Polynome in der Koordinatendarstellung von $\rho$ ist, wird
 man dann $m=N+2$ w"ahlen (denn ein Polynom vom Grad $N$ ist durch $N+2$ Werte
 "`"uberbestimmt"').

In der Praxis f"uhrt die skizzierte Heuristik in der Regel schnell zum Ziel.\\

Im folgenden Satz demonstrieren wir dieses Verfahren f"ur die Kozyklen f"ur
die Gruppen $\SL_{n},\GL_{n}$ und $\Sp_{n}$, die
wir gem"a"s Proposition \ref{kompl} und Korollar \ref{NagatasNoComplementCorollar}
bereits als nichttrivial erkannt haben. Dazu verwenden wir den Homomorphismus
\eqref{HomGa}. Zus"atzlich erhalten wir dabei das Ergebnis, dass der Kozyklus
auch bei Einschr"ankung auf viele endliche Untergruppen $G\subseteq \GL_{n}$ nichttrivial
bleibt. Dies liefert dann in der Regel kleinere Beispiele von $G$-Moduln $V$
mit $\cmdef K[V]^{G}\ge k-2$ als die Konstruktion gem"a"s Fall $2.$ in Satz
\ref{ExpliziteDarstellungen}, wo "uber die regul"are Darstellung von $G$
gegangen wird.

Um die entsprechenden Resultate f"ur die in $2.$ und $4.$ genannten Gruppen von Bemerkung
\ref{ExpliziteBspFuerKlassGruppen} zu erhalten (insbesondere also f"ur deren
endlichen Untergruppen), kann man die Homomorphismen
\eqref{HomSOn},~\eqref{HomOn3} und \eqref{HomOn4} verwenden.

\begin{Satz}[Kohls \cite{DAPaper}] \label{PKoz}
Sei $\chr K=p>0$, $n \ge 2$ und $G$ eine abgeschlossene Untergruppe von $\GL_{n}$ mit\\
\begin{enumerate}
\renewcommand{\labelenumi}{(\alph{enumi})}
\item  Falls $p=2$:  
$
\left(
\begin{array}{ccc}
 1 & a &\\
&1&\\
&&I_{n-2}
\end{array} \right)
\in G
$
      f"ur wenigstens drei verschiedene Werte von $a\in K$ ($a=0$ ist stets ein solcher).\\
\item Falls $p \ge 3$: 
$
\left(
\begin{array}{ccc}
1  & 1 &\\
0&1&\\
&&I_{n-2}
\end{array} \right)
\in G.
$
\end{enumerate}
Sei weiter $V=S^{p}(\langle X_{1},\ldots,X_{n} \rangle)$ die $p$-te
symmetrische Potenz der nat"urlichen Darstellung von $G$,
$
W:=F^{p}(\langle X_{1},\ldots,X_{n} \rangle) \subseteq V 
$ die $p$-te Frobenius-Potenz der nat"urlichen Darstellung
 und
$
U:=\Hom_{K}(V,W)_{0}
$.

Sei ferner $\iota \in \Hom_{K}(V,W)$ gegeben durch 
$
\iota |_{W}=\id_{W}
$
und 
$
\iota$ gleich $0$  auf allen Monomen, die nicht in $W$  liegen.
Dann ist durch
$
g:G\rightarrow U,\,\, \sigma \mapsto g_{\sigma}:=(\sigma-1) \iota 
$
ein nichttrivialer Kozyklus $g\in Z^{1}(G,U)$ gegeben.
\end{Satz}

\Bew Nach Proposition \ref{kompl} ist jedenfalls $g\in Z^{1}(G,U)$. Wir
zeigen  $g\not\in B^{1}(G,U)$. Wir w"ahlen dazu 
f"ur $V$ eine monomiale Basis $\mathcal{B}$, wobei wir die Reihen\-fol\-ge der ersten $n+1$ bzw. $n+2$ Monome in den F"allen (a) bzw. (b) vorge\-ben, und zwar\\
im Fall (a):
\[
\mathcal{B}=\{X_{1}^{2},\ldots,X_{n}^{2},X_{1}X_{2},\ldots\}
\]\\
im Fall (b):
\[
\mathcal{B}=\{X_{1}^{p},\ldots,X_{n}^{p},X_{1}^{p-1}X_{2},X_{1}^{p-2}X_{2}^{2},\ldots\}
\]
Als Basis $\mathcal{C}$ von $W$ dienen die ersten $n$ Eintr"age von
 $\mathcal{B}$. Sei $N:=|\mathcal{B}|= {n+p-1 \choose p}$.
 Wir verwenden die folgende {\bf Notation}: F"ur
$A \in K^{n \times N}$ sei $A_{\mathcal{C},\mathcal{B}}$ das Element
von $\operatorname{Hom}_{K}(V,W)$, dass $A$ als Darstellungsmatrix
 bez"uglich der Basen $\mathcal{B}$ von $V$ und $\mathcal{C}$ von $W$ hat.
Analog sei f"ur $x\in K^{N}$ durch $x_{\mathcal{B}}$ das Element
von $V$ gegeben, das $x$ als Koordinatenvektor bez"uglich der Basis 
$\mathcal{B}$ von V hat.

 Wir bezeichnen mit $f_{p}:\GL_{n}(K) \rightarrow \GL_{n}(K)$ den koeffizientenweisen Frobenius-Ho\-mo\-mor\-phis\-mus, also $f_{p}(a_{ij})=(a_{ij}^{p})$. Ist $A_{\sigma} \in K^{N \times N}$ die Dar\-stell\-ungs\-matrix von $\sigma \in G$ bzgl. der Basis $\mathcal{B}$, so hat diese die Form
\[
A_{\sigma}=\left( \begin{array}{cc}
f_{p}(\sigma) & *\\
0_{(N-n)\times n}&*\\
\end{array} \right)\in \GL_{N}(K).
\]
(Dabei schreiben wir $0_{k\times l}\in K^{k\times l}$ f"ur die Nullmatrix, analoges gilt sp"ater f"ur
 $0_{k}\in K^{k}$.)
Weiter haben wir bzgl. der Basen $\mathcal{C},\mathcal{B}$
\[
U = \left\{
\left( \begin{array}{cc}
0_{n \times n} & B\\
\end{array} \right)_{\mathcal{C},\mathcal{B}} \in K^{n \times N}:  B \in K^{n \times (N-n)} \right\},
\]
und f"ur $f =   
\left( \begin{array}{cc}
0_{n \times n} & B\\
\end{array} \right)_{\mathcal{C},\mathcal{B}}$  haben wir die Operation gegeben durch
\[
\sigma \cdot f=\sigma \circ f \circ \sigma^{-1} = \left( f_{p}(\sigma) \cdot
\left( \begin{array}{cc}
0_{n \times n} & B\\
\end{array} \right)
\cdot A_{\sigma^{-1}}\right)_{\mathcal{C},\mathcal{B}}.
\]
Die Darstellungsmatrix von $\iota$ bzgl. $\mathcal{C},\mathcal{B}$ ist gegeben durch
\[
J:=\left( \begin{array}{cc}
I_{n} & 0_{n\times (N-n)}\\
\end{array} \right) \in K^{n \times N}, \textrm{ also }
\iota =J_{\mathcal{C},\mathcal{B}}.
\]
Um $g \not\in
B^{1}(G,U)$ zu zeigen, gehen wir vor wie im Vorspann beschrieben, d.h. wir nehmen die Existenz eines
\[
 u = Z_{\mathcal{C},\mathcal{B}} \in U\quad \textrm{ mit }
Z=\left( \begin{array}{cc}
0_{n \times n} & \hat{Z}\\
\end{array} \right) \in K^{n \times N}, \,  \hat{Z}=(z_{ij}) \in K^{n \times (N-n)}
\]
an mit $g_{\sigma}=(\sigma -1) \iota =(\sigma -1)u \quad
\textrm{f"ur alle } \sigma \in G$, d.h.
\begin{equation} \label{zerf}
f_{p}(\sigma)J A_{\sigma^{-1}}-J = f_{p}(\sigma)Z
A_{\sigma^{-1}}-Z \quad \textrm{f"ur alle } \sigma \in G.
\end{equation}
Wir werden diese Gleichung nun in beiden F"allen zum Widerspruch f"uhren.\\

\noindent (a) Mit
\[
\sigma:=
\left(
\begin{array}{ccc}
 1 & a &\\
0&1&\\
&&I_{n-2}
\end{array} \right) =\sigma^{-1},\quad  a\in K \textrm{ so dass } \sigma \in G, 
\]
berechnen wir die $(n+1)$te Spalte von $A_{\sigma^{-1}}$:
\[
\begin{array}{rcl}
\sigma^{-1} \cdot X_{1}X_{2}&=&X_{1}(aX_{1}+X_{2})\\
&=&aX_{1}^{2}+X_{1}X_{2}\\
&=&(a,0_{n-1},1,0_{N-n-1})^{T}_{\mathcal{B}}.\\
\end{array}
\]
Nun vergleichen wir auf beiden Seiten von (\ref{zerf}) die Eintr"age in der
ersten Zeile  und
der $(n+1)$ten Spalte:\\
Linke Seite:
\[
(1,a^{2},0_{n-2})
\left( \begin{array}{cc}
I_{n} & 0_{n \times (N-n)}\\
\end{array} \right)
\left( \begin{array}{c}
a\\
0_{n-1}\\
1\\
0_{N-n-1}\\
\end{array} \right) =a
\]
Rechte Seite:
\[
\begin{array}{rcl}
(1,a^{2},0_{n-2})
\left( \begin{array}{cc}
0_{n \times n} & \hat{Z}\\
\end{array} \right)
\left( \begin{array}{c}
a\\
0_{n-1}\\
1\\
0_{N-n-1}\\
\end{array} \right) -z_{11}&=&z_{11}+a^{2}z_{21}-z_{11}\\
&=&a^{2}z_{21}\\
\end{array}
\]
Mit $c:=z_{21}$ und gleichsetzen beider Seiten erhalten wir, dass
\[
ca^{2}+a=0
\]
f"ur wenigstens drei Werte $a \in K$ gelten muss. Dies ist nicht m"oglich.\\

\noindent (b) Wir betrachten
\[
\sigma=
\left(
\begin{array}{ccc}
1  & 1 &\\
0&1&\\
&&I_{n-2}
\end{array} \right),
\sigma^{-1}=
\left(
\begin{array}{ccc}
1  & -1 &\\
0&1&\\
&&I_{n-2}
\end{array} \right)
\in G
\]
und berechnen die  $(n+1)$te und $(n+2)$te Spalte von $A_{\sigma^{-1}}$:\\
$(n+1)$te Spalte:
\[
\begin{array}{rcl}
\sigma^{-1} \cdot X_{1}^{p-1}X_{2} &=& X_{1}^{p-1}(-X_{1}+X_{2})\\
&=&-X_{1}^{p}+X_{1}^{p-1}X_{2}\\
&=&(-1,0_{n-1},1,0_{N-n-1})^{T}_{\mathcal{B}}
\end{array}
\]
$(n+2)$te Spalte:
\[
\begin{array}{rcl}
\sigma^{-1} \cdot X_{1}^{p-2}X_{2}^{2} &=& X_{1}^{p-2}(X_{1}^{2}-2X_{1}X_{2}+X_{2}^{2})\\
&=&X_{1}^{p}-2X_{1}^{p-1}X_{2}+X_{1}^{p-2}X_{2}^{2}\\
&=&(1,0_{n-1},-2,1,0_{N-n-2})^{T}_{\mathcal{B}}
\end{array}
\]
Wir vergleichen wieder beide Seiten von (\ref{zerf}):\\
(i) \underline{Erste Zeile, $(n+1)$te Spalte}\\
Linke Seite:
\[
\left( \begin{array}{ccc} 1 & 1 & 0_{n-2}\\ \end{array} \right)
\left( \begin{array}{cc} I_{n} & 0_{n \times (N-n)}\\ \end{array} \right)
\left( \begin{array}{c}
-1\\
0_{n-1}\\
1\\
0_{N-n-1}\\
\end{array} \right)=-1
\]
Rechte Seite:
\[
\left( \begin{array}{ccc} 1 & 1 & 0_{n-2}\\ \end{array} \right)
\left( \begin{array}{cc} 0_{n \times n} & \hat{Z}\\ \end{array} \right)
\left( \begin{array}{c}
-1\\
0_{n-1}\\
1\\
0_{N-n-1}\\
\end{array} \right)-z_{11}=z_{11}+z_{21}-z_{11}=z_{21}.
\]
Gleichsetzen beider Seiten liefert
\begin{equation} \label{wid1}
z_{21}=-1.
\end{equation}

\noindent (ii) \underline{Zweite Zeile, $(n+2)$te Spalte}\\
Linke Seite:
\[
\left( \begin{array}{ccc} 0 & 1 & 0_{n-2}\\ \end{array} \right)
\left( \begin{array}{cc} I_{n} & 0_{n \times (N-n)}\\ \end{array} \right)
\left( \begin{array}{c}
1\\
0_{n-1}\\
-2\\
1\\
0_{N-n-2}\\
\end{array} \right)=0
\]
Rechte Seite:
\[
\left( \begin{array}{ccc} 0 & 1 & 0_{n-2}\\ \end{array} \right)
\left( \begin{array}{cc} 0_{n \times n} & \hat{Z}\\ \end{array} \right)
\left( \begin{array}{c}
1\\
0_{n-1}\\
-2\\
1\\
0_{N-n-2}\\
\end{array} \right)-z_{22}=-2z_{21}+z_{22}-z_{22}=-2z_{21}
\]
Da $p \ge 3$ ist $2 \ne 0$, und gleichsetzen beider Seiten liefert 
\[
z_{21}=0,
\]
im Widerspruch zu (\ref{wid1}). \qed\\

Eine Matrix (bzw. die zugeh"orige lineare Abbildung) hei"st
\emph{Transvektion},\index{Transvektion} wenn sie die Matrix aus Fall (b) als Jordan-Normalform
hat. Offenbar ist $A\in K^{n\times n}$ genau dann eine Transvektion, wenn
$\rang(A-I_{n})=1$ und $(A-I_{n})^{2}=0$ ist ($I_{n}$ die Einheitsmatrix). Dies zeigt, dass $A$ genau dann
eine Transvektion ist, wenn es Spaltenvektoren $0\ne u,v\in K^{n}$ gibt mit
$A=I_{n}+uv^{T}$ und $v^{T}u=0$. 

Da die Voraussetzungen in Satz \ref{PKoz} nat"urlich nur
bis auf Konjugation in $\GL_{n}$ zu verstehen sind, liefert dieser also (f"ur
$p\ge 3$) f"ur jede abgeschlossene Untergruppe $G\subseteq \GL_{n}$, die eine
Transvektion enth"alt, einen
nichttrivialen Kozyklus. Die Nichttrivialit"at des Kozyklus ist nach
Proposition \ref{kompl} "aquivalent dazu, dass $W=F^{p}(X)$ (mit $X=\langle X_{1},\ldots,X_{n}\rangle$) kein Komplement in $V=S^{p}(X)$ hat.
Da Konjugation von $G$ in $\GL_{n}$ lediglich einem Basiswechsel von $X$
entspricht, $F^{p}(X)$ aber Basisunabh"angig ist (vgl. Definition \ref{FpV}), sehen wir, dass $F^{p}(X)$
kein Komplement in $S^{p}(X)$ hat, wenn $G$ eine Transvektion enth"alt ($p\ge 3$). Dies
erweitert Korollar \ref{NagatasNoComplementCorollar}.\\

\begin{Def}
Ist $G\subseteq \GL_{n}$ eine abgeschlossene Untergruppe, $\chr K=p>0$ und
$q=p^{m}>1$, so schreiben wir
\[
G({\mathbb F}_{q}):=G \cap \GL_{n}({\mathbb F}_{q}):=\{(a_{ij})\in G:
a_{ij}^{q}-a_{ij}=0 \,\,\,\myforall i,j=1,\ldots,n \},
\]
und dies ist dann ebenfalls eine abgeschlossene Untergruppe.
Insbesondere sind so $\SL_{n}({\mathbb F}_{q})$, $\GL_{n}({\mathbb
  F}_{q})$, $\Sp_{n}({\mathbb F}_{q})$, $\SO_{n}({\mathbb F}_{q})$, $\On_{n}({\mathbb F}_{q})$ "uber dem K"orper $K$ definierte lineare
  algebraische Gruppen.
\end{Def}

Die Voraussetzungen von Satz \ref{PKoz} sind offenbar f"ur die endlichen
  Gruppen $\SL_{n}({\mathbb F}_{q})$, $\GL_{n}({\mathbb
  F}_{q})$, $\Sp_{n}({\mathbb F}_{q})$ erf"ullt, falls $q\ge 3$ ist, und
  dieser liefert
  so einen nichttrivialen Kozyklus. Entsprechendes gilt auch f"ur die Gruppe
  $({\mathbb F}_{q},+)$ via der Einbettung \eqref{HomGa}. Korollar
  \ref{depthToInfty} ist also anwendbar und liefert

\begin{Korollar}
Sei $G\subseteq \GL_{n}$ eine reduktive Untergruppe, die die Voraussetzungen
von Satz \ref{PKoz} (bis auf Konjugation in $\GL_{n}$) erf"ullt, z.B. 
$G=\SL_{n}({\mathbb F}_{q})$, $\GL_{n}({\mathbb
  F}_{q})$, $\Sp_{n}({\mathbb F}_{q})$,$({\mathbb F}_{q},+)$ f"ur $n\ge 2,
q=p^{m}\ge 3$, oder $G$ enthalte eine Transvektion und $p\ge 3$, oder auch $G$
eine der Gruppen $\SL_{n}$, $\GL_{n}$ oder $\Sp_{n}$, oder $G$ eine reduktive
Gruppe, die eine dieser Gruppen enth"alt. Ist dann $U$ der Modul aus
Satz \ref{PKoz} und $\tilde{U}$ der zu dem nichttrivialen Kozyklus $g$
geh"orige erweiterte $G$-Modul, so ist
\[
\cmdef K\left[U^{*}\oplus\bigoplus_{i=1}^{k}\tilde{U}\right]^{G} \ge k-2.
\]
\end{Korollar}

Da hier jeweils $U\subseteq S^{2}(W\oplus (V/W)^{*})$ ist, ist also $W\oplus
(V/W)^{*}$ jeweils ein treuer $G$-Modul mit $H^{1}(G,K[W^{*}\oplus (V/W)])\ne
0$. F"ur endliches $G$ ist also auch Satz \ref{FiniteCMdef} mit $W^{*}\oplus
(V/W)$ (statt $V$)
und $r_{0}=1$ anwendbar.\\

Wie bereits bemerkt, kann man das Korollar unter geeigneten Voraussetzungen noch auf die Gruppen $\SO_{n}({\mathbb F}_{q})$ und
$\On_{n}({\mathbb F}_{q})$ (sowie deren reduktiven (z.B. endlichen) Obergruppen) ausdehnen, indem man einen zu Satz \ref{PKoz}
analogen Satz mit Hilfe des Homomorphismus \eqref{HomSOn} beweist. Da der
Beweis letztlich genauso l"auft wie im Fall (a) von Satz \ref{PKoz}, wollen wir
hier darauf verzichten. 

Jedenfalls l"asst sich Satz \ref{PKoz} so wie er dasteht
nicht auf orthogonale Gruppen anwenden. Man kann n"amlich zeigen, dass
orthogonale Gruppen in ungerader Charakteristik "uberhaupt keine Transvektionen enthalten,
 und orthogonale Gruppen in gerader Charakteristik keine zwei Transvektionen
 mit demselben Fixraum enthalten. Daf"ur geben wir im n"achsten Abschnitt eine
 alternative, elementare Konstruktion von nichttrivialen Kozyklen f"ur gewisse
 orthogonale Gruppen an, die einfacher, aber weniger allgemein ist.

\subsection{Weitere Beispiele f"ur einige orthogonale Gruppen}
Der wesentliche Schritt einen nichttrivialen Kozyklus zu konstruieren ist nach
Proposition \ref{kompl}, einen Untermodul ohne Komplement zu finden. Ist $G=\SO_{n}$
oder $\On_{n}$ und $V=K^{n}=\langle X_{1},\ldots,X_{n}\rangle$ die nat"urliche
Darstellung, so bietet
sich (abgesehen von $F^{p}(V)\subseteq S^{p}(V))$ der von der kanonischen
Invariante $X_{1}^{2}+\ldots+X_{n}^{2}$ erzeugte Untermodul in $S^{2}(V)$ an.

\begin{Satz}
Sei $\chr K=p\ne 2$, $G=\SO_{n}$
oder $\On_{n}$ und $\langle X_{1},\ldots,X_{n}\rangle$ die nat"urliche
Darstellung. Genau dann hat $K\cdot(X_{1}^{2}+\ldots+X_{n}^{2})$ ein
Komplement in $S^{2}(\langle X_{1},\ldots,X_{n}\rangle)$, wenn $p\not| n$ gilt.
\end{Satz}

\Bew Wir setzen $V:=S^{2}(\langle X_{1},\ldots,X_{n}\rangle)$ und
$e:=X_{1}^{2}+\ldots+X_{n}^{2}\in V^{G}$.\\
"`$\Leftarrow$"' Sei $p\not| n$. Wir behaupten, dass dann
\[
U:=\langle X_{i}X_{j}:\, i\ne j \rangle_{K} \oplus
\left\{a_{1}X_{1}^{2}+\ldots+a_{n}X_{n}^{2}: \,\, a_{i}\in K, \,\,a_{1}+\ldots+a_{n}=0 \right\} 
\]
ein $G$-invariantes Komplement zu $Ke$ ist. Offenbar ist $\dim U=\dim
V-1$. Ferner ist $e\not\in U$, sonst w"are $1+\ldots+1=n=0$, aber $p\not|
n$. Also ist $V=Ke\oplus U$. Es
ist also nur noch $GU\subseteq U$ zu zeigen. Sei $\sigma=(s_{ij})\in G$. Dann
ist
\[
\sigma \cdot X_{i}X_{j} = (\textrm{Terme in } X_{k}X_{l} \textrm{ mit } k\ne
l)+\sum_{k=1}^{n}s_{ki}s_{kj}X_{k}^{2}.
\]
F"ur $i\ne j$ ist
$\sum_{k=1}^{n}s_{ki}s_{kj}=(\sigma^{T}\sigma)_{ij}=\delta_{ij}=0$, und damit
$\sigma \cdot X_{i}X_{j}\in U$. Ist weiter $a_{1}+\ldots+a_{n}=0$ mit $a_{i}\in K$,
so gilt
\[
\sigma \cdot (a_{1}X_{1}^{2}+\ldots+a_{n}X_{n}^{2})=(\textrm{Terme in } X_{k}X_{l} \textrm{ mit } k\ne
l)+\sum_{i=1}^{n}\sum_{j=1}^{n}a_{j}s_{ij}^{2}X_{i}^{2}\in U,
\]
denn
$\sum_{j=1}^{n}a_{j}\underbrace{\sum_{i=1}^{n}s_{ij}^{2}}_{1}=\sum_{j=1}^{n}a_{j}=0$.
Also ist $\sigma U\subseteq U$.\\
"`$\Rightarrow$"' Sei $p|n$. Angenommen, es g"abe einen Untermodul $U$ mit
$V=Ke\oplus U$. Wir f"uhren die Annahme durch Induktion zum Widerspruch.

{\it 1. $p=n$:}  Sei $\sigma$ die lineare Fortsetzung der
Permutation\[
X_{1}\mapsto X_{2}\mapsto X_{3}\mapsto \ldots\mapsto X_{p}\mapsto
X_{1}.\]
Dann ist $\sigma \in \On_{p}$, und wegen $\det \sigma =\sgn (12\dots
p)=(-1)^{p+1}=1$ ist sogar $\sigma\in \SO_{p}$ (weil die entsprechende
Eigenschaft f"ur gerades n nicht gilt, ben"otigen wir den Schritt {\it 2.}). Wir betrachten die von
$\sigma$ erzeugte zyklische Untergruppe der Ordnung $p$, $Z_{p}:=\langle
\sigma \rangle \subseteq G$. Dann ist $V=Ke\oplus U$ erst recht eine Zerlegung
von 
$Z_{p}$-Moduln $(*)$. Auch $W:=\langle X_{1}^{2},\ldots,X_{p}^{2}\rangle$ ist ein
$Z_{p}$-Modul (aber kein $G$-Modul), der isomorph ist zur regul"aren
Darstellung von $Z_{p}$ ($X_{i}^{2}$ entspricht $e_{\sigma^{i}}$). Nach Korollar \ref{MaschkeNoComplement} hat dann $Ke$ kein
$Z_{p}$-Komplement in $W$, also nach Bemerkung \ref{kleinKompl} erst recht
nicht in $V$, im Widerspruch zu $(*)$.

{\it 2.  $p<n$:} Sei $n=k+m$ mit $p|k,m$ und $k,m>0$. Sei
$e_{1}:=X_{1}^{2}+\ldots+X_{k}^{2}$  und
$e_{2}=X_{k+1}^{2}+\ldots+X_{n}^{2}$. Falls $e_{1}\notin U$, so ist wegen $\dim
U=\dim V-1$ dann $V=Ke_{1}\oplus U$. Mittels $A\mapsto \diag(A,I_{m})$ 
($I_{m}$ die Einheitsmatrix) kann
man die Gruppe $\SO_{k}$ in $G$ einbetten, und dann ist $V=Ke_{1}\oplus U$
eine Zerlegung in $\SO_{k}$-Moduln $(*)$. Aber auch $W:=S^{2}\left(\langle
X_{1},\ldots,X_{k}\rangle\right)\subseteq V$ ist ein $\SO_{k}$-Untermodul, und wegen
$(*)$ und Bemerkung \ref{kleinKompl} h"atte $Ke_{1}$ dann auch ein
$\SO_{k}$-Komplement in $W$, im Widerspruch zur Induktionsvoraussetzung. Also
ist $e_{1}\in U$. Analog folgt $e_{2}\in U$ und damit $e=e_{1}+e_{2}\in U$, im
Widerspruch zu $Ke\cap U=\{0\}$.\qed\\

Fast genauso beweist man

\begin{Satz}
Sei $\chr K=2$, $2|n$, $\langle X_{1},\ldots,X_{2n}\rangle$ die nat"urliche
Darstellung der $\SO_{2n}$. Dann hat $K\cdot(X_{1}X_{2}+\ldots+X_{2n-1}X_{2n})$ kein
Komplement in $S^{2}(\langle X_{1},\ldots,X_{2n}\rangle)$.
\end{Satz}

\Bew Sei  $V:=S^{2}(\langle
X_{1},\ldots,X_{2n}\rangle)$, $e:=X_{1}X_{2}+\ldots+X_{2n-1}X_{2n}\in V$ die
$\SO_{2n}$ definierende quadratische Form, und $U:=\langle
X_{1}X_{2},X_{3}X_{4},\ldots,X_{2n-1}X_{2n}\rangle\subseteq V$.
Wir betrachten die lineare Fortsetzung $\sigma$ von $X_{i}\mapsto
X_{i+2}$ (zyklisch). Dann gilt $\sigma\cdot X_{2k-1}X_{2k}=X_{2k+1}X_{2k+2}$,
insbesondere $\sigma \cdot e=e$, also $\sigma \in \SO_{2n}$ und
 $U$ ist die regul"are Darstellung von $Z_{n}:=\langle \sigma \rangle$. Nach
 Korollar \ref{MaschkeNoComplement} hat dann also $Ke$ kein $Z_{n}$-Komplement in $U$,
 nach Bemerkung \ref{kleinKompl} also auch nicht in $V$. Dann gibt es aber
 erst recht kein $\SO_{2n}$-Komplement in $V$. \qed\\

Diese beiden S"atze liefern (mit Proposition \ref{kompl} bzw. der "`Quintessenz"',
S. \pageref{Quintessenz}) also auf elementarem Weg nichttriviale Kozyklen f"ur
einige orthogonale Gruppen. Au"serdem sind die zugeh"origen Moduln von
niedrigerer Dimension als diejenigen aus dem Satz von Nagata. F"ur die Gruppe
$\SO_{n}$ mit $p\ne 2, \,p|n$ und dem aus dem vorletzten Satz konstruierten Modul mit Kozyklus
$U$ gilt dann etwa f"ur die zugeh"orige Erweiterung $\tilde{U}$, dass
$\tilde{U}^{*}=S^{2}(\langle X_{1},\ldots,X_{n}\rangle)$, also $\dim
\tilde{U}={n+1 \choose 2}$. Bei der Konstruktion nach Nagata dagegen h"atte man
$\dim \tilde{U}=n\left({n+p-1 \choose p}-n \right)+1$, siehe
S. \pageref{DimOfM1}, was im Allgemeinen deutlich gr"o"ser ist. Insbesondere
haben dann nat"urlich auch die nach Korollar \ref{depthToInfty} konstruierten
$\SO_{n}$-Moduln $V$ mit $\cmdef K[V]^{\SO_{n}}\ge k-2$ eine entsprechend
geringere Dimension.

\subsection{Die Beispiele f"ur die $\SL_{2}$ in Charakteristik $2$ und $3$}
Wir wollen hier die nach Satz \ref{ExpliziteDarstellungen}
konstruierten $\SL_{2}$-Moduln $V$ mit $\cmdef K[V]^{\SL_{2}}\ge k-2$  mit
$p=\chr K\in\{2,3\}$ genauer untersuchen. Zum einen werden wir $V$ zun"achst
wie vom Satz geliefert, aber etwas expliziter angeben. Schliesslich "andern
wir das Konstruktionsverfahren f"ur $V$ etwas ab, um eine geringere Dimension
zu erhalten. Anstatt n"amlich wie in Korollar \ref{depthToInfty} einfach
$V^{*}=U\oplus \bigoplus_{i=1}^{k}\tilde{U}^{*}$ zu setzen, wobei dann jeder
Summand $\tilde{U}^{*}$ einen Annullator $a_{i}$ enth"alt, ersetzen wir hier
grob gesprochen $\tilde{U}^{*}$ durch Moduln kleinerer Dimension $W$, so dass
$\tilde{U}^{*}\subseteq S^{2}(W)$ gilt. Die Annullatoren $a_{i}$ des
nichttrivialen Kozyklus (mit Werten in $U$) liegen also
hier in der zweiten symmetrischen Potenz. Manchmal enth"alt $S^{2}(W)$ sogar
gleich mehrere Kopien von $\tilde{U}^{*}$ und damit auch mehrere
Annullatoren. Diese kann man aber nicht immer alle verwenden, denn um den
Hauptsatz \ref{BigMainTheorem} anzuwenden, m"ussen sie zus"atzlich noch ein phsop im Polynomring
bilden, d.h. man muss eine geschickte Auswahl treffen. Wie wir au"serdem
bereits gesehen haben (siehe Bemerkung \ref{smallCocycles}),
kann man auch den Kozyklus in die zweite symmetrische
Potenz bringen, um die Dimension von $V$ zu verringern.\\

F"ur die hier dargestellten Ergebnisse ben"otigen wir Resultate, die sich in meiner
Diplomarbeit \cite{Diplomarbeit} und in dem daraus
entstandenen Artikel \cite{DAPaper} finden. In diesen Arbeiten ging es darum,
nicht Cohen-Macaulay Invariantenringe zu konstruieren; Mit dem Hauptsatz
k"onnen wir nun auch noch den Cohen-Macaulay-Defekt gegen unendlich
treiben. Die in diesem Abschnitt dargestellten Ergebnisse versch"arfen also die Resultate aus
\cite{Diplomarbeit} und \cite{DAPaper}.
S"amtliche hier weggelassenen
Zwischenrechungen (die zwar alle einfach, aber teilweise doch umfangreich
sind), finden sich in \cite[Abschnitt 6]{Diplomarbeit}.\\

Wir beschr"anken uns hier der Einfachheit halber auf die Gruppe $\SL_{2}$. S"amtliche Ergebnisse gehen
nat"urlich auch f"ur reduktive Untergruppen der $\SL_{2}$ durch, die die
entsprechende Voraussetzung aus Satz \ref{PKoz} erf"ullen, z.B. endliche
$\SL_{2}({\mathbb F}_{q})$ wobei $q$ entsprechend eine $2$- oder $3$-Potenz
ist. In den Arbeiten  \cite{Diplomarbeit} bzw. \cite{DAPaper} habe ich
au"serdem gezeigt, wie man durch "`Tensorieren mit dem Inversen der
Determinante"' die betrachteten $\SL_{2}$-Moduln zu $\GL_{2}$-Moduln machen
kann, so dass auch hier die Ergebnisse mit etwas Modifikation durchgehen.\\

Wir erinnern nochmal an die auf S. \pageref{mynotation}  eingef"uhrte {\bf Notation.} Mit $\langle X,Y \rangle$
bezeichnen wir die nat"urliche Darstellung der $\SL_{2}$. Mit
$\sigma:=\left(\begin{array}{cc}a&b\\c&d\end{array}\right)\in \SL_{2}$ gilt dann also
\[
\sigma \cdot X=aX+cY, \quad \sigma \cdot Y=bX+dY.
\]
Wenn wir diese Operation entsprechend auf homogene Polynome in $X$ und $Y$
fortsetzen, erhalten wir die Operation auf den symmetrischen Potenzen
\[
S^{k}(\langle X,Y\rangle)=:\langle
X^{k},X^{k-1}Y,\ldots,XY^{k-1},Y^{k}\rangle.
\]
Um Matrizen bez"uglich gegebener Basen lineare Abbildungen zuzuordnen,
verwenden wir die Notation, die wir im Beweis von Satz \ref{PKoz} eingef"uhrt
haben. 

Die nat"urliche Darstellung der $\SL_{2}$ ist selbstdual (f"ur beliebige
Charakteristik), genauer haben $\langle X,Y \rangle^{*}=\langle X^{*},Y^{*}
\rangle$ (mit Dualbasis) und  $\langle Y,-X \rangle$ die gleiche Darstellung.

\subsubsection{Charakteristik 2}
Wir gehen die Konstruktion nach Satz \ref{ExpliziteDarstellungen} Schritt f"ur
Schritt durch. Da $\SL_{2}^{0}$ kein Torus ist, ben"otigen wir zun"achst einen
treuen $\SL_{2}$-Modul und w"ahlen $V=\langle X,Y \rangle$. Dann
ist
\begin{eqnarray*}
S^{2}(V)&=&\langle X^{2},Y^{2},XY \rangle,\\
F^{2}(V)&=&\langle X^{2},Y^{2}\rangle.
\end{eqnarray*}
Wir schreiben $\mathcal{B}$ und $\mathcal{C}$ f"ur die beiden gegebenen
Basen. Die Darstellungen der beiden Moduln bez"uglich dieser Basen sind dann gegeben durch
\[
\sigma \mapsto \left( \begin{array}{ccc}
a^{2} & b^{2} & ab\\
c^{2} & d^{2} & cd\\
0 & 0 & 1
\end{array} \right) \textrm{ und }\sigma \mapsto\left( \begin{array}{cc}
a^{2} & b^{2}\\
c^{2} &  d^{2}
\end{array} \right).
\]
Gem"a"s Satz \ref{ExpliziteDarstellungen} betrachten wir den Modul
$U:=\Hom_{K}(S^{2}(V),F^{2}(V))_{0}$. Dessen Elemente (gewisse lineare Abbildungen)
werden bez"uglich der Basen $\mathcal{B}$ und $\mathcal{C}$ durch $2\times 3$
Matrizen beschrieben, die nur in der letzten Spalte von Null verschiedene
Eintr"age haben. Mit dieser Beschreibung erh"alt man dann (bei naheliegender
Basiswahl), dass eine Darstellung von $U$ durch
$
\sigma \mapsto
\left( \begin{array}{cc}
a^{2} & b^{2}\\
c^{2} &  d^{2}
\end{array} \right)
$ gegeben ist. Damit ist 
$
U\cong \langle X^{2},Y^{2} \rangle
$, und wir identifizieren im Folgenden $\langle X^{2},Y^{2} \rangle$ mit
$U$. Als n"achstes brauchen wir ein wie in Satz \ref{ExpliziteDarstellungen}
gefordertes $\iota$, und wir w"ahlen nat"urlich $\iota =
\left( \begin{array}{ccc} 1&0&0\\ 0&1 &0 \end{array}
\right)_{\mathcal{C},\mathcal{B}}$. F"ur die Operation von $\SL_{2}$,
angewendet auf $\iota$ erhalten wir dann mit obiger Identifikation
 $\sigma \iota=abX^{2}+cdY^{2}+\iota$. 
Damit ist die Operation auf $\tilde{U}=U\oplus K\iota = \langle
X^{2},Y^{2},\iota\rangle$ gegeben durch $\sigma \mapsto \left( \begin{array}{ccc}
a^{2} & b^{2} & ab\\
c^{2} & d^{2} & cd\\
0 & 0 & 1
\end{array} \right)$, also 
$
\tilde{U} = \langle X^{2},Y^{2},XY \rangle
$,
wobei wir $XY$ mit $\iota$ identifiziert haben. Der nichttriviale Kozyklus in
$Z^{1}(\SL_{2},U)$ ist dann gegeben durch
$
g_{\sigma}=(\sigma-1)\iota=(\sigma-1)XY
$. Da $U=\langle X^{2},Y^{2}\rangle$ mit $\langle X,Y\rangle$ selbstdual ist, liefert Satz \ref{ExpliziteDarstellungen} sofort:
\begin{Bsp} \label{simpleBsp2}
Sei $\chr K=2$ und $\langle X,Y\rangle$ die nat"urliche Darstellung der
$\SL_{2}$. Dann gilt
\[
\cmdef K\left[\langle X^{2},Y^{2}\rangle\oplus\bigoplus_{i=1}^{k}\langle
  X^{2},Y^{2},XY \rangle \right]^{\SL_{2}}\ge k-2.
\]
\end{Bsp}
Die Dimension des Invariantenringes ist nach Bemerkung
\ref{exakteDimOfSL2FrobInvs} gegeben durch $3k-1$. (Der zugrundeliegende Modul ist weder selbstdual
noch vollst"andig reduzibel, da er $\tilde{U}$ als Summanden enth"alt).

Wir versuchen nun, $\tilde{U}^{*}$ in einer zweiten symmetrischen Potenz
wiederzufinden. Wir verwenden folgende Notation f"ur die zugeh"orige Dualbasis:
\[
\mu:=(X^{2})^{*}, \nu:=(Y^{2})^{*}, \pi:=(XY)^{*}, \textrm{ also }\tilde{U}^{*} =\langle \mu,\nu,\pi \rangle.
\]
Dann ist $\pi$ die den Kozyklus $g$ annullierende Invariante gem"a"s
Proposition \ref{AnnulatorProp}. Die Darstellung von $\tilde{U}^{*}$
bez"uglich dieser Basis ist gegeben durch
\begin{equation} \label{munupi}
\sigma \mapsto
\left( \begin{array}{ccc}
a^{2} & b^{2} & ab\\
c^{2} & d^{2} & cd\\
0 & 0 & 1
\end{array} \right)^{-T}=
\left( \begin{array}{ccc}
d^{2} & c^{2} & 0\\
b^{2} & a^{2} & 0\\
bd & ac & 1
\end{array} \right).
\end{equation}

Wir betrachten nun das Tensorprodukt der nat"urlichen Darstellung mit sich selbst, $\langle X,Y \rangle \otimes
\langle X,Y \rangle$, und berechnen seine Darstellung bez"uglich der Basis  $\{Y
\otimes Y,  X \otimes X,  X\otimes Y-Y\otimes X,Y\otimes X \}$ zu
\[
\sigma \mapsto
\left( \begin{array}{cccc}
d^{2} & c^{2} & 0 & cd\\
b^{2} & a^{2} & 0 & ab\\
bd & ac & 1 & bc\\
0 & 0 & 0 & 1\\
\end{array} \right).
\]
Da die linke obere $3 \times 3$ Block-Matrix aber die Darstellung von $\langle
\mu, \nu, \pi \rangle$ ist, siehe~\eqref{munupi}, haben wir also (bis auf Isomorphie)
\[
\langle \mu, \nu, \pi \rangle \subseteq  \langle X,Y \rangle \otimes
\langle X,Y \rangle \subseteq  S^{2} \left( \langle X_{1},Y_{1} \rangle \oplus
\langle X_{2},Y_{2} \rangle \right). 
\]
Dabei entspricht dann $X_{1}Y_{2}-X_{2}Y_{1}$ der annullierenden Invariante
$\pi$ des Kozyklus $g$. Hat man nun die $k$-fache direkte Summe der nat"urlichen
Darstellung, so liegen in ihrer zweiten Potenz offenbar ${k\choose 2}$ solche
Annullatoren. Wir werden in Lemma \ref{phsopXY} zeigen, dass $k-1$ von ihnen ein phsop im
Polynomring bilden. Damit haben wir nach dem Hauptsatz \ref{BigMainTheorem} (mit $k-1$ statt $k$)
\begin{Bsp} \label{HauptBspChar2}
Ist $\chr K=2$ und $\langle X,Y \rangle$ die nat"urliche Darstellung der
$\SL_{2}$, so gilt
\[
\cmdef K\left[\langle X^{2},Y^{2} \rangle \oplus \bigoplus_{i=1}^{k}\langle
  X,Y \rangle \right]^{\SL_{2}}\ge k-3.
\]
\end{Bsp}
\noindent Die Dimension des Invariantenringes ist $2k-1$ nach Korollar \ref{exakteDimOfSL2Invs}.

Bemerkenswert ist, dass der hier auftretende, treue $\SL_{2}$-Modul selbstdual
ist (als direkte Summe selbstdualer Moduln), und au"serdem noch vollst"andig
reduzibel (als Summe irreduzibler Moduln). K.N. Raghavan stellte in einer
E-Mail an Gregor Kemper die Frage, ob ein vollst"andig reduzibler $G$-Modul
$V$ mit nicht Cohen-Macaulay Invariantenring $K[V]^{G}$ existiert. Mit diesem
Beispiel lautet die Antwort also "`Ja"'. 

Mit ziemlichem Aufwand werden wir dieses Beispiel im n"achsten Abschnitt (mit fast v"ollig anderer
Methode) auf Charakteristik $p$ verallgemeinern.

Es fehlt noch das versprochene phsop im Polynomring:

\begin{Lemma} \label{phsopXY}
Es sei $\chr K$ beliebig und $R=K[X_{1},Y_{1},\ldots,X_{n},Y_{n}]$
der Polynomring in $2n$ Variablen. Mit\[
g_{ij}:=X_{i}Y_{j}-X_{j}Y_{i}
\]
ist $G:=\{g_{12},g_{23},g_{34},\ldots,g_{n-1,n}\}$ ein phsop in $R$ der L"ange $n-1$ .
\end{Lemma}

\Bew Wir zeigen zun"achst, dass $G$ eine Gr"obner-Basis des Ideals $I:=(G)$ ist,
und zwar bez"uglich graduierter lexikographischer Ordnung mit
\[
X_{1}>Y_{1}>X_{2}>Y_{2}>\ldots>X_{n}>Y_{n}.
\]
F"ur $f\in R$ bezeichnen wir mit $\LM(f)$ das \emph{Leitmonom}\index{Leitmonom} (normiert) von
$f$ bez"uglich dieser Ordnung. Offenbar ist
\[
\LM(g_{i,i+1})=X_{i}Y_{i+1}.
\]
Gem"a"s dem Buchberger-Kriterium (Eisenbud \cite[Theorem 15.8]{Eisenbud})
m"ussen wir zeigen, dass der Divisionsalgorithmus \cite[Division Algorithm
15.7]{Eisenbud}  bzgl. $G$, angewendet auf die \emph{s-Polynome}
\[
\spol(g_{i,i+1},g_{j,j+1})=\LM(g_{j,j+1})g_{i,i+1}-\LM(g_{i,i+1})g_{j,j+1}
\]
 f"ur $1\le i,j \le n-1$ jeweils ohne Rest aufgeht.
Aufgrund der Struktur unserer Monomordnung k"onnen wir O.E. $i=1, \,j>1$
annehmen. Es ist dann also
\begin{eqnarray*}
r_{0}:=\spol(g_{12},g_{j,j+1})&=&X_{j}Y_{j+1}(X_{1}Y_{2}-X_{2}Y_{1})-X_{1}Y_{2}(X_{j}Y_{j+1}-X_{j+1}Y_{j})\\
&=&X_{1}Y_{2}X_{j+1}Y_{j}-X_{2}Y_{1}X_{j}Y_{j+1}.
\end{eqnarray*}
Es wird $\LM(r_{0})=X_{1}Y_{2}X_{j+1}{Y_{j}}$ von $\LM(g_{12})$ geteilt, also
\[
r_{1}:=r_{0}-X_{j+1}Y_{j}(X_{1}Y_{2}-X_{2}Y_{1})=X_{2}Y_{1}X_{j+1}Y_{j}-X_{2}Y_{1}X_{j}Y_{j+1}.
\]
Schliesslich wird $\LM(r_{1})=X_{2}Y_{1}X_{j}Y_{j+1}$ von $\LM(g_{j,j+1})$
geteilt, und es ist
\[
r_{2}:=r_{1}+X_{2}Y_{1}(X_{j}Y_{j+1}-X_{j+1}Y_{j})=0,
\]
also ist $G$ eine Gr"obner-Basis.

F"ur eine Menge $M\subseteq\{X_{1},Y_{1},\ldots,X_{n},Y_{n}\}$ minimaler
M"achtigkeit mit der Eigenschaft, dass f"ur alle $g_{i,i+1}\in G$ das
Leitmonom $\LM(g_{i,i+1})=X_{i}Y_{i+1}$ wenigstens eine Variable aus $M$
enth"alt, gilt offenbar $|M|=n-1$ (z.B. $M=\{X_{1},\ldots,X_{n-1}\}$), denn
keine Variable kommt in zwei Leitmonomen $\LM(g_{i,i+1}),\, \LM(g_{j,j+1})$
mit $i\ne j$ gleichzeitig vor. Nach
\cite[Algorithm 1.2.4]{DerksenKemper} gilt dann $\dim I=2n-|M|$ und damit
$\height I=|M|=n-1$. Also ist $G$ ein phsop nach Lemma \ref{phsopHeight}. \qed\\

\subsubsection{Charakteristik 3}\label{SL2char3}
Wir gehen wieder Schritt f"ur Schritt Satz \ref{ExpliziteDarstellungen} durch,
 wobei wir wieder die nat"urliche Darstellung $V=\langle X,Y \rangle$ als
 treue Darstellung verwenden. 
Dann ist
\[
\begin{array}{rcll}
S^{3}(V)&=&S^{3}(\langle X,Y \rangle)=\langle X^{3},Y^{3},X^{2}Y,XY^{2} \rangle
& \textrm{mit Basis } \mathcal{B}\\
\textrm{und }F^{3}(V)&=&\langle X^{3},Y^{3} \rangle & \textrm{mit Basis } \mathcal{C}.
\end{array}
\]
Mit $\sigma
=\left( \begin{array}{cc} a&b\\ c&d \end{array} \right)\in
\operatorname{SL}_{2}(K)$, erhalten wir die Darstellung von  $S^{3}(V)$
bez"uglich der Basis~$\mathcal{B}$~zu
\[
\sigma \mapsto A_{\sigma}=
\left( \begin{array}{cccc}
a^{3}&b^{3}&a^{2}b&ab^{2}\\
c^{3}&d^{3}&c^{2}d&cd^{2}\\
0&0&a&-b\\
0&0&-c&d
\end{array} \right).
\]
Die linke obere Block-Matrix gibt dabei die Darstellung von $F^{3}(V)$  und
die rechte untere Block-Matrix diejenige von
$S^{3}(V)/F^{3}(V)$. Transponieren  und auswerten bei
$\sigma^{-1}$ liefert die Darstellung von $(S^{3}(V)/F^{3}(V))^{*}$, wir erhalten $\sigma \mapsto
\left( \begin{array}{cc} d&c\\ b&a \end{array} \right)$, und das ist die
Darstellung von  $\langle Y,X \rangle$. Es gilt also  $(S^{3}(V)/F^{3}(V))^{*}\cong
\langle X,Y \rangle$, und f"ur den Modul $U$ mit nichttrivialem Kozyklus aus Satz
\ref{ExpliziteDarstellungen} erhalten wir
\[
U:=\Hom_{K}(S^{3}(V),F^{3}(V))_{0} \cong F^{3}(V) \otimes (S^{3}(V)/F^{3}(V))^{*}=\langle X^{3},Y^{3} \rangle \otimes \langle X,Y \rangle.
\]
Gem"a"s  Bemerkung \ref{smallCocycles} werden wir sp"ater $U$ durch $\langle X^{3},Y^{3} \rangle \oplus
\langle X,Y \rangle$ ersetzen, so dass unser nichttrivialer Kozyklus dann in
Grad $2$ liegt.
Wir berechnen nun die Darstellungen von $U$ und $\tilde{U}$. Die Elemente von
$
U=\{ f \in \Hom_{K}(S^{3}(V),F^{3}(V)): f|_{F^{3}(V)}=0\}
$
werden bez"uglich der Basen
$\mathcal{B}$ und $\mathcal{C}$ durch Darstellungsmatrizen der Form $
\left( \begin{array}{cccc}
0 & 0 & x_{1} & x_{2}\\
0 & 0 & x_{3} & x_{4}
\end{array} \right)
$ beschrieben. Die zugeh"orige Operation von $G$ erhalten wir dann durch
\[
\sigma \cdot
\left( \begin{array}{cccc}
0 & 0 & x_{1} & x_{2}\\
0 & 0 & x_{3} & x_{4}
\end{array} \right)=
\left( \begin{array}{cc} a^{3}&b^{3}\\ c^{3}&d^{3} \end{array} \right)
\left( \begin{array}{cccc}
0 & 0 & x_{1} & x_{2}\\
0 & 0 & x_{3} & x_{4}
\end{array} \right) A_{\sigma^{-1}}.
\]
Bei dieser Operation spielt nur der rechte untere Block von $A_{\sigma^{-1}}$
eine Rolle, d.h. wir k"onnen die Operation auch durch
\[
\begin{array}{rcl}
\sigma \cdot
\left( \begin{array}{cc}
x_{1} & x_{2}\\
x_{3} & x_{4}
\end{array} \right)&=&
\left( \begin{array}{cc} a^{3}&b^{3}\\ c^{3}&d^{3} \end{array} \right)
\left( \begin{array}{cc}
x_{1} & x_{2}\\
x_{3} & x_{4}
\end{array} \right) 
\left( \begin{array}{cc} d&c\\ b&a \end{array} \right)^{T}
\end{array}
\]
beschreiben.
Aber dies ist eine Darstellung von  $\langle X^{3},Y^{3} \rangle
\otimes \langle Y,X \rangle$, wobei die Koordinaten
$(x_{1},x_{2},x_{3},x_{4})^{T}$ ihre Entsprechung in $x_{1}X^{3}\otimes
Y+x_{2}X^{3}\otimes X+x_{3}Y^{3}\otimes Y+x_{4}Y^{3}\otimes X$
haben. Bez"uglich dieser Koordinaten ist die Darstellung von $U$ gegeben durch
\begin{equation} \label{S4XY}
\sigma \mapsto
\left( \begin{array}{cc} a^{3}&b^{3}\\ c^{3}&d^{3} \end{array} \right)
\otimes
\left( \begin{array}{cc} d&c\\ b&a \end{array} \right)
=
\left( \begin{array}{cccc}
a^{3}d&a^{3}c&b^{3}d&b^{3}c\\
a^{3}b&a^{4}&b^{4}&ab^{3}\\
c^{3}d&c^{4}&d^{4}&cd^{3}\\
bc^{3}&ac^{3}&bd^{3}&ad^{3}
\end{array} \right).
\end{equation}
Dies ist zugleich die Darstellung des Untermoduls $\langle
X^{3}Y,X^{4},Y^{4},XY^{3} \rangle$ von $S^{4}(\langle X,Y \rangle)$, wie man
leicht anhand des durch  $X^{3} \otimes Y \mapsto X^{3}Y, \; X^{3} \otimes X
\mapsto X^{4}, \; Y^{3} \otimes Y \mapsto Y^{4}, \; Y^{3} \otimes X \mapsto
XY^{3}$ gegebenen Isomorphismus sieht.\\
Nun ben"otigen wir noch die Darstellung des Kozyklus, den wir gem"a"s Satz
\ref{ExpliziteDarstellungen} aus
\[
\begin{array}{rcl}
\sigma \cdot  
\left( \begin{array}{cccc}
1 & 0 & 0 & 0\\
0 & 1 & 0 & 0
\end{array} \right)&=&
\left( \begin{array}{cccc}
1 & 0 & -ab(ad+bc) & a^{2}b^{2}\\
0 & 1 & c^{2}d^{2} & -cd(ad+bc)
\end{array} \right)
\end{array}
\]
erhalten. 
Der rechte $2\times 2$ Block enth"alt die Koordinaten von $g_{\sigma}$, also
\[
g_{\sigma} = (-ab(ad+bc),a^{2}b^{2},c^{2}d^{2}, -cd(ad+bc))^{T}.
\]
Zusammen mit der Darstellung \eqref{S4XY} von $U$ erhalten wir so als
Darstellung von $\tilde{U}$:
\begin{equation} \label{US3}
\sigma \mapsto
\left( \begin{array}{ccccc}
a^{3}d&a^{3}c&b^{3}d&b^{3}c&-ab(ad+bc)\\
a^{3}b&a^{4}&b^{4}&ab^{3}&a^{2}b^{2}\\
c^{3}d&c^{4}&d^{4}&cd^{3}&c^{2}d^{2}\\
bc^{3}&ac^{3}&bd^{3}&ad^{3}&-cd(ad+bc)\\
0&0&0&0&1
\end{array} \right).
\end{equation}
Man kann $g$ als Kozyklus mit Werten in $\langle X^{3}Y,X^{4},Y^{4},XY^{3}
\rangle$ interpretieren. Dazu berechnen wir in $S^{4}(\langle X,Y \rangle)$
\[
\begin{array}{rcl}
\sigma \cdot X^{2}Y^{2}&=&-ab(ad+bc)X^{3}Y+a^{2}b^{2}X^{4}+c^{2}d^{2}Y^{4}-cd(ad+bc)XY^{3}+X^{2}Y^{2},
\end{array}
\]
d.h. $g$ entspricht dem Kozyklus $\sigma \mapsto (\sigma-1)X^{2}Y^{2}$. Wir
fassen zusammen:

\begin{Lemma}
Zu dem Untermodul $U:=\langle X^{4},X^{3}Y,XY^{3},Y^{4} \rangle$ von 
\[\tilde{U}=\langle X^{4},X^{3}Y,X^{2}Y^{2},XY^{3},Y^{4} \rangle
=S^{4}\left(\langle X,Y \rangle \right)\] existiert ein nichttrivialer
Kozyklus $g\in Z^{1}(\SL_{2},U)$,  der durch
$g_{\sigma}:=(\sigma-1)X^{2}Y^{2}$ gegeben ist. Es gilt \[
U \cong \langle X^{3},Y^{3} \rangle \otimes \langle X,Y \rangle \subseteq S^{2}\left(\langle X^{3},Y^{3} \rangle \oplus \langle X,Y \rangle \right).
\]
$U$ ist mit seinen Faktoren $\langle X^{3},Y^{3} \rangle$ und $\langle X,Y
\rangle$ selbstdual.
Eine Darstellung von  $\tilde{U}$ ist durch~(\ref{US3}) gegeben.
\end{Lemma}

\noindent Korollar \ref{depthToInfty} liefert nun sofort
\begin{Bsp}
Ist $\chr K=3$ und $\langle X,Y \rangle$ die nat"urliche Darstellung der
$\SL_{2}$, so gilt
\[
\cmdef K\left[\langle X^{3},Y^{3} \rangle \oplus \langle X,Y \rangle
\oplus \bigoplus_{i=1}^{k} S^{4}\left(\langle X,Y \rangle \right)\right]^{\SL_{2}} \ge k-2.
\]
\end{Bsp}
\noindent Die Dimension des Invariantenringes ist $5k+1$ nach Korollar
\ref{exakteDimOfSL2Invs}. (Der zugrundeliegende Modul ist weder selbstdual
noch vollst"andig reduzibel, da er $\tilde{U}$ als Summanden enth"alt).

Wir wollen dieses Beispiel noch etwas vereinfachen (insbesondere die
Dimension reduzieren), indem wir $\tilde{U}^{*}$ als Untermodul einer zweiten
symmetrischen Potenz zu finden versuchen. Dazu berechnen wir zun"achst die
Darstellung von
$\tilde{U}^{*}$, in dem wir (\ref{US3}) invertieren (d.h bei $\sigma^{-1}$
auswerten) und transponieren. Wir erhalten
\begin{equation} \label{USS3}
\sigma \mapsto
\left( \begin{array}{ccccc}
ad^{3} & -bd^{3} & -ac^{3} & bc^{3} & 0\\
-cd^{3} & d^{4} & c^{4} & -c^{3}d & 0\\
-ab^{3} & b^{4} & a^{4} & -a^{3}b & 0\\
b^{3}c & -b^{3}d & -a^{3}c & a^{3}d &0\\
bd(ad+bc) & b^{2}d^{2} & a^{2}c^{2} & ac(ad+bc) & 1
\end{array} \right).
\end{equation}
Diese Darstellung ist bez"uglich einer Basis gegeben, die als letztes Element
die annullierende Invariante $\pi$ gem"a"s Proposition \ref{AnnulatorProp} enth"alt.\\
Wir betrachten nun den Modul $M:=\langle X^{2},Y^{2},XY \rangle$. Eine kurze
Berechnung der Darstellung zeigt, dass $M$ selbstdual ist. Da $M$ keinen
eindimensionalen Untermodul enth"alt, ist $M$ dann auch irreduzibel. Nun
berechnen wir die Darstellung von
\[
S^{2}\left( M \right) = \langle (XY)Y^{2},-(Y^{2})^{2},-(X^{2})^{2},(XY)X^{2},X^{2}Y^{2}-(XY)^{2},(XY)^{2} \rangle
\]
bez"uglich der angegebenen Basis  (man beachte, dass man hier zwischen
 $(XY)^{2}$ und $X^{2}Y^{2}$ unterscheiden muss!). Eine etwas l"angliche, aber
 einfache Rechnung zeigt, dass die Darstellung von $S^{2}(M)$ durch
\begin{equation} \label{DarstS2M}
\sigma \mapsto
\left( \begin{array}{cccccc}
ad^{3} & -bd^{3} & -ac^{3} & bc^{3} & 0&-cd(ad+bc)\\
-cd^{3} & d^{4} & c^{4} & -c^{3}d & 0&-c^{2}d^{2}\\
-ab^{3} & b^{4} & a^{4} & -a^{3}b & 0&-a^{2}b^{2}\\
b^{3}c & -b^{3}d & -a^{3}c & a^{3}d &0&-ab(ad+bc)\\
bd(ad+bc) & b^{2}d^{2} & a^{2}c^{2} & ac(ad+bc) & 1 &-abcd\\
0&0&0&0&0&1
\end{array} \right)
\end{equation}
gegeben ist. Vergleichen wir diese mit der Darstellung (\ref{USS3}), so erkennen wir
$\tilde{U}^{*}$ als Untermodul von
$S^{2}(\langle X^{2},Y^{2},XY \rangle)$, wobei die annullierende Invariante
$\pi$ ihre Entsprechung in $X^{2}Y^{2}-(XY)^{2}$ hat. Die folgende Bemerkung
zeigt, dass $k$ Kopien von ihnen ein phsop im Polynomring liefern:

\begin{Bemerkung}
Seien $V_{1},\ldots,V_{k}$ $G$-Moduln, $V=V_{1}\oplus\ldots\oplus V_{k}$ und
$0 \ne f_{i}\in S^{d_{i}}(V_{i})\subseteq S(V),\,\,
i=1,\ldots,k$ (mit
$d_{i}\ge 1$). Dann bilden $f_{1},\ldots,f_{k}$ ein phsop in $S(V)$.
\end{Bemerkung}

\Bew Da $S(V_{i})$ ein Polynomring und $f_{i}\ne 0$ ist, ist
$\height(f_{i})_{S(V_{i})}=1$. Also ist $f_{i}$ ein phsop in $S(V_{i})$, und
man kann dies zu einem hsop $f_{i1}:=f_{i},f_{i2},\ldots,f_{in_{1}}$ (mit
$n_{i}=\dim V_{i}$) von $S(V_{i})$ erg"anzen. Dann ist
$f_{11},\ldots,f_{kn_{k}}$ ein hsop von $S(V)$. Denn es hat die M"achtigkeit
$n_{1}+\ldots+n_{k}=\dim S(V)$, und $S(V)$ ist ganz "uber
$K[f_{11},\ldots,f_{kn_{k}}]$, da dies offenbar f"ur die ($S(V)$ erzeugenden)
Elemente von Basen von $S^{1}(V_{i})$ gilt. \qed\\

 Der Hauptsatz
\ref{BigMainTheorem} liefert nun sofort
 
\begin{Bsp} \label{c3}
Ist $\chr K=3$ und $\langle
X,Y \rangle$ die nat"urliche Darstellung der $\SL_{2}$, so ist
\[
\cmdef K\left[\langle X,Y \rangle \oplus \langle X^{3},Y^{3} \rangle \oplus
  \bigoplus_{i=1}^{k} \langle X^{2},Y^{2},XY \rangle\right]^{\SL_{2}} \ge k-2.
\]
Der hier auftretende Modul ist selbstdual und vollst"andig reduzibel als
direkte Summe selbstdualer, irreduzibler Moduln.
\end{Bsp}
\noindent Die Dimension des Invariantenrings ist $3k+1$ nach Korollar \ref{exakteDimOfSL2Invs}.

\subsection{Beispiele f"ur $\SL_{2}$ und $\Ga$ in beliebiger positiver Charakteristik}
{\it In diesem Abschnitt sei stets $p=\chr K>0.$} \label{HauptBspFuerSL2}\\

Die von Satz \ref{ExpliziteDarstellungen} gelieferten $\SL_{2}$-Moduln $V$ mit
$\cmdef K[V]^{\SL_{2}}\ge k-2$ haben eine mit der Charakteristik $p$ schnell wachsende
Dimension, siehe Gleichung \eqref{DimOfM1} (S. \pageref{DimOfM1}). Auch im vorigen Abschnitt
hatten die Moduln f"ur Charakteristik $3$ gr"o"sere Dimension als f"ur
Charakteristik $2$. In diesem Abschnitt werden wir Beispiel \ref{HauptBspChar2}
auf beliebige Charakteristik $p>0$ verallgemeinern, und so Darstellungen
erhalten, deren Dimension nicht von $p$ abh"angt. Unser Vorgehen wird
folgendes sein: Wir geben einen nichttrivialen Kozyklus und Annullatoren f"ur
die additive Gruppe an. Die Urbilder der Annullatoren unter Roberts'
Isomorphismus bilden ein phsop im Polynomring, also auch im entsprechenden
$\SL_{2}$-Invariantenring. Nochmals Roberts' Isomorphismus angewandt zeigt,
dass die Annullatoren auch ein phsop im $\Ga$-Invariantenring bilden. Der
Hauptsatz in der allgemeineren Fassung liefert die Aussage "uber den
Cohen-Macaulay-Defekt f"ur den $\Ga$-Invariantenring. Ein letztes Mal Roberts'
Isomorphismus angewendet liefert die Aussage f"ur den $\SL_{2}$-Invariantenring.

Um diese auf den
ersten Blick vielleicht unn"otig kompliziert erscheinende Methode zu rechtfertigen, zeigen wir kurz die Probleme auf,
die beim Versuch der Verallgemeinerung etwa gem"a"s dem letzten Abschnitt
entstehen.

\subsubsection{Wo's hakt}
In Charakteristik $2$ hatten wir durch $g_{\sigma}=(\sigma-1)XY$ einen
nichttrivialen Kozyklus in $Z^{1}(\SL_{2},\langle X^{2},Y^{2} \rangle)$. In beliebiger
Charakteristik liegt jedoch $(\sigma-1)XY=abX^{2}+cdY^{2}+(ad+bc)XY$ im
Allgemeinen nicht mehr in $\langle X^{2},Y^{2} \rangle$. Eine
Verallgemeinerung auf einen Kozyklus in  $Z^{1}(\SL_{2},\langle X^{p},Y^{p}
\rangle)$ liegt nicht auf der Hand.

Ein weiterer naheliegender Ansatz w"are, von der Invariante $X\otimes
Y-Y\otimes X\in \langle X,Y \rangle \otimes \langle X,Y \rangle$
auszugehen. Diese hat jedoch f"ur $p\ne 2$ ein Komplement (d.h. der durch "Ubergang zum Dual
erhaltene Kozyklus ist trivial), wie folgende "Uberlegung zeigt: $\langle X,Y
\rangle$ ist selbstdual, und der Dual hat die Darstellung $\sigma\mapsto
\sigma^{-T}$. Bez"uglich einer Koordinatenmatrix $C\in K^{2\times 2}$, die ein
Element aus $\langle X,Y \rangle \otimes \langle Y,-X \rangle$ darstellt, ist
die Operation dann gegeben durch $\sigma \cdot C=\sigma C
(\sigma^{-T})^{T}=\sigma C\sigma^{-1}$, also durch Konjugation. Die
Einheitsmatrix $I_{2}$ ist invariant und entspricht genau dem Element $X\otimes
Y-Y\otimes X$. In Charakteristik $2$ gilt $\Spur I_{2}=0$. In Charakteristik
$p\ne 2$ dagegen ist $\Spur I_{2}\ne 0$, und die Matrizen mit Spur $0$ bilden
ein $\SL_{2}$-invariantes Komplement zu $KI_{2}$.

\subsubsection{Charakteristik-$p$-Relationen von Binomialkoeffizienten}
Wir beginnen mit einigen Charakteristik-$p$-Relationen von
Binomialkoeffizienten.

\begin{Lemma} \label{pRelations1}
Sei $\chr K=p>0$. Dann gilt f"ur $0\le i<p$ bzw. $0\le i \le j <p$ in $K$
\begin{enumerate}
\renewcommand{\labelenumi}{(\alph{enumi})}
\item $i!=(-1)^{i}\frac{(p-1)!}{(p-1-i)!}$.
\item ${j \choose i}=(-1)^{i+j}{p-1-i \choose p-1-j}$.
\item ${p-1 \choose i}=(-1)^{i}$.
\end{enumerate}
\end{Lemma}

\Bew (a) Es gilt
\begin{eqnarray*}
i!&=&i(i-1)(i-2)\cdot\ldots\cdot 3\cdot 2 \cdot
1\\
&=&(-1)^{i}(p-i)(p-i+1)(p-i+2)\cdot\ldots\cdot(p-3)(p-2)(p-1)\\
&=&(-1)^{i}\frac{(p-1)!}{(p-1-i)!}.
\end{eqnarray*}

(b) Mit (a) gilt
\begin{eqnarray*}
{j \choose
  i}&=&\frac{j!}{i!(j-i)!}=\frac{(-1)^{j}\frac{(p-1)!}{(p-1-j)!}}{(-1)^{i}\frac{(p-1)!}{(p-1-i)!}\cdot
  (j-i)!}\\&=&(-1)^{i+j}\frac{(p-1-i)!}{(p-1-j)!(j-i)!}=(-1)^{i+j}{p-1-i \choose p-1-j}.
\end{eqnarray*}

(c) Wir verwenden (b) mit $j=p-1$:
\[
{p-1 \choose i}=(-1)^{i+p-1}{p-1-i \choose p-1-(p-1)}=(-1)^{i}{p-1-i \choose 0}=(-1)^{i}.
\]
\qed\\

Mit den folgenden Regeln lassen sich auch Binomialkoeffizienten behandeln,
deren Eintr"age gr"o"ser als $p$ sind.

\begin{Lemma} \label{pRelations2}
Sei $\chr K=p>0$ und $n\ge 0$. Dann gilt 
\begin{enumerate}
\renewcommand{\labelenumi}{(\alph{enumi})}
\item Ist $0\le k<p$, so gilt ${n+p\choose k}={n\choose k}$.
\item Sei $0\le k \le n+p<2p$ (also $n<p$). Dann gilt
\[
{n+p\choose k}=\left\{\begin{array}{cl}
{n\choose k}&\textrm{ f"ur } k\le n\\
0&\textrm{ f"ur } n<k<p\\
{n\choose k-p}&\textrm{ f"ur } p\le k \le n+p.
\end{array}
\right.
\]
\end{enumerate}
\end{Lemma}

\Bew (a) Es ist\[
{n+p\choose
  k}=\frac{(n+p)(n+p-1)\cdot\ldots\cdot(n+p-k+1)}{k!}=\frac{n(n-1)\cdot\ldots\cdot(n-k+1)}{k!}={n\choose k}.\]
(b) In den ersten beiden F"allen ist die linke Seite nach (a) (wegen $n<p$ ist $k<p$)   gleich
${n\choose k}$, und im zweiten Fall ist dies gleich $0$.

Im dritten Fall $ p\le k \le n+p$ betrachten wir zun"achst ${n+p\choose k}={n+p\choose
  n+p-k}$. Da $n+p-k\le n<p$, reduziert sich dies nach dem 1. Fall zu ${n\choose
  n+p-k}={n\choose  k-p}$.
\qed\\

\subsubsection{Ein nichttrivialer Kozyklus}
Ab jetzt betrachten wir die nat"urliche Darstellung $\langle X,Y \rangle$ der
$\SL_{2}$ wieder als $\Ga$-Modul via des Homomorphismus \eqref{sigmat}
(S.\pageref{sigmat}), also 
\begin{eqnarray*}
t\cdot X&=&X\\
t\cdot Y&=&tX+Y \quad\quad \textrm{f"ur } t\in \Ga=(K,+).
\end{eqnarray*}

Wir berechnen als erstes die Darstellungen $t\mapsto A_{t}=(a_{ij}^{t})_{i,j=0,\ldots,k}\in
K^{(k+1)\times (k+1)}$ der symmetrischen Potenzen
$S^{k}(X,Y)=\langle X^{k},X^{k-1}Y,\ldots,XY^{k-1},Y^{k} \rangle$ bez"uglich
der gegebenen Basis. Aus
\[
t\cdot X^{k-j}Y^{j}=X^{k-j}(tX+Y)^{j}=X^{k-j}\sum_{i=0}^{j}{j\choose
  i}(tX)^{j-i}Y^{i}=\sum_{i=0}^{j}t^{j-i}{j\choose i}X^{k-i}Y^{i}
\]
folgt
\begin{equation} \label{DarstSkGa}
a_{ij}^{t}=\left\{\begin{array}{cl} 
t^{j-i}{j\choose i}\quad& \textrm{f"ur } 0\le i\le j \le k\\
0\quad& \textrm{f"ur } i>j.
\end{array}\right.
\end{equation}
Insbesondere ist $A_{t}$ eine unipotente obere Dreiecksmatrix. F"ur $k<l$ hat
man ein Einbettung
\[
S^{k}(\langle X,Y\rangle) \rightarrow S^{l}(\langle X,Y\rangle), \quad
f\mapsto X^{l-k}\cdot f,
\]
wobei man beide Moduln als Teilmengen von $S(\langle X,Y\rangle)$ auffasst
und daher mit $X^{l-k}\in S(\langle X,Y\rangle)$ multiplizieren darf, mit Bild
\[
\langle X^{l},X^{l-1}Y,\ldots X^{l-k}Y^{k}\rangle \cong S^{k}(\langle
X,Y\rangle).
\]
Entsprechend ist dann auf diesem Bild auch die Multiplikation mit $\frac{1}{X^{l-k}}$ als die
Umkehrung der Multiplikation mit $X^{l-k}$ definiert.

\begin{Lemma} \label{VorlKozyklus}
Sei $U:=S^{p-2}(\langle X,Y\rangle)$. Dann ist durch
\[
t\mapsto g_{t}:=\frac{1}{X}\left((t-1)\cdot Y^{p-1}\right), \quad t\in \Ga
\]
ein nichttrivialer Kozyklus $g\in Z^{1}(\Ga,U)$ gegeben und es ist
$\tilde{U}=S^{p-1}(\langle X,Y\rangle)$ der zugeh"orige erweiterte $\Ga$-Modul.
Weiter ist $\tilde{U}$ selbstdual, und die $g$ annullierende Invariante in $\tilde{U}$
gem"a"s Proposition \ref{AnnulatorProp} ist gegeben durch $X^{p-1}$.
\end{Lemma}

\Bew  Da
\[
g_{t}=\frac{1}{X}\left((tX+Y)^{p-1}-Y^{p-1}
\right)=\sum_{i=0}^{p-2}{p-1\choose i}t^{p-1-i}X^{(p-2-i)}Y^{i},
\]
gilt tats"achlich $g_{t}\in S^{p-2}(\langle X,Y\rangle) \,\myforall t\in \Ga$. Die
Kozyklus-Eigenschaft ist dann auch klar, da $X$ invariant ist. Der Koeffizient
von $Y^{p-2}$ in $g_{t}$ ist 
\begin{equation} \label{coeffyp2}
{p-1\choose p-2}t^{(p-1)-(p-2)}=-t
\end{equation}
 (Lemma
\ref{pRelations1} (c)). Dagegen ist der Koeffizient von $Y^{p-2}$ in 
\begin{equation}\label{coeffyp2zero}
(t-1)\cdot X^{p-2-j}Y^{j}=X^{p-2-j}\left((tX+Y)^{j}-Y^{j}\right)\quad
(j=0,\ldots,p-2)
\end{equation}
gleich $0$ f"ur alle $j=0,\ldots,p-2$, und damit ist der Koeffizient von $Y^{p-2}$ in
$(t-1)\cdot v$ gleich $0$ f"ur jedes $v\in U$ (denn $v$ ist Linearkombination
der $X^{p-2-j}Y^{j}, \, j=0,\ldots,p-2$). Also ist $g$ nichttrivial.

Als n"achstes zeigen wir, dass $\tilde{U}$ selbstdual ist, genauer
$\tilde{U}^{*}\cong\langle Y^{p-1},XY^{p-2},\ldots,X^{p-1} \rangle$. Da der Kozyklus in
$U\cong\langle X^{p-1},X^{p-2}Y,\ldots,XY^{p-2} \rangle \subseteq \tilde{U}$ durch
$(t-1)Y^{p-1}$ gegeben ist, insbesondere also $\tilde{U}\cong\langle
X^{p-1},X^{p-2}Y,\ldots,Y^{p-1} \rangle=S^{p-1}(\langle X,Y \rangle)$ gilt,
 zeigt dies dann gem"a"s Proposition
\ref{AnnulatorProp} auch die Aussage "uber den Annullator. Ist $t\mapsto
A_{t}$  die Darstellung von $\tilde{U}$ gem"a"s \eqref{DarstSkGa}, so hat
$\tilde{U}^{*}$ bez"uglich der entsprechenden Dualbasis die Darstellung
$t\mapsto B_{t}=(b_{ij}^{t})=A_{-t}^{T}$, also
$b_{ij}^{t}=a_{ji}^{-t}$. Drehen wir nun die Reihenfolge der Vektoren in der
Dualbasis um, und bezeichnen die zugeh"orige Darstellung mit $t\mapsto
C_{t}=(c_{ij}^{t})$, so gilt also

\begin{eqnarray*}
c_{ij}^{t}&=&b_{p-1-i,p-1-j}^{t}=a_{p-1-j,p-1-i}^{-t}\\
&\stackrel{~\eqref{DarstSkGa}}{=}&
\left\{\begin{array}{cl}
(-t)^{(p-1-i)-(p-1-j)} {p-1-i \choose p-1-j}&\textrm{f"ur }0\le p-1-j\le
p-1-i\le p-1\\
0&\textrm{sonst}
\end{array}
 \right.\\
&\stackrel{\textrm{Lemma } \ref{pRelations1}(b)}{=}&
\left\{\begin{array}{cl}
t^{j-i}{j\choose i} &\textrm{f"ur
} 0\le i \le j \le p-1\\
0&\textrm{sonst}\\
\end{array}
 \right.\\
&\stackrel{~\eqref{DarstSkGa}}{=}& a_{ij}^{t}.  
\end{eqnarray*}
Aus $A_{t}=C_{t}\,\myforall t\in \Ga$ folgt die behauptete Isomorphie.
\qed

\subsubsection{Der endg"ultige Kozyklus}
Das Lemma "uber den Kozyklus aus dem letzten Abschnitt haben wir nur
ben"otigt, um einfach an einen Annullator f"ur den folgenden Kozyklus
zu kommen. Er entsteht aus dem vorigen Kozyklus durch multiplizieren mit einer
Invariante und wird unser nichttrivialer Kozyklus zur Anwendung des
Hauptsatzes werden.

\begin{Lemma} \label{FinalerKozyklus}
F"ur $U:=\langle X^{p},Y^{p}\rangle \otimes S^{p-2}(\langle X,Y \rangle)$ ist
mit \[\Ga\rightarrow U, \, t\mapsto g_{t}:=X^{p}\otimes
\frac{1}{X}\left((t-1)\cdot Y^{p-1}\right)\] ein nichttrivialer Kozyklus
$g\in Z^{1}(\Ga,U)$ gegeben.
\end{Lemma}

\Bew Wir bezeichnen mit $h$ den nichttrivialen Kozyklus aus Lemma
\ref{VorlKozyklus}, also $h_{t}=\frac{1}{X}\left((t-1)\cdot Y^{p-1}\right)\in
S^{p-2}(\langle X,Y \rangle)$. Da $X^{p}\in\langle X^{p},Y^{p}\rangle$
invariant unter $\Ga$ ist, ist also auch $t\mapsto g_{t}=X^{p}\otimes h_{t}$ ein
Kozyklus. Wir m"ussen zeigen, dass $g$ nichttrivial ist, und rechnen hierzu
mit der Basis $\mathcal{B}:=\{X^{p}\otimes X^{p-2},\ldots,X^{p}\otimes Y^{p-2}, Y^{p}\otimes
X^{p-2},\ldots, Y^{p}\otimes Y^{p-2}\}$.

Der Koeffizient von $X^{p}\otimes Y^{p-2}$  f"ur $j=0,\ldots,p-2$ in
\[
(t-1)\left(X^{p}\otimes X^{p-2-j}Y^{j}\right)=X^{p}\otimes \left((t-1)
  X^{p-2-j}Y^{j}\right)
\]
ist gleich $0$ gem"a"s Gleichung \eqref{coeffyp2zero}.

Der Koeffizient von $X^{p}\otimes Y^{p-2}$ f"ur $j=0,\ldots,p-2$ in
\begin{eqnarray*}
&&(t-1)\left(Y^{p}\otimes
  X^{p-2-j}Y^{j}\right)\\
&=&\left(t^{p}X^{p}+Y^{p}\right)\otimes \left(X^{p-2-j}(tX+Y)^{j}\right)-Y^{p}\otimes  X^{p-2-j}Y^{j}
\end{eqnarray*}
ist jedenfalls gleich $0$ f"ur $j=0,\ldots,p-3$ (wegen des Faktors $X^{p-2-j}$
auf der rechten Seite aller auftretenden tensoriellen Produkte); F"ur $j=p-2$
ist der Koeffizient von $X^{p}\otimes Y^{p-2}$ in
\[
\left(t^{p}X^{p}+Y^{p}\right)\otimes
(tX+Y)^{p-2}-Y^{p}\otimes  Y^{p-2}
\]
dagegen gleich $t^{p}$.

F"ur ein beliebiges $v\in U$ (welches Linearkombination der Basiselemente
$\mathcal{B}$ ist), ist der Koeffizient von $X^{p}\otimes Y^{p-2}$ in
$(t-1)v$ also gleich \[\lambda \cdot t^{p},\] wobei $\lambda\in K$ der
Koeffizient von $Y^{p}\otimes Y^{p-2}$ in $v$ ist.  

Der Koeffizient von $X^{p}\otimes
Y^{p-2}$ in der Darstellung von $g_{t}=X^{p}\otimes h_{t}$ bez"uglich
$\mathcal{B}$ ist nach Gleichung 
\eqref{coeffyp2} dagegen gleich \[-t.\]
H"atten wir also
$g_{t}=(t-1)v\,\, \myforall t\in \Ga$, so w"are \[\lambda \cdot t^{p}=-t \quad \myforall
t\in K.\] Da $|K|=\infty$, ist dies ein Widerspruch. Also ist $g$ ein
nichttrivialer Kozyklus. \qed\\

\subsubsection{Annullatoren des Kozyklus}\label{AnnulatorenDesKoz}
Als n"achstes geben wir vier verschiedene Typen von Annullatoren des Kozyklus
aus dem letzten Lemma an. Dabei betten wir den Kozyklus gleich in den
Polynomring ein, auf den wir dann den Hauptsatz anwenden wollen. Seien dazu ab
jetzt $\langle X_{i},Y_{i} \rangle \cong \langle X,Y \rangle $ f"ur $i\ge 1$ Kopien der nat"urlichen
Darstellung und $\langle X_{0}, Y_{0} \rangle:\cong \langle X^{p},Y^{p} \rangle$,
also 
\begin{eqnarray*}
t\cdot X_{0}&=&X_{0}\\
t\cdot Y_{0}&=&t^{p}X_{0}+Y_{0} \quad\quad \textrm{f"ur } t\in \Ga=(K,+).
\end{eqnarray*}
(Entsprechend auch, falls wir $\langle X_{0}, Y_{0} \rangle$ als
$\SL_{2}$-Modul auffassen.)

Wir erinnern an folgende Gleichung: Ist $V=\bigoplus_{i=1}^{n}V_{i}$ (mit
$V_{i}$ jeweils $G$-Modul), so ist
\begin{equation} \label{SVOfDirektSum}
S(V)=\bigoplus_{i_{1},\ldots,i_{n}\ge 0}S^{i_{1}}(V_{1})\otimes\ldots\otimes S^{i_{n}}(V_{n})
\end{equation}
Dabei ist $S^{0}(V_{i})\cong K$, so dass man in den Summanden Faktoren mit $i_{k}=0$ weglassen kann.\\

{\it F"ur den Rest dieses Abschnitts behalten wir die Bezeichnungen des
  folgenden Korollars bei.}

\begin{Korollar} \label{KozyklusInSV}
Sei $V:=\langle X_{0}, Y_{0} \rangle \oplus \bigoplus_{i=1}^{k}\langle
X_{i},Y_{i} \rangle.$ Dann ist durch
\[
t\mapsto g_{t}:=X_{0}\cdot \frac{1}{X_{1}}\left((t-1)\cdot Y_{1}^{p-1}\right) 
\]
ein nichttrivialer Kozyklus $g\in Z^{1}(\Ga,S(V))$ gegeben.
\end{Korollar}

\Bew Aufgefasst als Kozyklus in $Z^{1}(\Ga,U)$ mit $U:=\langle X_{0},Y_{0}\rangle \otimes
S^{p-2}(\langle X_{1},Y_{1}\rangle)$ ist $g$ nichttrivial gem"a"s Lemma
\ref{FinalerKozyklus}. Da aber $U$ nach \eqref{SVOfDirektSum} ein
\emph{direkter} Summand von $S(V)$ ist, ist $g$ auch aufgefasst als Kozyklus
in $Z^{1}(\Ga,S(V))$ nichttrivial. \qed\\

Wir k"onnen nun die ersten beiden Typen von Annullatoren angeben:

\begin{Lemma} \label{AnnulatorTyp1}
Der Kozyklus $g$ wird annulliert von den Invarianten $X_{1},X_{i}^{p-1}\in
S(V)^{\Ga}$ mit $i\ge 2$, d.h. $X_{1}g=X_{i}^{p-1}g=0\in H^{1}(\Ga,S(V))$.
\end{Lemma}

\Bew Nach Definition von $g$ ist 
\[
X_{1}g_{t}=X_{0}\cdot\left((t-1)\cdot Y_{1}^{p-1}\right)=(t-1)\cdot X_{0}Y_{1}^{p-1}, \]
also $X_{1}g\in B^{1}(\Ga,S(V))$.

Wir betrachten nun den Kozyklus $t\mapsto
h_{t}:=\frac{1}{X_{1}}\left((t-1)\cdot Y_{1}^{p-1}\right)\in S^{p-2}(\langle
X_{1},Y_{1}\rangle)=:U$. Nach Lemma \ref{VorlKozyklus} gilt
$\tilde{U}^{*}\cong S^{p-1}(\langle X_{i},Y_{i} \rangle)$ und $X_{i}^{p-1}$
ist der zugeh"orige Annullator, also $X_{i}^{p-1}\otimes h = 0 \in H^{1}(\Ga,
U\otimes \tilde{U}^{*})$. Mit den offensichtlichen Einbettungen von $U$ und
$\tilde{U}^{*}$ in $S(V)$ ist dann auch $X_{i}^{p-1}h=0\in H^{1}(\Ga,S(V))$,
d.h. es gibt ein $v\in S(V)$ mit $X_{i}^{p-1}h_{t}=(t-1)v \,\myforall t\in
\Ga$. Dann ist aber auch
\[
X_{i}^{p-1}g_{t}=X_{i}^{p-1}X_{0}h_{t}=(t-1)(X_{0}v)\quad \myforall t\in \Ga,
\]
also $X_{i}^{p-1}g\in B^{1}(\Ga,S(V))$. \qed\\

Leider bilden $X_{1},X_{2}^{p-1},\ldots,X_{i}^{p-1}$ f"ur $i\ge 3$ kein phsop
in $S(V)^{\Ga}$. Daher ben"otigen wir weitere Typen von Annullatoren.

\begin{Lemma} \label{AnnullatorTyp2}
Der Kozyklus $g$ wird annulliert von den Invarianten $X_{1}Y_{i}-X_{i}Y_{1}\in
S(V)^{\Ga}$ mit $i\ge 2$. Genauer gilt
\[
\left(X_{1}Y_{i}-X_{i}Y_{1}\right)g_{t}=(t-1)\cdot \left(Y_{1}^{p-1}Y_{i}X_{0}-X_{1}^{p-1}X_{i}Y_{0}
\right) \quad \myforall t\in \Ga.
\]
\end{Lemma}

\Bew O.E. sei $i=2$. Die rechte Seite der behaupteten Gleichung ist
\begin{eqnarray*}
&&\left(tX_{1}+Y_{1}\right)^{p-1}\left(tX_{2}+Y_{2}\right)X_{0}-X_{1}^{p-1}X_{2}\left(t^{p}X_{0}+Y_{0}\right)-\left(Y_{1}^{p-1}Y_{2}X_{0}-X_{1}^{p-1}X_{2}Y_{0}
\right)\\
&=&\sum_{j=0}^{p-1}{p-1 \choose j}t^{j}X_{1}^{j}Y_{1}^{p-1-j}(tX_{2}+Y_{2})X_{0}-t^{p}X_{1}^{p-1}X_{2}X_{0}-Y_{1}^{p-1}Y_{2}X_{0},
\end{eqnarray*}
wobei sich der Term $X_{1}^{p-2}X_{2}Y_{0}$ weggehoben hat. Der zweite Term
$-t^{p}X_{1}^{p-1}X_{2}X_{0}$ hebt sich mit dem Term aus der Summe f"ur
$j=p-1$ und dem Faktor $tX_{2}$ (aus $(tX_{2}+Y_{2})$) weg. Der dritte Term
$-Y_{1}^{p-1}Y_{2}X_{0}$ hebt sich mit dem Term aus der Summe f"ur $j=0$ und
dem Faktor $Y_{2}$ weg. Zusammen mit Lemma \ref{pRelations1} (c) bleibt daher
\begin{eqnarray*}
&&X_{0}X_{2}Y_{1}\sum_{j=0}^{p-2}(-1)^{j}t^{j+1}X_{1}^{j}Y_{1}^{p-2-j}+X_{0}X_{1}Y_{2}\sum_{j=0}^{p-2}(-1)^{j+1}t^{j+1}X_{1}^{j}Y_{1}^{p-2-j}\\
&=&(X_{1}Y_{2}-X_{2}Y_{1})X_{0}\sum_{j=0}^{p-2}(-1)^{j+1}t^{j+1}X_{1}^{j}Y_{1}^{p-2-j}\\
&=&(X_{1}Y_{2}-X_{2}Y_{1})X_{0}\frac{1}{X_{1}}\left((tX_{1}+Y_{1})^{p-1}-Y_{1}^{p-1}\right)\\
&=&(X_{1}Y_{2}-X_{2}Y_{1})g_{t},
\end{eqnarray*}
also die linke Seite und damit die Behauptung. \qed\\

Da die Annullatoren dieses Lemmas alle im Summanden $\langle X_{1},Y_{1}
\rangle$ "`verankert"' sind, k"onnen wir wieder nur maximal zwei f"ur ein
phsop verwenden, und etwa in der Kombination
$X_{1},X_{2}^{p-1},X_{1}Y_{3}-X_{3}Y_{1}$ sogar nur einen. Man k"onnte
vermuten, dass auch die $X_{i}Y_{j}-X_{j}Y_{i}$ Annullatoren sind, doch leider
ist dem nicht so. Nach Erheben in die $p-1$-te Potenz sind sie es aber, und
dies liefert uns den letzten "`Typ"' von Annullatoren.

\begin{Lemma} \label{AnnullatorTyp3}
Der Kozyklus $g$ wird annulliert von den Invarianten $(X_{i}Y_{j}-X_{j}Y_{i})^{p-1}\in
S(V)^{\Ga}$ mit $i,j\ge 1$.
\end{Lemma}

\Bew Sei O.E. $i=2, j=3$, und $t\in \Ga$.
Wir nummerieren in diesem Beweis alle Zeilen und Spalten von Matrizen von $0$ beginnend,
es ist also z.B. $e_{0}=(1,0,\ldots)$, $e_{1}=(0,1,\ldots)$ usw. jeweils ein
$0$ter bzw. $1$ter Basisvektor. Wir setzen
\[
M_{2}:=\langle X_{2}^{p-1},X_{2}^{p-2}Y_{2},\ldots,Y_{2}^{p-1}\rangle
\]
und
\[
M_{3}:=\langle Y_{3}^{p-1},X_{3}Y_{3}^{p-2},\ldots,X_{3}^{p-1}\rangle.
\]
Die Darstellung von $M_{2}$ ist nach \eqref{DarstSkGa} gegeben durch $t\mapsto
A_{t}=(a_{ij}^{t})$ mit
\begin{equation} \label{nochmalDarstGa}
a_{ij}^{t}=\left\{\begin{array}{cl} 
t^{j-i}{j\choose i}\quad& \textrm{f"ur } 0\le i\le j \le p-1\\
0\quad& \textrm{f"ur } i>j.
\end{array}\right.
\end{equation}
Nach dem Beweis von Lemma \ref{VorlKozyklus} gilt
$M_{3}\cong M_{2}^{*}$, und zwar so, dass die Basis von $M_{3}$ der
dualen Basis zu $M_{2}^{*}$ entspricht. Also hat $M_{3}$ die Darstellung
$t\mapsto A_{-t}^{T}$. Wir betrachten das Tensorprodukt $M_{2}\otimes
M_{3}$. Wir identifizieren $M_{2}\otimes M_{3}$ mit dem entsprechenden
Untermodul von $S(V)$. Weiter identifizieren wir die Elemente von $M_{2}\otimes
M_{3}$ mit ihren zugeh"origen Koordinatenmatrizen $X\in K^{p\times
  p}$. Die Operation von $\Ga$ ist dann gegeben durch 
\[
t\cdot X=A_{t}XA_{-t}^{TT}=A_{t}XA_{-t}.
\]
Sei $\pi=I_{p\times p}\in K^{p\times p}$ die $p\times p$ Einheitsmatrix. Wir sehen sofort,
dass $\pi$ unter der Operation von $\Ga$ invariant ist. Es gilt
\begin{eqnarray*}
\pi&=&\sum_{i=0}^{p-1}X_{2}^{p-1-i}Y_{2}^{i}\otimes X_{3}^{i}Y_{3}^{p-1-i}\\
&\stackrel{\textrm{Lemma } ~\ref{pRelations1} (c)}{=}&\sum_{i=0}^{p-1} {p-1 \choose i}(-1)^{i}X_{2}^{p-1-i}Y_{2}^{i}X_{3}^{i}Y_{3}^{p-1-i}\\
&=&(X_{2}Y_{3}-X_{3}Y_{2})^{p-1}.
\end{eqnarray*}
Wir sehen, dass $\pi$ der Invariante entspricht, von der wir behaupten,
dass sie $g$ annulliert. Sei nun $U$ wie in Lemma \ref{FinalerKozyklus}.
Wenn wir nun
einen $\pi$ enthaltenden Untermodul $M$ von $M_{2}\otimes M_{3}$ angeben, der zu
$\tilde{U}^{*}$ isomorph ist, so dass $\pi$ dessen annullierender Invariante gem"a"s
Proposition \ref{AnnulatorProp} entspricht, so sind wir fertig. Wir setzen
\[
M:=\langle v_{0},v_{1},\ldots,v_{p-2},w_{0},\ldots,w_{p-2},\pi \rangle
\subseteq K^{p\times p}
\]
mit folgenden Basisvektoren: 
$v_{i}\in K^{p\times p}, i=0,\ldots,p-2$ seien
die Matrizen, die genau in der $i+1$-ten oberen Nebendiagonale Einsen haben, also
\[
v_{0}:=\left(\begin{array}{ccccc}
0&1&0&\dots&0\\
&\ddots&\ddots&\ddots&\vdots\\
&&\ddots&\ddots&0\\
&&&\ddots&1\\
&&&&0\\
\end{array}
\right),\ldots,
v_{p-2}:=\left(\begin{array}{ccccc}
0&\dots&\dots&0&1\\
&\ddots&\ddots&\ddots&0\\
&&\ddots&\ddots&\vdots\\
&&&\ddots&0\\
&&&&0\\
\end{array}
\right)\in K^{p\times p}.
\]
Entsprechend seien 
$w_{i}\in K^{p\times p}, i=0,\ldots,p-2$
die Matrizen, die genau in der $i+1$-ten unteren Nebendiagonale Einsen haben,
also $w_{i}=v_{i}^{T}$. Dann besteht $M$ genau aus den Matrizen mit
"`konstanten Diagonalen"'. Wir zeigen, dass $M$ ein $\Ga$-Modul ist
(d.h. $\Ga\cdot M \subseteq M$) mit $M\cong \tilde{U}^{*}$ so, dass $\pi$ der
annullierenden Invariante von $\tilde{U}^{*}$ entspricht.

Als erstes bestimmen wir die $\Ga$-Operation auf den $v_{k}$, und berechnen dazu die $i$-te Zeile und $j$-te Spalte von $t\cdot
v_{k}$, mit $i,j=0,\ldots,p-1$. Seien dazu $e_{0}=(1,0,\ldots,0)^{T}, \,e_{1}=(0,1,0,\ldots,0)^{T}, \ldots,\,
e_{p-1}=(0,\ldots,0,1)^{T}\in K^{p}$ die Spalteneinheitsvektoren. Dann ist
\begin{eqnarray*}
(t\cdot
v_{k})_{i,j}&=&(A_{t}v_{k}A_{-t})_{i,j}=e_{i}^{T}A_{t}v_{k}A_{-t}e_{j}\\
&=&(a_{i0}^{t},\ldots,a_{i,p-1}^{t})v_{k}\left(
\begin{array}{c}
a_{0j}^{-t}\\
a_{1j}^{-t}\\
\vdots\\
a_{p-1,j}^{-t}
\end{array}\right)\\
&=&(0,\ldots,0,a_{i0}^{t},\ldots,a_{i,p-2-k}^{t})\left(
\begin{array}{c}
a_{0j}^{-t}\\
a_{1j}^{-t}\\
\vdots\\
a_{p-1,j}^{-t}
\end{array}\right),
\end{eqnarray*}
wobei man bei dem Zeilenvektor $k+1$ f"uhrende Nullen hat.
Mit Gleichung \eqref{nochmalDarstGa} erhalten wir also
\begin{eqnarray*}
(t\cdot
v_{k})_{i,j}&=&\sum_{l=0}^{p-2-k}a_{i,l}^{t}a_{k+1+l,j}^{-t}=\sum_{l=i}^{j-k-1}a_{i,l}^{t}a_{k+1+l,j}^{-t}\\
&=&\sum_{l=i}^{j-k-1}t^{l-i}{l\choose i}(-t)^{j-k-1-l}{j \choose k+1+l}\\
&=&\sum_{l=i}^{j-k-1}t^{j-k-i-1}(-1)^{j-k-l-1}{l\choose i}{j \choose k+l+1}.
\end{eqnarray*}
Bei dem zweiten Gleichheitszeichen haben wir dabei f"ur die Summationsindizes
verwendet, dass $A_{t}$ eine obere Dreiecksmatrix ist. Die leere Summe (falls
$i>j-k-1$) ist hier wie "ublich als $0$ zu lesen. Ab jetzt sei daher
\begin{equation}\label{jkiGe0}
 0\le j-k-i-1.
\end{equation}
Wir verschieben nun $l$
um $-i$ und erhalten
\begin{eqnarray*}
(t\cdot
v_{k})_{i,j}&=&\sum_{l=0}^{j-k-i-1}t^{j-k-i-1}(-1)^{j-k-l-i-1}{l+i\choose i}{j
  \choose k+l+i+1}\\
&\stackrel{\textrm{Lemma }~\ref{pRelations1} (b)}{=}&\sum_{l=0}^{j-k-i-1}t^{j-k-i-1}(-1)^{j-k-i-1}{p-1-i\choose p-1-i-l}{j
  \choose k+l+i+1}.
\end{eqnarray*}
Dabei d"urfen wir Lemma \ref{pRelations1} (b) anwenden, denn $0\le i \le
l+i\le j-k-1 \le p-2<p$

Nun wenden wir auf beide Binomialkoeffizienten die Regel ${n\choose k}={n
  \choose n-k}$ an und erhalten
\begin{eqnarray*}
(t\cdot
v_{k})_{i,j}&=&\sum_{l=0}^{j-k-i-1}(-t)^{j-k-i-1}{p-1-i\choose l}{j
  \choose (j-k-i-1)-l}.
\end{eqnarray*}
Mit der Formel ${m+n\choose k}=\sum_{j=0}^{k}{m\choose j}{n \choose k-j}$
ergibt sich schliesslich
\begin{eqnarray*}
(t\cdot
v_{k})_{i,j}&=&(-t)^{j-k-i-1}{p+j-i-1\choose j-i-1-k}.
\end{eqnarray*}
Da $0 \stackrel{~\eqref{jkiGe0}}{\le} j-i-1-k\le j-i-1\le p-1$,  folgt mit Lemma~\ref{pRelations2}~(a) weiter
\begin{eqnarray*}
(t\cdot
v_{k})_{i,j}&=&(-t)^{j-k-i-1}{j-i-1\choose j-i-1-k}\\
&=&(-t)^{(j-i-1)-k}{j-i-1\choose k}.
\end{eqnarray*}
Zusammen mit $(t\cdot v_{k})_{i,j}=0$ f"ur $k>j-i-1$ erhalten wir damit durch
Vergleich mit \eqref{nochmalDarstGa}
\[
(t\cdot
v_{k})_{i,j}=\left\{\begin{array}{cl}
a_{k,j-i-1}^{-t}& \textrm{f"ur }j-i-1\ge 0\\
0&\textrm{sonst.}
\end{array}
\right.
\]
Dies bedeutet, dass die Matrix $t\cdot v_{k}\in K^{p\times p}$ eine nilpotente
obere Dreiecksmatrix ist und in der oberen Nebendiagonalen Nr. $j-i$ konstant der
Eintrag $a_{k,j-i-1}^{-t}$ steht. Also gilt
\[
t\cdot
v_{k}=\sum_{l=0}^{p-2}a_{k,l}^{-t}v_{l},\quad k=0,\ldots,p-2,
\]

oder {\bf Zusammengefasst}\label{zusfassdarstga}: Sei $B_{t}\in K^{(p-1)\times (p-1)}$ die linke obere $(p-1)\times (p-1)$
Teilmatrix von $A_{t}$. Dann ist $\langle v_{0},\ldots,v_{p-2} \rangle$ ein Untermodul
mit Darstellung $t\mapsto B_{-t}^{T}$.\\

Im n"achsten Schritt bestimmen wir die Operation auf den $w_{k}, \,
k=0,\ldots,p-2$. Es ist
\begin{eqnarray*}
(t\cdot
w_{k})_{i,j}&=&(A_{t}w_{k}A_{-t})_{i,j}=e_{i}^{T}A_{t}w_{k}A_{-t}e_{j}\\
&=&(a_{i0}^{t},\ldots,a_{i,p-1}^{t})w_{k}\left(
\begin{array}{c}
a_{0j}^{-t}\\
a_{1j}^{-t}\\
\vdots\\
a_{p-1,j}^{-t}
\end{array}\right)\\
&=&(a_{i,k+1}^{t},\ldots,a_{i,p-1}^{t},0,\ldots,0)\left(
\begin{array}{c}
a_{0j}^{-t}\\
a_{1j}^{-t}\\
\vdots\\
a_{p-1,j}^{-t}
\end{array}\right)\\
&=&\sum_{l=0}^{p-k-2}a_{i,k+1+l}^{t}a_{l,j}^{-t}.
\end{eqnarray*}
Mit Gleichung \eqref{nochmalDarstGa} erhalten wir also
\begin{eqnarray*}
(t\cdot
w_{k})_{i,j}&=&\sum_{l=\max(0,i-k-1)}^{\min(p-k-2,j)}t^{k+1+j-i}(-1)^{j-l}{k+1+l
\choose i}{j \choose l}.
\end{eqnarray*}
Da $i<p$, ist diese Summe genau dann leer, wenn $i-k-1>j$. Sei daher ab jetzt
\begin{equation} \label{ikj1Ge0}
0\le j-i+k+1 .
\end{equation}
Da $i\le k+1+l\le k+1+p-k-2=p-1$ erhalten wir mit Lemma \ref{pRelations1} (b)
\begin{eqnarray*}
(t\cdot
w_{k})_{i,j}&=&\sum_{l=\max(0,i-k-1)}^{\min(p-k-2,j)}t^{k+1+j-i}(-1)^{k+1+j-i}{p-1-i
  \choose p-2-k-l}{j \choose l}\\
&=&(-t)^{k+1+j-i} {p-1-i+j\choose p-2-k},
\end{eqnarray*}
wobei wir nochmals dir Formel ${m+n\choose k}=\sum_{j=0}^{k}{m\choose j}{n
  \choose k-j}=\sum_{j=\max(0,k-n)}^{\min(k,m)}{m\choose j}{n \choose k-j}$
  verwendet haben. Wir machen nun eine Fallunterscheidung:

\underline{1. Fall:}  $0\le i-j-1$. Dann ist nat"urlich $i-j-1\le p-2$, und wir erhalten
\begin{eqnarray*}
(t\cdot w_{k})_{i,j}&=&(-t)^{(p-2-(i-j-1))-(p-2-k)}{p-2-(i-j-1)\choose
  p-2-k}\\
&=&a_{p-2-k,p-2-(i-j-1)}^{-t}.
\end{eqnarray*}
Wir sehen, dass dies auch im Fall der leeren Summe, also wenn $i-j-1>k$,
richtig ist.

\underline{2. Fall:} $0=i-j$. Dann ist
\begin{eqnarray*}
(t\cdot w_{k})_{i,j}&=&(-t)^{k+1} {p-1\choose p-2-k}\\
&=&a_{p-2-k,p-1}^{-t}.
\end{eqnarray*}

\underline{3. Fall:} $0 \le j-i-1$. Jedenfalls ist
\begin{eqnarray*}
(t\cdot w_{k})_{i,j}&=&(-t)^{k+1+j-i}{p+(j-i-1)\choose p-2-k}.
\end{eqnarray*}

\underline{3.1. Fall:} $0\le j-i-1<p-2-k$. Nach Lemma \ref{pRelations2} (b, 2. Fall) ist dann
\begin{eqnarray*}
(t\cdot w_{k})_{i,j}&=&0.
\end{eqnarray*}

\underline{3.2. Fall:} $p-2-k\le j-i-1$. Dann ist nach Lemma \ref{pRelations2} (b, 1. Fall)
\begin{eqnarray*}
(t\cdot w_{k})_{i,j}&=&(-t)^{p}(-t)^{j-i-1-(p-2-k)}{j-i-1\choose p-2-k}\\
&=&(-t)^{p}a_{p-2-k,j-i-1}^{-t}.
\end{eqnarray*}
Da $A_{-t}$ eine obere Dreiecksmatrix ist, ist dieses Ergebnis auch im Fall 3.1
 richtig.

Wir fassen zusammen:
\begin{eqnarray*}
(t\cdot w_{k})_{i,j}&=&\left\{\begin{array}{cl}
a_{p-2-k,p-2-(i-j-1)}^{-t}&\textrm{ f"ur }0\le i-j-1\\
a_{p-2-k,p-1}^{-t}&\textrm{ f"ur } i=j\\
(-t)^{p}a_{p-2-k,j-i-1}^{-t}&\textrm{ f"ur } 0\le j-i-1
\end{array}\right.
\end{eqnarray*}
Wir sehen zun"achst, dass das Ergebnis nur von der Differenz $i-j$ abh"angt,
also dass in jeder Nebendiagonalen gleiche Eintr"age stehen. Insgesamt
erhalten wir
\[
t\cdot w_{k}=\sum_{j=0}^{p-2}a_{p-2-k,p-2-j}^{-t}w_{j}+\sum_{j=0}^{p-2}(-t)^{p}a_{p-2-k,j}^{-t}v_{j}+a_{p-2-k,p-1}^{-t}\pi.
\]
Wir setzen nun $\tilde{w}_{k}:=w_{p-2-k}$ f"ur $0\le k\le p-2$ und $\tilde{w}_{p-1}:=\pi$.
Dann ist f"ur $0\le k\le p-2$
\begin{eqnarray*}
t\cdot \tilde{w}_{k}&=&t\cdot
w_{p-2-k}\,=\,\sum_{j=0}^{p-2}a_{k,p-2-j}^{-t}w_{j}+\sum_{j=0}^{p-2}(-t)^{p}a_{k,j}^{-t}v_{j}+a_{k,p-1}^{-t}\pi\\
&=&\sum_{j=0}^{p-1}a_{k,j}^{-t}\tilde{w}_{j}+\sum_{j=0}^{p-2}(-t)^{p}a_{k,j}^{-t}v_{j}
\end{eqnarray*}
Weiter ist auch
\begin{eqnarray*}
t\cdot \tilde{w}_{p-1}&=&t\cdot
\pi\,=\,\pi\\
&=&\sum_{j=0}^{p-1}a_{p-1,j}^{-t}\tilde{w}_{j},
\end{eqnarray*}
denn $A_{-t}$ ist eine unipotente obere Dreiecksmatrix. Wir schreiben nun
\[
A_{t}=\left(\begin{array}{cc}
B_{t}&h_{t}\\
0_{1\times (p-1)}&1
\end{array}\right),
\]
d.h. $h_{t}\in K^{p-1}$ sind die ersten $p-1$ Zeilen der letzten Spalte von
$A$. Mit der
Zusammenfassung von S. \pageref{zusfassdarstga} erhalten wir nun:

Die Darstellung von $M=\langle
v_{0},\ldots,v_{p-2},\tilde{w}_{0},\ldots,\tilde{w}_{p-1} \rangle$ ist gegeben
durch
\begin{eqnarray*}
t\mapsto C_{t}&:=&\left(\begin{array}{c|c}
B_{-t}^{T}&\begin{array}{c|c}(-t)^{p}B_{-t}^{T}&0_{(p-1)\times 1}\end{array}\\
\hline
0_{p\times (p-1) }&A_{-t}^{T}
\end{array}\right)\in K^{(2p-1)\times (2p-1)}\\
&=&
\left(\begin{array}{ccc}
B_{-t}^{T}&(-t)^{p}B_{-t}^{T}&\\
&B_{-t}^{T}&\\
&h_{-t}^{T}&1
\end{array}\right).
\end{eqnarray*}
Wir interpretieren nun $M$ als ein $\tilde{W}^{*}$ und gehen Proposition
\ref{AnnulatorProp} "`r"uckw"arts"', um zu einem Modul $W$ mit einem Kozyklus
zu gelangen, der von $\pi=\tilde{w}_{p-1}$ annulliert wird. Dann hat $\tilde{W}=M^{*}$ die
Darstellung
\[
t\mapsto C_{-t}^{T}=\left(\begin{array}{ccc}
B_{t}&&\\
t^{p}B_{t}&B_{t}&h_{t}\\
&&1
\end{array}\right).
\]
Streichen wir hiervon die letzte Zeile und Spalte, so erhalten wir die Matrix
$\left(\begin{array}{cc} 1\\t^{p}&1 \end{array}\right) \otimes B_{t}$, was
gerade der Darstellung von $U=\langle Y^{p},X^{p} \rangle \otimes
S^{p-2}(\langle X,Y \rangle)$ entspricht. In der rechten unteren $2\times 2$
Blockmatrix von $C_{-t}^{T}$ steht
die Darstellung von $S^{p-1}(\langle X,Y\rangle)$. Damit entspricht der letzten Spalte von
$C_{-t}^{T}$ gerade die Erweiterung von $U$ zu $\tilde{U}$ durch den Kozyklus
$t\mapsto X^{p}\otimes \frac{1}{X}(t-1)Y^{p-1}$ (vgl. Gleichung
\eqref{darstVSchlange}, S. \pageref{darstVSchlange}). Es gilt also
$\tilde{W}\cong \tilde{U}$ oder $M\cong
\tilde{U}^{*}$, und $\pi=\tilde{w}_{p-1}\in M$ entspricht genau der annullierenden Invariante
in $\tilde{U}^{*}$. Dies wollten wir zeigen. \qed\\
\begin{BemRoman}
Mit der Formel \eqref{Annulator}, S.\pageref{Annulator}  folgt 
\begin{eqnarray*}
-\left(X_{2}Y_{3}-X_{3}Y_{2}\right)^{p-1}g_{t}&=&(t-1)\cdot
\Bigg(Y_{0}\sum_{i=0}^{p-2}\sum_{j=0}^{p-2-i}(X_{1}^{p-2-i}Y_{1}^{i})(X_{2}^{p-1-j}Y_{2}^{j})(X_{3}^{i+j+1}Y_{3}^{p-2-i-j})\\
&&+X_{0}\sum_{i=0}^{p-2}\sum_{j=0}^{p-2-i}(X_{1}^{i}Y_{1}^{p-2-i})(X_{2}^{j}Y_{2}^{p-1-j})(X_{3}^{p-2-i-j}Y_{3}^{i+j+1})
\Bigg).
\end{eqnarray*}
Ausgehend von dieser Formel, die ich durch Experimente mit \Magma geraten
habe, ist das Lemma und der Beweis entstanden.
\end{BemRoman}

\subsubsection{Ein phsop}
Wir w"ahlen nun aus den annullierenden Invarianten ein phsop aus:

\begin{Lemma} \label{AnnullatorPhsopOfGa}
Die annullierenden Invarianten
\[X_{1},X_{2}^{p-1},X_{1}Y_{3}-X_{3}Y_{1},(X_{i}Y_{i+1}-X_{i+1}Y_{i})^{p-1}\in S(V)^{\Ga}\] mit
  $i=3,\ldots,k-1$ bilden ein phsop der L"ange $k$ in $S(V)^{\Ga}$.
\end{Lemma}

\Bew Wir betrachten
\[
V:=\langle X_{0}, Y_{0} \rangle \oplus \bigoplus_{i=1}^{k}\langle
X_{i},Y_{i} \rangle
\]
und $V\oplus \langle X_{k+1},Y_{k+1}\rangle$ sowohl als $\Ga$ als auch als
$\SL_{2}$-Moduln. Nach Roberts' Isomorphismus, Korollar~\ref{RobertsSymmAlg}, gilt
$S(V\oplus \langle X_{k+1},Y_{k+1}\rangle)^{\SL_{2}}\cong S(V)^{\Ga}$, wobei
der Isomorphismus durch Einsetzen von $X_{k+1}=0$ und $Y_{k+1}=1$ gegeben
ist. Insbesondere werden die Invarianten
\begin{equation} \label{SL2phsopcollection}
X_{1}Y_{k+1}-X_{k+1}Y_{1},(X_{2}Y_{k+1}-X_{k+1}Y_{2})^{p-1},X_{1}Y_{3}-X_{3}Y_{1},(X_{i}Y_{i+1}-X_{i+1}Y_{i})^{p-1}\end{equation}
aus 
   $S(V\oplus \langle X_{k+1},Y_{k+1}\rangle)^{\SL_{2}}$ mit
  $i=3,\ldots,k-1$
in dieser Reihenfolge abgebildet auf die Invarianten
\[
X_{1},X_{2}^{p-1},X_{1}Y_{3}-X_{3}Y_{1},(X_{i}Y_{i+1}-X_{i+1}Y_{i})^{p-1}\] aus
  $S(V)^{\Ga}$  mit
  $i=3,\ldots,k-1$. Schreiben wir nun die bei den $\SL_{2}$-Invarianten \eqref{SL2phsopcollection}
  vorkommenden Indizes in der Reihenfolge
\[
2,k+1,1,3,4,5,\ldots,k
\]
auf, so sehen wir, dass f"ur genau zwei benachbarte $i,j$ jeweils eine Potenz
von $X_{i}Y_{j}-X_{j}Y_{i}$ in \eqref{SL2phsopcollection} vorkommt. Nach den
Lemmata \ref{phsopXY} (umnummerieren!) und \ref{phsopPower} ist also \eqref{SL2phsopcollection}
ein phsop im Polynomring $S(V\oplus \langle X_{k+1},Y_{k+1}\rangle)$. Da $\SL_{2}$ aber
reduktiv ist, bildet nun \eqref{SL2phsopcollection} nach Lemma \ref{redphsop}
auch ein phsop im Invariantenring $S(V\oplus \langle X_{k+1},Y_{k+1}\rangle)^{\SL_{2}}$. Die Bilder
von \eqref{SL2phsopcollection} unter Roberts' Isomorphismus bilden dann
nat"urlich ein phsop in $S(V)^{\Ga}$ - zun"achst nur bez"uglich der via
Roberts' Isomorphismus vererbten Graduierung, aber da die Bilder hier auch
bez"uglich der Standardgraduierung homogen sind, dann auch bez"uglich dieser
(vgl. Bemerkung \ref{phsopZweiGradus}). \qed\\

\subsubsection{Ernte}
Nun haben wir alles zusammen, um die Beispiele f"ur $\SL_{2}$ und $\Ga$ zu
formulieren. Sie stellen ein weiteres Hauptresultat dieser Arbeit dar.

\begin{Satz} \label{ZweitesHauptResultat}
Sei $\langle X,Y \rangle$ die nat"urliche Darstellung der $\SL_{2}$, $\chr K=p>0$. Wir
fassen jeden $\SL_{2}$-Modul via $t\mapsto \left(\begin{array}{cc}1 & t\\&1\\
  \end{array}\right)$ auch als $\Ga$-Modul auf. Sei 
\[
V:=\langle X^{p},Y^{p}\rangle \oplus \bigoplus_{i=1}^{k} \langle X,Y \rangle.
\]
Dann gilt
\[
\cmdef K[V]^{\Ga} \ge k-2 \quad \textrm{ und } \quad \cmdef K[V]^{\SL_{2}} \ge k-3.
\]
Als direkte Summe selbstdualer $\Ga$- bzw. $\SL_{2}$- Moduln ist $V$
selbstdual. Aufgefasst als $\SL_{2}$-Modul ist $V$ au"serdem vollst"andig
reduzibel als direkte Summe irreduzibler $\SL_{2}$-Moduln. Ferner gilt
\[
\dim K[V]^{\Ga}=2k+1 \quad \textrm{ und } \quad \dim K[V]^{\SL_{2}}=2k-1.
\]
Damit erhalten wir die Absch"atzungen
\[
\depth K[V]^{\Ga}\le k+3  \quad \textrm{ und } \quad \depth K[V]^{\SL_{2}}\le k+2.
\]
\end{Satz}

\Bew Die Aussagen "uber die Dimension folgen sofort aus den Korollaren
\ref{exakteDimOfGaInvs} und \ref{exakteDimOfSL2Invs}, und mit der Absch"atzung f"ur den
Cohen-Macaulay-Defekt folgt dann die Absch"atzung "uber die Tiefe.

Nach Roberts' Isomorphismus, Satz \ref{RobertsIsom}, gilt
$K[V]^{\Ga}\cong K[V\oplus \langle X,Y \rangle]^{\SL_{2}}$ (insbesondere ist
damit auch der erste Invariantenring endlich erzeugt), so dass es gen"ugt,
die Aussage f"ur die $\Ga$ zu beweisen (vgl. auch \eqref{cmdefOfRobertsIsom},
S. \pageref{cmdefOfRobertsIsom}). Da $V$ selbstdual ist, ist
$K[V]=S(V^{*})\cong S(V)$ (als $\Ga$-Algebren) und
$K[V]^{\Ga}=S(V^{*})^{\Ga}\cong S(V)^{\Ga}$. Nach Korollar \ref{KozyklusInSV}
existiert ein nichttrivialer Kozyklus $g\in Z^{1}(\Ga, S(V))$. Nach Lemma
\ref{AnnullatorPhsopOfGa} hat man ein phsop der L"ange $k$ in $S(V)^{\Ga}$,
welches nach den Lemmata \ref{AnnulatorTyp1},~\ref{AnnullatorTyp2}
und~\ref{AnnullatorTyp3} aus Annullatoren von $g$ besteht. Ferner sind die
ersten beiden phsop Elemente $X_{1},X_{2}^{p-1}$ aus Lemma
\ref{AnnullatorPhsopOfGa} offenbar teilerfremd in $S(V)$. Der Hauptsatz
\ref{BigMainTheorem} in seiner allgemeineren Formulierung (da $\Ga$ nicht
reduktiv ist), liefert nun sofort die Behauptung. 
\qed\\

Wir erinnern nochmals daran, dass wenn man in $V$ den Summanden $\langle
X^{p},Y^{p}\rangle$ durch $\langle X,Y\rangle$ ersetzt, die zugeh"origen
Invariantenringe dann Cohen-Macaulay sind - siehe die Diskussion auf S. \pageref{MotivationsBeispiel}.

\subsection{Additive und unipotente Gruppen}
Ist $V$ ein $\SL_{2}$-Modul, der sich schreiben l"asst als direkte Summe
$V=U\oplus \langle X,Y \rangle$, so ist nach Roberts' Isomorphismus
$K[V]^{\SL_{2}}\cong K[U]^{\Ga}$. Insbesondere folgt aus $\cmdef
K[V]^{\SL_{2}} \ge k-2$ sofort $\cmdef K[U]^{\Ga}\ge k-2$. L"asst sich $V$
nicht auf diese Weise schreiben, wurde aber $\cmdef
K[V]^{\SL_{2}} \ge k-2$ mit Hilfe des Hauptsatzes gefolgert, so gilt auch $\cmdef
K[V\oplus \langle X,Y \rangle ]^{\SL_{2}} \ge k-2$; Denn ein phsop in $K[V]$
 bleibt auch eines nach Einbetten in $K[V\oplus \langle X,Y \rangle]$, und
 genauso bleiben Kozyklen nichttrivial - der Hauptsatz kann also weiterhin
 (mit "`selbem"' phsop und Kozyklus) angewendet werden. Nach Roberts'
 Isomorphismus gilt dann also auch $\cmdef K[V]^{\Ga}=\cmdef
K[V\oplus \langle X,Y \rangle ]^{\SL_{2}} \ge k-2$. Somit lassen sich aus
$\SL_{2}$-Beispielen immer $\Ga$-Beispiele mit gro"sem Cohen-Macaulay-Defekt
konstruieren. Auf diese Weise erhalten wir etwa aus Beispiel \ref{simpleBsp2}

\begin{Bsp}
Sei $\chr K=2$ und $\langle X,Y\rangle$ die nat"urliche Darstellung der
additiven Gruppe~$\Ga$. Dann gilt
\[
\cmdef K\left[\langle X^{2},Y^{2}\rangle\oplus\bigoplus_{i=1}^{k}\langle
  X^{2},Y^{2},XY \rangle \right]^{\Ga}\ge k-2.
\]
Die Dimension des Invariantenringes ist $3k+1$ nach Korollar \ref{exakteDimOfGaInvs}.
\end{Bsp}

Analog erhalten wir aus Beispiel \ref{c3} durch
weglassen des Summanden $\langle X,Y \rangle$

\begin{Bsp}
Ist $\chr K=3$ und $\langle X,Y \rangle$ die nat"urliche Darstellung der
additiven Gruppe~ $\Ga$, so gilt
\[
\cmdef K\left[\langle X^{3},Y^{3} \rangle \oplus
  \bigoplus_{i=1}^{k} \langle X^{2},Y^{2},XY \rangle\right]^{\Ga} \ge k-2.
\]
Die Dimension des Invariantenringes ist $3k+1$ nach Korollar \ref{exakteDimOfGaInvs}.
\end{Bsp}

Sind $G$ und $H$ lineare algebraische Gruppen und ist $f: G\rightarrow H$ ein
surjektiver (algebraischer) Homomorphismus, so ist jeder $H$-Modul $V$
 via $f$ auch ein $G$-Modul. Aufgrund der Surjektivit"at von $f$ gilt dann
 auch $K[V]^{G}=K[V]^{H}$. Aus Beispielen f"ur $H$-Invariantenringe mit
 gro"sem Cohen-Macaulay-Defekt erh"alt man so Beispiele f"ur $G$. 
Das folgende Lemma, das sich etwa in
 \cite[kurz vor Abschnitt 3]{Bryant} findet, stellt einen solchen
 Homomorphismus f"ur nichttriviale zusammenh"angende unipotente Gruppen zur Verf"ugung. 

\begin{Lemma} \label{unipotenToGa}
Sei $G$ eine nichttriviale, zusammenh"angende unipotente lineare algebraische
Gruppe. Dann gibt es einen surjektiven algebraischen Homomorphismus
$G\rightarrow \Ga$.
\end{Lemma}

\Bew Als unipotente Gruppe ist $G$ nilpotent (\cite[Corollary
17.5]{Humphreys}), also aufl"osbar. Nach \cite[Theorem 19.3]{Humphreys}
enth"alt $G$ einen abgeschlossenen Normalteiler $N$ mit $\dim G-\dim
N=1$. Nach \cite[Theorem 11.5]{Humphreys} gibt es eine rationale Darstellung $\varphi:
G\rightarrow \GL(V)$ mit $\ker \varphi=N$. Nach \cite[Proposition 2.2.5 (ii)]{SpringerLin} ist
$\varphi(G)\subseteq \GL(V)$ abgeschlossen, au"serdem mit $G$ zusammenh"angend und
unipotent (\cite[Theorem 2.4.8]{SpringerLin}). Weiter gilt nach
\cite[Corollary 4.3.4]{SpringerLin} $\dim \varphi(G)=\dim G - \dim \ker \varphi = 1$. Also
ist $\varphi(G)$ eine zusammenh"angende, unipotente, eindimensionale lineare
algebraische Gruppe,
und damit isomorph zu $\Ga$ (\cite[Proposition 2.6.6]{SpringerLin}). \qed\\

\begin{Satz}
Sei $G$ eine nichttriviale, zusammenh"angende unipotente lineare algebraische
Gruppe in positiver Charakteristik. Dann gibt es f"ur jedes $k\in{\mathbb N}$ einen $2k+2$-dimensionalen
$G$-Modul $V$ mit
\[
\cmdef K[V]^{G}\ge k-2.
\] 
\end{Satz}

\Bew Wir machen den von Satz
\ref{ZweitesHauptResultat} gelieferten $2k+2$-dimensionalen $\Ga$-Modul $V$ mit $\cmdef K[V]^{\Ga}\ge k-2$ mittels des
surjektiven Homomorphismus $G\rightarrow \Ga$ aus Lemma \ref{unipotenToGa} zu
einem $G$-Modul. Aus $K[V]^{G}=K[V]^{\Ga}$ folgt sofort die Behauptung. \qed\\

Man kann ein entsprechendes Resultat auch f"ur nichtzusammenh"angende
unipotente Gruppen
angeben. Dann hat n"amlich jedes Element von $G/G^{0}$ $p$-Potenzordnung
(ist $G\subseteq \GL_{n}$ unipotent, so gibt es zu $A\in G$ ein $N\in{\mathbb
  N}$ mit $(A-I_{n})^{N}=0$. F"ur $k$ mit $p^{k}\ge N$ ist dann
$0=(A-I_{n})^{p^{k}}=A^{p^{k}}-I_{n}$, vgl. \cite[2.4]{SpringerLin}), so dass
$H:=G/G^{0}$ eine (endliche) $p$-Gruppe ist. Man muss
hier also auf entsprechende Resultate im modularen Fall, etwa Satz \ref{FiniteCMdef} zur"uckgreifen.\\

Man beachte, dass unendliche unipotente Gruppen $G$ nicht
reduktiv sind, denn dann ist $G^{0}\ne\{\iota\}$ ein zusammenh"angender,
abgeschlossener, nichttrivialer
unipotenter Normalteiler von $G$ (vgl. Satz bzw. Definition \ref{zushgkomp} und \ref{defred}).

\newpage
\section{Algorithmische Untersuchungen}\label{Algorithmik}
In diesem Abschnitt wollen wir im Wesentlichen die in Satz
\ref{ZweitesHauptResultat} beschriebenen Invariantenringe mit "`a posteriori"'
Methoden untersuchen, d.h. wir berechnen zun"achst Erzeuger f"ur $K[V]^{G}$
(als $K$-Algebra). Hierf"ur entwickeln wir ein speziell angepasstes Verfahren,
da der allgemeine Algorithmus \cite{KemperCompRed} auch in kleinen F"allen
nicht durchkommt (ich habe die Rechnung bei knapp 70 \emph{Giga}byte
Speicherbedarf abgebrochen). Danach berechnen wir das Relationenideal und rechnen so im
Restklassenring eines Polynomrings weiter, und berechnen so explizit den
Cohen-Macaulay-Defekt f"ur die F"alle $(p,k)\in\{(2,3),(2,4),(3,3)\}$ zu
$1,2,1$ (f"ur die $\Ga$-Invarianten). Genau genommen berechnen wir jeweils
eine  obere Schranke f"ur den
Cohen-Macaulay-Defekt, die zusammen mit der unteren Schranke aus Satz
\ref{ZweitesHauptResultat} das angegebene Ergebnis liefert.

\subsection{Berechnung von Frobenius-Invarianten}
Sei $G$ eine lineare algebraische Gruppe in Charakteristik $p>0$, und $U,V$
zwei $G$-Moduln, so dass $S(U\oplus V)^{G}$ endlich erzeugt ist. Mit $F^{p}$ bezeichnen wir die $p$-te Frobenius-Potenz,
vgl. Definition \ref{FpV}. Wir
untersuchen folgendes Problem:

{\it Wie kann man Generatoren f"ur $S(F^{p}(U)\oplus V)^{G}$ effizient
  berechnen, wenn Generatoren f"ur $S(U\oplus V)^{G}$ bekannt sind?}\\

Ohne formale Definition nennen wir dieses Problem im
Folgenden die "`Berechnung von Frobenius-Invarianten"'.

\subsubsection{Der Isomorphismus}\label{DerIsomorphismus}
Seien 
\[
U=\langle X_{1},\ldots,X_{n} \rangle,\quad V=\langle Y_{1},\ldots,Y_{m}\rangle,\quad
\textrm{und }\quad
F^{p}(U)=\langle X_{1}^{p},\ldots,X_{n}^{p} \rangle=:\langle
Z_{1},\ldots,Z_{n} \rangle\]
$G$-Moduln. Dann ist 
\[
P:=S(U\oplus V)=K[X_{1},\ldots,X_{n},Y_{1},\ldots,Y_{m}],\]
\[
Q:=S(F^{p}(U)\oplus V)=K[Z_{1},\ldots,Z_{n},Y_{1},\ldots,Y_{m}].
\]
Wir betrachten den Algebrenhomomorphismus 
\[
\phi: Q\rightarrow P \quad \textrm{ mit } Z_{i}\mapsto X_{i}^{p},\,
Y_{i}\mapsto Y_{i},
\]
welcher als Abbildung des Polynomrings $Q$ eindeutig und wohldefiniert ist durch
Angabe der Bilder der unabh"angigen Variablen.

\begin{Satz}\label{IsomOfFrobInvs}
Durch die Einschr"ankung $\varphi:=\phi|_{Q^{G}}$ ist ein Isomorphismus
\[
\varphi: Q^{G}\rightarrow P^{G}\cap
K[X_{1}^{p},\ldots,X_{n}^{p},Y_{1},\ldots,Y_{m}]
\]
gegeben.
\end{Satz}

\Bew Wir bezeichnen mit $G\rightarrow \GL_{n}, \sigma \mapsto
A_{\sigma}=(a_{ij,\sigma})$ die Darstellung von $U$ und analog mit $G\rightarrow \GL_{m}, \sigma \mapsto
B_{\sigma}=(b_{ij,\sigma})$ die Darstellung von $V$. Die Darstellung von
$F^{p}(U)$ ist damit gegeben durch $\sigma\mapsto (a_{ij,\sigma}^{p})$. Es
gilt also
\[
\sigma \cdot X_{j}=\sum_{i=1}^{n}a_{ij,\sigma}X_{i},\quad \sigma \cdot
Z_{j}=\sum_{i=1}^{n}a_{ij,\sigma}^{p}Z_{i}, \quad \sigma \cdot
Y_{j}=\sum_{i=1}^{m}b_{ij,\sigma}Y_{i}.
\]

Es ist klar, dass $\varphi: Q^{G}\rightarrow P$ mit $\phi$ ein
Algebren-Homomorphismus ist. Da $X_{1}^{p},\ldots,X_{n}^{p}$, $Y_{1},\ldots,Y_{m}$
offenbar algebraisch unabh"angig sind, ist $\phi$ und damit auch die
Einschr"ankung $\varphi$ injektiv. Wir m"ussen nun nur noch $\varphi(Q^{G})=P^{G}\cap
K[X_{1}^{p},\ldots,X_{n}^{p},Y_{1},\ldots,Y_{m}]$ zeigen.

"`$\subseteq$"'. Die Inklusion $\varphi(Q^{G})\subseteq
K[X_{1}^{p},\ldots,X_{n}^{p},Y_{1},\ldots,Y_{m}]$ ist klar. Sei nun
\[f=f(Z_{1},\ldots,Z_{n},Y_{1},\ldots,Y_{m})\in Q^{G},\]
d.h.
\begin{equation} \label{fisInvariant}
\sigma \cdot
f\stackrel{(*)}{=}f\left(\sum_{i=1}^{n}a_{i1,\sigma}^{p}Z_{i},\ldots,\sum_{i=1}^{n}a_{in,\sigma}^{p}Z_{i},\sum_{i=1}^{m}b_{i1,\sigma}Y_{i},\ldots,\sum_{i=1}^{m}b_{im,\sigma}Y_{i}
\right)=f.
\end{equation}
Dann ist
\[\varphi(f)=f\left(X_{1}^{p},\ldots,X_{n}^{p},Y_{1},\ldots,Y_{m} \right),\]
und damit f"ur $\sigma \in G$
\begin{eqnarray*}
\sigma \cdot \varphi(f)&=&f\left(\left(\sum_{i=1}^{n}a_{i1,\sigma}X_{i}\right)^{p},\ldots,\left(\sum_{i=1}^{n}a_{in,\sigma}X_{i}\right)^{p},\sum_{i=1}^{m}b_{i1,\sigma}Y_{i},\ldots,\sum_{i=1}^{m}b_{im,\sigma}Y_{i}
\right)\\
&=&f\left(\sum_{i=1}^{n}a_{i1,\sigma}^{p}X_{i}^{p},\ldots,\sum_{i=1}^{n}a_{in,\sigma}^{p}X_{i}^{p},\sum_{i=1}^{m}b_{i1,\sigma}Y_{i},\ldots,\sum_{i=1}^{m}b_{im,\sigma}Y_{i}
\right)\\
&\stackrel{(*)}{=}&\varphi(\sigma \cdot f)\\
&\stackrel{~\eqref{fisInvariant}}{=}&\varphi(f),
\end{eqnarray*}
also $\varphi(f)\in P^{G}$. (Man kann auch argumentieren, dass $\varphi$
$G$-"aquivariant ist).

"`$\supseteq$"'. Sei nun umgekehrt
\[F=F(X_{1},\ldots,X_{n},Y_{1},\ldots,Y_{m})\in P^{G}\cap
K[X_{1}^{p},\ldots,X_{n}^{p},Y_{1},\ldots,Y_{m}].\]
Dann gibt es jedenfalls ein $f=f(Z_{1},\ldots,Z_{n},Y_{1},\ldots,Y_{m})\in Q$
mit 
\[
F=f(X_{1}^{p},\ldots,X_{n}^{p},Y_{1},\ldots,Y_{m})=\phi(f).\]
Wir m"ussen daher nur noch $f\in Q^{G}$ zeigen. F"ur $\sigma\in G$ ist
\begin{eqnarray*}
\sigma \cdot F&=&\sigma\cdot
f(X_{1}^{p},\ldots,X_{n}^{p},Y_{1},\ldots,Y_{m})\\
&=&f\left(\sum_{i=1}^{n}a_{i1,\sigma}^{p}X_{i}^{p},\ldots,\sum_{i=1}^{n}a_{in,\sigma}^{p}X_{i}^{p},\sum_{i=1}^{m}b_{i1,\sigma}Y_{i},\ldots,\sum_{i=1}^{m}b_{im,\sigma}Y_{i}\right)\\
&\stackrel{(*)}{=}&f(X_{1}^{p},\ldots,X_{n}^{p},Y_{1},\ldots,Y_{m})= F,
\end{eqnarray*}
da $F\in P^{G}$. Die Gleichung $(*)$ beschreibt eine Gleichung in
$K[X_{1}^{p},\ldots,X_{n}^{p},Y_{1},\ldots,Y_{m}]$. Da dieser Ring isomorph zu einem
Polynomring ist, darf man in $(*)$ die "`unabh"angige Variable"'  $X_{i}^{p}$ durch $Z_{i}$
ersetzen. Dies liefert
\begin{eqnarray*}
\sigma \cdot f&=&\sigma \cdot f(Z_{1},\ldots,Z_{n},Y_{1},\ldots,Y_{m})\\
&=&f\left(\sum_{i=1}^{n}a_{i1,\sigma}^{p}Z_{i},\ldots,\sum_{i=1}^{n}a_{in,\sigma}^{p}Z_{i},\sum_{i=1}^{m}b_{i1,\sigma}Y_{i},\ldots,\sum_{i=1}^{m}b_{im,\sigma}Y_{i}\right)\\
&\stackrel{(*)}{=}&f(Z_{1},\ldots,Z_{n},Y_{1},\ldots,Y_{m})=f,
\end{eqnarray*}
also $f\in Q^{G}$.\qed\\

\subsubsection{Zwischenspiel: Der Kern einer linearen Abbildung von Moduln}
Sei $K[X]=K[X_{1},\ldots,X_{n}]$ ein Polynomring und $A=K[f_{1},\ldots,f_{k}]$ mit $f_{i}\in
K[X] \,\,\myforall i=1,\ldots,k$ eine endlich erzeugte Unteralgebra. Sei ferner
$B=\sum_{i=1}^{r}At_{i}$ mit $t_{i}\in K[X]\,\,\myforall i=1,\ldots,r$ ein
endlich erzeugter $A$-Modul sowie $D: B\rightarrow K[X]^{m}$ eine $A$-lineare
Abbildung. Der folgende Algorithmus (der eine
Verallgemeinerung von Kemper \cite[Algorithm 4.2]{KemperCompRed} ist), liefert ein endliches Erzeugendensystem von $\ker D$
als $A$-Modul.

\begin{Algorithmus}\label{ComputeKernel}
Berechnung des Kerns einer linearen Abbildung von Moduln.\\
{\bf Eingabe:}
\begin{itemize}
\item Ein Polynomring $K[X]=K[X_{1},\ldots,X_{n}]$.
\item Generatoren $f_{1},\ldots,f_{k}\in K[X]$ von $A:=K[f_{1},\ldots,f_{k}]$.
\item Modulerzeuger $t_{1},\ldots,t_{r}\in K[X]$ von $B:=\sum_{i=1}^{r}At_{i}$.
\item Die Bilder $D(t_{1}),\ldots,D(t_{r})\in K[X]^{m}$ 
einer $A$-linearen Abbildung $D: B\rightarrow K[X]^{m}$. (Der Anwender
  hat sicherzustellen, dass die Bilder tats"achlich von einer $A$-linearen
  Abbildung $D$ abstammen!)
\end{itemize}
{\bf Ausgabe:} $A$-Modul Erzeuger von $\ker D$.\\
{\bf BEGIN}
\begin{enumerate}
\item Berechne Generatoren des Syzygien-Moduls
\[
M:=\left\{\left(a_{1},\ldots,a_{r}\right)\in K[X]^{r}:
a_{1}D(t_{1})
+\ldots+
a_{r}D(t_{r})=0
\right\}
\]
(als $K[X]$-Modul) mit einem der "ublichen Standard-Verfahren, siehe etwa Eisenbud
\cite[Chapter 15.5]{Eisenbud}.
\item Berechne Generatoren $b_{1},\ldots,b_{s}$ von $M\cap A^{r}$ als $A$-Modul mittels Kemper
  \cite[Algorithm 4.5]{KemperCompRed} (siehe auch Kemper \cite[Algorithm
  7]{KemperFiniteGroups}). Dieser Algorithmus verlangt als Eingabe die
  Generatoren von $A$ und die in Schritt 1. berechneten Generatoren von $M$
  als $K[X]$-Modul.
\item Setze 
\[
c_{i}:=\sum_{\mu=1}^{r}(b_{i})_{\mu}t_{\mu} \quad \textrm{ f"ur } i=1,\ldots,s
\]
mit $b_{i}=((b_{i})_{1},\ldots,(b_{i})_{r})\in A^{r}$.
\item {\bf RETURN} $c_{1},\ldots,c_{s}$. 
\end{enumerate}
{\bf END}
\end{Algorithmus}

\begin{Satz}\label{CompKernerlIsKorrekt}
Algorithmus \ref{ComputeKernel} ist korrekt, d.h.
\[
\ker D=\sum_{i=1}^{s}Ac_{i}.
\]
\end{Satz}

\Bew F"ur $f\in K[X]$ gilt $f\in \ker D$ genau dann, wenn es $a_{1},\ldots,a_{r}\in
A$ gibt mit
\begin{equation}\label{fIsInKerD}
f=\sum_{\mu=1}^{r}a_{\mu}t_{\mu} \quad \textrm{und} \quad D(f)=\sum_{\mu=1}^{r}a_{\mu}D(t_{\mu})=0,
\end{equation}
denn die $t_{\mu}$ erzeugen $B$ als $A$-Modul, und $D$ ist $A$-linear.
Wir zeigen nun beide Inklusionen der Behauptung.

"`$\supseteq$"'.  Da $b_{i}=((b_{i})_{1},\ldots,(b_{i})_{r})\in M\cap A^{r}$ (Schritt~2.), gilt nach
Definition von $M$ (Schritt~1.) jedenfalls $\sum_{\mu=1}^{r} (b_{i})_{\mu}D(t_{\mu})=0$, und damit (Schritt~3.)
$c_{i}=\sum_{\mu=1}^{r}(b_{i})_{\mu}t_{\mu}\in \ker D$ nach  \eqref{fIsInKerD}.

"`$\subseteq$"'. Sei $f\in \ker D$, d.h. $f$ habe eine Darstellung
wie in \eqref{fIsInKerD}. Dann gilt
$a=(a_{1},\ldots,a_{r})\in M\cap A^{r}$ nach Definition von $M$. Nach Schritt~2. gibt es dann
$p_{1},\ldots,p_{s}\in A$ mit 
\[
a=\sum_{i=1}^{s}p_{i}b_{i}, \,\,\textrm{ also } a_{\mu}=\sum_{i=1}^{s}p_{i}(b_{i})_{\mu}.\]
Es folgt
\[
f=\sum_{\mu=1}^{r}a_{\mu}t_{\mu}=\sum_{\mu=1}^{r}\sum_{i=1}^{s}p_{i}(b_{i})_{\mu}t_{\mu}=\sum_{i=1}^{s}p_{i}\underbrace{\sum_{\mu=1}^{r}(b_{i})_{\mu}t_{\mu}}_{c_{i}},
\]
und damit gilt  $f\in\sum_{i=1}^{s}Ac_{i}$. \qed\\

\subsubsection{Schnitt einer Algebra mit $K[X^{p},Y]$}
Wir verwenden wieder die Bezeichnungen aus Abschnitt
\ref{DerIsomorphismus}. Au"serdem k"urzen wir mit $X^{p}$ die Variablen
$X_{1}^{p},\ldots,X_{n}^{p}$ ab. 
Von dem Isomorphismus $\varphi$ kann man auch leicht die
Umkehrung durchf"uhren, indem man in einem Polynom $f$ aus dem Bild einfach
in jedem Monom $X_{i}^{p}$ durch $Z_{i}$ ersetzt. Da ein Algebrenisomorphismus
Generatoren auf Generatoren abbildet, k"onnen wir zur Berechnung von $Q^{G}$
(was unser Ziel ist) zun"achst
 Generatoren f"ur das Bild von
$\varphi$ berechnen, um durch Anwendung von $\varphi^{-1}$ Generatoren von
$Q^{G}$ zu erhalten. Wir versuchen dazu, $\varphi(Q^{G})=P^{G}\cap K[X^{p},Y]$ als Kern einer
geeigneten $A$-linearen Abbildung zu erhalten, den wir mit Algorithmus
\ref{ComputeKernel} berechnen k"onnen.\\

Sei also allgemeiner $B\subseteq K[X,Y]$ eine endlich erzeugte Algebra (z.B. $B=P^{G}$),
\[
B:=K[f_{1},\ldots,f_{k},g_{1},\ldots,g_{l}] \quad \textrm{ mit } f_{i}\in K[X,Y],\,\, g_{i}\in
K[Y] \,\,\myforall i.
\]
Gesucht ist der Schnitt $B\cap K[X^{p},Y]$. Wir setzen dazu
\[
A:=K[f_{1}^{p},\ldots,f_{k}^{p},g_{1},\ldots,g_{l}]\subseteq B\cap K[X^{p},Y]
\]
und bilden 
\[
\{ t_{1},\ldots,t_{r}\} :=\left\{ f_{1}^{e_{1}}\cdot\ldots\cdot f_{k}^{e_{k}}:
  0\le e_{i}<p \right\}.
\]
Dann ist offenbar $B=\sum_{i=1}^{r}At_{i}$. F"ur $j=1,\ldots,n$ ist die
Abbildung
\[
\frac{\partial}{\partial X_{j}}: B \rightarrow K[X,Y]
\]
$A$-linear, denn f"ur $a\in A, b\in B$ gilt wegen $A\subseteq K[X^{p},Y]$
\[
\frac{\partial (ab)}{\partial X_{j}}=a\frac{\partial b}{\partial
  X_{j}}+b\underbrace{\frac{\partial a}{\partial X_{j}}}_{=0}. 
\]
Weiter liegt ein $f\in B$ genau dann in $K[X^{p},Y]$, wenn $\frac{\partial
  f}{\partial X_{j}}=0$ f"ur $j=1,\ldots,n$ gilt, also wenn $f$ im Kern der
  linearen Abbildung 
\begin{equation}\label{Derivationdxj}
D:B\rightarrow K[X,Y]^{n}, f \mapsto D(f):=\left(\frac{\partial
  f}{\partial X_{1}},\ldots,\frac{\partial
  f}{\partial X_{n}}\right)
\end{equation}
liegt. Damit erhalten wir folgenden Algorithmus, um Algebra-Erzeuger f"ur
$B\cap K[X^{p},Y]=\ker D$ zu erhalten:

\begin{Algorithmus}\label{AlgIntersectXpY}
Schnitt einer Algebra $B$ mit $K[X^{p},Y]$, wobei $p=\chr K>0$.\\
{\bf Eingabe:} \begin{itemize}\item Ein Polynomring
$K[X_{1},\ldots,X_{n},Y_{1},\ldots,Y_{m}]=:K[X,Y]$.

\item $B:=K[f_{1},\ldots,f_{k},g_{1},\ldots,g_{l}]$ mit $f_{i}\in K[X,Y]$, $g_{i}\in
K[Y]\,\,\myforall i$.
\end{itemize}
{\bf Ausgabe:} Generatoren von $B\cap
K[X_{1}^{p},\ldots,X_{n}^{p},Y_{1},\ldots,Y_{m}]=:B\cap K[X^{p},Y]$.\\
{\bf BEGIN}
\begin{enumerate}
\item Setze  $A:=K[f_{1}^{p},\ldots,f_{k}^{p},g_{1},\ldots,g_{l}]$ und bilde
\[
T:=\{ t_{1},\ldots,t_{r}\} :=\left\{ f_{1}^{e_{1}}\cdot\ldots\cdot f_{k}^{e_{k}}:
  0\le e_{i}<p \right\}.
\]
\item "Ubergib die Daten $A,T$ und $D(t_{1}),\ldots,D(t_{r})$ mit $D$ wie in
  \eqref{Derivationdxj} an Algorithmus~\ref{ComputeKernel}. Erhalte Generatoren
  $c_{1},\ldots,c_{s}$ von $\ker D=B\cap K[X^{p},Y]$ als $A$-Modul.
\item {\bf RETURN} $f_{1}^{p},\ldots,f_{k}^{p},g_{1},\ldots,g_{l},c_{1},\ldots,c_{s}$. 
\end{enumerate}
{\bf END}
\end{Algorithmus}

\begin{Satz}
Der Algorithmus ist korrekt, d.h.
\[
B\cap
K[X^{p},Y]=K[f_{1}^{p},\ldots,f_{k}^{p},g_{1},\ldots,g_{l},c_{1},\ldots,c_{s}].\]
\end{Satz}

\Bew Da $B\cap
K[X^{p},Y]=\ker D=\sum_{i=1}^{s}Ac_{i}$ nach Satz \ref{CompKernerlIsKorrekt}, gilt erst recht $B\cap
K[X^{p},Y]=A[c_{1},\ldots,c_{s}]=K[f_{1}^{p},\ldots,f_{k}^{p},g_{1},\ldots,g_{l},c_{1},\ldots,c_{s}]$.\qed\\

\subsubsection{Der Algorithmus}
Wir erhalten nun den folgenden Algorithmus zur Berechnung von
$S(F^{p}(U)\oplus V)^{G}$:

\begin{Algorithmus}\label{ComputeFrobeniusInvs} Berechnen von $S(F^{p}(U)\oplus V)^{G}$, falls $S(U\oplus
  V)^{G}$ bekannt ist.\\
{\bf Eingabe:} Generatoren von $S(U\oplus V)^{G}$.\\
{\bf Ausgabe:} Generatoren von $S(F^{p}(U)\oplus V)^{G}$.\\
\begin{samepage}
{\bf BEGIN}
\begin{enumerate}
\item Schreibe $S(U\oplus V)=K[X,Y]$, $S(F^{p}(U)\oplus V)=K[Z,Y]$.
\item Sei  $S(U\oplus V)^{G}=K[f_{1},\ldots,f_{k},g_{1},\ldots,g_{l}]$ mit $f_{i}\in K[X,Y]$, $g_{i}\in
K[Y]$.
\item Berechne Generatoren $H_{1},\ldots,H_{s}$ von $S(U\oplus V)^{G}\cap
  K[X^{p},Y]$ mit Hilfe von Algorithmus~\ref{AlgIntersectXpY}.
\item Ersetze f"ur $i=1,\ldots,s$ in $H_{i}$ jedes $X_{j}^{p}$ durch $Z_{j}$
  und schreibe f"ur das Ergebnis $h_{i}$.
\item {\bf RETURN} $h_{1},\ldots,h_{s}$. 
\end{enumerate}
{\bf END}
\end{samepage}
\end{Algorithmus}

Die Korrektheit von Algorithmus~\ref{ComputeFrobeniusInvs}, also $S(U\oplus
  V)^{G}=K[h_{1},\ldots,h_{s}]$, folgt sofort aus Satz \ref{IsomOfFrobInvs}.

\subsubsection{Beispiele}
Algorithmus \ref{ComputeFrobeniusInvs} l"asst sich einfach in \Magma
implementieren. Wir wenden ihn an f"ur $G=\SL_{2}$ und $U=\langle
X_{0},Y_{0}\rangle$ (nat"urliche Darstellung) und
$V=\bigoplus_{i=1}^{k}\langle X_{i},Y_{i}\rangle$ ($k$-fache direkte Summe der
nat"urlichen Darstellung). Nach de Concini und Procesi \cite{ConciniProcesi} wird $S(U\oplus V)^{\SL_{2}}$ auch
in positiver Charakteristik erzeugt von
\[
\{X_{i}Y_{j}-X_{j}Y_{i}: 0\le i<j\le k \}.
\]
Damit erh"alt man mittels des Frobenius-Homomorphismus auch leicht Erzeuger
f"ur $S(U\oplus V)^{\Ga}$, indem man in obiger Menge zun"achst $k+1$ statt $k$
w"ahlt und dann die Ersetzung $X_{k+1}=0,Y_{k+1}=1$ vornimmt. Damit wird
$S(U\oplus V)^{\Ga}$ erzeugt von
\[
\{X_{i}Y_{j}-X_{j}Y_{i}: 0\le i<j\le k \} \cup \{X_{i}: 0\le i \le k\}.
\]
Auch wenn man an $S(F^{p}(U)\oplus V)^{\SL_{2}}$ interessiert ist, empfiehlt sich
dringend, stattdessen $S(F^{p}(U)\oplus V)^{\Ga}$ (mit $k-1$ statt k) zu
berechnen und dann darauf die Umkehrung von Roberts' Isomorphismus (Korollar
\ref{UmkehrungRoberts}) anzuwenden. So dauert etwa die Rechnung f"ur das Tupel
$(G=\SL_{2},p=3,k=4)$ etwa $250s$, dagegen die "aquivalente Rechnung f"ur
$(G=\Ga,p=3,k=3)$ nur $73s$.

Wir haben den Schnitt
\[S(U\oplus
V)^{\Ga}\cap K[X_{0}^{p},Y_{0}^{p},X_{1},Y_{1},\ldots,X_{k},Y_{k}]\]
f"ur einige weitere Werte von $p$ und $k$ mit Hilfe einer Implementierung von
Algorithmus \ref{ComputeFrobeniusInvs} in \Magma berechnet. Die folgenden Tabelle gibt die Entwicklung der Rechenzeit (t) sowie der Anzahl (\#)
der Generatoren von $S(F^{p}(U)\oplus V)^{\Ga}$ f"ur diese Werte wieder. 
Die L"ucken in der Tabelle resultieren daher, dass ich die Rechnung f"ur die
entsprechenden Werte nach etwa einem Monat ohne Ergebnis abgebrochen habe.

\begin{center}
\begin{tabular}{c|cccccc}
p$\backslash$k&2&3&4&5\\
\hline
2&t=0.1s, \#6&t=1.750s, \#11&t=53.56s, \#20&t=6.25h, \#37\\
3&t=1.1s, \#6&t=73.16s, \#14&t=20.4h, \#46\\
5&t=50.41s, \#6&t=32.15h, \#30\\
\end{tabular}\\

{\small Rechenzeiten (t) von Algorithmus \ref{ComputeFrobeniusInvs} und Anzahl (\#)
der berechneten Generatoren.}
\end{center}

Der Speicherbedarf in den aufw"andigeren F"allen betrug bis zu 23GB. Zur Abrundung wollen wir f"ur den kleinsten interessanten Fall $p=2, k=3$ die
berechneten Generatoren des Schnitts explizit angeben:

\begin{Satz}[\Magma und Algorithmus \ref{ComputeFrobeniusInvs}]
Sei $p=\chr K=2$ sowie $\langle
X_{i},Y_{i}\rangle$ f"ur $i=0,\ldots,3$ jeweils die nat"urliche Darstellung
der additiven Gruppe $\Ga$. Es wird 
\begin{eqnarray*}
&&S\left(\langle
X_{0},Y_{0}\rangle \oplus \bigoplus_{i=1}^{3}\langle
X_{i},Y_{i}\rangle\right)^{\Ga}\cap
K[X_{0}^{2},Y_{0}^{2},X_{1},Y_{1},X_{2},Y_{2},X_{3},Y_{3}]\\
&\cong&S\left(\langle
X_{0}^{2},Y_{0}^{2}\rangle \oplus \bigoplus_{i=1}^{3}\langle
X_{i},Y_{i}\rangle\right)^{\Ga}
\end{eqnarray*} erzeugt von den
$11$ Generatoren
\[
    X_{0}^2,\quad
    X_{1},\quad
    X_{2},\quad
    X_{3},
\]
\[ 
   X_{1}Y_{2} + Y_{1}X_{2},\quad
    X_{1}Y_{3} + Y_{1}X_{3},\quad
    X_{2}Y_{3} + Y_{2}X_{3},
\]
\[
    X_{1}^2 Y_{0}^2 + Y_{1}^2 X_{0}^2,\quad
    X_{2}^2 Y_{0}^2 + Y_{2}^2 X_{0}^2,\quad
    X_{3}^2 Y_{0}^2 + Y_{3}^2 X_{0}^2,
\]
\[
    X_{1} X_{2} X_{3} Y_{0}^2 + X_{1} Y_{2} Y_{3} X_{0}^2 + Y_{1} X_{2} Y_{3} X_{0}^2 + Y_{1} Y_{2} X_{3} X_{0}^2.
\]
\end{Satz}

\subsection{Untere Schranken f"ur die Tiefe}
W"ahrend der Hauptsatz \ref{BigMainTheorem} "`a priori"' eine obere Schranke
f"ur die Tiefe gibt, geben wir hier eine einfache Heuristik zur "`a posteriori"' Bestimmung
unterer Schranken f"ur die Tiefe. Wir gehen dabei von der Situation aus, dass
der Ring $R$ (bis auf Isomorphie) gegeben ist als Restklassenring eines Polynomrings
$P=K[X_{1},\ldots,X_{n}]$ modulo eines Ideals $I\unlhd P$, welches durch ein
endliches Erzeugendensystem gegeben ist, also $R\cong P/I$. Im Fall
$R=K[V]^{G}$ m"ussen dazu zun"achst Generatoren $f_{1},\ldots,f_{n}\in K[V]$
von $K[V]^{G}$ berechnet werden, und dann $I$ als Kern des
Einsetzungshomomorphismus $P\rightarrow K[V]^{G}, X_{i}\mapsto f_{i}$ (siehe
etwa \cite[Proposition 15.30]{Eisenbud}). Beide Schritte sind rechenintensiv
und k"onnen dem Vorhaben bereits ein Ende setzen.

Um nun zu pr"ufen, ob die Restklassen gegebener $g_{1},\ldots,g_{m}\in P$ eine
regul"are Sequenz in $R$ bilden, muss man testen, ob jeweils $g_{i}$ ein
Nichtnullteiler in $R/(g_{1},\ldots,g_{i-1})_{R}\cong
P/I+(g_{1},\ldots,g_{i-1})_{P}$ ist (Details hierzu auf
S. \pageref{MagmaTestAufRegSeq}). (Falls die Restklassen der $g_{i}$ in $R$
nicht homogen und positiven Grades
sind bzw. $R$ nicht graduiert ist, muss man noch $(g_{1},\ldots,g_{m})\ne R$ testen, damit man wirklich eine
regul"are Sequenz hat.)

Unsere Heuristik besteht nun einfach daraus, dass der Anwender eine Folge von
Elementen $g_{1},\ldots,g_{m}\in P$ vorgibt. Falls im $i$-ten Schritt dann
$g_{i}$ ein Nichtnullteiler modulo $I$ ist, so merken wir uns $i$ und setzen
$I:=I+(g_{i})$. Die gemerkten $g_{i}$ geben dann einer regul"are Sequenz und
ihre Anzahl eine untere Schranke f"ur die Tiefe. Formal aufgeschrieben:

\begin{samepage}
\begin{Algorithmus}\label{ScanReg} (Heuristik zum finden regul"arer Sequenzen)\\
{\bf Eingabe:}\begin{itemize}
\item Polynomring $P=K[X_{1},\ldots,X_{n}]$.
\item Ein Ideal $I\unlhd P$, gegeben durch endlich viele Generatoren.
\item Eine Testsequenz $g_{1},\ldots,g_{m}\in P$.
\end{itemize}
{\bf Ausgabe}: Eine in $R:=P/I$ regul"are Teilfolge der Testsequenz
$g_{i_{1}},\ldots,g_{i_{k}}$, im Fall eines graduierten Rings $R$ und alle
$g_{i}+I\in R_{+}$ insbesondere also eine untere Schranke f"ur die
Tiefe,
\[
k\le \depth(R).
\]
{\bf BEGIN}
\begin{enumerate}
\item Setze $k:=0, J:=I$.
\item F"ur $i=1,\ldots, m$
\begin{itemize}
\item FALLS $g_{i}$ ein Nichtnullteiler modulo $J$ ist, so setze 
\[
k:=k+1,\quad i_{k}:=i,\quad J:=J+(g_{i})_{P}.
\]
\end{itemize}
\item Falls die Restklassen der $g_{1},\ldots,g_{m}$ modulo $I$ nicht im
  maximalen homogenen Ideal $R_{+}$ eines
  graduierten Rings $R=P/I$ lagen, so erniedrige evtl. $k$, damit
  $(g_{i_{1}},\ldots,g_{i_{k}})+I\ne P$.
\item {\bf RETURN} Regul"are Sequenz $g_{i_{1}},\ldots,g_{i_{k}}$, untere
  Schranke f"ur die Tiefe $k$.
\end{enumerate}
{\bf END}
\end{Algorithmus}
\end{samepage}

\subsection{Ein hsop f"ur $S\left(\bigoplus_{i=1}^{n}\langle X,Y \rangle\right)^{\SL_{2}}$}
F"ur die Anwendung von Algorithmus \ref{ScanReg} ist die Angabe einer
geeigneten Testsequenz essentiell. Zum einen sollte man von ihr erwarten
k"onnen, dass sie eine m"oglichst lange regul"are Sequenz enth"alt, zum
anderen m"ussen ihre Elemente einfach genug sein, dass der Test auf  Nullteiler in
Schritt 2. noch schnell genug durchgef"uhrt werden kann. Hieran scheiterten
(bei meinen Versuchen) alle Testsequenzen, die durch Erg"anzen des phsop aus
Lemma \ref{AnnullatorPhsopOfGa} entstanden sind. Zum Erfolg f"uhrte dagegen
eine Testsequenz, die aus den $f_{k}$ (mit geeigneten Exponenten $e_{ij}$)
des folgenden Satzes bestand.

\begin{Satz}\label{zeroDimIdealsVonSL2}
Sei $n\ge 2$ und $R:=S\left(\bigoplus_{i=1}^{n}\langle X_{i},Y_{i} \rangle\right)^{\SL_{2}}$. Dann
 gilt $\dim R=2n-3$. Sei 
$g_{ij}:=X_{i}Y_{j}-X_{j}Y_{i}\in R$, und f"ur $k=3,\ldots,2n-1$ sei
\[
f_{k}:=\sum_{{i+j=k \atop i<j}}g_{ij}^{e_{ij}}\quad \textrm{mit }e_{ij}\ge 1.
\]
Dann ist $\dim (f_{3},\ldots,f_{2n-1})_{R}=0$. Insbesondere bilden die $2n-3$
Elemente $f_{3},\ldots,f_{2n-1}$ bei geeigneter Wahl der $e_{ij}$ und einer
Graduierung von $R$ (z.B. Standardgraduierung und alle $e_{ij}=1$) ein hsop
von $R$. 
\end{Satz}

Im Beweis verwenden wir
\begin{Lemma}\label{x1x2xnEqZero}
 Sind $x_{1},\ldots,x_{n}\in
K$ mit 
\[
x_{1}^{e_{1}}+\ldots+x_{n}^{e_{n}}=0 \textrm{ mit }e_{i}\ge 1 \myforall i
\quad\textrm{und}\quad x_{i}x_{j}=0 \,\,\textrm{f"ur alle } 1\le i<j\le n,
\] 
so gilt
\[
x_{1}=\ldots=x_{n}=0.
\]
\end{Lemma}

\Bew F"ur $n=1$ folgt die Behauptung aus $x_{1}^{e_{1}}=0$. Wegen
$x_{i}x_{n}=0$ f"ur alle $i<n$ k"onnen wir
O.E. $x_{n}=0$ voraussetzen. Die restlichen Gleichungen stellen genau die
Voraussetzung an $x_{1},\ldots,x_{n-1}$ im Fall $n-1$ dar und liefern durch
Induktion $x_{1}=\ldots=x_{n-1}=0$. \qed\\

{\it\noindent Beweis von Satz \ref{zeroDimIdealsVonSL2}.} {\it 1. Schritt.}
Wir zeigen zun"achst $\dim (f_{3},\ldots,f_{2n-1})_{R}=0$.

Sei $P:=K[G_{ij}: 1\le i< j\le n]$ der Polynomring mit ${n \choose 2}$
unabh"angigen Variablen $G_{ij}$. Wir schreiben zus"atzlich $G_{ji}:=-G_{ij}$ f"ur $i<j$ 
Nach de Concini und Procesi \cite{ConciniProcesi} ist auch in
positiver Charakteristik der Einsetzungshomomorphismus
\[
\phi: P\rightarrow R,\quad  G_{ij}\mapsto g_{ij}
\]
surjektiv, und $\ker \phi=:J$ wird erzeugt von den Pl"ucker-Relationen
\[
G_{ij}G_{kl}-G_{ik}G_{jl}+G_{il}G_{jk} \quad \textrm {mit }i,j,k,l
\textrm{ paarweise verschieden}.
\]
Wir schreiben $F_{k}:=\sum_{i+j=k, i<j}G_{ij}^{e_{ij}}$ f"ur
$k=3,\ldots,2n-3$, so dass $\phi(F_{k})=f_{k}$, und setzen
\[
I:=J+(F_{3},\ldots,F_{2n-1})_{P}.
\]
Dann ist $P/I\cong R/(f_{3},\ldots,f_{2n-1})_{R}$. Wir zeigen, dass die Nullstellenmenge $\mathcal{V}(I)$ in $K^{{n\choose 2}}$
nur aus der $0$ besteht. Dann ist $0=\dim P/I=\dim
R/(f_{3},\ldots,f_{2n-1})_{R}$. 

Offenbar ist $0\in \mathcal{V}(I)$. Sei nun
umgekehrt $x=(x_{12},x_{13},\ldots,x_{n-1,n})\in \mathcal{V}(I)$. Da $F_{3}=G_{12}^{e_{12}},
F_{4}=G_{13}^{e_{13}}\in I$, folgt sofort $x_{12}=x_{13}=0$. Dies ist der
Induktionsanfang ($k=5$) der folgenden

{\bf Behauptung:} Es gilt
\begin{equation}\label{BehauptXij}
x_{ij}=0\quad \textrm{f"ur }1\le i<j \le n, \,\, i+j\le k-1.
\end{equation}
Wir zeigen, dass diese Aussage dann auch f"ur $k+1$ statt $k$ gilt. Dazu
gen"ugt es, $x_{1,k-1}=x_{2,k-2}=\ldots=0$ zu zeigen, und hierf"ur gen"ugt es
zu zeigen, dass die involvierten Variablen die Voraussetzung von Lemma
\ref{x1x2xnEqZero} erf"ullen. Die erste ben"otigte Gleichung folgt
sofort aus $F_{k}(x)=0$. Seien $i_{1}+j_{1}=k=i_{2}+j_{2}$ mit $i_{1}<j_{1},
\,\, i_{2}<j_{2}$ und $i_{1}<i_{2}$. Wir m"ussen
$x_{i_{1},j_{1}}x_{i_{2},j_{2}}=0$ zeigen. Die Pl"ucker-Relation liefert
\begin{equation}\label{Pluckerxij}
x_{i_{1},j_{1}}x_{i_{2},j_{2}}-x_{i_{1},i_{2}}x_{j_{1},j_{2}}+x_{i_{1},j_{2}}x_{j_{1},i_{2}}=0,
\end{equation}
wobei wir auch hier $x_{ji}:=-x_{ij}$ f"ur $i<j$ setzen. Aus $i_{1}<i_{2}$
folgt $j_{1}>j_{2}$, und damit $i_{1}+j_{2}<i_{1}+j_{1}=k$. Damit ist nach
\eqref{BehauptXij} $x_{i_{1},j_{2}}=0$, und damit ist auch der dritte Summand in
\eqref{Pluckerxij} gleich $0$. 

Falls $i_{1}+i_{2}<k$ oder $j_{1}+j_{2}<k$, so
ist wieder nach \eqref{BehauptXij} der zweite Summand in \eqref{Pluckerxij}
ebenfalls gleich $0$, und damit auch der erste, was zu zeigen ist. 

Sei daher
jetzt $i_{1}+i_{2}\ge k$ und $j_{1}+j_{2}\ge k$. Da aber
$i_{1}+j_{1}+i_{2}+j_{2}=2k$, gilt beidemale Gleichheit. Dann ist also
$i_{1}+i_{2}=k=i_{1}+j_{1}$, also $i_{2}=j_{1}$ und genauso
$i_{1}=j_{2}$. Dann folgt aber $i_{1}=j_{2}>i_{2}=j_{1}$, im Widerspruch zu
$i_{1}<j_{1}$. 

Der letzte Fall tritt also nicht auf, und wir haben
$x_{i_{1},j_{1}}x_{i_{2},j_{2}}=0$ gezeigt. Aus Lemma \ref{x1x2xnEqZero} folgt
dann die Behauptung \eqref{BehauptXij} f"ur $k+1$ statt $k$, und damit insgesamt
$\mathcal{V}(I)=0$. Dies zeigt dann
 \[\dim (f_{3},\ldots,f_{2n-1})_{R}=0.\]
{\it 2. Schritt.} Da $R$ ein Integrit"atsring ist, folgt aus der
Nulldimensionalit"at des Ideals  \[\dim R=\height
(f_{3},\ldots,f_{2n-1})_{R}.\] Nach Korollar \ref{exakteDimOfSL2Invs}
gilt $\dim R= 2n-3$. F"ur $e_{ij}=1$
(bzw. geeignete Graduierung von $R$) bilden daher die $2n-3$ Elemente
$f_{3},\ldots,f_{2n-1}$ nach Lemma \ref{phsopHeight} ein hsop in $R$. \qed\\

\subsection{Anwendung auf $S(\langle X^{p}, Y^{p}\rangle \oplus \bigoplus_{i=1}^{k} \langle X, Y\rangle)^{\SL_{2}}$}
Wir setzen nun\[
\langle \tilde{X},\tilde{Y} \rangle :=\langle X^{p},Y^{p}\rangle\]
und 
\[
V:=\langle \tilde{X}, \tilde{Y}\rangle \oplus \bigoplus_{i=1}^{k} \langle
X_{i}, Y_{i}\rangle\cong\langle X^{p}, Y^{p}\rangle \oplus \bigoplus_{i=1}^{k} \langle X, Y\rangle.
\]
Im Polynomring $S(V)=K[\tilde{X},\tilde{Y},X_{1},Y_{1},\ldots,X_{n},Y_{n}]$
sollen alle unabh"angigen Variablen den Grad $1$ haben.
Weiter betrachten wir
\[
U:=\langle X_{0}, Y_{0}\rangle \oplus \bigoplus_{i=1}^{k} \langle
X_{i}, Y_{i}\rangle,
\]
wobei alle $\langle X_{i},Y_{i}\rangle$ f"ur $i=0,\ldots,k$ Kopien der
nat"urlichen Darstellung sind. In
$S(U)=K[X_{0},Y_{0},X_{1},Y_{1},\ldots,X_{n},Y_{n}]$ w"ahlen wir die
Graduierung
\[
\deg X_{0}:=\deg Y_{0}:=\frac{1}{p} \quad \textrm{und}\quad \deg X_{i}:=\deg
Y_{i}=1 \,\,\textrm{ f"ur }i\ge 1.
\]
Wer sich bei den gebrochen rationalen Graden unwohl f"uhlt, kann immer alle Grade
mit $p$ durchmultiplizieren.

Nach Satz \ref{IsomOfFrobInvs} haben wir die Isomorphie\[
S(V)^{\SL_{2}}\cong S(U)^{\SL_{2}} \cap
K[X_{0}^{p},Y_{0}^{p},X_{1},Y_{1},\ldots,X_{n},Y_{n}]=:R,
\]
gegeben durch Ersetzen von $\tilde{X}\mapsto X_{0}^{p},\tilde{Y}\mapsto
Y_{0}^{p}$ und Beibehalten der restlichen Variablen $X_{i},Y_{i}$. Mit der
gew"ahlten Graduierung ist dieser Isomorphismus sogar graderhaltend. Wir
w"ahlen nun f"ur $S(U)^{\SL_{2}}$ das hsop aus Satz \ref{zeroDimIdealsVonSL2}, wobei $n=k+1$ und wir
$X_{0},Y_{0}$ f"ur $X_{k+1},Y_{k+1}$ einsetzen. Mit unserer Graduierung ist
dann
\begin{eqnarray*}
f_{3}&=&X_{1}Y_{2}-X_{2}Y_{1}\\
f_{4}&=&X_{1}Y_{3}-X_{3}Y_{1}\\
f_{5}&=&(X_{1}Y_{4}-X_{4}Y_{1})+(X_{2}Y_{3}-X_{3}Y_{2})\\
&\vdots&\\
f_{k+1}&=&(X_{1}Y_{k}-X_{k}Y_{1})+(X_{2}Y_{k-1}-X_{k-1}Y_{2})+(X_{3}Y_{k-2}-X_{k-2}Y_{3})+\ldots\\
f_{k+2}&=&(X_{0}^{p}Y_{1}^{p}-X_{1}^{p}Y_{0}^{p})^{2}+(X_{1}Y_{k+1}-X_{k+1}Y_{1})^{p+1}+(X_{2}Y_{k-1}-X_{k-1}Y_{2})^{p+1}+\ldots\\
f_{k+3}&=&(X_{0}^{p}Y_{2}^{p}-X_{2}^{p}Y_{0}^{p})^{2}+(X_{2}Y_{k+1}-X_{k+1}Y_{2})^{p+1}+(X_{3}Y_{k-1}-X_{k-1}Y_{3})^{p+1}+\ldots\\
&\vdots&\\
f_{2k+1}&=&(X_{0}^{p}Y_{k}^{p}-X_{k}^{p}Y_{0}^{p})
\end{eqnarray*}
ein hsop f"ur $S(U)^{\SL_{2}}$ mit $\deg f_{i}=2$ f"ur $i=3,\ldots,k+1$,  $\deg f_{i}=2(p+1)$ f"ur
$i=k+2,\ldots,2k-1$ und $\deg f_{2k}=\deg f_{2k+1}=p+1$. Dann ist $S(U)^{\SL_{2}}$ also ein endlich erzeugter
$A:=K[f_{3},\ldots,f_{2k+1}]$-Modul. Da $A\subseteq R\subseteq S(U)^{\SL_{2}}$ und $A$
ein noetherscher Ring ist, ist also auch der $A$-Untermodul $R$ von $S(U)^{\SL_{2}}$
endlich erzeugt. Insbesondere ist also $f_{3},\ldots,f_{2k+1}$ ein hsop von
$R$. Ersetzt man nun jeweils $X_{0}^{p}$ durch $\tilde{X}$ und $Y_{0}^{p}$ durch
$\tilde{Y}$, so erh"alt man ein hsop f"ur $S(V)^{\SL_{2}}$. Um den Test auf
Regularit"at schneller zu machen, war es in der Praxis n"utzlich, die Elemente
$f_{k+2},\ldots,f_{2k-1}$ nicht wie angegeben zum Grad $2(p+1)$ zu homogenisieren (also in den
geklammerten Termen die Exponenten $2,p+1,\ldots,p+1$ wegzulassen) - man kann die
Tiefe ja auch mit nicht homogenen regul"aren Sequenzen messen, solange die
Testsequenz nur im maximalen homogenen Ideal $R_{+}$ liegt. 

Wir f"uhren nun
also folgende Schritte durch:

\begin{Algorithmus}\label{AlgorithmDepthXpYp}
 Untere Schranke f"ur die Tiefe von $K[\langle X^{p}, Y^{p}\rangle \oplus \bigoplus_{i=1}^{k} \langle X, Y\rangle]^{\SL_{2}}$.
\begin{enumerate}
\item Berechne Generatoren von \[
K[\{X_{i}Y_{j}-X_{j}Y_{i}: 0\le i<j\le n\}]\cap
K[X_{0}^{p},Y_{0}^{p},X_{1},\ldots,Y_{k}]
\]
mit Algorithmus \ref{AlgIntersectXpY}. Ersetze $X_{0}^{p}$ und $Y_{0}^{p}$
durch $X_{0}$ und $Y_{0}$, um den n"achsten Schritt zu beschleunigen.
\item Berechne des Relationenideal $I\unlhd P$ der Generatoren aus Schritt 1,
  wobei der Polynomring $P$ soviele Variablen hat, wie es Generatoren gibt.
\item Bilde in $P$ eine Testsequenz, die obigen $f_{3},\ldots,f_{2k+1}$
  entspricht (evtl. ab $f_{k+2}$ nicht mehr homogenisiert). 
\item Untersuche diese Testsequenz mit Algorithmus \ref{ScanReg}, und erhalte
  eine regul"are Sequenz und eine untere Schranke f"ur die Tiefe.
\end{enumerate}
\end{Algorithmus}

Wir haben hier das Verfahren f"ur die $\SL_{2}$ angegeben. In der Praxis wird
man vorher noch Roberts' Isomorphismus anwenden, um sich im ersten
Schritt zwei Variablen zu sparen.

\subsection{Ergebnisse der Untersuchung}
Wir geben hier die Ergebnisse wieder, die wir durch Anwenden der Methode aus
dem letzten Abschnitt mit Hilfe von \Magma erhalten haben. In den drei F"allen,
in denen $K[V]^{\SL_{2}}$ dann nicht Cohen-Macaulay war, war jeweils
$f_{3},f_{4},f_{k+2},\ldots,f_{2k+1}$ eine regul"are Sequenz der L"ange $k+2$,
was dann eine untere Schranke f"ur die Tiefe ist. 
Insbesondere war damit die obere Schranke f"ur die Tiefe aus Satz
\ref{ZweitesHauptResultat} scharf, und damit auch die Schranke f"ur den Cohen-Macaulay-Defekt. Dies l"asst vermuten, dass sie es im Fall
$k\ge 3$ immer ist.

Man beachte, dass gen"ugende
gro"se Potenzen der $f_{3},\ldots,f_{k+1}$ unter Roberts' Isomorphismus aus
Summen von Annullatoren eines Kozyklus bestehen (siehe Abschnitt
\ref{AnnulatorenDesKoz}, insbesondere Lemma \ref{AnnullatorTyp3}), also selbst Annullatoren sind; Nach dem Hauptsatz
\ref{BigMainTheorem} l"asst sich also $f_{3},f_{4}$ durch kein Element aus
$f_{5},\ldots,f_{k+1}$ zu einer regul"aren Sequenz erg"anzen (eine Folge von
homogenen Elementen ist genau dann regul"ar, wenn ihre beliebig potenzierten
Folgeglieder regul"ar sind, siehe etwa \cite[Lemma 2.30]{Diplomarbeit}). Damit
ist klar, dass die  in den untersuchten F"allen  ausgew"ahlte regul"are
Sequenz $f_{3},f_{4},f_{k+2},\ldots,f_{2k+1}$  die "`maximal m"ogliche"'
 ist.

 In \Magma haben wir immer mit $\Ga$ gerechnet.
Die so erhaltenen exakten Werte f"ur den Cohen-Macaulay-Defekt
finden sich in der folgenden Tabelle:\\
\begin{center}
\begin{tabular}{c|ccc}
p$\backslash$k&2&3&4\\
\hline
2&0&1&2\\
3&0&1&-\\
5&0&-&-
\end{tabular}\\
{Werte von $\cmdef K[\langle X^{p}, Y^{p}\rangle \oplus \bigoplus_{i=1}^{k} \langle X, Y\rangle]^{\Ga}$.}
\end{center}

\noindent Hier die ben"otigten Gesamtrechenzeiten f"ur Algorithmus
\ref{AlgorithmDepthXpYp}, wobei die $f_{k+2},\ldots,f_{2k-1}$ nicht
homogenisiert waren; Mit Homogenisierung hatte die Rechenzeit etwa
im Fall $(p,k)=(2,4)$ 11 Tage (!) betragen, also deutlich l"anger als ohne die Homogenisierung.
\begin{center}
\begin{tabular}{c|ccc}
p$\backslash$k&2&3&4\\
\hline
2&0.11s&2.29s&4.23h\\
3&1.13s&3.73h&-\\
5&50.4s&-&-\\
\end{tabular}\\
{Gesamtlaufzeiten von Algorithmus \ref{AlgorithmDepthXpYp} (f"ur $\Ga$).}
\end{center}
Der Gro"steil der Rechenzeit wurde dabei f"ur die Berechnung des
Relationenideals ben"otigt. Der Gesamtspeicherbedarf betrug etwa 30MB.
F"ur andere Parameter habe ich die Rechnung nach zu gro"ser Rechenzeit ($>1$
Monat) abgebrochen.\\

\noindent Die folgende Tabelle gibt die mit Roberts' Isomorphismus
"ubersetzten Werte an:
\begin{center}
\begin{tabular}{c|ccc}
p$\backslash$k&3&4&5\\
\hline
2&0&1&2\\
3&0&1&-\\
5&0&-&-
\end{tabular}\\
{Werte von $\cmdef K[\langle X^{p}, Y^{p}\rangle \oplus \bigoplus_{i=1}^{k} \langle X, Y\rangle]^{\SL_{2}}$.}
\end{center}

\subsection{Bemerkung zur Buchsbaum-Eigenschaft}
Ein graduierter Ring $R$ hei"st \emph{Buchsbaum}, \index{Buchsbaum} wenn jedes phsop eine
\emph{schwach regul"are Sequenz}\index{regul\"are Sequenz!schwach} ist. Dabei hei"st eine Folge
$a_{1},\ldots,a_{k}\in R$ schwach regul"ar, wenn f"ur alle $i\le k$ und $m\in
R$ mit $ma_{i}\in(a_{1},\ldots,a_{i-1})$ auch $mR_{+}\in
(a_{1},\ldots,a_{i-1})$ gilt. (Regularit"at w"urde dagegen wegen
$m\in(a_{1},\ldots,a_{i-1})$ sogar  $mR\in
(a_{1},\ldots,a_{i-1})$ implizieren).
In Buchsbaum-Ringen gilt die Eigenschaft, dass jedes hsop die Tiefe
misst, d.h. wenn $a_{1},\ldots,a_{n}\in R$ ein hsop ist und $k=\depth R$, so
ist $a_{1},\ldots,a_{k}$ eine regul"are Sequenz (siehe etwa \cite[section 5]{CampEtAl}). Wenn $R=K[V]^{G}$ die
Voraussetzungen des Hauptsatzes \ref{BigMainTheorem} mit $k\ge 3$ erf"ullt und
wenn zus"atzlich $\depth R>2$ gilt, so kann damit $R$ nicht Buchsbaum sein,
denn dann ist das phsop $a_{1},a_{2},a_{3}$ keine regul"are Sequenz - kein
hsop, dass mit $a_{1},a_{2},a_{3}$ beginnt, kann also die Tiefe messen. Mit
dem Hauptsatz d"urfte es also extrem schwer fallen, einen nicht-Cohen-Macaulay
Invariantenring zu konstruieren, der Buchsbaum ist. F"ur Invariantenringe
endlicher Gruppen gilt sogar, dass beide Eigenschaften "aquivalent
sind (vermutet in \cite[Conjecture 27]{CampEtAl}, bewiesen in Kemper
\cite[Theorem 3.4]{KemperLoci}).

Da die
im letzten Abschnitt untersuchten nicht Cohen-Macaulay Ringe alle eine Tiefe
gr"o"ser als $2$ hatten, ist jedenfalls keiner von ihnen Buchsbaum.

\newpage
\addcontentsline{toc}{section}{Literatur}
\bibliographystyle{plain}
\bibliography{Diss}
\newpage
\addsec{Notation}
\thispagestyle{empty}
Die Eintr"age sind thematisch geordnet.\\

\begin{tabular}{cl}
$\cmdef$&Cohen-Macaulay-Defekt\\
$\depth$&Tiefe\\
$\dim$&(Krull-)Dimension\\
$\height$&H"ohe\\
$K$&algebraisch abgeschlossener K"orper\\
$p$& Charakteristik von $K$, $p=\chr K$\\
$R\subseteq S$&affine (graduierte) Algebren\\
$\wp\lhd R,\mathcal{P}\lhd S$&Primideale\\
$R_{+}$&maximales homogenes Ideal\\
$M$&$R$-Modul\\
$\Ass_{R}M$&Menge der assoziierten Primideale von $M$ in $R$\\
$G$&(lineare algebraische) Gruppe\\
$KG$& Gruppenring\\
$G^{0}$& Zusammenhangskomponente des Einselements von $G$\\
$\Ga$&additive Gruppe $(K,+)$\\
$\Gm$&multiplikative Gruppe $(K\setminus\{0\},\cdot)$\\
$\sigma\in G$&Gruppenelement\\
$\iota\in G$&neutrales Element von $G$\\
$V,W$&$G$-Modul\\
$\Hom_{K}(V,W)_{0}$&S. \pageref{HomKVW0}\\
$\Hom_{G}(V,W)$& S. \pageref{homgvwdef}\\
$V=\langle X_{1},\ldots,X_{n}\rangle$&$G$-Modul mit geordneter Basis $\{X_{1},\ldots,X_{n}\}$\\
$\langle X,Y\rangle$& nat"urliche Darstellung der $\SL_{2},\GL_{2}$  oder $\Ga$\\
$A_{\sigma}$&Darstellungsmatrix eines $G$-Moduls bzgl. der angegebenen Basis\\
$X$&$G$-Variet"at\\
$K[V]$&Polynomring, Ring der Polynomfunktionen auf dem $G$-Modul $V$\\
$K[X]$&Koordinatenring, Ring der Polynomfunktionen auf der ($G$)-Variet"at
$X$\\
&oder: $K[X]:=K[X_{1},\ldots,X_{n}]$ Polynomring mit $n$
unabh"angigen Variablen\\
$K[V]^{G}, K[X]^{G}$&Invariantenring\\
$K(V)$&K"orper der rationalen Funktionen auf $V$\\
$K(V)^{G}$&Invariantenk"orper\\
$F^{p}(V)$&$p$-te Frobenius Potenz von $V$, S. \pageref{FpV}\\
$S^{p}(V)$&$p$-te symmetrische Potenz von $V$\\
$C^{n}(G,V)$&Gruppe der $n$-Koketten mit Werten in $V$\\
$B^{n}(G,V)$&Gruppe der $n$-Kor"ander mit Werten in $V$\\
$Z^{n}(G,V)$&Gruppe der $n$-Kozyklen mit Werten in $V$\\
$H^{n}(G,V)$&$n$-te Kohomologiegruppe mit Werten in $V$\\
$P_{n}$&Freie Moduln der bar resolution\\
$d_{n}: P_{n}\rightarrow P_{n-1}$&Differentiale der bar resolution\\
$g\in Z^{n}(G,V)$&$n$-Kozyklus\\
$\partial_{n}^{V}$&S. \pageref{DefVonPartnV}\\
$\tilde{V}$& Erweiterung des $G$-Moduls $V$ durch einen $1$-Kozyklus $g\in
Z^{1}(G,V)$, S. \pageref{firstCohom}
\end{tabular}

\newpage
\addcontentsline{toc}{section}{Index}
\printindex
\end{document}